\documentclass{article}
\usepackage{amsmath,amssymb,bm} 
\usepackage[all]{xy}
\usepackage{caption,graphicx}
\font\ibf=cmbxti10

\title{Geometric simple connectivity and finitely presented groups}
\author{Valentin {\sc Po\'enaru}\footnote{Professor Emeritus, Universit\'e Paris-Sud, UMR 8628 du CNRS, Math\'ematiques, B\^atiment 425, 91405 Orsay Cedex, France. e-mail: valpoe@hotmail.com}}
\date{\sl (31 mars 2014)}

\begin{document}

\maketitle

\vglue 1cm

\setcounter{section}{-1}
\section{Introduction}\label{sec0}
\setcounter{equation}{0}

This is the second one of our series of papers giving a complete proof that {\ibf all finitely presented groups $\Gamma$ have the property QSF} of S.~Brick, M.~Mihalik and J.~Stallings \cite{3}, \cite{35}. The first one in the series is \cite{29} and then together with \cite{27}, where the result was first announced they should serve as introduction and background for the whole series (which will eventually consist of three papers). The main step in all this story is the {\ibf GSC-theorem} stated below and proved in the present paper. We will rely here heavily on \cite{29} and also on \cite{31}. That last paper was really written in the context of smooth 3-manifolds and their fundamental groups, but it so happens that it is taylor-made for extension to arbitrary groups $\Gamma$.

\smallskip

In order to make this present introduction independently readable, some concepts from \cite{27}, \cite {29} will have to be reviewed now, trying to keep the repetitions at minimum length.

\smallskip

Any group $\Gamma$ is, from now on, finitely presented. Our favourite presentation for $\Gamma$ will be like in section~II of \cite{29}, namely $\Gamma = \pi_1 \, M (\Gamma)$, where $M (\Gamma)$ is a compact singular $3$-manifold, with singularities ${\rm Sing} \, M (\Gamma) \subset M (\Gamma)$ which, by definition, are {\ibf immortal}. We do work with an arbitrary finitely presented group, but then we will be very choosy when it comes to the presentations of $\Gamma$, our $M(\Gamma)$'s.

\smallskip

The basic notion of our approach is that of REPRESENTATION for $\Gamma$ or for $\widetilde M (\Gamma)$ which, remember, up to quasi-isometry is the same thing as $\Gamma$ itself. Our representations are, by definition, nondegenerate simplicial maps
\begin{equation}
\label{eq0.1}
X \overset{f}{\longrightarrow} \widetilde M (\Gamma) \, .
\end{equation}

By definition, the set of points ${\rm Sing} \, (f) \subset X$ where $f$ fails to be an immersion, are the, so-called {\ibf mortal}, singularities of $f$. The following three features will be required for the representation (\ref{eq0.1}) above.

\medskip

I) $X$ is a not necessarily locally finite simplicial complex (or cell-complex) of dimension two or three, which is {\ibf geometrically simply connected} (GSC). This means that the $1$-cells and $2$-cells are in cancelling position. The GSC notion includes, as a special case, arborescence, i.e. being gettable from a point by Whitehead dilatations. But, in general, our $X$'s come with $\pi_2 \ne 0$.

\medskip

II) We start by considering now two equivalence relations on $X$, namely the
$$
\Psi (f) \subset \Phi (f) \subset X \times X \, ,
$$
where $\Phi (f) \ni (x,y)$ simply means $fx=fy$ while, by definition, $\Psi (f)$ is the smallest equivalence relation compatible with $f$, such that the natural map $X / \Psi (f) \to \widetilde M (\Gamma)$ is immersive. With the $\Phi (f) \supset \Psi (f)$ so defined, in the case of a representation (\ref{eq0.1}) we insist that 
$$
\Psi (f) = \Phi (f) \, .
$$

\smallskip

III) The map $f$ is ``essentially surjective'', a notion for which the exact definition is given in \cite{29}. Anyway, it should be clear that $\dim X = 2$ or $3$ and this $\dim X$ is, by definition, the {\ibf dimension} of the representation. In the $3$-dimensional case, our essential surjectivity implies that $\overline{fX} = \widetilde M (\Gamma)$ and in the $2$-dimensional case it implies that $\widetilde M (\Gamma)$ can be gotten from the closure $\overline{fX}$, by adding cells of dimensions $\lambda = 2$ and $\lambda = 3$.

\medskip

Concerning the general features of $\Psi / \Phi$ one should see \cite{22}. There they are discussed in the special case when $\Gamma = \pi_1 \, M^3$ with $M^3$ a smooth 3-manifold and $\widetilde M (\Gamma) = \widetilde M^3$, but that little theory generalizes immediately for our arbitrary $\Gamma$'s. See here \cite{29} too.

\smallskip

There is also an alternative definition for $\Psi (f)$. This $\Psi (f)$ is the smallest equivalent relation on $X$, compatible with $f$, which kills all the mortal singularities. As it is explained in \cite{29}, the quotient-space $X / \Psi (f)$ is realized by {\ibf zipping}, which means closing in, starting from the mortal singularities of $f$, all the continuous paths inside the following set
$$
\widehat M^2 (f) \underset{\rm def}{=} M^2 (f) \cup {\rm Sing} \, (f) \subset X \times X \, ,
$$
where $M^2 (f) \subset X \times X$ is by definition the set of the $(x,y)$, $x \ne y$, $fx = fy$, and when by ``${\rm Sing} \, (f)$'' we really mean here ${\rm Diag} \, ({\rm Sing} \, (f))$. A typical zipping path is drawn in Figure~1.1 from \cite{29}, reappearing in this paper as figure~5.3. But catching any specific $(x,y) \in M^2 (f)$ by the zipping process, is really a matter of {\ibf strategy}, not to be discussed now. Also, we have introduced above the double point set $M^2 (f) \subset X \times X$, but we will also need its other incarnation $M_2 (f) \subset X$, consisting of those $x \in X$ such that ${\rm card} \, f^{-1} \, fx > 1$.

\smallskip

While the immortal singularities concern the source $X$, the mortal ones concern the map $f$, i.e. rather the target.

\smallskip

Notice that our present representations (\ref{eq0.1}), with target $\widetilde M (\Gamma) \sim \Gamma$ (quasi-isometry) take a dual form with respect to the usual representations $\Gamma \to \{$some other group$\}$, since the map in (\ref{eq0.1}) goes rather like $\to \Gamma$. But then, this opens the door for the interesting possibility of enforcing (if and when possible) a free action of $\Gamma$ on $X$, and then asking for the map $f$ to be equivariant. This is how group theory really enters the scene here. A priori, representations can be defined (and have been used, for intance in \cite{18}, \cite{19}, \cite{32}) for targets other than our present $\widetilde M (\Gamma)$. The conditions I) to III) actually force, automatically, the target of the representation to be simply-connected.

\smallskip

In \cite{29}, first paper of this series, the main result was a $3^{\rm d}$ representation theorem for arbitrary $\Gamma$. The first result of the present paper will be an analogous $2^{\rm d}$ representation theorem, essentially gotten by taking a very dense skeleton of the $3^{\rm d}$ representation space from \cite{29}, plus some other necessary refinements. This is the following

\bigskip

\noindent 2-DIMENSIONAL REPRESENTATION THEOREM. {\it For any finitely presented group $\Gamma$ there is a $2$-dimensional representation
\begin{equation}
\label{eq0.2}
X^2 \overset{f}{\longrightarrow} \widetilde M (\Gamma) \, ,
\end{equation}
with the following additional features.}

\medskip

1) {\it (First finiteness condition.) The representation space $X^2$ is {\ibf locally finite} and, inside it, ${\rm Sing} \, (f) \subset X^2$ is discrete, i.e. it has no accumulation points at finite distance.}

\medskip

2) {\it There is a free action $\Gamma \times X^2 \to X^2$ such that $f$ is {\ibf equivariant}}
$$
gf(x) = f(gx) \, , \qquad \forall \, g \in \Gamma \, , \ x \in X^2 \, .
$$

3) {\it When appropriate strategies are being chosen, then there is a {\ibf uniformly bounded zipping length}. More explicitly, there is a uniform bound $0 < M < \infty$, such that for any double point $(x,y) \in M^2 (f)$, we can find a zipping path $\lambda (x,y) \subset \widehat M^2 (f)$, with
$$
{\rm length} \, \lambda (x,y) \leq M \, .
$$
[The length in question is well-defined up to quasi-isometry.]}

\medskip

4) {\it (Second finiteness condition.) We move now from the double points $M^2 (f) \subset X^2 \times X^2$ to the double points in $X^2$, $M_2 (f) \subset X^2$. For any tight compact transversal $\Lambda$ to $M_2 (f) \subset X^2$ we have then}
$$
{\rm card} \, (\lim (\Lambda \cap M_2 (f))) < \infty \, .
$$

\smallskip

5) {\it The closed subset below, where the union is over all compact tight transversals,
\begin{equation}
\label{eq0.3}
{\rm LIM} \, M_2 (f) \underset{\rm def}{=} \ \bigcup_{\Lambda} \ \lim (\Lambda \cap M_2 (f)) \subset X^2
\end{equation}
is a locally finite graph and, moreover, the
$$
f \, {\rm LIM} \, M_2 (f) \subset f \, X^2
$$
is also a closed subset.}

\bigskip

The first section of the present paper, strongly relying on \cite{29} gives the full proof of this theorem. Concerning our $2^{\rm d}$ representation theorem, we CANNOT assume, a priori, that $M_2 (f) \subset X^2$ is a {\ibf closed} subset and, once  this is so, our second finiteness condition above is the next best thing which could happen. One should be aware that, in the absence of special measures (see here the ``decantorianization process'' in \cite {29} and in \cite{31}), the accumulation pattern of $\Lambda \cap M_2 (f)$ may become chaotic, like in \cite{32}, where Julia sets pop up.

\smallskip

The vilain in our whole story is the set ${\rm LIM} \, M_2 (f)$. There is an equivalence
\begin{equation}
\label{eq0.4}
{\rm LIM} \, M_2 (f) = \phi \Longleftrightarrow M_2 (f) \subset X^2 \ \mbox{is closed,}
\end{equation}
and we will call {\ibf ``easy''}, any $\Gamma$ admitting a locally finite $2^{\rm d}$ representation which verifies (\ref{eq0.4}), i.e. which comes with a closed $M_2 (f)$ and which, in addition to this, is s.t. $fX^2$ is a closed subset of $\widetilde M (\Gamma)$. No additional conditions like equivariance and uniformly bounded zipping length are required here and I suspect that, short of weakening the GSC condition I) in the definition of REPRESENTATIONS, given above, generically, we cannot enforce them without creating some ${\rm LIM} \, M_2 (f) \ne \emptyset$. These conditions of equivariance and bounded zipping length, are essential for proving the $\forall \, \Gamma \in {\rm QSF}$, and then we have to live with the possibility of a non closed $M_2 (f)$. But as papers like \cite{24}, \cite{25} abundantly show, representations with ${\rm LIM} \, M_2 (f) = \emptyset$, even without any other embellishment, do quite naturally occur and can be very useful. Groups $\Gamma$ for which any representation exhibits a ${\rm LIM} \ne \emptyset$, will be called {\ibf difficult}. One should not give to our definitions of difficult and easy any other connotation beyong their purely technical meaning.

\smallskip

One should notice that, once $M_2 (f)$ is not closed, neither the $fX^2$ (contrary to $X^2$) nor the $f \, {\rm LIM} \, M_2 (f)$, will be locally finite, sources of much headache. Also, what generates a ${\rm LIM} \, M_2 (f) \ne \emptyset$ is the following basic phenomenon. Consider some generic, arbitrary representation (\ref{eq0.1}). The map $f$ wanders randomly through $\widetilde M (\Gamma)$, inspecting every nook and hook; in general, any given compact fundamental domain of $\widetilde M (\Gamma)$ will be hit infinitely many times by it. Thinking of the classical Whitehead manifold and possibly of the Casson Handles too \cite{28}, \cite{4}, \cite{5}, \cite{10}, I have called this kind of phenomenon, which does generate a ${\rm LIM} \, M_2 (f) \ne \emptyset$, the {\ibf Whitehead nightmare}. Losely speaking, the easy groups are those which manage to avoid the nightmare in question.

\smallskip

To the best of my knowledge, for every known group $\Gamma$ there is, one way or another, a proof that the $\Gamma$ in question is easy. But that might be a very difficult matter. For instance, with the present state of the art, the only way to prove it for a $\Gamma = \pi_1 \, M^3$, when $M^3$ is a closed 3-manifold, has to make use of the full Thurston geometrization of $3$-manifolds, i.e. use the very spectacular work of G.~Perelman \cite{15}, \cite{16}, \cite{17}, \cite{1}, \cite{2}, \cite{11}, \cite{12}, \cite{13}, \cite{38}. Strongly related to this, just to show that all $\pi_1 \, M^3$'s are QSF, besides the proof presented in this series of three papers and which eventually proves that {\ibf all} finitely presented groups are QSF, the only other known way is to invoke Perelman's work.

\smallskip

All this having been said, I do believe that making use of the results and techniques of the present series of papers, plus some other ingredients which I will not go into here, we can prove that all $\Gamma$'s are indeed easy. I hope to come back on this issue, another time soon. For our very present purposes I want to put down on record here, that there is a very easy implication (and ``easy'' is to be  taken here with its usual meaning)
\begin{equation}
\label{eq0.5}
\{\Gamma \ \mbox{is easy}\}  \Longrightarrow \{\Gamma \in \mbox{QSF}\} \, .
\end{equation}

This is actually explicitly proved in the joint paper with Daniele Otera \cite{36}. But I thought it would be still appropriate to sketch the proof of this very simple implication (\ref{eq0.5}) here. This should be a good preparation for the arguments developed in the present paper. It is also a good illustration of how, in our arguments, a mixture of dimension two, dimension three and high dimension, has to occur.

\smallskip

So, we are given now, for our supposedly easy $\Gamma$, a $2^{\rm d}$ representation (\ref{eq0.2}) where $M_2 (f) \subset X^2$ is a {\ibf closed} subset and where, moreover, the subset $fX^2 \subset \widetilde M (\Gamma)$ is also closed. All the other embellishments from the $2^{\rm d}$ representation theorem are now irrelevant and certainly not assumed to be there, (except when they are obviously with us, like the two finiteness conditions). Because $\Phi (f) = \Psi (f)$, we can break the big quotient-space projection, where the target is now a bona fide locally finite simplicial complex
$$
X^2 \longrightarrow fX^2 = X^2 / \Phi (f) = X^2 / \Psi (f) \, ,
$$
into an infinite sequence of elementary compact zipping operations, each of which is either a very simple-minded simple homotopy equivalence or, homotopically speaking, an equally simple-minded addition of a $2$-cell
\begin{equation}
\label{eq0.6}
X^2 \equiv X_0 \to X_1 \to X_2 \to \ldots \to fX^2 \, .
\end{equation}
These are, actually, the elementary moves $O(i)$ from \cite{8}, \cite{18}, \cite{19}.

\smallskip

One can show that, after thickening the $X_i$'s into smooth $n$-manifolds, where some $n \geq 5$ is fixed, the following happens (and then, below, concerning these matters some additional explanation will be given).

\medskip

A) The various arrows in (\ref{eq0.6}) become smooth embeddings, each of which is either a compact smooth Whitehead dilatation, or the addition of a $2$-handle
\begin{equation}
\label{eq0.7}
\Theta^n (X^2) = \Theta^n (X_0) \hookrightarrow \Theta^n (X_1) \hookrightarrow \Theta^n (X_2) \hookrightarrow \ldots
\end{equation}

\medskip

B) Since $M_2 (f) \subset X^2$ is closed, one can put together the (\ref{eq0.7}) and assemble it into a smooth $n$-manifold with large non-empty boundary $\underset{i=0}{\overset{\infty}{\bigcup}} \, \Theta^n (X_0)$ (with ``$\Theta^n$'' standing for $n$-dimensional thickening, but see then below too). Because $X^2 = X_0 \in {\rm GSC}$, so is $\underset{i=0}{\overset{\infty}{\bigcup}} \, \Theta^n (X_0)$ (the GSC being meant now in the DIFF category) {\ibf and} because $M_2 (f)$ is closed in $X^2$ and $fX^2$ is closed in $\widetilde M (\Gamma)$, we also have
\begin{equation}
\label{eq0.8}
\bigcup_{i=0}^{\infty} \ \Theta^n (X_{\rm I}) \underset{\rm DIFF}{=} \, \Theta^n (fX^2) \, .
\end{equation}

\medskip

C) Point III in the general definition of representations makes that $\Theta^n (\widetilde M (\Gamma))$ $= \Theta^n (fX^2) + \{$handles of index $2$ and $3\}$, and hence we also have $\Theta^n (\widetilde M (\Gamma)) \in {\rm GSC}$.

\smallskip

Our representation was not assumed $\Gamma$-equivariant, but the $\Theta^n (\widetilde M (\Gamma))$ on which $\Gamma$ acts freely, certainly is so. Moreover the fundamental domain $\Theta^n (\widetilde M (\Gamma)) / \Gamma$ is compact. It is now a standard fact that, if $Y$ is a locally finite simplicial complex which is GSC, comes with a free $\Gamma$-action such that $Y/\Gamma$ is compact, with $\pi_1 \, Y/\Gamma = \Gamma$, then $\Gamma \in {\rm QSF}$. End of the argument.

\medskip

Before giving the additional explanations concerning the A), B), C) above, here are some comments concerning the very trivial implication GSC $\Rightarrow$ QSF. To begin with, instead of GSC one can also use here the weaker notion WGSC ($=$ weak geometric simple connectivity, i.e. exhaustibility by simply-connected compacta) which has been thoroughly studied by L.~Funar and D.~Otera \cite{6}, \cite{7}, \cite{14}. Then, Funar and Otera also proved a sort of converse to GSC $\Rightarrow$ QSF, in the group-theoretical context: If $\Gamma \in {\rm QSF}$ then there exists a smooth compact manifold $N$ such that $\pi_1 N = \Gamma$ and $\widetilde N \in {\rm GSC}$. But then, while QSF is presentation independent, GSC is not.  I thing that the correct standpoint here is to think of QSF as being the presentation independent way of saying that $\Gamma \in {\rm GSC}$. And the GSC {\ibf is} a central concept.

\smallskip

Here are now some additional explanations concerning our little argument for proving (\ref{eq0.5}). To begin with here, the $n \geq 5$ needs to be explained and then, the whole line of proof clearly asks for a {\ibf canonical} $\Theta^n (X)$ too. This brings us to the following little digression which I would like to go into now. The little argument for the proof of (\ref{eq0.5}) which was sketched above is essentially inspired by the very initial, very easy part of my own approach to the $3^{\rm d}$ Poincar\'e Conjecture and I have here in mind both my papers \cite{18}, \cite{19}, living very up-stream in the approach in question, and also David Gabai's review paper \cite{8} which presents things as they stood about twenty years ago (for an up-date see here \cite{20}, \cite{21}). Anyway, in terms of \cite{18}, \cite{19}, \cite{8}, all the mortal singularities involved in (\ref{eq0.6}) are of the ``undrawable'' kind, the acyclic elementary  now are $O(i)$'s, with $0 \leq i \leq 2$, while the ones which homotopically add $2$-cells are $O(3)$'s.

\smallskip

But here comes a little subtelty which we need to discuss now. My whole approach to the $3^{\rm d}$ Poincar\'e Conjecture needs dimension four and there, in order to thicken singular objects like our $X$'s above, a {\ibf desingularization} $R$ has to be specified too. So one ends up with smooth 4-manifolds $\Theta^4 (X,R)$ and, in that context, when sequences like (\ref{eq0.6}) are changed  into $4$-dimensional versions of (\ref{eq0.7}) then there is a very serious obstruction for getting embeddings. In particular, only for the so-called COHERENT zippings, and this is now an $R$-dependent notion, are the $4$-dimensional thickenings of $O(3)$ moves embeddings, actually additions of $2$-handles. This coherence issue is very central to my $4$-dimensional approach to the Poincar\'e Conjecture. The paper \cite{8} should be a very clear exposition concerning this question. I have finally managed to kill the obstruction to coherence in 2006 (see here \cite{21}, which as said above, together with \cite{20} should give a good view of my whole program for the Poincar\'e Conjecture). Now, I certainly need dimension four in my Poincar\'e program but then, as soon as one takes the product with $[0,1]$ or with $[0,1] \times [0,1] , \ldots$ and goes to dimensions $n \geq 5$, then this washes away the $R$-dependence, and hence for $p \geq 1$ the $\Theta^{p+4} (X) = \Theta^4 (X,R) \times B^p$ are now all, {\ibf canonical}. The same stabilization kills the obstruction for coherence and it is these things which make our simple-minded proof of (\ref{eq0.5}) sketched above, work as soon as $n \geq 5$. This really ends the discussion of the implication (\ref{eq0.5}).

\smallskip

We turn back now to our general finitely presented group $\Gamma$, which we no longer can assume to be easy. We still have our $2^{\rm d}$ representation theorem for $\Gamma$ coming now with ${\rm LIM} \, M_2 (f) \ne \emptyset$ and with the nasty $2^{\rm d}$ object $fX^2 \subset \widetilde M (\Gamma)$. We would like now to put up, starting from our $2^{\rm d}$ representation theorem a proof that $\Gamma \in {\rm QSF}$, for an arbitrary finitely presented $\Gamma$, following very roughly speaking the same general plan as in the proof of (\ref{eq0.5}) just sketched. This, of course, can be so only in a very first approximation, once ${\rm LIM} \, M_2 (f) \ne \emptyset$.

\smallskip

The problem here is that, in this very general set-up, the $fX^2$ is now a very pathological space, which is non locally finite. This means that, once we drop the simplifying hypothesis that ``$\Gamma$ is easy'' and face a real-life ${\rm LIM} \, M_2 (f) \ne \emptyset$, then we get an infinitely rougher ride than in the previous proof of the toy-case (\ref{eq0.5}).

\smallskip

Since, contrary to the $fX^2$ from the easy case discussed above, the real-life $fX^2$ is not locally finite, there are no smooth regular neighbourhoods for it. There still is a redeeming feature. Together with the second finiteness condition in the $2^{\rm d}$ REPRESENTATION THEOREM stated above, comes the fact that the set of points where $fX^2$ fails to be locally finite, and I will denote now these points, generically, by $p_{\infty\infty} (\Gamma)$, like in the main body of this paper, is such that the subset $\{ p_{\infty\infty} (\Gamma)\} \subset fX^2$ is {\ibf discrete}. Starting from this, one can build up a very high-dimensional {\it sort of} thickening of $fX^2$ which, like in the paper, will be denoted by $S_u \, \widetilde M (\Gamma)$. The $S_u \, \widetilde M (\Gamma)$ is only a locally-finite cell-complex and NOT a smooth manifold. Of course, a cell-complex can always be thickened further into a manifold. But the singular, non-manifold points which $S_u \, \widetilde M (\Gamma)$ possesses are essential for the technology which we will develop. Of course, $S_u \, \widetilde M (\Gamma)$ is a smooth cell-complex, in the sense that it is put together from smooth cells, glued together by smooth maps, but topologically, it is not locally like $R^n$. Also, for technical reasons, $S_u \, \widetilde M (\Gamma)$ is not gotten directly from $fX^2$. We must go first to a sort of locally finite $3^{\rm d}$ thickening $\Theta^3 (fX^2)$ and next, like in the discussion which has come with the proof of (\ref{eq0.5}) above, to a $\Theta^4 (\Theta^3 (fX^2),R)$. Finally more additional dimensions are thrown in, to our $\Theta^4 (\Theta^3 (fX^2),R)$ in order to get $S_u \, \widetilde M (\Gamma)$.

\smallskip

And it is in those {\ibf additional dimensions} that the action of our proof takes place. I think it is useful to make here the distinction between more ``high dimensions'' and ``additional dimensions''.

\smallskip

So, we need high-dimensionality, became part of our constructions will take place inside the {\ibf supplementary dimensions}, i.e. those in addition to the four of $\Theta^4 (\Theta^3 (fX^2),R)$. In view of the high dimensions involved, the $R$-dependence is eventually washed away too.

\smallskip

Dimension two is necessary for us because it is there that the zipping is most transparent and high dimensions are necessary because there one proves the all-important GSC feature. But one goes from $d=2$ to high $d$ via the intermediary 3-dimensional and 4-dimensional steps. And then, in the third and last paper of this series when will go from the very high dimensional, infinitely foamy and violently non-co-compact $S_u \, \widetilde M (\Gamma)$ to $\Gamma \in {\rm QSF}$ one will need to trace back our steps, from high $d$ to $d=4$ and then to $d=3$.

\smallskip

Section~II of the present object will present the definition of this sort of thickening of $fX^2$, which we call $S_u \, \widetilde M (\Gamma)$. I want to review here the kind of features which this $S_u \, \widetilde M (\Gamma)$ will need to have and which this and the next paper in this series, will establish.

\medskip

1) $S_u \, \widetilde M (\Gamma)$ should come equipped with a free $\Gamma$-action and $\pi_1 (S_u \, \widetilde M (\Gamma) / \Gamma)$ $= \Gamma$.

\medskip

2) Actually, we will want the ``$S_u$'' occurring in $S_u \, \widetilde M (\Gamma)$ to be a functor of sorts, with good localization and glueing properties.  So ``$S_u$'' can be applied to other things than $\widetilde M (\Gamma)$, too. In particular, $S_u \, M (\Gamma)$ will have to make sense too, and we will very much need things like the functorial property
\begin{equation}
\label{eq0.9}
S_u \, \widetilde M (\Gamma) = (S_u \, M (\Gamma))^{\sim} \, .
\end{equation}

\medskip

3) Here comes now the item which is the hardest to get, in our whole list. We will want that $S_u \, \widetilde M (\Gamma) \in {\rm GSC}$.

\medskip

4) As a consequence of 1) and 2) we also have, of course, that
$$
S_u  \, \widetilde M (\Gamma) / \Gamma = S_u \, M (\Gamma) \, ,
$$
but it is important here that the $S_u \, M (\Gamma)$ can also be defined directly.

\smallskip

Once we have gotten this far, notice that if it would be the case that $S_u \, M (\Gamma)$ is compact, then $\Gamma \in {\rm QSF}$ would follow automatically. But $S_u \, M (\Gamma)$ is NOT compact. The best one might say here is that if one starts from $\widetilde M (\Gamma)$ as a union of compact fundamental domains, then $S_u \, \widetilde M (\Gamma)$ is gotten by changing these into {\it infinitely foamy, high dimensional, non-compact objects}. So here comes a fifth requirement to be added to our list. It will be proved in the third paper of this series.

\medskip

5) Notwithstanding the lack of a compact fundamental domain, our $S_u \, \widetilde M (\Gamma)$ should be good enough for the following implication to work
\begin{equation}
\label{eq0.10}
S_u \, \widetilde M (\Gamma) \in {\rm GSC} \Longrightarrow \Gamma \in {\rm QSF} \, .
\end{equation}
Let us elaborate a bit more on this sixth requirement, (\ref{eq0.10}). The QSF property of S.~Brick, M.~Mihalik and J.~Stallings \cite{3}, \cite{35} has among its ancestors the concept of {\ibf Dehn-exhaustibility} which, in the smooth category is the following. 

\smallskip

\noindent The manifold $V$ is, by definition, Dehn-exhaustible if for any compact $K \subset V$ there is a simply connected compact manifold $M$ with $\dim M = \dim V$ and a commutative diagram

\begin{equation}
\label{eq0.11}
\xymatrix{
K \ar[rr]_{j} \ar[dr]_-i &&M \ar[dl]^-{g}  \\ 
&V
}
\end{equation}
where $i$ is the standard embedding, $j$ also injects, $g$ is a smooth {\ibf immersion}, and where the following Dehn property holds: $M_2 (g) \cap jK = \emptyset$. Here $M$ is ``abstract'' and not some submanifold of $V$.

\smallskip

Dehn-exhaustibility occurred first in my old papers \cite{23} to \cite{25} as well as in Casson's work (see \cite{9}).

\smallskip

If one tries to introduce this notion for groups $\Gamma$ then, it is not presentation-independent, unlike the QSF. But just like for GSC and QSF there is also a weak equivalence between QSF and Dehn-exhaustibility \cite{14}.

\smallskip

At least when $\dim V = 3$ and $\partial V = \emptyset$, one can extract the following fact from my old paper \cite{23}, namely
\begin{equation}
\label{eq0.12}
\mbox{If $V \times B^p \in {\rm GSC}$ for some $p \geq 0$, then $V$ is Dehn-exhaustible.}
\end{equation}
Now, the proof of (\ref{eq0.10}), which occupies the third and last paper in this series, will require, among others, arguments like those used in \cite{23} for proving (\ref{eq0.12}). Never mind that now $\dim V$ is arbitrary, that $\partial V \ne \emptyset$ and that the issue of $\pi_1^{\infty} \, V$ being zero, which was looming big in \cite{23} appears now rather as a red herring. But it is essential to recognize here that in the context of (\ref{eq0.12}) it is very important that $B^p$ be a compact $p$-ball, and not something like ${\rm int} \, B^p$ or $B^p - \{$some boundary points$\}$. This is a key point, since any simply-connected $V$ can be rendered GSC by multiplying it with $R^p$ but not by multiplying it  with $B^p$. We abstract from this discussion the following {\ibf transversal compactness} requirement for our $S_u \, \widetilde M (\Gamma)$, namely that it should be a regular neighbourhood (with compact fiber) of some appropriate low(er) dimensional spine. 

\smallskip

Finally, we are able to state the main result which this present paper proves, except that the proof of the implication (\ref{eq0.10}) will still need another subsequent paper. But that one will be easier, several orders of magnitude easier than the proof of the theorem below, something like three months versus thirty years.

\smallskip

So here is the main result of the present paper.

\bigskip

\noindent THE GSC THEOREM. {\it For any finitely presented group $\Gamma$ one can construct a smooth locally finite cell-complex $S_u \, \widetilde M (\Gamma)$ with the properties $1)$ to $5)$ listed above, in particular such that $S_u \, \widetilde M (\Gamma)$ is} GSC.

\bigskip

\noindent {\ibf An overview of the proof of the GSC theorem, in a nutshell.} The general idea is the following. Start with the gigantic quotient space projection, which can be assumed to be an equivariant zipping
\begin{equation}
\label{eq0.13}
X^2 \overset{f}{\longrightarrow} fX^2 \subset \widetilde M (\Gamma) \, ,
\end{equation}
followed by the sort of high dimensional thickening which goes from $fX^2$ to $S_u \, \widetilde M (\Gamma)$. Since this is, a priori, not very promising as far as GSC is concerned, we try then the following road. Start now with a high dimensional thickening of $X^2$, which is automatically GSC. Then, like in the process of changing (\ref{eq0.6}) into (\ref{eq0.7}), we try to imitate the quotient-space projections by a sequence of smooth, GSC-preserving embeddings. Became ${\rm LIM} \, M_2 (f) \ne \emptyset$, putting things together is now an infinitely trickier affair than in the toy-model case of (\ref{eq0.7}). We will still manage to get a second functor (of sorts), call it now $S_b$. There are now fewer morphisms than in the case of $S_u$, but the analogue of the functorial formula (\ref{eq0.9}) is valid for $S_b$ too. What makes such an $S_b$ possible at all is that the zipping is equivariant. Actually equivariance is needed everywhere in our whole approach, and this is how the group structure of $\Gamma$ comes in. Similarly, compactness, whenever it may occur, reflects finite presentation for $\Gamma$.

\smallskip

Finally, not without a lot of fatigue, one can show that the smooth locally finite cell-complex $S_b \, \widetilde M (\Gamma)$, which is of the same dimension as $S_u \, \widetilde M (\Gamma)$ is actually GSC. Incidentally, the subscripts ``$u$'', ``$b$'' stand for ``usual'' and ``bizarre''. The big step is to go now from the known $S_b \, \widetilde M (\Gamma) \in {\rm GSC}$ to the desired $S_u \, \widetilde M (\Gamma) \in {\rm GSC}$, which seems unreachable by any direct assault, since no direct assault seems to be able to connect the two objects (let us say by a diffeomorphism).

\smallskip

Actually two more functors are needed here, $S'_u$ and $S'_b$ where the ``$'$'' signalizes an important technical switch allowing us to deal with some really very nasty pathologies involved in (\ref{eq0.13}). The $S_u \, \widetilde M (\Gamma)$, $S_b \, \widetilde M (\Gamma)$, $S'_u \, \widetilde M (\Gamma)$, $S'_b \, \widetilde M (\Gamma)$ are smooth cell-complexes, all of the same dimension. Also there are two transformations
\begin{equation}
\label{eq0.14}
S'_u \, \widetilde M (\Gamma) \Longrightarrow S_u \, \widetilde M (\Gamma) \, , \quad S'_b \, \widetilde M (\Gamma) \Longrightarrow S_b \, \widetilde M (\Gamma)
\end{equation}
which are isomorphic transformations, except that their two sources are not a priori known to be so.

\smallskip

At this level we move downstairs to $M (\Gamma)$ which {\ibf is} compact and to the two $S'_u \, M (\Gamma)$, $S'_b \, M (\Gamma)$ which are not. Making use of the compactness of $M (\Gamma)$ and of the uniform bound for the zipping length one can show that $S'_u \, M (\Gamma) =  S'_b \, M (\Gamma)$ and from there, by functoriality (see (\ref{eq0.9})) it follows that $S'_u \, \widetilde M (\Gamma) = S'_b \, \widetilde M (\Gamma)$ too. The equivariance of the whole construction is essential here, of course.

\smallskip

Now we finally use the isomorphism between the two transformations in (\ref{eq0.14}) and we hence manage to produce the desired diffeomorphism
\begin{equation}
\label{eq0.15}
S_u \, \widetilde M (\Gamma) \underset{\rm DIFF}{=} S_b \, \widetilde M (\Gamma) \, ,
\end{equation}
which in turn implies that $S_u \, \widetilde M^3 (\Gamma) \in {\rm GSC}$. The big diagram below schematizes the whole story with the three equalities (``$=$'') following in logical order from the bottom line to the top one
$$
\begin{matrix}
&\underset{\mbox{$\vert$} \qquad \qquad \qquad \qquad \mbox{=} \qquad \qquad \qquad \qquad \mbox{$\downarrow$}}{\!\!\!-\!\!\!-\!\!\!-\!\!\!-\!\!\!-\!\!\!-\!\!\!-\!\!\!-\!\!\!-\!\!\!-\!\!\!-\!\!\!-\!\!\!-\!\!\!-\!\!\!-\!\!\!-\!\!\!-\!\!\!-\!\!\!-\!\!\!-\!\!\!-\!\!\!-\!\!\!-\!\!\!-\!\!\!-\!\!\!-\!\!\!-\!\!\!-\!\!\!-\!\!\!-\!\!\!-\!\!\!-\!\!\!-\!\!\!-\!\!\!} \\
{ \ } \\ 
&S_u \, \widetilde M (\Gamma) \Longleftarrow S'_u \, \widetilde M (\Gamma) = S'_b \, \widetilde M (\Gamma) \Longrightarrow S_b \, \widetilde M (\Gamma) \\
{ \ } \\
&\downarrow \ \  \mbox{covering maps} \ \  \downarrow \\ 
{ \ } \\
&S'_u \, M (\Gamma) = S'_b \, M (\Gamma)
\end{matrix}
$$
This ends our overview of the general plan of the proof of the GSC theorem.

\smallskip

Here are, finally, some general kind of remarks concerning the technologies used in this paper. During our constructions, some boundary points have to be deleted, or rather sent to infinity as ``punctures''. This is a general kind of phenomenon occurring in the construction of $S_u$. But then, transversal compactness puts a limit on how much punctures we are allowed to use, without loosing our desired (\ref{eq0.10}). I call this the ``Stallings' barrier''. Remember that, according to classical results of John Stallings, it is very easy to get GSC by multiplying with $R^p$, $p$ large, which of course makes havoc of transversal compactness. Multiplying with $B^p$, like in (\ref{eq0.12}), i.e. very different affair. We have already mentioned these things in connection with (\ref{eq0.12}).
But then there is also a second ``nonmetrizability barrier'', to be respected, and it pushes in the opposite direction with respect to the Stallings' barrier. We want to stay on the good side of both of them. Here is how this second barrier occurs. During our construction of $S_b$, we need to drill some ditches in those additional dimensions, and then, later on, fill them in again with various material, but differently. Unless this is done with very great care, it leads to spaces which are not metrizable, hence useless for our purposes.

\bigskip

Many thanks are due to Louis Funar, David Gabai and Daniele Otera, for useful conversations.

\smallskip

Without the friendly help of the IHES, this paper could not have seen the light of day. Last, but not least, many thanks are due to C\'ecile Gourgues for the typing and Marie-Claude Vergne for the drawings.

\newpage

\section{The 2-dimensional representation theorem}\label{sec1}
\setcounter{equation}{0}

The REPRESENTATIONS of some completely arbitrary finitely presented group $\Gamma$ have been defined in \cite{29}, a paper of which the present one is a direct conti\-nuation. We will freely refer to the whole content of \cite{29}, notations included.

\smallskip

Here is the main result of the present section.

\bigskip

\noindent {\bf Theorem 1.1.} (2-DIMENSIONAL REPRESENTATION THEOREM) 

\smallskip

\noindent {\it For any finitely presented group $\Gamma$ there is a $2$-dimensional  REPRESENTATION
\begin{equation}
\label{eq1.1}
X^2 \overset{f}{\longrightarrow}  \widetilde M (\Gamma) \, ,
\end{equation}
with the following features}

\medskip

1) {\it (First finiteness condition.) The $2$-dimensional cell-complex $X^2$ is locally finite.}

\medskip

2) {\it (Equivariance.) There is a free action $\Gamma \times X^2 \to X^2$ s.t. for all $x \in X^2$, $\gamma \in \Gamma$ we have $f(\gamma x) = \gamma f(x)$.}

\medskip

3) {\it (The second finiteness condition.) For any tight compact transversal $\Lambda$ to $M_2 (f) \subset X^2$ we have
\begin{equation}
\label{eq1.2}
{\rm card} \, (\lim (\Lambda \cap M_2 (f))) < \infty \, .
\end{equation}
}

\medskip

4) {\it The closed subset, with $\underset{\Lambda}{\bigcup}$ running over the tight $\Lambda$'s above
\begin{equation}
\label{eq1.3}
{\rm LIM} \, M_2 (f) \underset{\rm def}{=} \, \bigcup_{\Lambda} \, \lim (\Lambda \cap M_2(f)) \subset X^2 \, ,
\end{equation}
which is the only place where $M_2(f) \subset X^2$ can accumulate, is a locally finite graph and $f \, {\rm LIM} \, M_2(f) \subset fX^2$ is also a {\ibf closed} subset. We also have the following feature. Let $\Lambda^*$ run over all tight transversals to ${\rm LIM} \, M_2 (f)$, then we have}
\begin{equation}
\label{eq1.4}
\bigcup_{\Lambda^*} \ (\Lambda^* \cap M_2(f)) = M_2 (f) \, .
\end{equation}

\medskip

5) {\it Let $(x,y) \in M^2 (f) \subset X^2 \times X^2$ be a double point of $f$ and
$$
\lambda (x,y) \subset \widehat M^2 (f) \equiv M^2 (f) \cup {\rm Diag} ({\rm Sing} (f)) \subset X^2 \times X^2
$$
be a {\ibf zipping path} for $(x,y)$. There exists a {\ibf uniform bound} $K > 0$ s.t. if $\Vert \lambda (x,y) \Vert \equiv \{$The length of $\lambda (x,y) \}$ (which is well-defined, up to quasi-isometry), then
\begin{equation}
\label{eq1.5}
\inf_{\lambda} \Vert \lambda (x,y) \Vert < K \, ,
\end{equation}
where $\lambda$ runs over all zipping paths for $(x,y)$.}

\bigskip

In this statement, exactly like in \cite{29}, $\widetilde M (\Gamma)$ is the universal covering space of a compact {\ibf singular} $3$-manifold $M(\Gamma)$ which is a {\ibf presentation} of the group $\Gamma$, i.e. $\pi_1 \, M(\Gamma) = \Gamma$. Actually our $M(\Gamma)$ is a gotten by starting with a smooth $3^{\rm d}$ handlebody of some high genus $g$, and then adding to it $2$-handles. Each of these 2-handles is attached along a smooth embedding $S^1_i \times I \to \partial H$, but the global map $\underset{i}{\sum} \, S_i^1 \times I \overset{s}{\longrightarrow} \partial H$ is no longer injective, it is only an immersion coming with a double point set ${\rm Sing} \, M(\Gamma) \subset \partial H \subset M(\Gamma)$, the points of which are called the {\ibf immortal singularities} of $M(\Gamma)$; these are the points where $M(\Gamma)$ fails to be a manifold. Here is a precise description of the local structure of $M(\Gamma)$ around a connected component $\bar S \subset {\rm Sing} \, M(\Gamma)$. [The immersion $s$ above is assumed generic and $\bar S$ is a like square.] We are given a copy of $R_+^3$, with $\partial R_+^3 = R^2$ and with two copies of $R_+^2 \times [0,1]$ parametrized as $R_{\varepsilon} \times R_+ \times [0,1]_{\varepsilon}$, with $\varepsilon = 1,2$. The $R_{\varepsilon} \times [0,1]_{\varepsilon} = R_{\varepsilon} \times [0,1]_{\varepsilon} \times \partial R_+ \subset \partial (R_{\varepsilon} \times R_+ \times [0,1]_{\varepsilon})$ come with two embeddings
$$
R_1 \times [0,1]_1 \longrightarrow (R^2 = \partial R_+^3) \longleftarrow R_2 \times [0,1]_2 \, ,
$$
cutting through each other transversally along the square $\bar S = [0,1]_1 \times [0,1]_2 \subset R^2$. The singular 3-space $N^3 (\bar S)$, gotten by glueing each of the two $R_{\varepsilon} \times R_+ \times [0,1]_{\varepsilon}$'s along the corresponding $R_{\varepsilon} \times [0,1]_{\varepsilon}$ to $R_+^3$ is an open neighbourhood of our $\bar S \subset M(\Gamma)$. One can notice that the $N^3 (\bar S)$ is exactly like the {\ibf undrawable singularities} in \cite{8}, \cite{18}, \cite{19} in their 3-dimensional version, without any attached map of which they would be the source. One can require that our $\overline S$'s should always be contained in the lateral surface of the $0$-handles of $M (\Gamma)$ and that a given $2$-handles should see at most one $\overline S$.

\bigskip

\noindent {\bf Proof of Theorem 1.1.} The proof of theorem~1.1 starts from a $3$-DIMENSIO\-NAL REPRESENTATION THEOREM (see theorem~1.2 below) proved in \cite{29}, proceeding afterwards on the general lines of \cite{31}, a paper which, although focusing on smooth 3-manifolds was like taylor-made for our present needs, i.e. for \cite{29} and for the present paper. But there will be now important technical details where we will proceed differently from  \cite{31} and, normally, we will signalize this at the appropriate time.

\smallskip

Exactly like in \cite{29} and \cite{31}, the $M(\Gamma)$ and its $\widetilde M(\Gamma)$ are unions of handles of index $\lambda \leq 2$, denoted generically by $h_i^{\lambda}$, while the $3^{\rm d}$ object $Y(\infty)$ to the considered next and introduced already in \cite{29}, consists of {\ibf bicollared handles of index $\lambda \leq 2$.} These handles are denoted, generically, by $H_i^{\lambda} (\gamma)$, and each of these $H_i^{\lambda} (\gamma)$'s corresponds to some usual $\lambda$-handle $h_i^{\lambda} \subset \{ M(\Gamma) \ \mbox{or} \ \widetilde M (\Gamma)\}$, the case $\widetilde M(\Gamma)$ being of real interest for us. As explained in \cite{29}, the index $\gamma$ belongs to a countable set which, in principle, is $(\lambda, i)$-dependent. Each bicollared $H_i^{\lambda}$ comes with its decomposition into usual handles, $H_i^{\lambda} = \overset{\infty}{\underset{n=1}{\bigcup}} \, H_{i,n}^{\lambda}$. More details concerning the bicollared handles are to be found in \cite{29} and \cite{31}.

\smallskip

Theorem~1.2 which is stated below is, essentially, the juxtaposition of the main result of \cite{29} and of the lemma~4.3 in the same paper. Like always in this paper, $\Gamma$ is a generic, arbitrary, finitely presented group.

\bigskip

\noindent {\bf Theorem 1.2.} (3-DIMENSIONAL REPRESENTATION THEOREM) 

\smallskip

\noindent {\it For any $\Gamma$, there is a $3$-dimensional REPRESENTATION
\begin{equation}
\label{eq1.6}
Y(\infty) \overset{g(\infty)}{-\!\!\!-\!\!\!-\!\!\!\longrightarrow} \widetilde M (\Gamma) \, ,
\end{equation}
with the following features.}

\medskip

1) {\it The space $Y(\infty)$ is a $3^{\rm d}$ {\ibf locally finite} cell complex which is a union of bicollared handles $Y(\infty) = \underset{\footnotesize\overbrace{\lambda,i,\gamma}}{\bigcup} \, H_i^{\lambda} (\gamma)$, where $\lambda \leq 2$, and where the system of indices $i$ corresponds to a given handle-decomposition $\widetilde M (\Gamma) = \underset{i,\lambda}{\bigcup} \, h_i^{\lambda}$.}

\medskip

2) {\it There is a free action $\Gamma \times Y(\infty) \to Y(\infty)$, respecting the bicollared handlebody structure, with respect to which the non-degenerate map $g(\infty)$ is {\ibf equivariant}.}

\medskip

3) {\it For each $H_i^{\lambda} (\gamma)$ we have the non-compact attaching zone $\partial H_i^{\lambda} (\gamma) \subset H_i^{\lambda} (\gamma)$ and the compact lateral surface $\delta H_i^{\lambda} (\gamma)$ lives at the infinity of $H_i^{\lambda} (\gamma)$. The embedding $g(\infty) \mid H_i^{\lambda} (\gamma)$ extends continuously to an embedding of $\widehat H_i^{\lambda} (\gamma) \equiv H_i^{\lambda} (\gamma) \, \cup \, \delta H_i^{\lambda} (\gamma)$, coming with a strict equality of sets $g(\infty)(\delta H_i^{\lambda} (\gamma)) = \delta h_i^{\lambda}$. Moreover, inside $\widetilde M(\Gamma)$, the $g(\infty) \, \widehat H_i^{\lambda} (\gamma)$ occupies, roughly, the position of $h_i^{\lambda} \subset \widetilde M(\Gamma)$.}

\medskip

4) {\it The various $\varepsilon$-skeleta $Y^{(\varepsilon)} \subset Y(\infty)$ contain canonical outgoing collars, such that each $H_i^{\lambda} (\gamma)$ is attached along $\partial H_i^{\lambda} (\gamma)$ to $Y^{(\lambda - 1)}$ in a collar-respecting manner at some {\ibf level} $k(i,\gamma) \in Z_+$ s.t.
\begin{equation}
\label{eq1.7}
\mbox{If $i$ is fixed and} \ \gamma_n \to \infty \ \mbox{then} \ \lim k(i,\gamma_n) = \infty \, .
\end{equation}
}

Before going on, let us notice that, topologically speaking, $H_i^{\lambda} (\gamma)$ is a non-compact 3-manifold with $\partial H_i^{\lambda} (\gamma)$ being, at the same time, the boundary and the $\lambda$-handle attaching zone.

\bigskip

\noindent {\bf Complement 1.3 to the $3^{\rm d}$ representation theorem.}

\smallskip

1) {\it For each $\varepsilon$-skeleton $Y^{(\varepsilon)}$ of $Y(\infty)$ we introduce the ideal boundary, living at infinity
\begin{equation}
\label{eq1.8}
\delta \, Y^{(\varepsilon)} \equiv \bigcup_{i,\gamma , \lambda \leq \varepsilon} \delta H_i^{\lambda} (\gamma) \quad \mbox{and} \quad \widehat Y^{(\varepsilon)} \equiv Y^{(\varepsilon)} \cup \delta \, Y^{(\varepsilon)} \, .
\end{equation}
As a consequence of} (\ref{eq1.7}), {\it inside $\widehat Y^{(\lambda - 1)}$ we find that}

\medskip

\noindent (1.9) \quad {\it When $\gamma_n \to \infty$, then $\{ \partial H_i^{\lambda} (\gamma_n)\}$ accumulates on $\delta \, Y^{(\lambda - 1)}$, and it is this which implies the local finiteness of $Y(\infty)$.}

\medskip

2) {\it We also introduce, now at the target the set}

\setcounter{equation}{9}
\begin{equation}
\label{eq1.10}
\Sigma_1 (\infty) \equiv \bigcup_{i,\gamma , \lambda} G (\delta H_i^{\lambda} (\gamma)) = \bigcup_{i,\lambda} \delta h_i^{\lambda} \subset \widetilde M (\Gamma)
\end{equation}
{\it and, with this, inside $\widetilde M (\Gamma)$, the $G\delta H_{i,m}^{\lambda} (\gamma_n)$ accumulate on
$$
G \, \delta \, Y^{(\lambda)} \subset \Sigma_1 (\infty) \, .
$$
[While the {\rm (1.9)} above is behind the first finiteness condition in the theorem~$1.1$, the present item is one of the essential ingredients behind the second finiteness condition from the same theorem. For the same purpose of getting the second finiteness condition in the theorem~$1.1$, an additional $2^{\rm d}$ condition will have to be imposed. For given $(\lambda , \gamma)$, consider the innermost compact wall (see {\rm (1.12)} below)
$$
W_i (\lambda) \subset X^2 \mid H_i^{\lambda} (\gamma) \, .
$$
Then, inside $X^2$, we have $\underset{i = \infty}{\lim} \ W_i (\lambda) = \infty$.]}

\medskip

3) (Complement to (\ref{eq1.7}) and (1.9).) {\it There are {\rm PROPER} individual embeddings $\partial H_i^{\lambda} (\gamma) \subset Y^{(\lambda - 1)}$ and the following global map is} PROPER
$$
\sum_{i,\gamma,\lambda} \partial H_i^{\lambda} (\gamma) \overset{j}{\longrightarrow} Y(\infty) \, .
$$
{\it We have ${\rm Im} \, j = {\rm Sing} \, (g(\infty)) \equiv \{$the points $x \in Y(\infty)$ where $g(\infty)$ fails to be immersive, i.e. the {\ibf mortal singularities} of $g(\infty)\}$.}

\medskip

\noindent {\it [COMMENT. While it is the vocation of the equivalence relations $\Psi (f)$ which come with our REPRESENTATIONS, to kill the {\ibf mortal} singularities, and these are always singularities of {\ibf maps}, nothing, in particular none of our maps, will ever kill the {\ibf immortal} singularities, which are singularities of {\ibf spaces}.]}

\medskip

4) {\it Our $Y(\infty)$ is a $3^{\rm d}$ train-track, which fails to be smooth exactly along the set ${\rm Im} \, j$ above. Finally, we have}
$$
g(\infty) ({\rm Sing} \, (g(\infty))) \cap {\rm Sing} \, \widetilde M (\Gamma) = \emptyset \, .
$$

\bigskip

Following, essentially, \cite{31} but with some specific modification to be explicitly developed below, we will show now how to deduce our theorem~1.1 from the results (\ref{eq1.2}) $+$ (\ref{eq1.3}) above, which have already been proved in \cite{29}.

\smallskip

Each individual $H_i^{\lambda} (\gamma)$ will be endowed, like in \cite{31} with three partial foliations ${\mathcal F}$ (COLOUR), our three colours being blue $(\lambda = 0)$, red $(\lambda = 1)$, black $(\lambda = 2)$. The foliations are invariant when the action of $\Gamma$ permutes the $H^{\lambda}$'s and the handle attachments will respect them. Unlike what we did in \cite{31} we do not ask now that the ${\mathcal F}$ (COLOUR) should be globally defined throughout $\widetilde M (\Gamma)$, we only insist that, for each individual colour, $g(\infty)({\mathcal F} \mid H_1)$ and $g(\infty)({\mathcal F} \mid H_2)$ should agree on $g(\infty) \, H_1 \cap g(\infty) \, H_2$. BUT, when $h_k^0 , h_i^2 , h_j^2$ participate in some immortal $S \subset {\rm Sing} \, \widetilde M (\Gamma)$, $S \subset \delta h_k^0$ then, inside the $\gamma$-independent $g(\infty) \, H_k^0 (\gamma)$, the $g(\infty)({\mathcal F} ({\rm BLACK}) \mid H_i^2 (\gamma_1))$ and $g(\infty)({\mathcal F} ({\rm BLACK}) \mid H_j^2 (\gamma_2))$ have to cut through each other transversally, for each choice of $\gamma_1 , \gamma_2$. This {\ibf is} a novelty with respect to \cite{31}.

\smallskip

The representation space of the $2^{\rm d}$ representation, the $X^2$ is, in a first approximation, a {\ibf very dense} $2$-skeleton of $Y(\infty)$, the $3^{\rm d}$ representation space, put together out of spare parts, which are always pieces of some leaf of one of the ${\mathcal F}$ (COLOUR)'s. These are either compact walls $W$ or ``non-compact'' {\ibf security walls} $W_{\infty}({\rm BLACK})$. We have put the non-compact between quotation marks since it may mean simply that $\partial W_{\infty} ({\rm BLACK})$ contains free boundary pieces (which we may decide to delete). The notation $W_{(\infty)}({\rm BLACK})$ will mean ``$W({\rm BLACK})$ or $W_{\infty}({\rm BLACK})$''. The $X^2 \mid H_i^{\lambda}$, the detailed structure of which is to be explained later, is to fulfill the following requirements:

\medskip

\noindent (1.11.1) \quad Each individual $X^2 \mid H_i^{\lambda}$ is GSC.

\medskip

\noindent (1.11.2) \quad Let the $H_{j_1}^{\lambda+1} , H_{j_2}^{\lambda+1} , \ldots$ be adjacent to $H_i^{\lambda}$ inside $Y(\infty)$, where $H_j^{\lambda+1}$ and $H_i^{\lambda}$ are glued together exactly along $\partial H_j^{\lambda+1} \cap H_i^{\lambda}$. We denote $H_i^{\lambda} \cup \underset{j}{\sum} \, H_j^{\lambda + 1} \equiv H_i^{\lambda} \cup H_{j_1}^{\lambda + 1} \cup H_{j_2}^{\lambda + 1} \cup \ldots$ and, with this, we will ask that
$$
\Psi \left( f \mid \left( X^2 \mid H_i^{\lambda} \cup \sum_j H_j^{\lambda + 1} \right)\right) = \Phi \left( f \mid \left( X^2 \mid H_i^{\lambda} \cup \sum_j H_j^{\lambda + 1} \right)\right) .
$$
The two features (1.11.1) $+$ (1.11.2) will make that, for the $2^{\rm d}$ REPRESENTATION from (\ref{eq1.1}), the conditions
$$
X^2 \in {\rm GSC} \quad \mbox{and} \quad \Psi (\psi) = \Phi (\psi)
$$
will be automatically fulfilled, once the analogous conditions are verified already of (\ref{eq1.6}). With this we give now the

\medskip

\noindent (1.12) \quad ({\ibf The structure of $X^2 \mid H_i^{\lambda}$, when $\lambda \leq 1$}). With some important differences, to be signalized in due time and also with some additions to be explained only later on, this will be very much like in the \cite{31} and the figures 2.8 to 2.20 from that paper should still be useful now. With this, modulo the additions mentioned above, for $\lambda \leq 1$ one will have $X^2 \mid H_i^{\lambda} = \{$the attaching zone $\partial H_i^{\lambda} \} \cup \{$infinitely many compact walls $W$ (natural COLOUR $\lambda$), glued in the case $\lambda = 1$ via their $\partial W ({\rm RED})$ to $\partial H_i^{\lambda}$ ($\lambda = 1$) and converging, in both cases $\lambda=0$ and $\lambda=1$ to the ideal $\delta H_i^{\lambda}$ which lives at infinity$\} \cup \{$whenever, at the level of $Y(\infty)$ we have $\partial H_j^{\mu > \lambda} \cap H_i^{\lambda} \ne \emptyset$, then the corresponding piece of the attaching zone $\partial H_j^{\mu}$ is already contained in $X^2 \mid H_i^{\lambda} \} \cup \{$infinitely many security walls $W_{\infty} ({\rm BLACK})_{H_i^{\lambda}} \}$. Neither the $W(\lambda)$'s nor the $W_{\infty} ({\rm BLACK})$'s accumulate at finite distance.

\smallskip

Figures 1.2, 1.3 display the $W_{\infty} ({\rm BLACK})$'s and after the EXPLANATIONS for these figures, it will also be shown how each of the individual $W_{\infty} ({\rm BLACK})$'s is attached to the rest of $X^2 \mid H_i^{\lambda}$. BUT, careful, our present $W_{\infty} ({\rm BLACK})$'s are larger than the ones from \cite{31}. Then another novelty with respect to \cite{31} is that, while in \cite{31} we had $\Gamma = \pi_1 M^3 (\mbox{\it smooth})$, our present $M(\Gamma)$, which replaces the $M^3$, has the immortal singularities $\bar S$ which we find now at the infinity of the $H_i^0$'s and which comes with a host of complications to be discussed in due time. This ENDS the item (1.12). \hfill $\Box$

\bigskip

\noindent (1.13) \quad ({\ibf The structure of $X^2 \mid H_i^2$.}) We will find now that $X^2 \mid H_i^2 = \{$the attaching zone $\partial H_i^2\} \cup \{$infinitely many $W({\rm BLACK})$'s glued to the preceeding piece along their $\partial W({\rm BLACK})$'s, and with no other glueings at the source $X^2$ of $f\}$.

\smallskip

For every given $H_i^2$ there is a {\ibf unique} $W({\rm BLACK \, complete})$ which is a 2-cell, actually a $2p$-gone for some $p > 1$, let us say a hexagon, like in the figure~1.1. All the other $W({\rm BLACK})$ of $H_i^2$ are annuli, wide enough to see all the double lines signalized in figure~1.1. Let us say that they are punched by a little open cell, a BLACK hole $H$, with $\partial H$ a piece of free boundary for $X^2$, not glued to anything at the source.

\smallskip

This very simple stucture for $X^2 \mid H_i^2$ corresponds to an earlier suggestion made by Dave Gabai, in a different but related context. \hfill $\Box$

\bigskip

When, inside $\widetilde M (\Gamma)$ we look at the totality of the $f$ (compact walls $W$) these have {\ibf limit positions}, the limit walls and, with the $\Sigma_1 (\infty)$ introduced in (\ref{eq1.10}), we find that there is the exact equality
\setcounter{equation}{13}
\begin{eqnarray}
\label{eq1.14}
\Sigma_1 (\infty) &= &\{\mbox{the union of the limit positions of the compact walls}\} \nonumber \\
&= &\sum \, S_{\infty}^2 (\mbox{BLUE}) \cup \sum \, (S^1 \times I)_{\infty} (\mbox{RED}) \cup \sum \, {\rm Hex}_{\infty} (\mbox{BLACK}) \subset \widetilde M (\Gamma) \, . \nonumber \\
\end{eqnarray}

\medskip

\noindent VERY IMPORTANT REMARKS.

\smallskip

A) If it would not be for the immortal singularities $\bar S \subset {\rm Sing} \, \widetilde M (\Gamma)$, the $\Sigma_1 (\infty)$ would be a simple-minded $\{$smooth surface ($=$ $2$-manifold) with {\ibf ramification} points$\}$ locally embeddable in $R^3$. But with $\bar S$'s present, we also get {\ibf branching} points. For the ramification points, the local models for $\Sigma_1 (\infty)$ are the following generic configurations
$$
(x=0) \cup [(y=0) \cap (x \geq 0)]
$$
and
$$
(x=0) \cup [(y=0) \cup (x \geq 0)] \cup [(z=0) \cap (x \geq 0) \cap (y \geq 0)] \, ,
$$
which are clearly embeddable in $R^3$.

\smallskip

For the branching points, the local model is
$$
\{ (x=0) \cup [(y=0) \cap ( x \geq 0)]\} + \{ (x=0) \cap [(z=0) \cap (x \geq 0)]\} \, ,
$$
with the two pieces being glued only along $(x=0)$; the $p_{\infty\infty} (\infty)$'s in figure 1.5.(c) are such points. They are very much like an undrawable singularity for $\Sigma_1 (\infty)$.

\medskip

B) The BLUE and RED limit walls ($S_{\infty}^2$ and $(S^1 \times I)_{\infty}$) are generated already at the level of the individual bicollared handles $H_i^0 (\gamma)$, $H_j^1 (\gamma)$ and, for $\lambda = 0$ or $\lambda = 1$, the various $H_i^{\lambda} (\gamma_1) , H_i^{\lambda} (\gamma_2) , \ldots$, with ``$i$'' meaning ``$i$ or $j$'', correspond to the same individual $S_{\infty}^2$ or $(S^1 \times I)_{\infty}$. By contrast with this, it takes the whole infinite collection $H_k^2 (\gamma_1) , H_k^2 (\gamma_2) , \ldots$ to generate an individual pair of walls ${\rm Hex}_{\infty} ({\rm BLACK})$.

\medskip

C) The BLACK holes mentioned above are God-given, i.e. they are part of the structure of $X^2$, from the beginning and they are not to be mixed up with the later, artificial i.e. man-made Holes $H$, introduced in sections III, IV. These latter ones are an indispensable tool for making sense of the geometric realization of the zipping process, which is the core of section IV below. \hfill $\Box$

\bigskip

\noindent EXPLANATIONS CONCERNING THE FIGURE 1.1. The drawing (A) gives a complete view of the unique $W \, ({\rm BLACK \, complete}) \subset H_i^2 (\gamma) \subset Y(\infty)$. The doubly collared structure of $H_i^2 (\gamma)$ imposes the telescopic system of hexagons, which is actually infinite. In real life, this system may be more irregular than suggested here (see, for instance the figure 1.4.III.C in \cite{31}); but that comes without any harm.

\smallskip

In the present context, figure (A) replaces the figure 2.20 from the paper \cite{31}, with which it should be compared. But it also differs in several important respects from that 2.20. One will notice that now, the central lines from 2.20 have been deleted. Most importantly still, while in \cite{31} we had $\Gamma = \pi_1 M^3$ with a smooth $M^3$, now $\widetilde M (\gamma)$ is singular and $\bar S \subset {\rm Sing} \, \widetilde M (\Gamma)$ is one of the immortal singularities.

\smallskip

In our drawing (A), only the double points coming from the $H^0$'s and $H^1$'s adjacent to $H_i^2 (\gamma) \supset W \, ({\rm BLACK})$ are explicitly displayed; a figure like (A) will be said to be {\ibf at the source}. Actually, the $W \, ({\rm BLACK})$ is glued to the $H^0$'s along $M,K$ and to $H^1$'s along $L,N$.

\vglue 15mm

$$
\includegraphics[width=165mm]{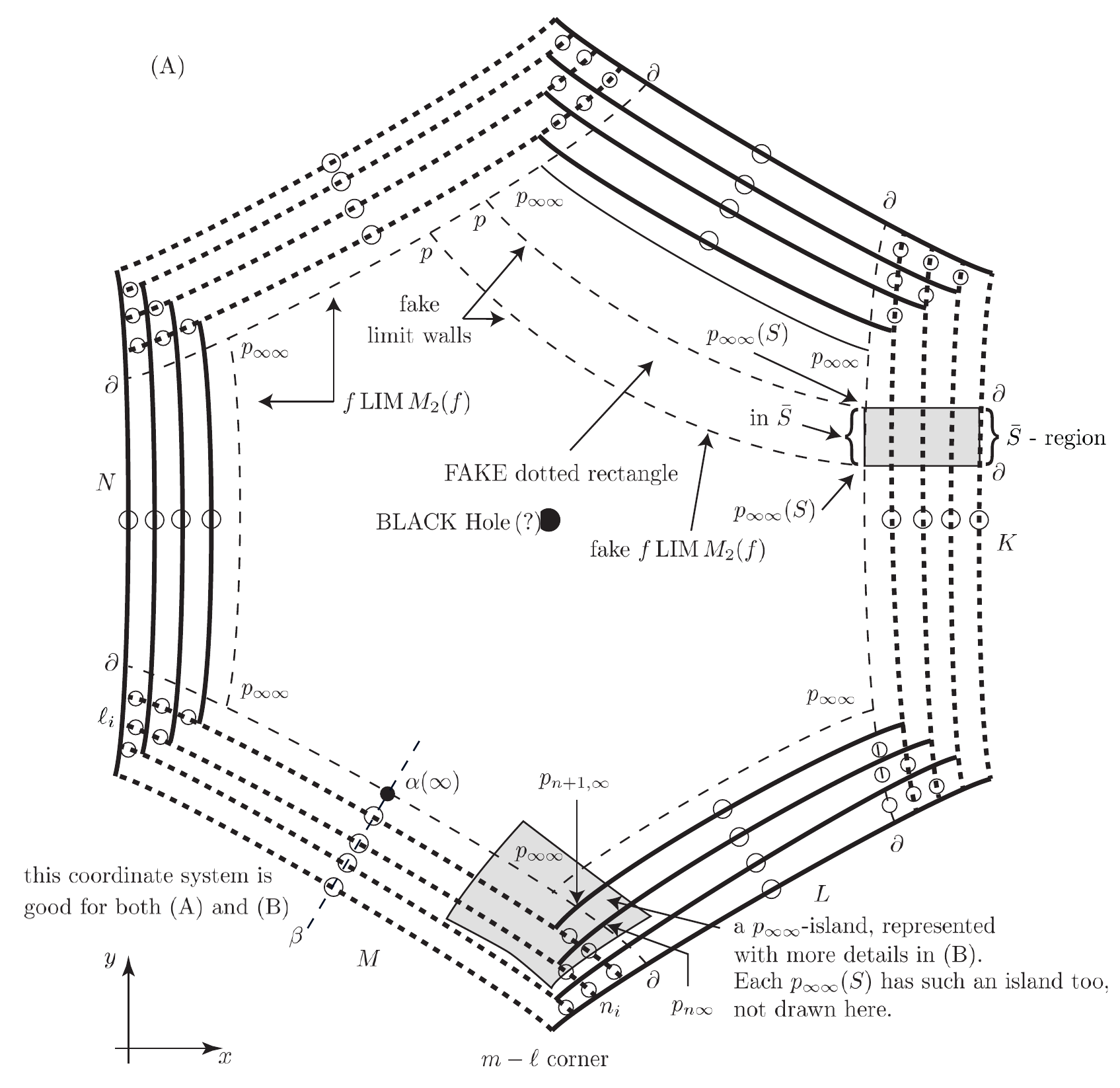}
$$
\vglue -15mm
$$
\includegraphics[width=13cm]{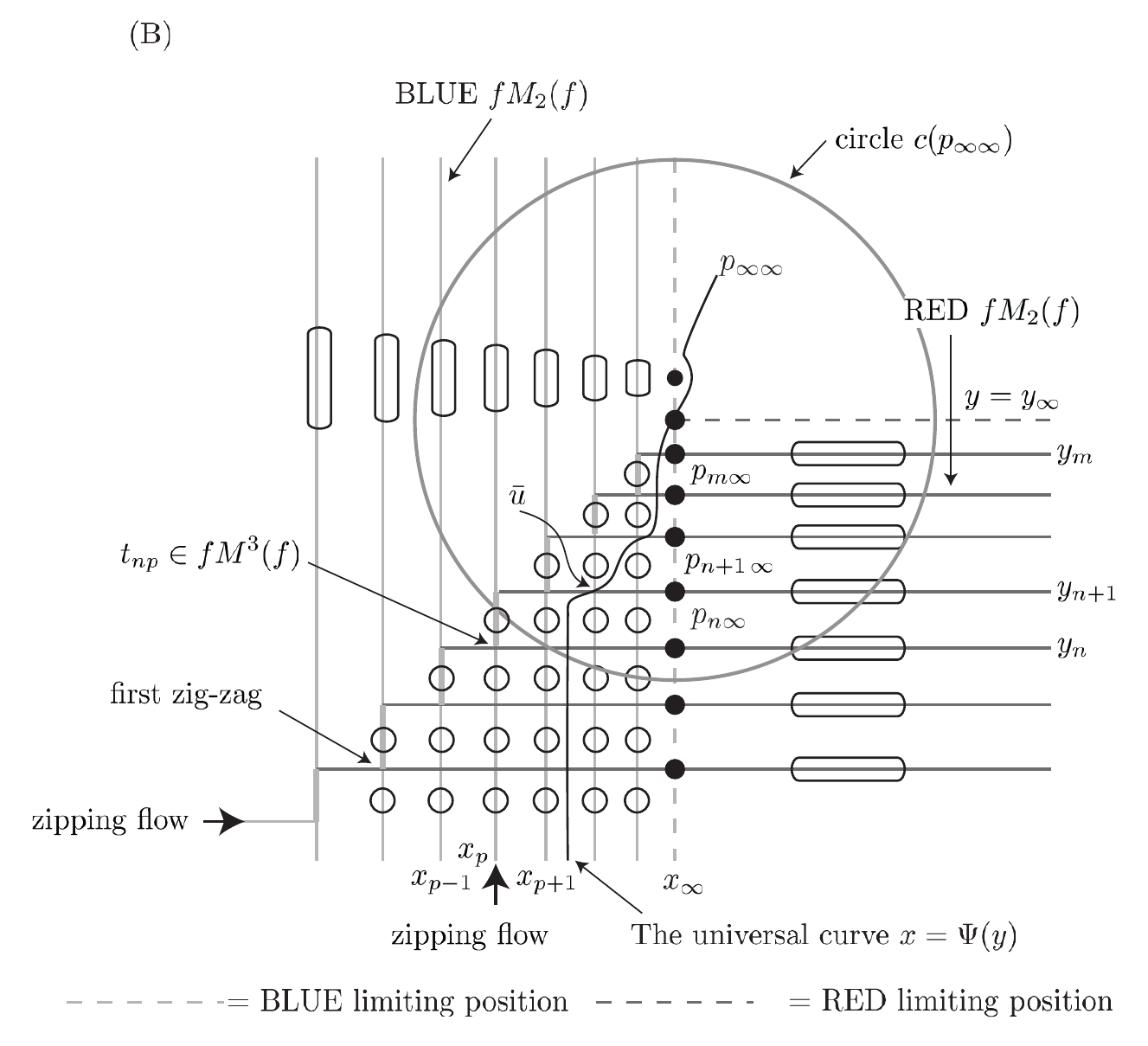}
$$
$$
\includegraphics[width=12cm]{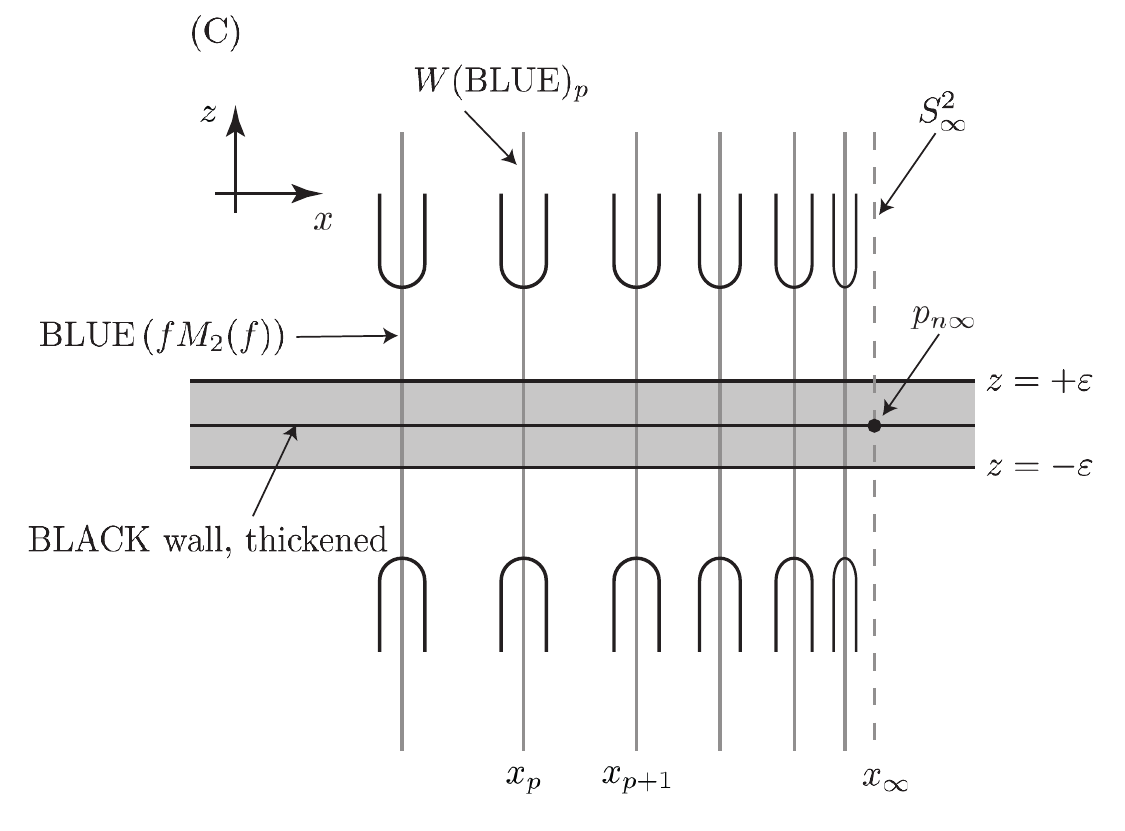}
$$
\label{fig1.1}

\centerline {\bf Figure 1.1.} 

\smallskip

\begin{quote}
In (A) one sees a $W({\rm BLACK})$. In (B) one sees with more details the $p_{\infty\infty}$-island.

\medskip

\noindent LEGEND: 

\smallskip

\noindent $-----$ (thin dotted line) $=$ trace of (limit wall) $\cap$ $W$(COLOUR)

\smallskip

\noindent {\Huge{$_-\!\!\!_-\!\!\!_-\!\!_-\!\!_-\!\!_-$}} (thick black line) $=$ (RED $\cap$ BLACK) in $fM_2 (f)$

\smallskip

\noindent {\Huge{$_{- - - }$}} (thick dotted line) $=$ (BLUE $\cap$ BLACK), in $fM_2 (f)$.

\smallskip

\noindent $\bigcirc$ $=$  Such a circle, riding on top of a double line, in the drawings (A), (B), signalizes the presence of two holes, at level $X^2 - H \subset X^2$, inside the corresponding $W \, ({\rm COLOUR}) \pitchfork W \, ({\rm BLACK})$, one over $W \, ({\rm BLACK})$, the other one under; see here drawing (C) too. For a discrete set of values $y=y_n$, figure (C) lives in a RED wall.
\end{quote}

\bigskip

There is also another kind of figure, {\ibf at the target}, where all the lines $M_2 (f) \cap W$ are drawn. This is a much denser version of our present figure which is at the source, but at the level of our coarse typography, it is graphically undistinguishable from it. The $W \, ({\rm BLACK \, complete}) \subset X^2$ are without the central BLACK hole, while each $W \, ({\rm BLACK \, reduced}) \subset X^2$ has one, permanently deleted.

\smallskip

When later on in this paper (see section IV below) we will move from $X^2$ to $X^2 - H$ ($=$ the $X^2$ with Holes), then from {\ibf some} of the $W \, ({\rm BLACK \, complete})$ the central BLACK Hole (occurring here with a question mark) is to be deleted too (when at level $X^2 - H$). In this last situation, the BLACK Hole $\subset W \, ({\rm BLACK}$ ${\rm complete})$ will be among our normal Holes $H$ from (4.6). Much more concerning $X^2 - H$ will be said later on.

\smallskip

The drawing (B) gives more details concerning the $p_{\infty\infty}$-island, and the not explicitly drawn $p_{\infty\infty} (S)$-island is completely similar; it is concentrated around the corresponding $p_{\infty\infty} (S)$ corner of the doubly shaded $\bar S$-region. The (C) which lives in the plane $y=y_n$ of $W({\rm RED})_n$ accompanies the (B), which itself lives in the $W({\rm BLACK})$ from (A).

\smallskip

Here is how the $\bar S$-region arises. At the level $\widetilde M (\Gamma)$, $\bar S \subset \delta h_k^0$ is an immortal singularity, generated by $h_i^2 , h_j^2$. At level $Y(\infty)$ to $h_i^2 , h_j^2$ correspond our $H_i^2 (\gamma)$ and also a dual $H_j^2 (\gamma')$, both glued to the $H_k^0$ (side $K$, figure~(A)), coming with $g(\infty) \, H_i^2(\gamma)$ and $g(\infty) \, H_j^2 (\gamma')$ which cut transversally through each other. The shaded $S$-region in (A) stands for $\{\mbox{our} \, W({\rm BLACK}) \cap H_j^2 (\gamma') \subset \widetilde M (\Gamma)\}$. When we go to (\ref{eq1.1}), then the $H_j^2 (\gamma'_1) , H_j^2 (\gamma'_2) , \ldots$ (all possible $\gamma'_1 , \gamma'_2$), generate infinitely many $W \, ({\rm BLACK \, complete})^*$'s. The dotted line $[p_{\infty\infty} (S),P]$ in (A) corresponds to what the trace of the $\{$limiting position of our $W^*$'s$\} \cap W({\rm BLACK})$ {\ibf would be} IF $W$ and $W^*$ would cut through each other physically, for $x > 0$. As it is $[p_{\infty\infty} (S),P]$ is not physical, so we call it FAKE. (Think here of the train-track structure coming with $\bar S$.) The fake $\Sigma_1 (\infty)$ will play no role in our story.

\smallskip

In the drawings (A) $+$ (B), which we consider now at the source, the outer zigzag line, going from the corner marked $m-\ell$ (see (A)) to the corresponding point $p_{\infty\infty}$, comes from the following intersection, consisting essentially entirely of double points of the $f$ from (\ref{eq1.1}):
$$
W \, ({\rm BLACK}) \cap \{\mbox{a pair} \, (X^2 \mid H^0) \underset{{\overbrace{\partial H^1}}}{\cup} (X^2 \mid H^1), \ \mbox{ glued together at the level of $X^2$}\} \, . \leqno (*)
$$

At the source $X^2$, each $W({\rm BLACK})$ is glued to the rest exactly along its $\partial W({\rm BLACK})$; as we shall see, for $W_{\infty}({\rm BLACK})$, things are more complicated, forced by the needs of (1.11.1).

\smallskip

In a figure at the source, a $p_{\infty\infty}$-island has one zigzag line, generating everything, via the BLUE and RED half-lines which come out of it, but it is only the RED ones which are glued to the zigzag at the level of $X^2$. At the target, there are infinitely many such arborescent drawings all superposed, creating an infinite checkerboard. Coming back to our zigzag at the source, in terms of (1.12) it is of the form $W({\rm BLACK}) \, \cap \, \partial H^1$ and this specific $H^1$, occurring in figure~1.1, is called $H_L^1 (\gamma^0)$. We have $g(\infty) \, \widehat H_L^1 (\gamma^0) \approx h_L^1 \subset \widetilde M (\Gamma)$.

\smallskip

As long as we only consider the figure at the source, all the red lines shoot out of the zigzag above, and all the BLUE lines, cutting transversally through them, stretch from the RED border of $\partial W({\rm BLACK})$ (let us say the $L$ in (A)), to our zigzag and beyond. When we consider now the figure~1.1 {\ibf at the target}, then there are now infinitely many zigzags, produced by all the $\partial H_L^1 (\gamma)$'s, one of these indices $\gamma$ being {\it our} $\gamma_0$. Each $\partial H_L^1 (\gamma)$ comes, like in 4) from theorem~1.2, with its attaching level $k(L,\gamma)$. Finitely many $\gamma$'s come with $k(L,\gamma) < k(L,\gamma_0)$ and their zigzags enter the figure at the target, through the BLUE side $M$. For the other, infinitely many of them, coming with $k(L,\gamma) > k(L,\gamma_0)$, their zigzags enter the figure through the RED side L, to the right of $(*)$, closer and closer to the BLUE limit wall, while the preceeding ones were, in terms of the figure~1.1.(A), to the left of $(*)$.

\smallskip

A given $W({\rm BLACK})$ has at most one $\bar S$-region, coming with its two $p_{\infty\infty} (S)$'s, and each $p_{\infty\infty} (S)$-island is treated just like the ordinary $p_{\infty\infty}$-islands. Each $W({\rm BLACK})$, complete or reduced, has a uniquely chosen side $M$, see here figure~1.1, with its line $(\alpha (\infty),\beta) (W({\rm BLACK}))$ (where $\alpha (\infty) \in f \, {\rm LIM} \, M_2 (f)$), cutting through all the infinitely many BLUE double lines, like $[\ell_i , n_i]$. This $(\alpha (\infty) , \beta)$ is not part of $M_2 (f) \cup {\rm LIM} \, M_2 (f)$.

\smallskip

We can easily arrange the zipping flow for (\ref{eq1.1}) so that
$$
\mbox{each spot} \ \{ [\ell_i , n_i] \ \mbox{for each of the various} \ W_n ({\rm BLUE}), n \to \infty\} \cap (\alpha (\infty) , \beta) \leqno (*)
$$
should correspond to an $O(3)$ move, defined like in \cite{8}, \cite{18}.

\smallskip

Notice that, modulo some embellishments, the drawing (B) is essentially figure~14.III.(A) from \cite{31}, what we called there the local Model~III. Of course, models I, II are present again now, but we will have also a new local model IV connected with the ${\rm Sing} \, M(\Gamma)$, which were absent in \cite{31}. This ``model IV'' is displayed in figure~1.5. Coming back to (B), we can see the triple points
$$
t_{np} = t_{np} (x = x_p , y = y_n) \in f M^3 (f) \, ,
$$
which come with the following accumulation pattern, inside $fX^2$
$$
\lim_{p=\infty} t_{np} = p_{n\infty} \in fM^2 (f) \cap f \, {\rm LIM} \, M_2 (f) \, .
$$
This ENDS our explanations concerning the figure~1.1. In the figures~1.2 and 1.3 we have displayed security walls $W_{\infty} ({\rm BLACK})$, and we give now EXPLANATIONS CONCERNING THESE FIGURES 1.2, 1.3.

\smallskip

Each of our present $W_{\infty}$ is larger than the corresponding security wall in \cite{31}. It {\ibf overflows}, like the $W_{\infty} ({\rm BLACK})_{H^0} (n)$ in figure~1.2, which enters deep inside the zone of the $1$-handle, it has free boundary and, most importantly, for each $W_{\infty}$ there is a unique $p_{\infty\infty} = p_{\infty\infty} (W_{\infty}) \in {\rm int} \, W_{\infty}$. For technical reasons to be clarified later, we have to allow now for $W_{\infty} - Y(\infty) \ne \emptyset$ and, moreover, in figure~1.3 the $W_{\infty} ({\rm BLACK})_{H^1} (n)$ {\ibf overflows} a bit beyond $W({\rm RED})_n$, towards the core of $H^1$. In the figure~1.2, $W_{\infty} ({\rm BLACK})$ continues beyond [CD] to the other $\partial H_j^2 (\gamma_n)$; at level $\widetilde M (\Gamma)$ we have here the $h_i^1 , h_j^1$ attached to $h^0$. At the level $X^2$, the $W_{\infty}$ is attached, in figure~1.2 to $\partial H_i^1 (\gamma_n) \cup W_n ({\rm BLUE}) \cup W_{n+1} ({\rm BLUE}) \cup \partial H_j^1 (\gamma_n)$ and similarly, in figure~1.3 the $W_{\infty}$ {\ibf is attached to} $\partial H^1 \cup W({\rm RED})_n \cup W({\rm RED})_{n+1}$. At least for $W_{\infty} ({\rm BLACK})_{H^0} (n)$ other attachments at the source $X^2 \mid H^0$ will be described later on. In particular, let us say that for $X^2 \mid H^0$ correspond not only the already mentioned $\partial H^1 = \partial H_i^1 (\gamma_n) \subset X^2 \mid H^0$, but also higher
$$
\partial H_i^1 (\gamma_{n+1}) + \partial H_i^1 (\gamma_{n+2}) + \ldots \subset X^2 \mid H^0 \, ,
$$
and our $W_{\infty} ({\rm BLACK})_{H^0} (n)$ is certainly to be attached to $\partial H_i^1 (\gamma_n)$, but not to the higher $\partial H_i^1$'s. In the figure~1.2, the $a_1 , b_1 , a_2 , s'_1 , s_1 , s'_2 , s_2 , \ldots $ are all mortal singularities of $f$. Actually, at $s_m^1$ there are effectively two singularities, involving $W_{\infty}$ and $W({\rm RED}) + W({\rm BLUE})$. At $s_m$ there is only one such. Very importantly, as things stand right now, with a simple-minded $W_{\infty}$, we would find that
$$
\lim_{m=\infty} (s_m , s'_m) = p_{\infty\infty} (W_{\infty}) \leqno (*)
$$
and similar things go for the figure~1.3 too. With this, if in (1.12) we take the $X^2 \mid H^{\lambda}$, $\lambda \leq 1$ to be the simple-minded union of smooth walls $W,W_{\infty}$, then $X^2$ WOULD FAIL to be {\ibf locally finite} at the points $p_{\infty\infty} (W_{\infty})$ and then, connected to this, the set ${\rm Sing} \, (f) \subset X^2$ would have unwanted accumulations at finite distance too, i.e. it would not be discrete. In order to take care of this difficulty, we will use a modified structure for the $W_{\infty} ({\rm BLACK})$'s, at their $p_{\infty\infty}$'s and then, for uniformity's reason (and also in anticipation of things to come), we will use a similar structure for the $W ({\rm BLACK})$'s, at their various points $p_{\infty\infty}$, $p_{\infty\infty}(S)$. I will explain here this modification for the $W_{\infty} ({\rm BLACK})$, where we start by considering a large enough circle $C(p_{\infty\infty}) \subset W_{\infty} ({\rm BLACK})$, centered at $p_{\infty\infty}$; next, we perform the following modification, deleting $p_{\infty\infty}$ and adding a compensating 2-handle instead,

\medskip

\noindent (1.15) \quad $X^2 \supset \{$the simple-minded, smooth $W_{\infty} ({\rm BLACK})\}  \Longrightarrow  \{ W_{\infty} ({\rm BLACK})  -  \{ p_{\infty\infty} \} \}$ 

\smallskip

\hglue1cm$\displaystyle{\underset{\footnotesize\overbrace{C(p_{\infty\infty})}}{\cup}} D^2 (H (p_{\infty\infty})),$

\medskip

\noindent where $D^2 (H(p_{\infty\infty}))$ (which is sometimes denoted $D^2 (p_{\infty\infty})$), is a round disk of boundary $C(p_{\infty\infty})$. The $C(p_{\infty\infty})$ becomes now a new mortal singularity of (\ref{eq1.1}), i.e. it consists of non-immersive points, but the image $fX^2$ does not feel the modification (1.15). The local finiteness of $X^2$ has now been restored. When it comes to the points of type $p_{\infty\infty}$, like we may see in the figures~1.1, 1.2, 1.3, 1.5, they come in two kinds
$$
\includegraphics[width=17cm]{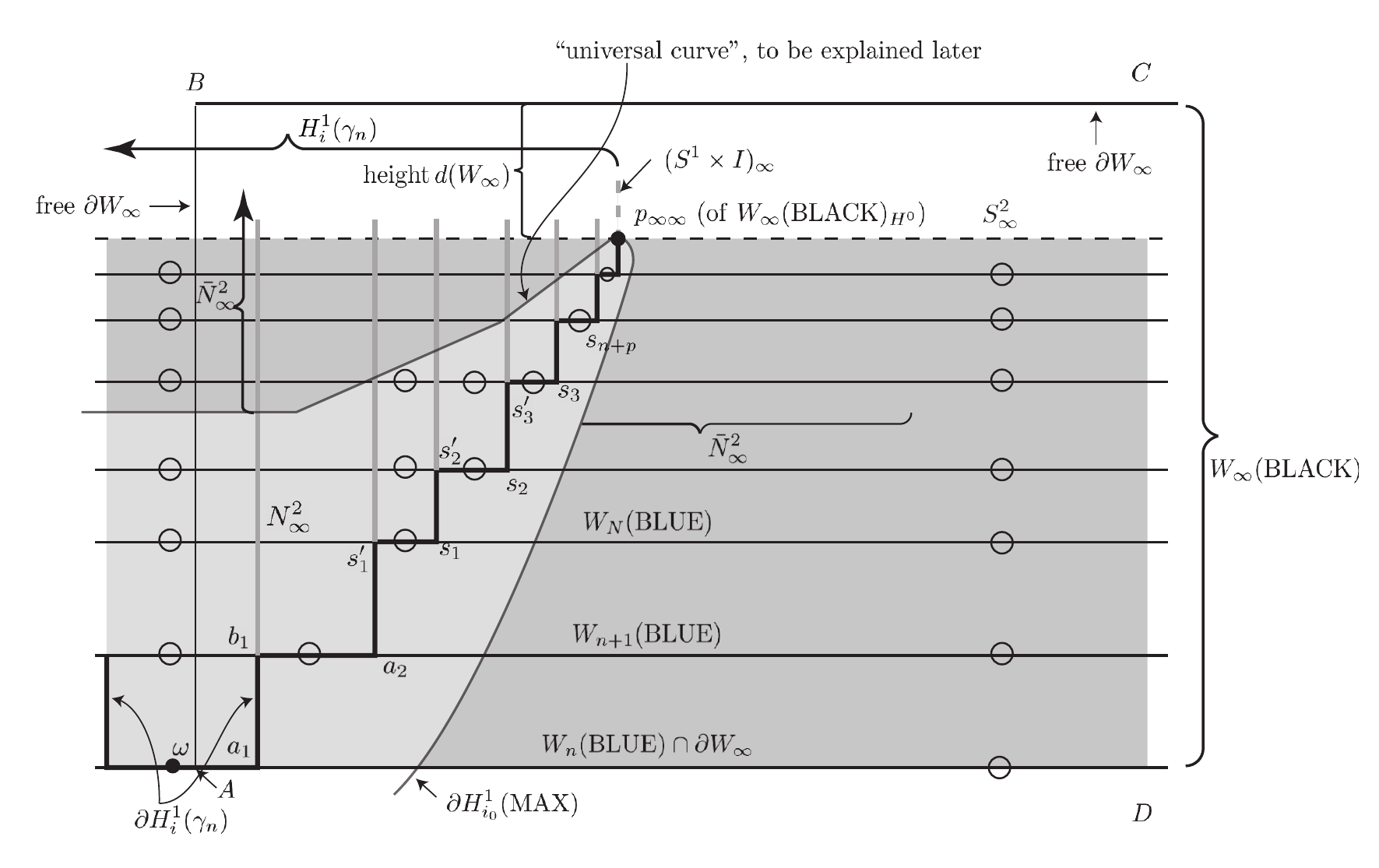}
$$
\label{fig1.2}
\centerline {\bf Figure 1.2.} 

\smallskip

\begin{quote} 
We see here one half of the security wall $W_{\infty} ({\rm BLACK})_{H^0} (n) \, \subset$ $X^2 \mid H^0$, which continues to the right, beyond the line [CD].
\end{quote}
$$
\includegraphics[width=16cm]{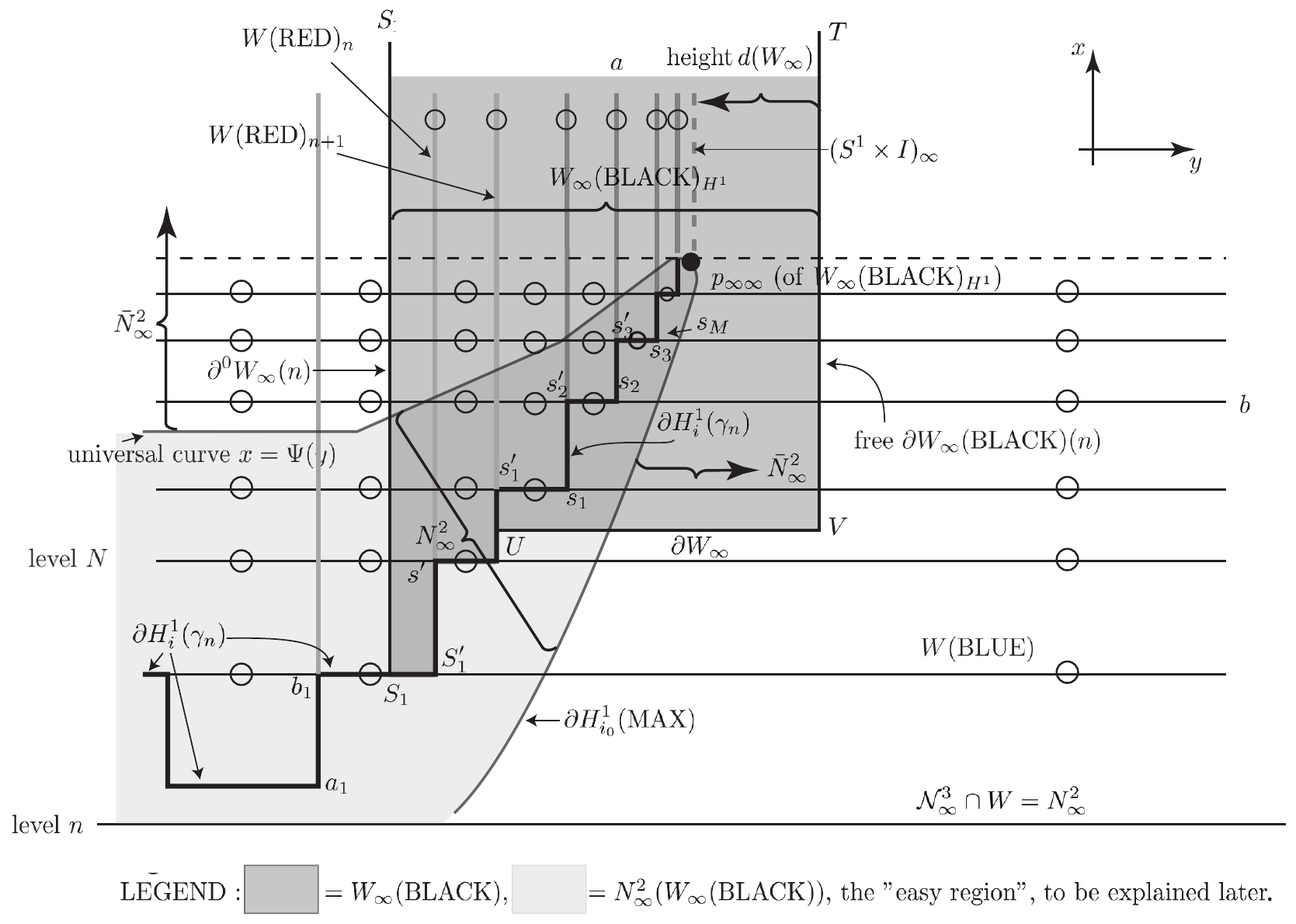}
$$
\label{fig1.3}
\centerline {\bf Figure 1.3.} 

\smallskip

\begin{quote}
We display here one half of the $W_{\infty} ({\rm BLACK})_{H^1} (n) \subset X^2 \mid H^1$, which continues beyond [ST], to the other end of $H^1$. The line [SS$_1$] is supposed to be very close to $W({\rm RED})_n$, on its left side.
\end{quote}

\bigskip

\noindent (1.15.0) \quad $X^2 \supset \{ p_{\infty\infty} ({\rm all}) \} = \{$the $p_{\infty\infty}$ (proper) which come with $fX^2$ locally embeddable in $R^3$; they are smooth points for $X^2$, but ramification points for ${\rm LIM} \, M_2 (f)$. They appear as corners of the main dotted hexagon in figure~1.1.(A), and each $W_{\infty} ({\rm BLACK})$, figures~1.2, 1.3 has one such point too$\} + \{$the $p_{\infty\infty} (S)$, which are created via the zipping by the immortal singularities $\bar S \subset {\rm Sing} \, \widetilde M (\Gamma)$, like in the formula (1.22) below. At the $p_{\infty\infty} (S)$'s which occur in figure~1.5 the zipping of two $W ({\rm BLACK})$'s cutting through each other at the target steps. The $p_{\infty\infty} (S)$'s are immortal singularities of $fX^2\}$.

\medskip

For reason of smoothness of the exposition and for further needs too, the transformation (1.15) will be performed now, not only at the $p_{\infty\infty} ({\rm proper}) \in W_{\infty} ({\rm BLACK})$, but at all $p_{\infty\infty}$'s in (1.15.0). This {\ibf new revised definition} of $X^2$ comes without technical problems of its own, because at the simple-minded level the (1.15.0) is a discrete subset of smooth points (of the simple-minded $X^2$).

\smallskip

This also means that now, for each individual $p_{\infty\infty}$'s there is an additional folding map in any zipping strategy for
$$
X^2 \longrightarrow \widetilde M (\Gamma) \, .
$$

The various free boundaries $\partial W_{\infty} ({\rm BLACK})_{H^0}$ or $\partial W_{\infty} ({\rm BLACK})_{H^1}$ accumulate, respectively, on $S_{\infty}^2$ or on $(S^1 \times I)_{\infty}$; see here the figures~1.2 and 1.3. 

\smallskip

Besides the zigzag of $\partial H_i^1 (\gamma_n)$ in figure~1.2, staying still at the source, there are infinitely many zigzags $\partial H_i^1 (\gamma_{m > n})$ to the left of it and finitely many to the right of it, $\partial H_i^1 (\gamma_{\ell < n})$'s. {\ibf All} of them are part of $X^2 \mid H^0$.

\smallskip

When we go to a figure~1.2 at the target, then more zigzags will appear. The extreme right one will be denoted $\partial H_{i_0}^1 ({\rm MAX})$ (but for typographical reasons it is drawn as a curved line, rather than a zigzag). In the figure at the target, between our $\partial H_i^1 (\gamma_n)$ and the $\partial H_{i_0}^1 ({\rm MAX})$, only finitely many new zigzags will appear.
\vglue -8mm
$$
\includegraphics[width=15cm]{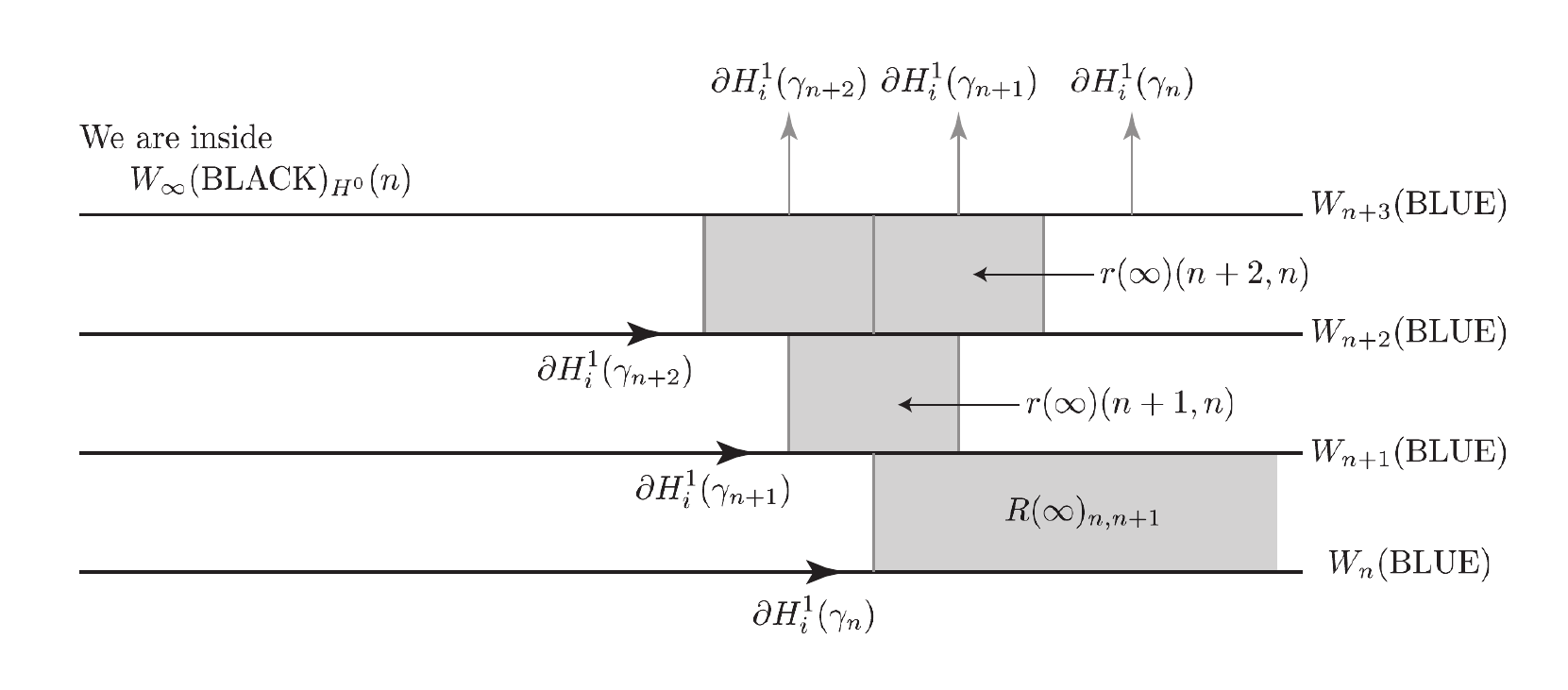}
$$
\vglue -5mm
\label{fig1.3.1}
\centerline {\bf Figure 1.3.1.} 
\begin{quote}
The $W_{\infty} ({\rm BLACK})_{H^0} (n)$, which contains the plane of the figure, is glued, at level $X^2 \mid H^0$ to $W_n ({\rm BLUE}) + W_{n+1} ({\rm BLUE}) + \partial H_i^1 (\gamma_n)$, and to the $R(\infty)_{n,n+1} \subset W_{\infty} ({\rm BLACK})_{H^0} (n)$. The $X^2 \mid H^0$ contains not only $\partial H_i^1 (\gamma_n)$, but also the higher $\partial H_i^1 (\gamma_{N > n})$, but our $W_{\infty} (n)$ is {\ibf not glued} to them. None of the rectangles $R(\infty)$ or the smaller $r(\infty)$ have any free face, they are all glued to the rest of the $X^2 \mid H^0$. But, while $R(\infty)_{n,n+1} \subset W_{\infty} ({\rm BLACK})_{H^0} (n)$, the infinitely many $\overset{\circ}{r} (\infty)$'s are, for right now just empty spots inside
$$
X^2 \mid H^0 - \{\mbox{the} \ W_{\infty} ({\rm BLACK})_{H^0}\mbox{'s}\} \, ;
$$
we will have to deal with them, later on, see $(*_2)$ below.
\end{quote}

\smallskip

In the figures 1.1, 1.2, 1.3 appear also, ``universal curves'', concerning which much more will be said later. For right now, it will suffice to know that $W_{\infty} ({\rm BLACK})$ is divided by its universal curve into an ``easy region'' $N_{\infty}^2$, between the universal curve and $\partial H_{i_0}^1 ({\rm MAX})$ (including the $\partial H_{i_0}^1 ({\rm MAX})$ too), and also a ``difficult region'' $\overline N_{\infty}^2$ containing two components, one higher than the universal curve, the other one lower than $\partial H_{i_0}^1 ({\rm MAX})$.

\smallskip

We will decide that $N^2 (W_{(\infty)})$ is closed and contains $p_{\infty\infty}$, so that
$$
W_{(\infty)} ({\rm BLACK}) \cap {\rm LIM} \, M_2 (f) - \{ p_{\infty\infty} (W_{(\infty)}) \} \subset \overline N^2 (W_{(\infty)}) \, .
$$
The $W_{\infty}$'s will be chosen far from the $W({\rm BLACK})$'s and the $\bar S$'s. This ENDS our explanations for the figures~1.2, 1.3.

\smallskip

Inside $H^{\lambda}$, when $\lambda = 0$, let $W_1, W_2 , W_3 , \ldots$ be the succession of walls of natural BLUE colour, all in $X^2 \mid H^0$. Inside the same $X^2 \mid H^0$ we will also have $\partial H^1 (\gamma_n) \cap H^0$'s, and for the good match between the indices of the  $W_n ({\rm BLUE})$ and of $\partial H^1 (\gamma_n)$ see the item (1.15.1) below. For each $W_n$ there is, inside $X^2 \mid H^0$ a $W_{\infty} (n)_{H^0} \equiv W_{\infty} ({\rm BLACK})_{H^0} (n)$ resting on it, and glued both to $W_n$ and to $W_{n+1}$. The $W_{\infty} (n)_{H^0}$, half of which is seeable in figure~1.2 goes from $H_i^1 (\gamma_n)$ to $H_j^1 (\gamma_n)$ and, at level $X^2 \mid H^0$ is glued to
$$
\partial H_i^1 (\gamma_n) \cup W_n \cup W_{n+1} \cup \partial H_j^1 (\gamma_n) \, , \leqno (*_1)
$$
and continues beyond, to its free boundary. The $(*_1)$ determines a rectangle $R(\infty)_{n,n+1} \subset W_{\infty} (n)_{H^0}$, see the figure~1.3.1. [The $(*_1)$ with its $R(\infty)$ is necessary for keeping $X^2 \mid H^0$ be GSC, as (1.11.1) demands. Explanations: The $W({\rm BLUE})_n$'s are glued to all the $\partial H^1 \cap H^0$'s to begin with. With this, the $\partial H^1 \cap H^0 \cup W_n \cup W_{n+1}$ creates a surface of genus $g \geq 1$ and this is rendered harmless by the $R(\infty)$'s. Then, the next $g \geq 1$ is taken care of by the $W_{\infty} (n+1)_{H^0}$, a.s.o.] So far, we discussed $H^{\lambda}$ when $\lambda =0$. But mutatis-mutandis, all this, including the $R(\infty)$'s goes for $\lambda = 1$ too; see here the discussion which comes with the figure~1.3.2. IF $W_1 ({\rm RED}) , W_2 ({\rm RED}),\ldots$ are the successive natural $W$'s of $H_i^1 (\gamma)$, then a {\ibf pair} of $W_{\infty} ({\rm BLACK})_{H^1}$'s will rest on $W_1$; (and get a bit beyond it, towards the core of $H^1$, see the line $[S,S_1]$ in the figure~1.3), then another pair of $W_{\infty} ({\rm BLACK})_{H^1}$'s will rest in $W_2$, a.s.o. The whole structure should be suggested by figure~1.4. But this structure will be needed only much later, in the third paper of our QSF trilogy.

\smallskip

We turn back now to $H^0$. For any pair $(m,p)$ when $m > p$, the
$$
W_m ({\rm BLUE}) + W_{m+1} ({\rm BLUE}) + \partial H_i^1 (\gamma_p) + \partial H_i^1 (\gamma_{p+1}) \leqno (*_2)
$$
creates an embedded torus inside $X^2 \mid H^0$, which we have to renders harmless (for the sake of $X^2 \mid H^0 \in {\rm GSC}$) by adding a small doubly shaded rectangle $r(\infty)(m,p)$, like in figure~1.3.1. In that figure only the cases $p \geq n$ are visible. But the $r(\infty)$'s, which have to be {\ibf added} into the definition of $X^2 \mid H^0$, completing the (1.12), are {\ibf not} inside any of the $W_{\infty} ({\rm BLACK})$'s. Except for rendering the $X^2$ be GSC, the $r(\infty)$'s will not play any other role in this paper.

\smallskip

Contrary to the $R(\infty)$'s the reason for filling the $r(\infty) (m,p) \subset X^2 \mid H^0$ vanishes when one moves from $X^2 \mid H^0$ to $X^2 \mid H^1$. Here is the reason why. When moving from $X^2 \mid H^0$ to $X^2 \mid H^1$, the $(*_2)$ which has led to the creation of our $r(\infty)(m,p)$ is to be replaced by
$$
W_m ({\rm RED}) + W_{m+1} ({\rm RED}) + \partial H_i^2 (\gamma_p) + \partial H_i^2 (\gamma_{p+1}) \, , \leqno (*)_3
$$
which cut by a plane orthogonal to the core line of $H^1$, would create something like in the figure~1.3.2.
\vglue -5mm
$$
\includegraphics[width=100mm]{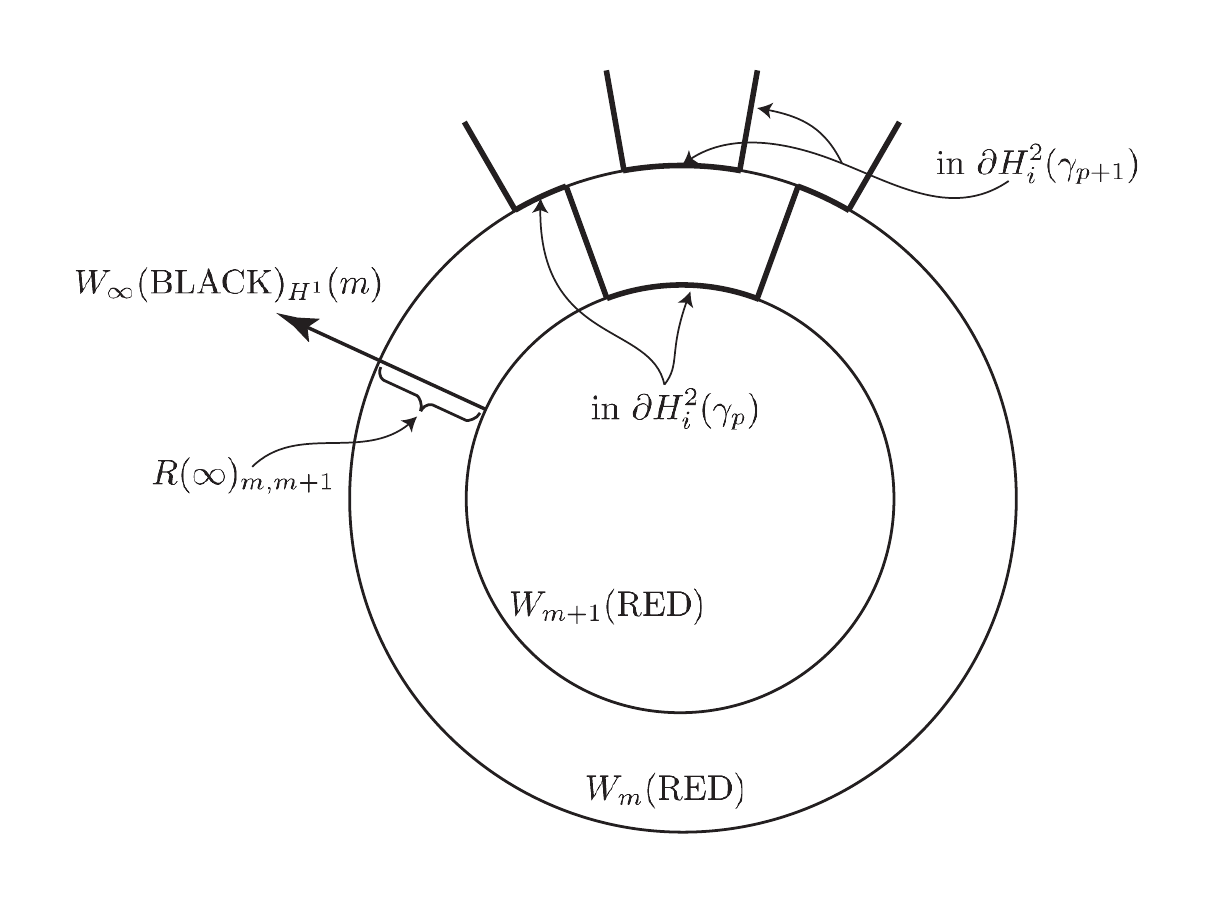}
$$
\label{fig1.3.2}
\vglue -10mm
\centerline {\bf Figure 1.3.2.} 
\begin{quote}
\vglue -3mm
A section through $H^1$. Here there are no longer unwanted tori of the type which have needed the doubly shaded $r(\infty)$'s in figure~1.3.1. The $X^2 \mid H^1$ does not need $r(\infty)$'s.
\end{quote}

\bigskip

When we consider the various $r(\infty) \subset X^2 \mid H^0$, created, like in figure~1.3.1 by $H^0$ and the adjacent $H_j^1 (\gamma_n)$'s (fixed $j$), then we have, as the only possible accumulation of the rectangles $r(\infty)$, that $\underset{m=\infty}{\lim} \, r(\infty) (m,p) = \{$some point on the limit circle $S_{\infty}^1 = S_{\infty}^2 \cap (S^1 \times I)_{\infty}$, coming with $H^0 , H^1\}$. No $r(\infty)$'s are necessary for $H^1$, as we just saw.

\smallskip

For the $H^0 \subset Y(\infty)$ some additional adjustments will be necessary at the level $X^2 \mid H^0$. Let us say that $H^0$ corresponds to $h^0 \subset \widetilde M (\Gamma)$, to which the $h_1^1 , h_s^1 , \ldots , h_{\alpha}^1$ are attached. For each of these $h_i^1$'s, there is an infinite family of 1-handles of $Y(\infty) : H_i^1 (\gamma_1) , H_i^1 (\gamma_2) , \ldots$ all attached to $H^0$, coming closer and closer to $S_{\infty}^2 = \delta H_0$.

\medskip

\noindent (1.15.1) \quad For any given $H^0$ and any level $n$, all the $\partial H_i^1 (\gamma_n)$'s, where $i=1,2,\ldots,\alpha$, should reach the same $W_n ({\rm BLUE})$ and, moreover, the $\partial H_i^1 (\gamma_1)$ should reach the innermost $W_1 ({\rm BLUE})$ of $H^0$.

\medskip

We can achieve these conditions by artificially adding, when necessary, to the $\partial H_i^1 (\gamma)$'s annuli $W({\rm RED})$ glued at level $X^2$ to the $W_n ({\rm BLUE}) + W_{n+1} ({\rm BLUE}) + \ldots + W_N ({\rm BLUE})$. These pieces are to be added to the zigzags $ZZ_i (j) \equiv \partial H_i^1 (\gamma_j) \subset X^2 \mid H^0$. The $[a_1 , b_1 , a_2 , s'_1]$ in figure~1.2 could be such an addition.

\medskip

\noindent (1.16) \quad For any given pair $(H^0 , n)$, all the $\partial H_i^1 (\gamma_n)$'s for $i=1,2,3,\ldots , \alpha$ should be joined {\ibf arborescently} by the $W_{\infty} ({\rm BLACK})_{H^0} (n)$'s which rest on $W_n ({\rm BLACK})$ and which are glued to them. So, there are exactly $\alpha-1$ such $W_{\infty}$'s, and the global picture should be readable on figure~1.4.

\medskip

In principle, the figures~1.2, 1.3 correspond to unique individual handles $H^0 , H^1$. But $H^0 , H^1$ correspond to $h^0 , h^1 \subset \widetilde M (\Gamma)$ and, when all the $H^0 , H^1$'s corresponding to these $h^0 , h^1$ are taken into account, then we get at the target, {\ibf complete figures}, similar to 1.2, 1.3, with the same accumulation pattern, but much denser. For a given $W_{\infty}$, when moving from figures~1.2, 1.3 ``at the source'', as drawn, to the complete figures ``at the target'', then there are no additional glueings, similar to $(*_1)$, prior to any zipping. In the figures~1.2, 1.3 we have defined heights $d(W_{\infty})$ for the $W_{\infty}$'s. We also have quantities $k (W_{\infty}) \equiv \{$the level at which $W_{\infty}$ is attached$\}$. With this, we will have
\setcounter{equation}{16}
\begin{equation}
\label{eq1.17}
\lim_{n=\infty} d(W_{\infty} ({\rm BLACK})(n)) = 0 \, , \quad \lim_{n=\infty} k(W_{\infty} ({\rm BLACK})(n)) = \infty \, ,
\end{equation}
in the complete figures. So, at the level of $\widetilde M (\Gamma)$ the $W_{\infty}$'s can only accumulates on $S_{\infty}^2 \cup (S^1 \times I)_{\infty}$, without generating their own limit walls. For given $W_{\infty}$, the $M_2(f) \cap W_{\infty}$ accumulates on tree-shaped figures
\begin{equation}
\label{eq1.18}
{\rm LIM} \equiv (S_{\infty}^2 \cup (S^1 \times I)_{\infty}) \cap W_{\infty} \subset {\rm LIM} \, M_2 (f) \subset X^2 \, .
\end{equation}
The limiting position of the LIM's inside $S_{\infty}^2 \cup (S^2 \times I)_{\infty}$ is denoted by $\lim \, {\rm LIM}$ and it follows from (\ref{eq1.17}) that $\lim \, {\rm LIM}$ is reduced to a collection of independent arcs, living inside $S_{\infty}^2$ OR inside $(S^1 \times I)_{\infty}$; see figure~1.4.

\smallskip

With all our various modulations just described, in particular the (1.15.1), we have the final $X^2 \mid H^{\lambda}$'s out of which we put together the global $X^2$. This object has only mortal singularities ($=$ non-immersion points), which are either of the undrawable type (\`a la \cite{18}) or the circles $C(p_{\infty\infty})$.

\smallskip

So ${\rm Sing} \, (f) = \{$a discrete set$\} \, \cup \, \sum C(p_{\infty\infty})$, $f \mid {\rm Sing} \, (f)$ injects and $f \, {\rm Sing} \, (f)$ $\cap \ {\rm Sing} \, \widetilde M (\Gamma) = \emptyset$. The $fM^2 (f) \cap {\rm Sing} \, \widetilde M (\Gamma) \ne \emptyset$ which is a big novelty with respect to \cite{31} will be discussed later, in great detail.

$$
\includegraphics[width=175mm]{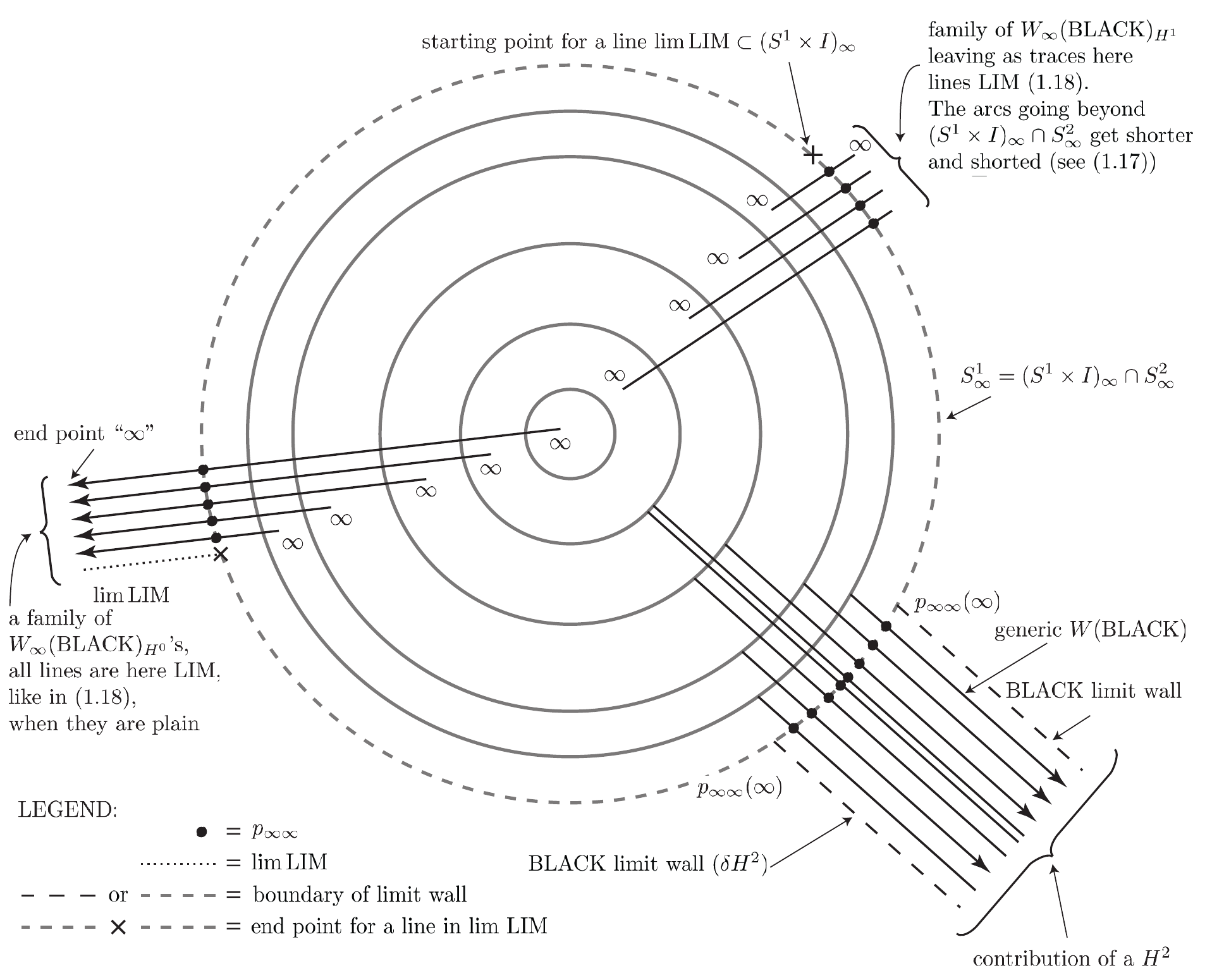}
$$\label{fig1.4}
\centerline {\bf Figure 1.4.} 

\smallskip

\begin{quote}
This is a piece of a BLUE limit wall $S_{\infty}^2$ cut by walls $W({\rm RED})$, $W({\rm BLACK})$ and $W_{\infty}({\rm BLACK})_{H^{\varepsilon}}$. This figure is at the target. The RED and BLACK limit walls contribute only with dotted lines, where they rest on $S_{\infty}^2$. The plain lines are traces of transversal contacts ($W_{\infty}$ or $W({\rm COLOUR})$) $\cap \, S_{\infty}^2$. When the lines in question come from $W_{\infty}$'s then their limiting position is $\lim \, {\rm LIM}$. None of the various dotted lines have counterparts in $X^2$. The points marked ``$\infty$'' are in the free part of $\partial W_{\infty} ({\rm BLACK})_{H^{\varepsilon}}$ and/or live at infi\-nity. The packages of $W_{\infty} ({\rm BLACK})$, $W ({\rm BLACK})$ which come with arrows continue to some other $S_{\infty}^1$ inside our $S_{\infty}^2$.
\end{quote}

\bigskip

At the infinity of $X^2$ we have the following basic graph
\begin{equation}
\label{eq1.19}
\sigma_1 (\infty) \equiv \bigcup \left( \bigcap \, \{\mbox{limit walls}\} \right) \cup \bigcup \lim \, {\rm LIM}
\end{equation}
and where
$$
\bigcup \left( \bigcap \, \{\mbox{limit walls}\} \right) = \sum S_{\infty}^1 (=S_{\infty}^2 \cap (S^1 \times I)_{\infty}) 
\cup \sum \left[{\rm Hex}_{\infty} ({\rm BLACK}) \cap (S_{\infty}^2 \cup (S^1 \times I)_{\infty})\right] \subset \Sigma_1 (\infty) \, . 
$$

By now enough has been said about (\ref{eq1.1}) so that the implication
$$
\{\mbox{Theorem 1.3 and its complement 1.3}\} \Longrightarrow \mbox{Theorem 1.1},
$$
should be clear. For instance, the second finiteness condition (\ref{eq1.2}) follows from the complement 1.3. We will take now a closer look at the (\ref{eq1.3}). Together with the ideal $\Sigma_1 (\infty)$ comes the useful set
\begin{equation}
\label{eq1.20}
\Sigma_2 (\infty) = g(\infty) \, Y(\infty) \cap \Sigma_1 (\infty) = \bigcup \, \{\mbox{int of the attaching zones} \ \partial h^1 , \partial h^2 \} \, ,
\end{equation}

\medskip

\noindent and also the following kind of objects, which are all {\ibf closed} subsets of $\Sigma_1 (\infty)$, respectively of $X^2$, respectively of $fX^2$,

\medskip

\noindent (1.21.1) \quad ${\rm Sing} \, \widetilde M(\Gamma) = \sum$ squares $\bar S \subset \bigcup S_{\infty}^2 \subset \Sigma_1 (\infty)$,

\medskip

\noindent (1.21.2) \quad ${\rm LIM} \, M_2 (f) = f^{-1} (fX^2 \cap \Sigma_1 (\infty)) \subset X^2$, same as in (\ref{eq1.3}),

\medskip

\noindent (1.21.3) \quad $f \, {\rm LIM} \, M_2 (f) = fX^2 \cap \Sigma_1 (\infty) \subset fX^2$.

\bigskip

\noindent COMMENTS. When it comes to $f \, {\rm LIM} \, M_2 (f)$ being closed (a feature which is {\ibf not} a purely mechanical consequence of the inputs) then, when we are outside the (1.21.1), this was proved in \cite{31}, while for each individual $\bar S \subset {\rm Sing} \, \widetilde M(\Gamma)$ the set $\bar S \cap f \, {\rm LIM} \, M_2 (f) = {\rm int} \, \bar S \cap f \, {\rm LIM} \, M_2 (f)$ will be soon described explicitly, and this should make our claimed closedness obvious.

\smallskip

We will call FAKE LIM $M_2 (f)$ the total contributions of the dotted lines $[p_{\infty\infty} (S) , P] \subset W$ (BLACK complete or reduced) in all the figures 1.1. This is {\ibf not} included in (1.21.2), (1.21.3), and it will never play any active role in our discussions.

\smallskip

In the formulae (1.21.2), (1.21.3) it is only the BLUE and RED limit walls of $\Sigma_1(\infty)$ which play an active role, generating the ${\rm LIM} \, M_2 (f)$ in the formulae in question. Inside $\widetilde M (\Gamma)$, the BLACK limit walls are disjoined from $fX^2$, and hence they are mute in the contexts of the same formulae.

\smallskip

We discuss now the impact of the immortal square $\bar S \subset {\rm Sing} \, \widetilde M (\Gamma)$. Contrary to what has happened in \cite{31} the $fX^2$ has now immortal singularities of its own. By definition, these are the points where $fX^2$ is NOT locally embeddable in $R^3$. Their set ${\rm Sing} \, fX^2 \subset fX^2$ is generated as follows
\setcounter{equation}{21}
\begin{equation}
\label{eq1.22}
{\rm Sing} \, fX^2 = \sum_{\bar S} {\rm int} \, \bar S \cap fM^2 (f) \subset f \, {\rm LIM} \, M_2 (f) \, .
\end{equation}
What is more serious, the $fX^2$ is {\ibf not locally finite} at the points $\partial \bar S \cap fX^2$, generically denoted $p_{\infty\infty} (S)$ (and see here the figure~1.1 too), where (\ref{eq1.22}) accumulates. Typically the arcs $[\omega,y]$, $[\omega_- , y_-]$ in figure~1.5 are such points $x \in \partial S \cap fX^2$ blown into arcs. We shall show later on how this serious difficulty, a complete novelty with respect to \cite{31}, is to be handled; it will require to thicken our $fX^2$ into dimension three, to begin with, but then it will require also much more than that.

\smallskip

The geometry of (\ref{eq1.22}), or what we have loosely called above the local model~IV, is explained in the figure~1.5, and here are some comments concerning it. Both (A) and (B) live in the plane of some $W({\rm BLACK \, complete})$ and, in the neighbourhood of each of their two $p_{\infty\infty} (S)$'s, except for some changes in colour, they are like figure~1.1.(B). We see here, in (A), (B) the unique $W({\rm BLACK \, complete})$'s of $X^2 \mid H_i$, $X^2 \mid H_j$. When we go to the complete figures, i.e. the figures at the target, with all the $f \, \underset{\gamma}{\sum} \, X^2 \mid H_i^2 (\gamma)$, $f \, \underset{\gamma}{\sum} \, X^2 \mid H_j^2 (\gamma)$ thrown in, then the (A), (B) get enriched with infinitely more walls of type $W({\rm BLACK \, complete})$, converging to the dotted limit walls.
\newpage
\vglue -15mm
$$
\includegraphics[width=175mm]{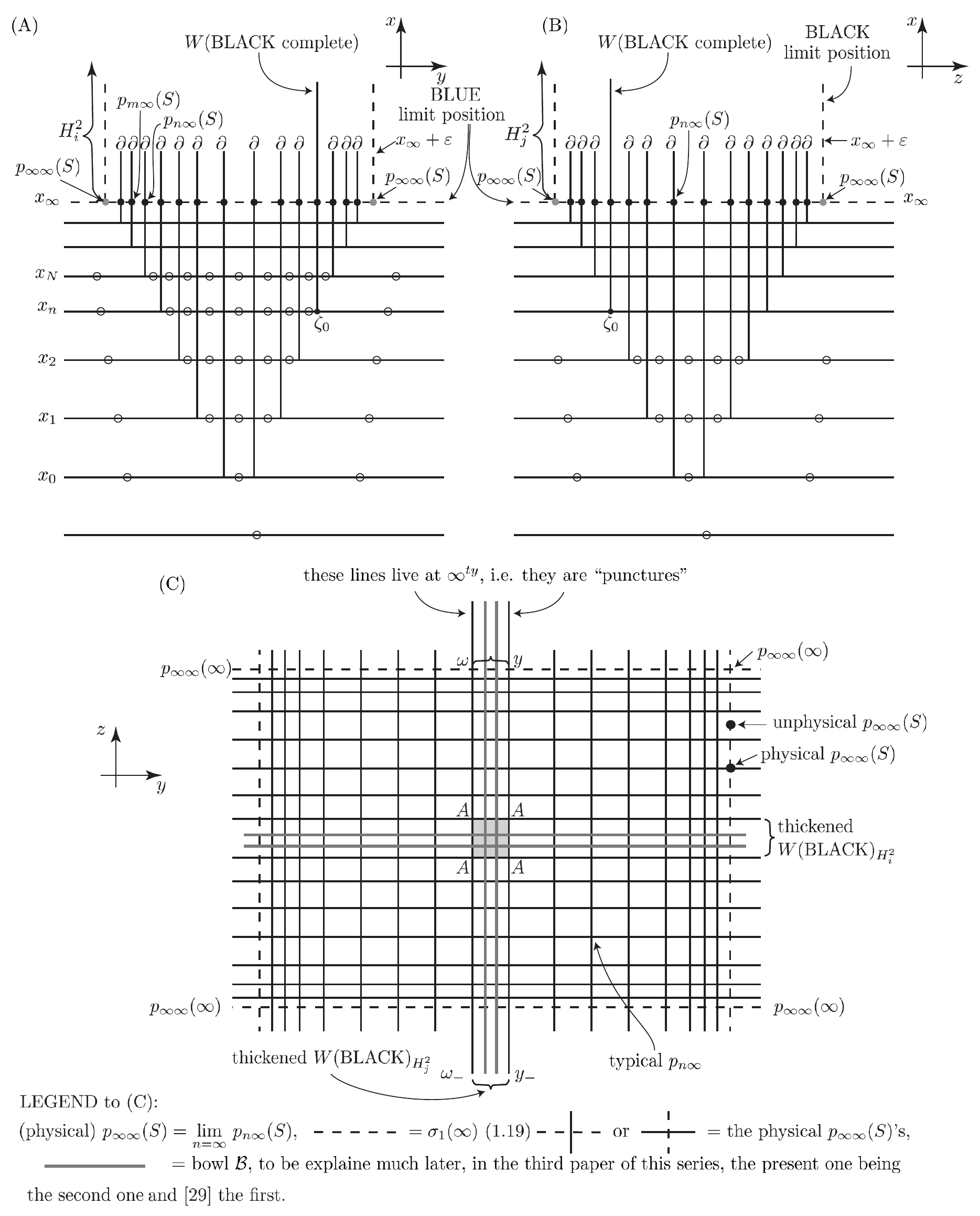}
$$
\label{fig1.5}
\vglue -7mm
\centerline {\bf Figure 1.5.} 

\smallskip

\begin{quote}
We see here the ``local'' model IV, attached to an immortal singularity $\bar S \subset \widetilde M^3 (\Gamma)$. The initial large $\bar S$ splits into a double infinity of smaller squares $S \subset \{ fX^2 \ \mbox{thickened into} \ 3^{\rm d}\}$. Never mind, right now what such a thickening should mean. Both (A) and (B) live inside the same $H_k^0 (\gamma)$. The (C), which lives at $x=x_{\infty}$, represents a small piece of the ideal boundary $\delta H_k^0 (\gamma)$.

\smallskip

The large square is the initial $\bar S \subset \widetilde M^3 (\Gamma)$. The (A) is in the plane $z = z_p$ of $W({\rm BLACK})_{H_j^2}$, and (B) in the plane $y=y_q$ of $W({\rm BLACK})_{H_i^2}$. The physical intersection $W({\rm BLACK})_{H_i^2} \ \cap$ $W$ ({\rm BLACK}$)_{H_j^2}$ stops at $x=x_{\infty}$. All the $W({\rm BLACK})$'s in the doubly infinite collection are to be thickened into things like $W({\rm BLACK}) \times [-\varepsilon , - \varepsilon]$, when we will go to the not yet defined $\Theta^3 (fX^2)$, in the next section. We only suggested this thickening for the $W({\rm BLACK})_{H_i^2}$, $W({\rm BLACK})_{H_j^2}$ in our drawing (C). With this thickening in mind, the $(w,y),(w_- , y_-)$ are arcs $p_{\infty\infty} (S) \times (-\varepsilon,\varepsilon)$.

\smallskip

Any intersection point of plain lines in the doubly infinite lattice in (C) is a $p_{n\infty} (S)$ or $p_{\infty n} (S)$, and these are the immortal singularity for $fX^2$. When we will go to $\Theta^3 (fX^2)$ these will be thickened into small disjoined immortal squares called generically $S$. So, every initial immortal square $\bar S \subset \widetilde M (\Gamma)$ breaks into a double infinity of smaller immortal $S \subset \Theta^3 (fX^2)$. The [AAAA] in (C) stands for such an $S$.
\end{quote}

\bigskip

As a big novelty with respect to \cite{31}, we have now double lines $\mbox{BLACK} \cap \mbox{BLACK} \subset M^2 (f)$, stopping short and abruptly when they reached the immortal singularity $\bar S \subset {\rm Sing} \, \widetilde M (\Gamma)$. One should also think that, in the context of a figure like 1.5.(C), the $x$-coordinate is now a $1^{\rm d}$ {\ibf train-track} manifold with two branches $x(H_i^2)$ and $x(H_j^2)$. With this, we can explain now the FAKE LIM $M_2 (f)$ in a fashion similar to the formulae (1.21.2) and (1.21.3). When it comes to the points where $\sum_1 (\infty)$ fails to be a 2-manifold (see the IMPORTANT REMARK after formula (\ref{eq1.14})), then the branching points of $\sum_1 (\infty)$, the $p_{\infty\infty} (\infty)$ never touch $fX^2$ while all the $p_{\infty\infty} (S)$'s are ramification ($=$ non-2-manifold) points for $\sum_1 (\infty)$, all touching $fX^2$.

\smallskip

Look in particular at the $p_{\infty\infty} (S)$ thickened into $[\omega,y]$, $[\omega_- , y_-]$ in figure~1.5.(C). Through $[p_{\infty\infty} (\infty)$, $p_{\infty\infty} (S) \approx [\omega,y] , p_{\infty\infty} (\infty)]$ goes a BLACK branch of $\sum_1 (\infty)$, in the direction $x (H_j^2) > 0$. Then, we also have our $W({\rm BLACK})$ $\subset H_i^2$, figure~1.1, shooting in the direction $x(H_i^2) > 0$. If we artificially decree that $x(H_i^2) = x(H_j^2)$, thereby collapsing the train-track $x$-coordinate line into a usual line, then the $[p_{\infty\infty} (S),P] \subset \{ W({\rm BLACK})$ in the figure~1.1.(A)$\}$ is defined by a formula like the (1.21.2), let us say by

\medskip

\noindent (1.22.2) \quad FAKE LIM $M_2 (f) =  f^{-1}$ (FAKE intersection $fX^2 \, \cap \, \sum_1 (\infty)) \subset X^2$, but unlike what had happened in (1.21.2), it is now the BLACK $\sum_1 (\infty)$ which will contribute actively.

\medskip

\noindent But the (1.22.2) will normally be mute, throughout this paper; there will be no role for it.

\smallskip

Just like the $p_{\infty\infty} (\infty)$'s, the unphysical $p_{\infty\infty} (S)$'s in the figure~1.5 will never play any role, they are outside our universe and irrelevant. This ends our discussion of figure~1.5, for the time being. The proof of theorem~1.1 is by now completed too, we will just add the following complement to it.

\bigskip

\noindent {\bf Lemma 1.4.} {\it Consider an infinite sequence $x_1 , x_2 , \ldots \in f \widehat M^2 (f) \subset fX^2$. At least one of the following three conditions is then satisfied for our sequence, and it is not claimed that they are mutually exclusive}
\begin{enumerate}
\item[1)] {\it The $x_1 , x_2 , \ldots$ accumulates on a compact subset of $f \widehat M^2 (f)$.}
\item[2)] {\it The $x_1 , x_2 , \ldots$ remains confined at finite distance, i.e. there is a finite union of handles of $\widetilde M (\Gamma)$ which contains it. For the part of $x_1 , x_2 , \ldots$ which lives inside any specific $f X^2 \mid h_k^{\lambda}$, there is a subsequence which accumulates on $(f \, {\rm LIM} \, M_2 (f) \cup \sigma_1 (\infty)) \mid h_{k}^{\lambda}$, with $\sigma_1 (\infty)$ like in {\rm (\ref{eq1.19})}.}
\item[3)] {\it Our sequence is not confined at finite distance. There exists then a subsequence $x_{i_1} , x_{i_2} , x_{i_3} , \ldots$ which is s.t
$$
x_{i_n} \in fX^2 \mid h_{k(n)}^{\lambda (n)} \, , \quad h_{k(n)}^{\lambda (n)} \subset \widetilde M (\Gamma) \quad \mbox{and} \quad \lim_{n=\infty} k(n) = \infty \, .
$$
}
\end{enumerate}

\bigskip

\noindent EQUIVARIANT ZIPPING. Since our (\ref{eq1.1}) is a REPRESENTATION, we have $\Psi (f) = \Phi (f)$ and then according to Lemma~2.4 in \cite{29} our $X^2 \overset{f}{\longrightarrow} \widetilde M (\Gamma)$ is zippable, i.e. the map $f$ can be exhausted by a sequence of folding maps, more precisely of $O(i)$ moves with $0 \leq i \leq 3$ defined like in \cite{8}, \cite{18}, \cite{31}
\setcounter{equation}{22}
\begin{equation}
\label{eq1.23}
X^2 = X_0^2 \overset{P_1}{-\!\!\!\longrightarrow} X_1^2 \overset{P_2}{-\!\!\!\longrightarrow} X_2^2 \overset{P_3}{-\!\!\!\longrightarrow} \ldots\ldots\ldots fX^2 = X^2 / \Phi (f) \, .
\end{equation}
For what follows next we need to allow the $O(0)$ to be sometimes trivial too, i.e. not meeting triple points at all. Moreover, we will allow for each $P_j$ in (\ref{eq1.23}) to be not just one single $O(i)$-move but to take the form $P_j = \underset{k}{\sum} \, P_{jk}$ where each individual $P_{jk}$ is such an elementary move.

\smallskip

Let us also choose an atlas $M(\Gamma) = \underset{p}{\bigcup} \, U_p$ such that, when we go to $\widetilde M (\Gamma) \overset{\pi}{\longrightarrow} M(\Gamma)$, then each $\pi^{-1} \, U_p$ is product, with $\pi^{-1} \, U_p = \underset{k}{\sum} \, U_{pk}$.

\bigskip

\noindent {\bf Lemma 1.5.}
\begin{enumerate}
\item[1)] {\it For the $2^{\rm d}$ equivariant REPRESENTATION $X^2 \overset{f}{\longrightarrow} \widetilde M(\Gamma)$ from theorem~{\rm 1.1}, there exists $\Gamma${\ibf -equiva\-riant zippings}, by which it is meant that there is a zipping {\rm (\ref{eq1.23})}, with $P_j = \underset{k}{\sum} \ P_{jk}$ s.t. for each $P_j$ there exists a neighbourhood $U_{p(j)}$ coming with $\pi^{-1} \, U_{p(j)} = \underset{k}{\sum} \ U_{pk}$, where $P_{jk}$ takes place inside $U_{pk}$, and where all the moves $\pi \, P_{jk} \subset U_{p(j)}$ are isomorphic; here $j$ is given (i.e. $p(j)$ too) and $k$ is arbitrary. We will write $P_j = \pi^{-1} \pi P_{jk}$.}
\item[2)] {\it Actually, for any atlas of $M(\Gamma)$ which locally trivializes $\pi$, we can find a zipping like in} 1).
\item[3)] {\it When we go to the diagram}
\begin{equation}
\label{eq1.24}
\xymatrix{
X^2 \ar[rr]_-{f} \ar[d]^-{\pi} &&\widetilde M (\Gamma) \ar[d]^-{\pi}  \\ 
\pi X^2 = X^2 / \Gamma \ar[rr]_-{f/\Gamma} &&M(\Gamma) = \widetilde M (\Gamma) / \Gamma \, ,
}
\end{equation}
{\it then the equivariant zipping upstairs, for $f$, induces a zipping downstairs for $f/\Gamma$.}
\end{enumerate}

\bigskip

The upper line in (\ref{eq1.24}) {\ibf is} a REPRESENTATION for $\Gamma$ but not necessarily the truncation $X_{i \geq1}^2 \to \widetilde M (\Gamma)$, since $X_i^2$ may not be GSC. As far as $X^2 / \Gamma$ is concerned, this is NOT GSC and so the lower line in (\ref{eq1.24}) is NOT a REPRESENTATION of anything.

\smallskip In the rest of this paper we will constantly work with the equivariant zippings for (\ref{eq1.1}).

\bigskip

\noindent FINAL COMMENTS. A) The $f \, {\rm Sing} (f) \subset fX^2$ is a discrete subset, without accumulations at finite distance. Also, it is disjoined from ${\rm Sing} (fX^2) \subset fX^2$ which does accumulate at finite distance, the $\{ p_{\infty\infty} (S)\}$. The lack of local finiteness of the $fX^2$, occurring exactly along $\{ p_{\infty\infty} \}$, will be handled in the next section. It is, indeed, a very serious issue.

\medskip

B) The (\ref{eq1.4}) in theorem~1.1 is made possible by the overflowing of the walls $W_{\infty}$. In the same vein, the set $M_2 (f) \cup {\rm LIM} \, M_2 (f)$ has its own holonomy, just like a foliation or a lamination, and this holonomy may be non-trivial.

\medskip

C) Generally speaking, ${\rm LIM} \, M_2 (f) \ne \emptyset$ and we have the obvious equivalence
$$
{\rm LIM} \, M_2 (f) \ne \emptyset \Longleftrightarrow M_2 (f) \subset X^2 \quad  \mbox{is NOT {\ibf closed}}.
$$
This is one of the basic difficulties we will have to fight with, in this paper.

\medskip

D) In the context of lemma~1.4, situation 2) concerns the local, ``man-made'' infinity of $H_k^{\lambda} (\gamma) \subset Y(\infty)$, while 3) concerns the much more mysterious INFINITY of $fX^2$ and/or $\widetilde M (\Gamma)$, which is the actual topic of this series of papers. We will manage to tame or handle this INFINITY by making use of equivariance; see the last section of the present paper.

\medskip

E) In everything which comes next, when the contrary is not explicitly stated, $p_{\infty\infty}$ may mean $p_{\infty\infty} ({\rm proper})$ or $p_{\infty\infty} (S)$. We will always have
\begin{equation}
\label{eq1.25}
p_{\infty\infty} = \lim_{n =\infty} (p_{n\infty} \ \mbox{or} \ p_{\infty n}), \ \mbox{where} \ p_{n\infty} \in M_2 (f) \cap {\rm LIM} \, M_2 (f) \subset X^2 \, .
\end{equation}

\newpage

\section{Where we go back to dimension three}\label{sec2}
\setcounter{equation}{0}

Remember that we have started with the 3-dimensional REPRESENTATION $Y(\infty) \overset{g(\infty)}{-\!\!\!-\!\!\!-\!\!\!\longrightarrow} \widetilde M (\Gamma) \approx \Gamma$ (quasi-isometry), where the $Y(\infty)$ was put together from bicollared handles of dimension three. This had allowed us to get local finiteness, equivariance and uniformly bounded zipping length, features which we will by all means try to stick to in what will follow from now on. But, in order to get a more transparent structure for the double points we moved to a 2-dimensional very dense skeleton of $Y(\infty)$, i.e. to the 2-dimensional representation $X^2 \overset{f}{\longrightarrow} \widetilde M (\Gamma)$, without losing our three desirable features just stated. What we have gained is a transparent double point structure
$$
X^2 \times X^2 \supset M^2 (f) \to M_2 (f) \subset X^2
$$
the rich complexity of which, schematically represented in the figures of the last section (see, in particular, 1.1 to 1.3), will be one of our tools.

\smallskip

On the 2-complex $X^2$ there are two, not everywhere well-defined flows, namely the {\ibf collapsing flow} and the {\ibf zipping flow}. Generally speaking, there will be transversal intersection of the two kind of trajectories, creating closed loops. Together with the non-metrizability barrier and with the Stallings barrier, this is one of the potential problems which the present paper will have to overcome. High dimensions will be needed here. To begin with, in order to take care of the Stallings barrier we will use transversally compact objects i.e. objects of the form
$$
\{ \mbox{$3^{\rm d}$, or intermediary $4^{\rm d}$, objects} \} \times B^{N({\rm high})} \, ,
$$
and it is the supplementary $B^N$ which will take care of the non-metrizability barrier. Our high dimensional objects will be GSC and will carry a free $\Gamma$-action. But this action will fail to be co-compact, so the technology from our old papers \cite{23}, \cite{24}, \cite{25} will here be necessary, before we can get to $\Gamma \in {\rm QSF}$.%

\smallskip

As already explained, we will never use the simple-minded $X^2$, which is not locally finite, but use instead its more sophisticated version $\{ X^2$ from (1.15)$\}$ where all the
$$
\{ p_{\infty\infty}\} = \{ p_{\infty\infty} \, (\mbox{proper}) , \ \mbox{fig. 1.1.(B)}\} + \{ p_{\infty\infty} (S), \ \mbox{fig. 1.5}\}
$$
are deleted, and compensating 2-handles are added, instead. But then, $fX^2$ is {\ibf not locally finite} at the $p_{\infty\infty}$'s. At this point, notice that, except for the change of colour RED $\Rightarrow$ BLACK, figure~1.5.(B) looks, deceptively, quite like 1.1.(B).

\medskip

\noindent (2.1) \quad But the two figures are quite different, although in both we see the standard $\underset{n=\infty}{\lim} \ p_{n\infty} = p_{\infty\infty}$. At the level of figure~1.1 $fX^2 \mid p_{\infty\infty}$ naturally embeddable in $R^3$, while in figure~1.5.(B), the $p_{n\infty}$, $p_{\infty n}$ are undrawable singularities, like in \cite{8}, \cite{18}, not embeddable in $R^3$, without destroying their train-track structure. The most serious lack of local finiteness for $fX^2$ is the accumulation of the undrawable singularities above at the $p_{\infty\infty} (S)$. A workable definition of $fX^2$ requires the deletion of the $p_{\infty\infty} (S)$'s and their replacement by compensating disks. End of (2.1).

\medskip

Two, both locally finite ``thickenings'' for $fX^2$ will  be introduced now, namely the $\Theta^3 (fX^2)$ and the $\Theta^3 (fX^2)'$, and when we will write $\Theta^3 (fX^2)^{(')}$, then this will mean one, or the other, or both. Before even going to the explicit definitions, some features of these objects will be reviewed already.

\smallskip

The $\Theta^3 (fX^2)^{(')}$ certainly is {\ibf locally finite}. There will be a natural free action $\Gamma \times \Theta^3 (fX^2)^{(')} \to \Theta^3 (fX^2)^{(')}$, and in the context of the lower line of (\ref{eq1.24}) we will define, downstairs, a $\Theta^3 (f/\Gamma (X^2/\Gamma))^{(')}$, coming with a map $\Theta^3 (f/\Gamma (X^2 /\Gamma))^{(')} \to M (\Gamma)$, and with the functorial properties
\setcounter{equation}{1}
\begin{equation}
\label{eq2.2}
\Theta^3 (f/\Gamma (X^2/\Gamma))^{(')} = \Theta^3 (fX^2)^{(')} / \Gamma \, ,
\end{equation}
and
$$
(\Theta^3 (f/\Gamma (X^2/\Gamma))^{(')})^{\sim} = \Theta^3 (fX^2)^{(')} \, ,
$$
accompanied by $\pi_1 \, \Theta^3 (f/\Gamma (X^2/\Gamma))^{(')} = \Gamma$. It should be stressed here that the $\Theta^3 (f/\Gamma (X^2/\Gamma))^{(')}$ can also be defined directly downstairs, strictly from ``first principles'', i.e. using the same local recipees as for $\Theta^3 (fX^2)^{(')}$ upstairs; then the formulae (\ref{eq2.2}) occur as a posteriori facts, consequences of this. An alternative would be to start with the first line in (\ref{eq2.2}) and take it as a definition. We have, from (\ref{eq1.1}), the simple-minded inclusion $fX^2 \subset \widetilde M (\Gamma)$ from which, after the appropriate treatment of the $p_{\infty\infty}$'s and the passage to regular neighbourhoods, we get things like
$$
\Theta^3 (fX^2)^{(')} \longrightarrow \widetilde M (\Gamma) \, ,
$$
or its reflex dowstairs, already mentioned.

\smallskip

Both $\Theta^3 (fX^2)$ and $\Theta^3 (fX^2)'$ are cell-complexes (and not 3-manifolds), the second one slightly more singular than the first one.

\smallskip

We will show now how to construct the $\Theta^3 (fX^2)$, and this will be done in several successive steps.

\bigskip

\noindent {\bf Step I.} Start with a decomposition $fX^2 = \underset{1}{\overset{\infty}{\bigcup}} \, U_i$ into small smooth, embedded 2-dimensional sheets. When $U_i$ enters the canonical open neighbourhood of some $\bar S \subset \widetilde M (\Gamma)$, touching the singularity $\bar S$, then it may always be assumed (see figure~1.5) that $U_i$ is contained inside one of the two $R_+^3 \cup (R_{\varepsilon} \times R_+ \times [0,1]_{\varepsilon})$'s. This allows us to define unambiguous thickening $U_i \times [-\varepsilon_i , \varepsilon_i] \subset \widetilde M (\Gamma)$ when $\varepsilon_1 > \varepsilon_2 > \ldots 0$ converges very fast to zero. Putting together these objects we get a first, very coarse thickening $\underset{1}{\overset{\infty}{\bigcup}} \, U_i \times [-\varepsilon_i , \varepsilon_i]$ for $fX^2$. With the $\sum_1 (\infty)$ which was defined in (\ref{eq1.10}), (\ref{eq1.14}) we go now to a first {\ibf provisional} definition
\begin{equation}
\label{eq2.3}
\sum (\infty) \equiv \left\{ \Sigma_1 (\infty) \cap \ \bigcup_1^{\infty} \ U_i \times [-\varepsilon_i , \varepsilon_i] \right\} \, ,
\end{equation}
coming with
$$
\partial \, \sum (\infty) \equiv \sum (\infty) \cap \partial \left( \bigcup_1^{\infty} \ U_i \times [-\varepsilon_i , \varepsilon_i] \right) \quad \mbox{and} \quad {\rm int} \, \sum (\infty) \equiv \sum (\infty) - \partial \, \sum (\infty) \, .
$$
Formulae (\ref{eq2.3}) are appropriate for the present Step~I, but they will need to get sharpened and amplified, by the time we will get to the final Step~III of our definition for $\Theta^3 (fX^2)$.

\smallskip

The ${\rm int} \, \sum (\infty)$ in (\ref{eq2.3}), which can be thought of as a $2^{\rm d}$ thickening of the ${\rm LIM} \, M_2 (f)$ from (\ref{eq2.3}), is an open surface with a ramification ($=$ non-manifold) locus, having the local form $\{$figure $Y\} \times R$. These ramifications occur, exactly, along the open arcs $p_{\infty\infty} ({\rm proper}) \times (-\varepsilon,\varepsilon)$.

\smallskip

We leave for the time being the $p_{\infty\infty} (S)$'s in the dark, and only in the context of Step~III will we finally make sense of the
\begin{equation}
\label{eq2.4}
\Theta^3 (fX^2) \mid \{ p_{\infty\infty} (S) \} \, ?
\end{equation}

At the points $p_{\infty\infty} ({\rm proper}) \times \{ \pm \, \varepsilon \}$, the $\partial \sum (\infty)$ is not locally finite and, with all these things, we introduce the following first and very provisional definition
\begin{equation}
\label{eq2.5}
\Theta^3 (fX^2)({\rm I}) \equiv \bigcup_1^{\infty} \ U_i \times [-\varepsilon_i , \varepsilon_i] - \partial \sum (\infty) \, .
\end{equation}
With this definition, at the level of $\Theta^3 (fX^2) ({\rm I})$, the lack of local finiteness coming with the $p_{\infty\infty} ({\rm proper})$, and which the figure~1.1.(B) might illustrate has been eliminated. Outside of the $p_{\infty\infty} (S)$'s, the $\Theta^3 (fX^2)({\rm I})$ is now a $3^{\rm d}$ manifold with boundary. At the level of (\ref{eq2.5}), we got local finiteness at the level of $p_{\infty\infty} ({\rm proper})$ by deleting just $p_{\infty\infty} \times \{ \pm \, \varepsilon \}$, with the $p_{\infty\infty} ({\rm proper}) \times (-\varepsilon , \varepsilon)$ left in place.

\smallskip

When it comes to $p_{\infty\infty} (S) \times [-\varepsilon , \varepsilon] \subset W \times [-\varepsilon , \varepsilon] \times \{ x = x_{\infty} \}$, then in the  plane $\{ x = x_{\infty} \} \cap W \times [-\varepsilon \leq z \leq \varepsilon]$ live infinitely many thinner and thinner immortal singularities of $\Theta^3 (fX^2)({\rm I})$ which accumulate on $p_{\infty\infty} (S) \times [-\varepsilon , \varepsilon]$ and deleting just $p_{\infty\infty} (S) \times \{ \pm \, \varepsilon \}$ will not be enough; see here (\ref{eq2.13}) below.

\bigskip

\noindent Notational Remark: Our present ${\rm int} \, \sum (\infty)$ used to be called ``$\sum (\infty)$'' in \cite{31}.

\medskip

\noindent (2.6) \quad Outside of the $p_{\infty\infty} (S) \times [-\varepsilon , \varepsilon]$, the $\Theta^3 (fX^2)({\rm I})$ is actually a smooth non-compact 3-manifold with non-compact boundary, inside which the $p_{\infty\infty}$ (pro\-per) $\times \, (-\varepsilon , \varepsilon)$ are PROPERLY embedded.

\smallskip

But, at the level of the figure~1.5.(C), the $p_{n\infty}$, $p_{\infty n}$ stand for $3^{\rm d}$ undrawable singularities (see \cite{8}, \cite{18}) for $\Theta^3 (fX^2)({\rm I})$ which accumulate at the $p_{\infty\infty} (S)$, where $\Theta^3 (fX^2)({\rm I})$ {\ibf fails to be locally finite}.

\bigskip

\noindent {\bf Step II.} We start now by codifying the structure of $fM_2 (f) \subset X^2$ as presented in the figures~1.1 to 1.3 and 1.5. Keeping the $p_{\infty\infty}$'s out of focus, here are our typical {\it local models for} $fM^2 (f)$. There is always a main sheet $V = \{ W \ \mbox{or} \ W_{(\infty)}\}$, cut transversally by $\sum_1 (\infty)$ along a smooth line $L_{\infty} = V \cap f \, {\rm LIM} \, M_2 (f)$, which carries, possibly, a unique $p_{\infty\infty} (V)$. Next, there are infinitely many double lines $L_j = W_j \cap V$, with $\underset{j=\infty}{\lim} \ L_j = L_{\infty}$, too. Finally, there are also lines $T_i = W^i \cap V$, which cut transversally the $(L_j , L_{\infty})$'s at points
\setcounter{equation}{6}
\begin{equation}
\label{eq2.7}
t_{ij} = T_i \cap L_j \subset W^i \cap W_j \cap V \, , \ \mbox{accumulating on $L_{\infty}$ when} \ j \to \infty \, .
\end{equation}
This $t_{ij}$ may be a ramification point of $X^2$ or a triple point, hence a ramification point somewhere higher on the zipping road from $X^2$ to $fX^2$. In our figures~1.1 to 1.3 or 1.5 we see that
$$
\lim_{j = \infty} \, t_{ij} = p_{n\infty} \in fM^2 (f) \cap f \, {\rm LIM} \, M_2 (f) = \{\mbox{the set of $p_{n\infty} , p_{\infty n}$'s}\} \, .
$$
We will distinguish three kinds of local models, and only the first two contain a $p_{\infty\infty}$.

\medskip

\noindent (2.8.i) \quad The generic case $V = W_{(\infty)} ({\rm BLACK}) \ni p_{\infty\infty}$, $W^i = W^i ({\rm RED})$, $W_j = W_j ({\rm BLUE})$; the figures are here 1.1, 1.2, 1.3. Also, in the next two cases, we will again have $W_j = W_j ({\rm BLUE})$. It should be understood here that, in terms of figures like 1.1, 1.2, 1.3, for any individual line $T \cap W_{\infty}$ ($=$ red line), all the triple points except for finitely many of them live inside the difficult region $T \cap \overline N_{\infty}^2$; this will be introduced formally, later on.

\medskip

\noindent (2.8.ii) \quad The immortal case, when $V = W({\rm BLACK}) \ni p_{\infty\infty}$ and $W^i =\{$transversal $W^* ({\rm BLACK})\}$, see here figure~1.5. This model is imposed by the existence of the $\bar S \subset {\rm Sing} \, \widetilde M (\Gamma) \subset \widetilde M (\Gamma)$.

\medskip

\noindent (2.8.iii) \quad The RED case, which had no $p_{\infty\infty}$, $V = W ({\rm RED} \cap H^0) \supset \{$some smooth line ${\rm LIM} \, M_2 (f)\}$, $W^i = W_{(\infty)} ({\rm BLACK})$. Here $V = \{$cyclinder $S^1 \times [0,x_{\infty} + \zeta]\}$, $\zeta > 0$, with $S^1 \times \{ 0 \} \subset W({\rm BLUE})$ and $S^1 \times \{ x_{\infty}\} = V \cap \sum_1 (\infty)$. Figures~2.17 to 2.19 in \cite{31} may illustrate this case. At this point I will introduce now a localized system of {\ibf dilatations}, which will change $\Theta^3 (fX^2)({\rm I})$ into a $\Theta^3 (fX^2)({\rm II})$. These dilatations which will always grow out of the $W_{(\infty}) ({\rm BLACK})$'s are called {\ibf fins}. The fins are not introduced neither for purposes of $3^{\rm d}$ local finiteness nor of smoothness but for the needs of the geometric realization of the zipping, in high dimensions. There, in high $d$, they will help keeping our variously needed objects locally finite. Very loosely speaking, during the geometric realization of the zipping, in high $d$, various undesirable accumulations will be corralled at infinity, along $\partial \sum (\infty)$. The fins will come with {\ibf rims}, also living at infinity (see formula (\ref{eq2.15}) below) and these, added to $\partial \sum (\infty)$ will take care of the gaps, which things like the triple points create for $\partial \sum (\infty)$, making thus the corralling above possible. This, admittedly vague description will be make precise, later on.

\smallskip

Here is the GEOMETRY OF THE FINS. We will need fins for handling all the three local models from (2.8), and I will start now with the paradigmatic case i). To make our discussion smooth, we will use now $\sum (\infty) \subset \underset{1}{\overset{\infty}{\bigcup}} \ U_i \times [-\varepsilon_i , \varepsilon_i]$. In the context of (2.8.i) we have a local coordinate system $(x,y,z)$ with $V = (z=0)$, $W_j = (x = x_j)$, $W^i = (y=y_i)$, to which we will add a fourth dimension, with the coordinate $u$. We have also a square, at the level of which the next contribution will be localized
$$
\sum (\infty) \mid T_i = (x=x_{\infty}) \cap \{ (- \varepsilon \leq z \leq \varepsilon) \cap (y_i - \varepsilon_i \leq y \leq y_i + \varepsilon_i)\} \, .
$$
Let

\medskip

\noindent (2.9) \quad $\frac12 \, D_{\pm}^2 (p_{i\infty}) \equiv \{$the half-disk of diameter ($x=x_{\infty}$, $y_i - \varepsilon_i \leq y \leq y_i + \varepsilon_i$, $z = \pm \, \varepsilon$), contained in the plane $(y,u)$. Here $\frac12 \, D^2$ in a dilatation added to $V \times [-\varepsilon \leq z \leq \varepsilon]$, at $z = \pm \, \varepsilon\}$, and $u$ is an additional dimension. We define
$$
{\rm rim} \, \left( {\rm of} \ \frac12 \, D^2 \right) \equiv \partial \, \frac12 \, D^2 - \{\mbox{the diameter}\} \, .
$$
With this we will want the addition in (2.9) to be part of a larger $3^{\rm d}$ dilatation based on $V \times \{ \pm \, \varepsilon \}$, denoted by $F_{\pm}$ (with ``$F$'' like ``fin''), and which does not touch otherwise to $W^i \times [-\varepsilon_i , \varepsilon_i] - V \times [-\varepsilon , \varepsilon]$. In the description of this dilatation, we have to make use of a fourth dimension $u$, and figure~2.1 should suggest what one can talking about.
$$
\includegraphics[width=13cm]{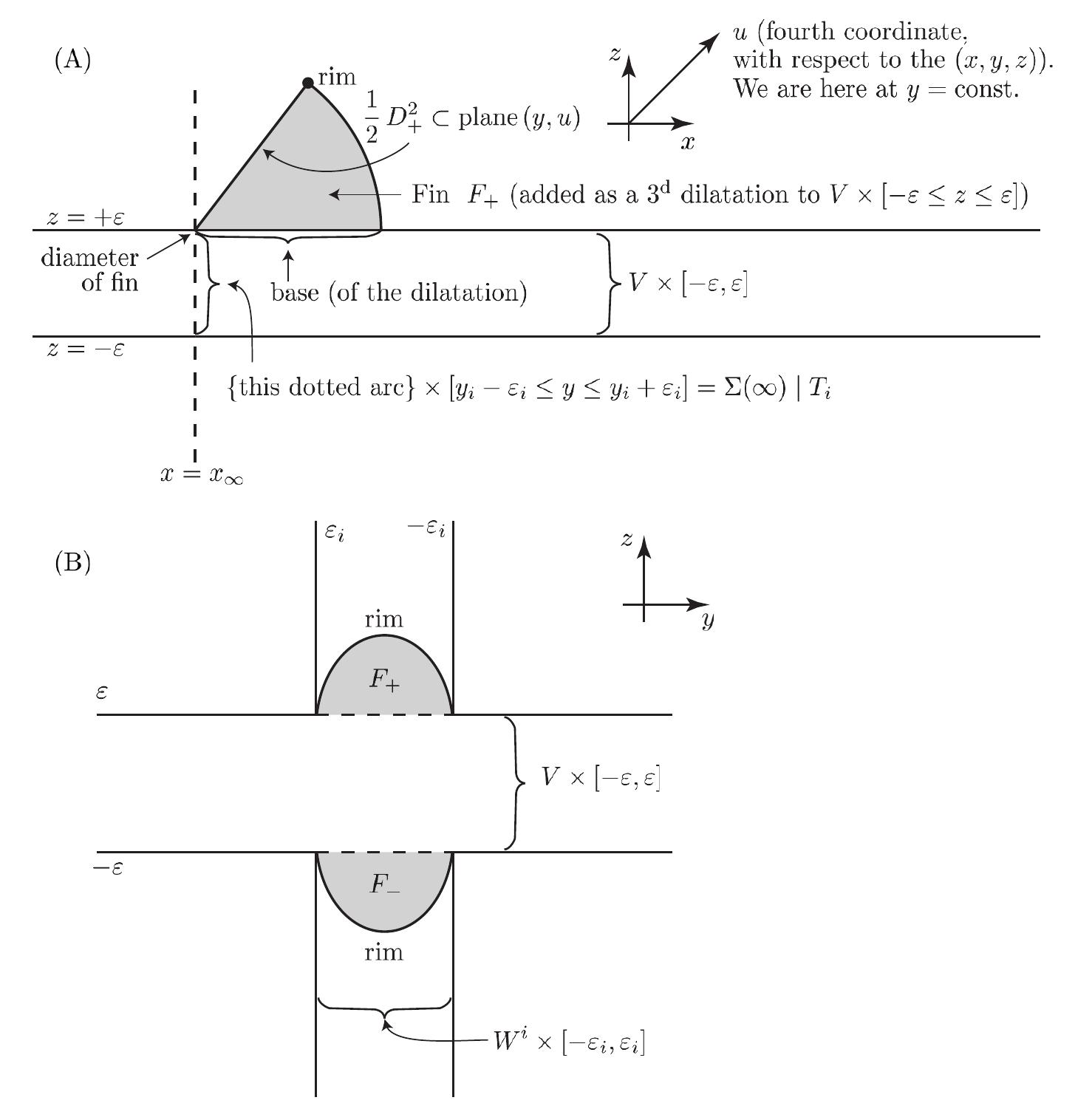}
$$
\label{fig2.1}
\centerline {\bf Figure 2.1.} 

\smallskip

\begin{quote}
The fin $F_+$ in (A) lives inside ($x,y,z = +\varepsilon$, $u \geq 0$) and outside of its base, which lives at ($z = +\varepsilon$, $y \in [y_i - \varepsilon_i , y_i + \varepsilon_i]$, $x \geq x_{\infty}$, $u=0$), it avoids $W^i \times (y_i - \varepsilon_i \leq y \leq y_i + \varepsilon_i) \subset (x,y,z,u=0)$. In (A) we are at some fixed value $y \in [y_i  - \varepsilon_i , y + \varepsilon_i]$, and there is also another fin $F_-$, added at $z=-\varepsilon$. Figure (B) presents a projection of $F_{\pm}$ on $(y,z)$. Both the boundary vertical or horizontal ones, and the rims, are living at infinity.
\end{quote}

\bigskip

The (2.9) tells we which fins to add in the situation (2.8.i) and we move now to (2.8.ii). Here, at the level (\ref{eq1.1}), for each $\bar S \subset {\rm Sing} \, \widetilde M (\Gamma)$ we find the situation from figure 1.5.(C), i.e. a doubly infinite checkerboard of walls $W({\rm BLACK})$, $(W(n) , W^* (m))$ for $(n,m) \in Z \times Z$, and we will use the notation $W^{**} = W$. At each site $(n,m)$ we will make a {\ibf choice}. A label {\ibf ``overflows''} or {\ibf ``subdued''} will be attached to each of the two terms $W,W^*$, for each fixed $(n,m)$ so that the following RULES should be satisfied for $(W (n) , W^* (m))$.

\medskip

\noindent (2.10) \quad We certainly want that $W$ overflows $\Longleftrightarrow$ $W^*$ is subdued; at $(n,m)$.

\medskip

\noindent (2.10.1) \quad  For each $W(n)$, when $n$ is fixed, there should be a finite central zone where for all the indices $m\in Z$ we should have
$$
\{\mbox{At $(n,m)$ $W(n)$ is subdued, hence $W^* (m)$ is overflowing}\} \, ,
$$
leaving us with two infinite ends of $W(n)$, each of these connected, where $W(n)$ overflows.

\smallskip

The same kind of thing should be valid when $m$ is fixed. END of (2.10.1).

\medskip

Each $p_{\infty\infty} (S)$ is an endpoint of some $\lim W \cap (x = x_{\infty})$, with $W$ like above. We will denote by $W_+ (p_{\infty\infty} (S)) \subset W = \{ W(m)$ OR $W(n)\}$, the half line corresponding to the end in question and which is such that at each site $S \in W_+ (p_{\infty\infty} (S)) \subset W$, it is exactly $W$ which overflows (and hence $W^*$ is subdued). Here comes now our second RULE, complementing (2.10.1).

\medskip

\noindent (2.10.2) \quad The $\sum W_+ (p_{\infty\infty} (S))$ induces a disjoined partition of $\{$the set of all the sites $(n,m) \approx S(n,m)\}$.

\medskip

The figure 5.2 in section V below, displays one possible recipee for implementing the rules (2.10.1) $+$ (2.10.2), via the diagonal $0(3)$-lines; any half line going from the diagonal $0(3)$-lines to infinity, is a $W_+$. On the unique central vertical $W$ which does not touch the $0(3)$-lines (see the figure~5.2), we have made the arbitrary choice that $W$ overflows everywhere. But the precise way in which our rules (2.10.1) $+$ (2.10.2) are implemented, is immaterial, our figure~5.2 is rather a mere illustration. On the other hand, once our choice of overflow versus subdued, has been made at each site $S(m,n)$, we can talk about how to place the fins too.

\medskip

\noindent (2.11) \quad For each $W_+ (p_{\infty\infty} (S))$ in (2.10.2), at each site $(n,m) \in W_+ (p_{\infty\infty} (S)) \subset \sum (\infty)$, the site is standing now for a triple point of the form $fM_2 (f) \cap S_{\infty}^2 ({\rm BLUE})$, {\ibf two fins $F_{\pm}$ will be attached} just like in (2.9), {\ibf to the $W_+$, which overflows} at $(n,m)$. Figure~5.2 illustrates this rule too. All this takes care of~(2.8.ii). The same kind of procedure as for (2.8.i) can be implemented for the case (2.8.iii) too; see here the figure~5.4.

\medskip

With all this, we define now the next
\setcounter{equation}{11}
\begin{equation}
\label{eq2.12}
\Theta^3 (fX^2) ({\rm II}) \equiv \Bigl\{ \Theta^3 (fX^2)({\rm I}) \cup \sum_{\rm all \, fins} F_{\pm} \ \mbox{with both the}
\end{equation}
$$
\partial \, \sum (\infty) \ \mbox{{\ibf and} the rims of fins, which rest on} \ \partial \, \sum (\infty), \, \mbox{deleted} \Bigl\} \, .
$$
Clearly, we have a dilatation $\Theta^3 (fX^2)({\rm I}) \nearrow \Theta^3 (fX^2)({\rm II})$.

\bigskip

\noindent {\bf Step III.} So far, the (\ref{eq2.12}) is still not locally finite at the $p_{\infty\infty} (S)$'s. To take care of this, we introduce the following cell-complex, which {\ibf is} locally-finite, and which is also the final $\Theta^3 (fX^2)$, now correctly defined

\newpage

\begin{equation}
\label{eq2.13}
\Theta^3 (fX^2) \equiv \Theta^3 (fX^2)({\rm II}) \ \mbox{(from (2.12))} \ - \sum_{p_{\infty\infty} (S)} p_{\infty\infty} (S) \times [-\varepsilon , \varepsilon] 
\end{equation}
$$
+ \sum_{p_{\infty\infty} (S)} D^2 (H (p_{\infty\infty} (S))) \times \left[ -\frac{\varepsilon}4 , \frac{\varepsilon}4 \right] \, ,
$$
with the last batch of $2$-handles added along the $\sum \, C(p_{\infty\infty}) \times \left[ -\frac{\varepsilon}4 , \frac{\varepsilon}4 \right]$. Notice that the $p_{\infty\infty} (S) \times \{ \pm \, \varepsilon \}$ was already deleted at the level of $\Theta^3 (fX^2)$ (I and II), now we delete the rest of $p_{\infty\infty} (S) \times [-\varepsilon,\varepsilon]$ too. The notation ``$H(p_{\infty\infty} (S))$'', going with the synonimous ``$- p_{\infty\infty} (S) \times (-\varepsilon , \varepsilon)$'' means that the deleted arc in question is a Hole, on par with the other Holes to be deleted in the section~IV. But, contrary to the normal Holes (\ref{eq4.6}), which are open sets, the $H(p_{\infty\infty} (S)) \times (-\varepsilon , \varepsilon)$ (as the larger $\{ H (p_{\infty\infty} ({\rm all})) \times (-\varepsilon , \varepsilon)\}$) is closed.

\smallskip

Together with (\ref{eq2.13}) comes the following definition, which supersedes (\ref{eq2.3})

\medskip

\noindent (2.13.1) \quad $\sum (\infty)$ (now correctly defined) $\equiv \, \{$the $\sum (\infty)$ from (\ref{eq2.3})$\} \, - \underset{p_{\infty\infty}(S)}{\sum}$ $p_{\infty\infty}(S) \times [-\varepsilon , \varepsilon]$, coming with ${\rm int} \, \sum (\infty) \equiv \{$the $\sum (\infty)$ as just defined above$\} - \partial \, \sum (\infty)$ (from (\ref{eq2.3})). This definition will be accompanied by the formula
$$
\partial \sum (\infty) \, (\mbox{now correctly defined}) 
\equiv \left\{ \partial \sum (\infty) \, (\mbox{from (\ref{eq2.3})}) \cup \sum_{p_{\infty\infty} (S)} p_{\infty\infty} (S) \times [-\varepsilon , \varepsilon] \right\} \, .
$$
END of (2.13.1).
$$
\includegraphics[width=12cm]{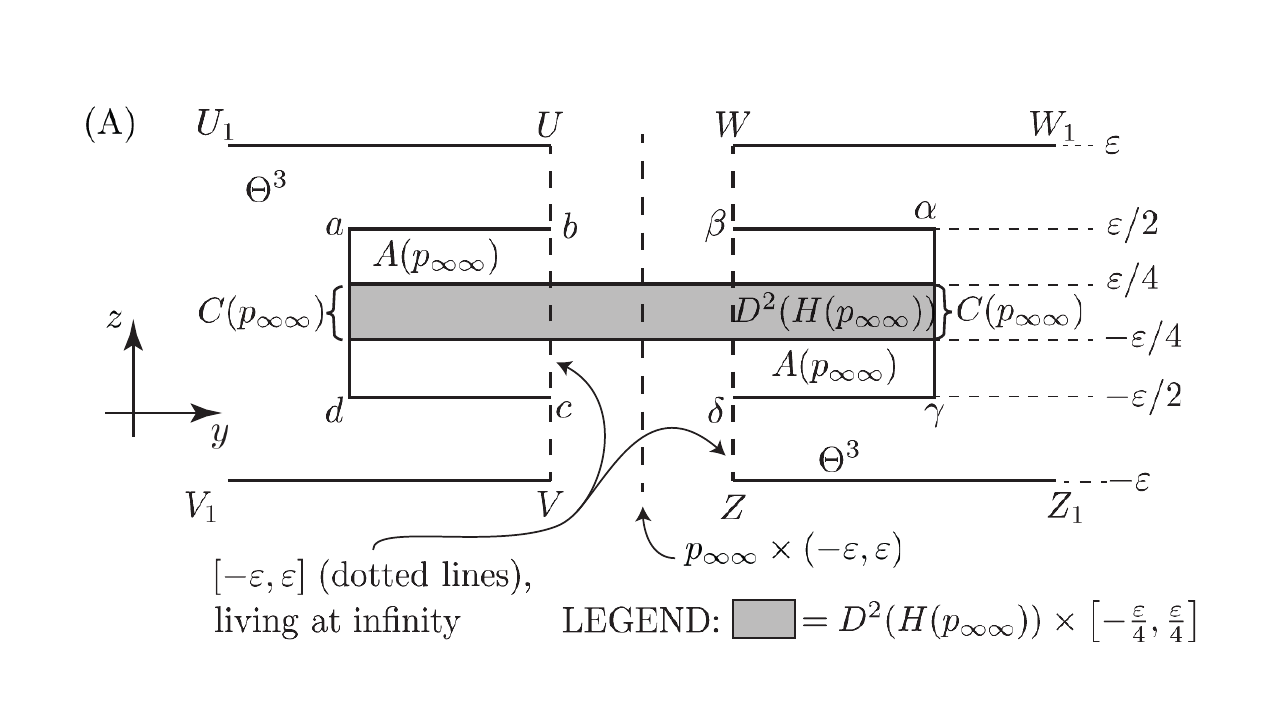}
$$
$$
\includegraphics[width=14cm]{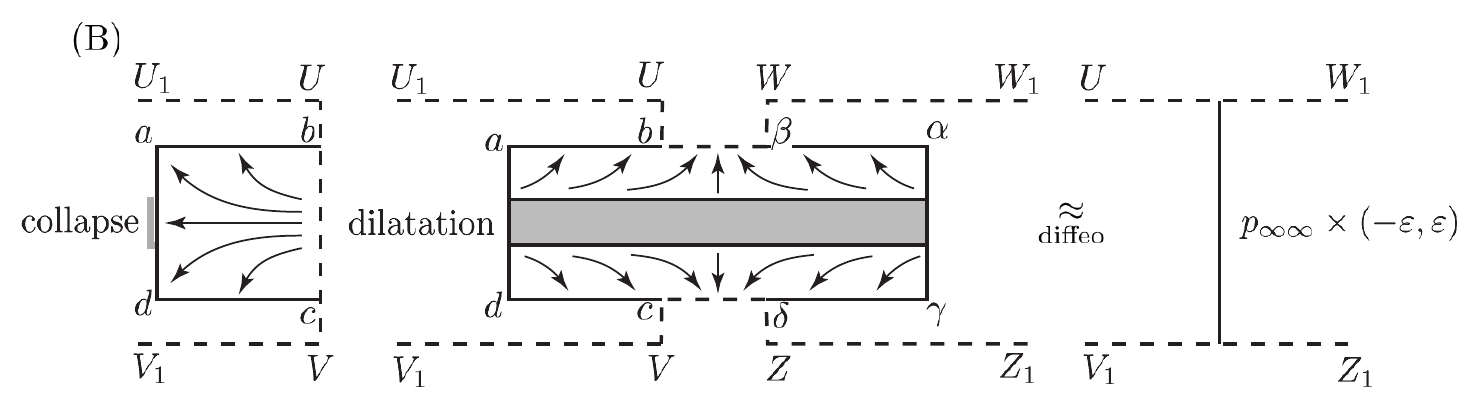}
$$
\label{fig2.2}

\centerline {\bf Figure 2.2.} 

\smallskip

\begin{quote}
Figure (A) is a section $x = {\rm const}$, through the space $\Theta^3 = \Theta^3 (fX^2)$, which is singular here exactly along $C(p_{\infty\infty}) \times \left[-\frac{\varepsilon}4 , \frac{\varepsilon}4 \right]$ where the $D^2 (H (p_{\infty\infty}))$ $\times \left[-\frac{\varepsilon}4 , \frac{\varepsilon}4 \right]$ is glued to the rest. The $[a,b,c,d] + [\alpha , \beta , \gamma , \delta]$ is bounding the rotationally invariant region $A(p_{\infty\infty}) \subset \Theta^3 - p_{\infty\infty} \times [-\varepsilon , \varepsilon]$. If our section $x = {\rm const}$ is actually $x = x_{\infty}$, then not only the $[UV] + [WZ]$ are automatically at infinity, but as well the whole $[U_1 U] + [V_1 V] + [W W_1] + [Z Z_1]$. When our present $p_{\infty\infty}$ is actually of type $p_{\infty\infty} (S)$, then the immortal singularities certainly touch $\Theta^3$, cutting through $A(p_{\infty\infty} (S))$, but they stay away from the 2-handle $D^2 (p_{\infty\infty} (S)) \times \left[-\frac{\varepsilon}4 , \frac{\varepsilon}4 \right]$.

\smallskip

The collapsing part in formula (\ref{eq2.23}) below corresponds, at the present $3^{\rm d}$ level to what one sees in the drawing (B), while the dilatation part corresponds to (C). Both (B) and (C) live at $x=x_{\infty}$ and they concern {\ibf exclusively} the $p_{\infty\infty} \, ({\rm proper})$'s. Notice how, at the level of (C), the exceptional Hole $H (p_{\infty\infty} \, ({\rm proper})) = p_{\infty\infty} \, ({\rm proper}) \times [-\varepsilon , \varepsilon]$ is finally healed, i.e. filled, in a {\ibf smooth} manner. The holes $H (p_{\infty\infty} (S))$ {\ibf cannot} be healed this way, since the infinitely many immortal $S$'s cutting through $A(p_{\infty\infty})$, accumulates now on $p_{\infty\infty} (S) \times [-\varepsilon , \varepsilon]$. In all this discussion, one should think of $\Theta^3 (fX^2)'$ as being the primary object, while
$$
\Theta^3 (fX^2) \equiv \{\mbox{the} \ \Theta^3 (fX^2)' \ \mbox{with the hole} \ H (p_{\infty\infty} \, ({\rm proper})) \ \mbox{healed smoothly}\}.
$$
\end{quote}

\bigskip

Starting from the final $\sum (\infty)$, as just defined in (\ref{eq2.13}), we also introduce the bigger
\begin{equation}
\label{eq2.14}
\sum (\infty)^{\wedge} \equiv \sum (\infty) \cup \sum_{{\rm Fin} \, F} \frac12 \, D^2 (F) \, ,
\end{equation}
the union being made along the diameters, and this comes with the following set, living at infinity,
\begin{equation}
\label{eq2.15}
\partial \sum (\infty)^{\wedge} = \left\{ \partial \sum (\infty) \ \mbox{(from (2.13.1))} \cup \sum_F {\rm rim} \, (F) \right\} \supset \{\mbox{all the} \ p_{\infty\infty} \times \{ \pm \, \varepsilon \} \} \, .
\end{equation}
This also comes with the following PROPER inclusions
$$
\partial \sum (\infty)^{\wedge} - \{{\rm the} \ p_{\infty\infty} (S) \times [-\varepsilon , \varepsilon] \} \subset \sum (\infty)^{\wedge} \, , \eqno (2.16.1)
$$
and
$$
{\rm int} \sum (\infty)^{\wedge} \equiv \sum (\infty)^{\wedge} - \partial \sum (\infty)^{\wedge} \subset \Theta^3 (fX^2) \, . \eqno (2.16.2)
$$
With $p_{\infty\infty} = p_{\infty\infty} (S)$, the figure 2.2.(A) (which was a priori drawn for $p_{\infty\infty}$ $({\rm proper})$), should serve as an illustration for the formula (\ref{eq2.13}).

\smallskip

On par with the $\Theta^3 (fX^2)$, we will introduce now another locally finite cell-complex, namely the following, and for it the figure~2.2 can also serve as an illustration
\setcounter{equation}{16}
\begin{equation}
\label{eq2.17}
\Theta^3 (fX^2)' \equiv \Theta^3 (fX^2)({\rm II}) - \sum_{p_{\infty\infty} (\bm{all})} p_{\infty\infty} \times (-\varepsilon , \varepsilon)  
+ \, \sum_{p_{\infty\infty} (\bm{all})} D^2 (H(p_{\infty\infty})) \times \left[-\frac{\varepsilon}4 , \frac{\varepsilon}4 \right] \, ,
\end{equation}
 with the $2$-handle $D^2 \times \left[-\frac{\varepsilon}4 , \frac{\varepsilon}4 \right]$ always added along the $C(p_{\infty\infty}) \times \left[-\frac{\varepsilon}4 , \frac{\varepsilon}4 \right]$.

\smallskip

The figures~2.2.(A,B,C) can serve as an illustration for the (\ref{eq2.17}) too. The (\ref{eq2.17}) is also accompanied by (\ref{eq2.15}). We will not bother now to develop the analogues of (2.13.1), (\ref{eq2.14}) for the $\Theta^3 (fX^2)'$ (\ref{eq2.17}), and we will not need them either. But this whole kind of issue will be pricked up again with more detail, in the section~VI. Clearly both (\ref{eq2.13}) and (\ref{eq2.17}) make use of the same idea as in (1.15), and for good measure we also introduce the following thickened version of $X^2$ (1.15), in $3^{\rm d}$, namely

\newpage

$$
\Theta^3 (X^2) \equiv \Bigl\{\mbox{The smooth regular neighbourhood of thickeness} \ [-\varepsilon , \varepsilon], \eqno (2.17.1)
$$
$$
\Theta^3 \Bigl(\{ X^2 \ \mbox{simple-minded}\} \, - \sum_{p_{\infty\infty} (\bm{all})} p_{\infty\infty} \Bigl)\Bigl\} \, + \sum_{p_{\infty\infty} (\bm{all})} D^2 (H(p_{\infty\infty})) \times \left[-\frac{\varepsilon}4 , \frac{\varepsilon}4 \right] \, .
$$
For $fX^2$ we have introduced distinct $\Theta^3 (fX^2)$ and $\Theta^3 (fX^2)'$, a distinction which will not be made, neither for $X^2$ (1.15) nor for $\Theta^3 (X^2)$ (2.17.1), both objects being rather like $\Theta^3 (fX^2)'$, anyway, from the beginning.

\smallskip

Here are some additional explanations concerning the figure~2.2, which is the background behind all the three formulae (\ref{eq2.13}), (\ref{eq2.17}), (2.17.1). The action is concentrated in a region
$$
p_{\infty\infty} \times (-\varepsilon , \varepsilon) \underset{\rm PROPER \ embedding}{-\!\!\!-\!\!\!-\!\!\!-\!\!\!-\!\!\!-\!\!\!-\!\!\!-\!\!\!-\!\!\!-\!\!\!-\!\!\!-\!\!\!-\!\!\!-\!\!\!-\!\!\!-\!\!\!\longrightarrow} W_{\infty} \times [-\varepsilon \leq z \leq \varepsilon] \subset \Theta^3 (fX^2)^{(')} 
$$
((\ref{eq2.13}) or (\ref{eq2.17})). What we actually see in figure~2.2.(A) is a small region, PROPERLY embedded inside $\Theta^3 (fX^2)^{(')}$, and spanning from $z=-\frac\varepsilon2$ to $z=\frac\varepsilon2$
\begin{equation}
\label{eq2.18}
A(p_{\infty\infty}) \underset{\overbrace{\mbox{\footnotesize$C(p_{\infty\infty}) \times \left[-\frac{\varepsilon}4 , \frac{\varepsilon}4 \right]$}}}{\cup} \left[ D^2 (H (p_{\infty\infty})) \times \left[-\frac{\varepsilon}4 , \frac{\varepsilon}4 \right] \right] \longrightarrow \Theta^3 (fX^2)^{(')} \, .
\end{equation}
The $A(p_{\infty\infty})$ which is rotationally symmetric (the axis of symmetry being the thickly dotted $p_{\infty\infty} \times [-\varepsilon , \varepsilon]$, which is certainly {\ibf NOT} part of $\Theta^3 (fX^2)^{(')}$), is a copy of $S^1 \times \left[-\frac{\varepsilon}2 , \frac{\varepsilon}2 \right] \times [0,\infty)$. The 2-handle $D^2 (H(p_{\infty\infty})) \times \left[-\frac{\varepsilon}4 , \frac{\varepsilon}4 \right]$ is glued to the rest of $\Theta^3 (fX^2)^{(')}$ along the $C(p_{\infty\infty}) \times \left[-\frac{\varepsilon}4 , \frac{\varepsilon}4 \right] \underset{\rm TOP}{=} S^1 \times \left[-\frac{\varepsilon}4 , \frac{\varepsilon}4 \right]$. We have

\medskip

\noindent (2.18.1) \quad ${\rm Sing} \, \Theta^3 (fX^2)^{(')} =\{$the singularities $p_{n\infty}$, $p_{\infty m}$ from figure~1.5 (connected to $\bar S \subset {\rm Sing} \, \widetilde M (\Gamma)$) plus, disjoined from them, $\sum C (p_{\infty\infty} ({\rm all})) \times \left[-\frac{\varepsilon}4 , \frac{\varepsilon}4 \right]$ for the case $\Theta^3 (fX^2)'$, respectively, $\sum C (p_{\infty\infty} (S)) \times \left[-\frac{\varepsilon}4 , \frac{\varepsilon}4 \right]$, for $\Theta^3 (fX^2)$ (\ref{eq2.13})$\}$.

\medskip

Desingularizations ${\mathcal R}$ \`a la \cite{8}, \cite{18}, \cite{19}, can be defined for all of the immortal $S = p_{n\infty}$ or $p_{\infty n}$, AND also for the $C(p_{\infty\infty}) \times \left[-\frac{\varepsilon}4 , \frac{\varepsilon}4 \right]$ occurring in the formula (2.18.1). With this, we may introduce the $4^{\rm d}$ cell complexes $\Theta^4 (\Theta^3 (fX^2)^{(')} , {\mathcal R})$, like in the definition below.

\bigskip

\noindent {\bf Definition (2.19).} We define, to begin with, the cell-complex
$$
\Theta^4 (\Theta^3 (fX^2) , {\mathcal R}) \equiv \Theta^4 ((\Theta^3 (fX^2) , \, \mbox{with the contribution} \ p_{\infty\infty} (S) \, \mbox{deleted}), \  {\mathcal R}) 
$$
$$
\cup \ \sum_{p_{\infty\infty} (S)} \left( D^2 (H (p_{\infty\infty})) \times \left[-\frac{\varepsilon}4 , \frac{\varepsilon}4 \right] \times I \right) \, ,
$$
a formula which, instead of being a simple-minded 4-dimensional thickening of (\ref{eq2.13}), is rather a $4^{\rm d}$ analogue of it. Its RHS requires some explanations. One uses here the $4^{\rm d}$ thickening functor from \cite{8}, \cite{18}, \cite{19}, $\Theta^4 (\ldots , {\mathcal R})$. This has good localization and glueing properties and, for any open set $U \subset \Theta^3 -\{{\rm singularities}\}$, $\Theta^4 (U,{\mathcal R}) = U \times I$. Without any loss of generality, for {\ibf any} $p_{\infty\infty}$ we have (see here also figure~2.2) that $\left( C(p_{\infty\infty}) \times  \left[-\frac{\varepsilon}4 , \frac{\varepsilon}4 \right] \right) \cap \{$the undrawable singularities $S = p_{n\infty} , p_{\infty n}$ of $\Theta^3 (fX^2)$, generated by the figure~1.5$\} = \emptyset$, a fact which gives us an embedding
$$
\sum_{p_{\infty\infty} (S)} C(p_{\infty\infty}) \times \left[-\frac{\varepsilon}4 , \frac{\varepsilon}4 \right] \times I \subset \Theta^4 ((\Theta^3 (fX^2) \
\mbox{with the contribution} \ p_{\infty\infty} (S) \ {\rm deleted}), {\mathcal R}) \, .
$$
We clearly have also another embedding
$$
\sum_{p_{\infty\infty} (S)} C(p_{\infty\infty} (S)) \times \left[-\frac{\varepsilon}4 , \frac{\varepsilon}4 \right] \times I \subset \sum_{p_{\infty\infty} (S)} D^2 (H (p_{\infty\infty} (S))) \times \left[-\frac{\varepsilon}4 , \frac{\varepsilon}4 \right] \times I \, .
$$

One can glue now the two terms in the RHS of the formula for $\Theta^4 (\Theta^3 (fX^2),{\mathcal R})$ along the two embeddings just listed. This completes the definition of $\Theta^4 (\Theta^3$ $(fX^2),{\mathcal R})$ and there is a completely similar definition for $\Theta^4 (\Theta^3 (fX^2)',{\mathcal R})$; this is gotten by proceeding exactly like for $\Theta^4 (\Theta^3 (fX^2),{\mathcal R})$, with the change $\underset{p_{\infty\infty} (S)}{\sum} \Longrightarrow \underset{p_{\infty\infty} ({\rm all})}{\sum}$. This ends definition (2.19).

\bigskip

\noindent {\bf Definition (2.19.bis).} We will define now the $(N+4)$-dimensional cell-complexes $S_u^{(')} \widetilde M (\Gamma)$, where $N$ is very large. We start by introducing the following pieces of $S_u \, \widetilde M (\Gamma)$:
$$
S_u \, \widetilde M (\Gamma) \mid \Biggl(\Theta^3 (fX^2) \Bigl(- \sum_{p_{\infty\infty} (S)} p_{\infty\infty} (S) \times [-\varepsilon , \varepsilon) \Bigl)\Biggl) \equiv \Biggl\{\mbox{the {\ibf smooth} manifold}
$$
$$
\Biggl(\Theta^4 \Biggl(\Theta^3 (fX^2) \Bigl( -\sum_{p_{\infty\infty} (S)} p_{\infty\infty} \times (-\varepsilon , \varepsilon) , {\mathcal R} \Bigl) \Biggl) \times B^n \Biggl)\Biggl\} \ 
$$
and then also the next
$
\mbox{piece} \ S_u \, \widetilde M (\Gamma) \mid D^2 (H (p_{\infty\infty} (S))) \equiv D^2 (H (p_{\infty\infty} (S))) \times \left[-\frac{\varepsilon}4 , \frac{\varepsilon}4 \right] \times \frac12 \, B^{N+1}, \mbox{where}
$
$
B^{N+1} = B^N \times \{\mbox{the canonical factor $I$ which occurs in}
$
$
\Theta^4 (U_{\rm smooth} \, {\mathcal R}) \equiv U \times I \subset \Theta^4 (\ldots , {\mathcal R}), \ 
$ see the definition (2.19)$\}$.

\bigskip

With all these things, comes an embedding
$$
\sum_{p_{\infty\infty} (S)} C(p_{\infty\infty} (S)) \times \left[-\frac{\varepsilon}4 , \frac{\varepsilon}4 \right] \times B^{N+1}
\subset S_u \, \widetilde M (\Gamma) \mid \Biggl( \Theta^3 (fX^2) \Bigl( -\sum_{p_{\infty\infty} (S)} p_{\infty\infty} (S) \times (-\varepsilon, \varepsilon) \Bigl) \Biggl) \, ,
$$
and the two spare parts for $S_u \, \widetilde M (\Gamma)$ which we have introduced above are then to be glued along the $\underset{p_{\infty\infty} (S)}{\sum} C(p_{\infty\infty} (S)) \times \left[-\frac{\varepsilon}4 , \frac{\varepsilon}4 \right] \times \frac12 \, B^{N+1}$, which lives inside both. This ends the definition of $S_u \, \widetilde M (\Gamma)$.

\bigskip

There is, again, an analogous definition for $S'_u \, \widetilde M(\Gamma)$, with the same change $\underset{p_{\infty\infty} (S)}{\sum} \Longrightarrow \underset{p_{\infty\infty} ({\rm all})}{\sum}$ as at the end of (2.19). This ends the definition (2.19.bis). 

\hfill $\Box$

\bigskip

\noindent COMMENTS. Both objects which we have just defined $S_u$ and $S'_u$ are singular, they are only cell-complexes and NOT smooth manifolds. The singular sets are $\underset{p_{\infty\infty} (S)}{\sum} C(p_{\infty\infty} (S)) \times \left[-\frac{\varepsilon}4 , \frac{\varepsilon}4 \right] \times \frac12 \, B^{N+1}$ for $S_u$, respectively the larger $\underset{p_{\infty\infty} ({\rm all})}{\sum} C(p_{\infty\infty})\times \left[-\frac{\varepsilon}4 , \frac{\varepsilon}4 \right] \times \frac12 \, B^{N+1}$, for $S'_u$. The $S_u \, \widetilde M(\Gamma)$ is the best thing we can offer as ``high dimensional regular neighbourhood of $fX^2$''. Now, the $\Theta^4 (\Theta^3 (fX^2)^{(')} , {\mathcal R})$ from the definition~(2.19) is certainly ${\mathcal R}$-dependent. But once we start taking a  cartesian product with any $B^{p \geq 1}$, then the ${\mathcal R}$-dependence gets washed away. This makes that both $S_u^{(')} \widetilde M (\Gamma)$'s are {\ibf canonical}; they have as good functorial properties as $\Theta^3$ in the context of the (2.2) above; we actually have the following

\bigskip

\noindent {\bf Lemma~2.1.} 1) {\it There is a natural free action}
\setcounter{equation}{19}
\begin{equation}
\label{eq2.20}
\Gamma \times S_u^{(')} \widetilde M (\Gamma) \longrightarrow S_u^{(')} \widetilde M (\Gamma) \, .
\end{equation}

\noindent 2) {\it Even better, the $S_u^{(')}$'s have good localization and glueing properties, making that one can define {\ibf directly}, downstairs the $S_u^{(')} M(\Gamma)$, proceeding like for $S_u^{(')} \widetilde M(\Gamma)$ upstairs; with this, we have that
\begin{equation}
\label{eq2.21}
S_u^{(')} M(\Gamma)^{\sim} = S_u^{(')} \widetilde M (\Gamma) \, .
\end{equation}
Moreover, for the directly defined $S_u^{(')} M(\Gamma)$ we find that
$$
S_u^{(')} M(\Gamma) = (S_u^{(')} \widetilde M (\Gamma))/\Gamma \, , \eqno (2.21.1)
$$
in terms of the action from} (\ref{eq2.20}) {\it above.}

\bigskip

In (2.17.1) we have introduced the $\Theta^3 (X^2)$, which has the $3^{\rm d}$ version of the undrawable singularities ${\rm Sing} \, (f) \subset X^2$, and of course, also, the $C(p_{\infty\infty})$'s, another kind of singularity.

\bigskip

\noindent {\bf Definition (2.22).} We will introduce now, on the same lines as in (2.19.bis) the $\Theta^{N+4} (X^2)$, which completes (2.17.1). With $p_{\infty\infty} = p_{\infty\infty} ({\rm all})$, we will take $\Theta^{N+4} (X^2) \mid \{$outside of $D^2 (H(p_{\infty\infty}))\} \equiv \{($the already smooth $\Theta^4 (\Theta^3 (X^2) , {\mathcal R}) \mid\{$outside of $D^2 (H(p_{\infty\infty}))\}\} \times B^N$, and $\Theta^{N+4} (X^2) \mid D^2 (H(p_{\infty\infty})) \equiv D^2 (H$ $(p_{\infty\infty})) \times \left[-\frac{\varepsilon}4 \leq z \leq \frac{\varepsilon}4 \right] \times \frac12 \, B^{N+1}$, the two pieces being put together the obvious way.

\bigskip

\noindent {\bf Lemma~2.2.} {\it There is a transformation leading from $S'_u \, \widetilde M (\Gamma)$ to $S_u \, M(\Gamma)$, which proceeds as follows. This transformation is localized in the neighbourhood of the $p_{\infty\infty} \, ({\rm proper})$, and {\ibf no other} $p_{\infty\infty}$'s will be talked about now. We start by an infinite PROPER collapse, destroying the singularity along $C(p_{\infty\infty}) \times \bigl[-\frac{\varepsilon}4 \leq z$ $\leq \frac{\varepsilon}4 \bigl] \times \frac12 \, B^{N+1}$, followed by a PROPER smooth infinite dilatation, the whole process being}
\setcounter{equation}{22}
\begin{equation}
\label{eq2.23}
S'_u (\widetilde M (\Gamma)) \underset{\rm collapse \, and \, dilatation}{-\!\!\!-\!\!\!-\!\!\!-\!\!\!-\!\!\!-\!\!\!-\!\!\!-\!\!\!-\!\!\!-\!\!\!-\!\!\!-\!\!\!-\!\!\!-\!\!\!-\!\!\!\longrightarrow} S_u (\widetilde M (\Gamma)) \, .
\end{equation}

\bigskip

\noindent {\bf Proof.} With a ``small'' $\varepsilon_0 > 0$, to be more precise later on, we consider the PROPER embedding
\begin{equation}
\label{eq2.24}
\sum_{p_{\infty\infty} ({\rm proper})} (A(p_{\infty\infty}) - \partial \, A (p_{\infty\infty})) \times \left( \frac12 + \varepsilon_0 \right) B^{N+1} \overset{j}{\longrightarrow} S'_u \, \widetilde M (\Gamma) \, .
\end{equation}
Here $A(p_{\infty\infty})$ is like in figure~2.2 and the embedding $j$ is disjoined from the $(N+4)$-dimensional 2-handle $D^2 (H(p_{\infty\infty})) \times \left[-\frac{\varepsilon}4 , \frac{\varepsilon}4 \right] \times \frac12 \, B^{N+1} \} - \{$its attaching zone$\}$. The collapsing part of (2.22) removes the set
$$
\sum_{p_{\infty\infty}} (A(p_{\infty\infty}) - \partial \, A(p_{\infty\infty})) \times \left( \frac12 + \varepsilon_0 \right) B^{N+1} \, ,
$$
leaving us with the 2-handle above still attached to the rest, along $C(p_{\infty\infty}) \times \left[-\frac{\varepsilon}4 , \frac{\varepsilon}4 \right] \times \frac12 \, B^{N+1}$, which by now is no longer singular.

\smallskip

So, we are left with (and below, $\Theta^3$ is like in the figure 2.2) the space
\begin{equation}
\label{eq2.25}
\left[\Theta^3 \times B^{N+1} - \left([A(p_{\infty\infty}) - \partial \, A(p_{\infty\infty})] \times \left( \frac12 + \varepsilon_0 \right) B^{N+1} \right) \right] \cup
\end{equation}
$$
D^2 (H(p_{\infty\infty})) \times \left[-\frac{\varepsilon}4 , \frac{\varepsilon}4 \right] \times \frac12 \, B^{N+1} , \ \mbox{which is no longer singular along} \ 
C(p_{\infty\infty}) \times \left[-\frac{\varepsilon}4 , \frac{\varepsilon}4 \right] \times \frac12 \, B^{N+1} \, .
$$
In order to describe the dilatation part of (\ref{eq2.23}) some notations will be introduced with $D^2 \equiv D^2 (H(p_{\infty\infty}))$, we have now the following decomposition, easily readable from the figure~2.2: To begin with, in $2^{\rm d}$ we have
$$
D^2 = (D^2 \cap A) \cup (D^2 - A) \, ,
$$
and here we should really read $H(p_{\infty\infty}) \underset{\rm TOP}{=} d^2 \times I$, thickened as the figure~2.2 really suggests, so that $D^2 - A \approx d^2$ and then next, in dimension three too
$$
A \mid \left( -\frac{\varepsilon}2 , \frac{\varepsilon}2 \right) = \left( (A \cap D^2) \mid \left( -\frac{\varepsilon}4 , \frac{\varepsilon}4 \right) \right) \cup  \left( (A \cap D^2) \mid \left( \pm \, \frac{\varepsilon}4 , \pm \, \frac{\varepsilon}2 \right) \right) \, .
$$
With these things, what the dilatation part of (\ref{eq2.23}) adds, to $\{ S'_u \, \widetilde M (\Gamma)$ after the collapse above has been done$\}$, is the following set
$$
\Biggl\{ A \mid \left( -\frac{\varepsilon}2 , \frac{\varepsilon}2 \right) \times \left( \frac12 + \varepsilon_0 \right) B^{N+1} -  \left(\left( (A \cap D^2) \mid \left( -\frac{\varepsilon}4 , \frac{\varepsilon}4 \right) \right) \times \frac12 \, B^{N+1} \right) \begin{pmatrix} \mbox{a piece which} \\ \mbox{was never} \\ \mbox{removed} \end{pmatrix} \Biggl\}
$$
$$
\cup \left\{ (D^2-A) \times \left[ -\frac{\varepsilon}4 , \frac{\varepsilon}4 \right] \times \left( \left( \frac12 + \varepsilon_0 \right) B^{N+1} - \frac12 \, B^{N+1} \right) \right\}
\cup \left\{ (D^2-A) \times \left[ \pm \, \frac{\varepsilon}4 , \pm \, \frac{\varepsilon}2 \right] \times \left( \frac12 + \varepsilon_0 \right) B^{N+1} \right\} \, .
$$

\medskip

To understand this formula, one should notice that the $2$-handle $D^2 \times \left[ -\frac{\varepsilon}4 , \frac{\varepsilon}4 \right] \times \frac12 \, B^{N+1}$ can be expressed as
$$
(A \cap D^2) \times \left[ -\frac{\varepsilon}4 , \frac{\varepsilon}4 \right] \times \frac12 \, B^{N+1} \cup (D^2 - A) \times \left[ -\frac{\varepsilon}4 , \frac{\varepsilon}4 \right] \times \frac12 \, B^{N+1} \, ,
$$
and that, as figure~2.2 tells us, the $A$ stretches, as part of $\Theta^3$, in the $z$-direction, all the way from $-\frac{\varepsilon}2$ to $+\frac{\varepsilon}2$.

\smallskip

If we forget the $(N+1)$-dimensional multiplying factors, $\frac12 \, B^{N+1}$, or $\left( \frac12 + \varepsilon_0 \right)$ $B^{N+1}$ or $B^{N+1}$, then in dimension three this corresponds to what one sees in figure~2.2 and, so far, everything is OK. But what we said so far is only rough, in the supplementary $N+1$ dimension, where the size of the multiplying ball $B^{N+1}$ {\ibf jumps}. So, here is what we will actually do, in real life.

\smallskip

Notice that, in the neighbourhood of the $p_{\infty\infty} \, ({\rm proper})$ considered, our BLACK wall $W_{(\infty)}$ is reduced to $D^2 = D^2 (H(p_{\infty\infty})) \subset W_{(\infty)}$. With $q \in D^2$ being the genric point in $W_{(\infty)} - \{p_{\infty\infty}\}$, we will let the size of the multiplying factor $B^{N+1}$ occurring in $(W_{(\infty)} - \{p_{\infty\infty}\}) \times [-\varepsilon , \varepsilon] \times B^{N+1}$ to depend on $q$.

\bigskip

\noindent {\bf Claim (2.26).} The size of $\varepsilon_0$ can be made $q$-dependent, $\varepsilon_0 = \varepsilon_0 (q)$ so that, in the context of (\ref{eq2.24}), as developed in the proof below, and in a manner which will be compatible with the geometric realization of the zipping in high dimensions, something which is still to come later, the whole process can proceed smoothly, without jumps in the size of the $(N+1)$-factor.

\bigskip

The proof of the CLAIM will be given in section~IV, where its content should become clearer too.

\smallskip

The next statement is the main result of the present paper and it is also the main step in our proof that all finitely presented groups are QSF. The whole proof requires three successive papers, of which \cite{29} is the first and the present the middle one.

\bigskip

\noindent {\bf The GSC theorem 2.3.} 

\smallskip

\noindent 1) {\it The $(S'_u \, \widetilde M (\Gamma))_{\rm II}$ is GSC.}

\medskip

\noindent 2) {\it As a consequence of this, $(S_u \, \widetilde M (\Gamma))_{\rm II}$ is also GSC.}

\bigskip

The proof occupies the rest of the present paper. Here are some COMMENTS concerning the GSC theorem.

\medskip

A) To begin with, I should explain the subscript II occurring in the statement above. There are, actually, two variants for the construction of our functor $S_u^{(')}$, the variant I which is what we just did, and then also a variant II too, to be described in detail later on in this paper. The subscript II refers to it. With this, the statement of theorem 2.3 should be taken as an indication of what we are after.

\medskip

B) If the $S_u \, \widetilde M (\Gamma)$ would be compact, which is certainly NOT, then our GSC theorem~2.3 would already imply that $\Gamma$ is QSF. As things stand, still another paper will be necessary for deducing that result from the GSC theorem~2.3 and its proof. With \cite{29} and the present one, that makes a total of these papers for proving that all $\Gamma$'s are QSF. Of the three, the longest and most difficult is certainly the present one.

\medskip

C) The GSC theorem is at the heart of the proof that $\forall \, \Gamma \in {\rm QSF}$, which certainly needs that $(S_u \, \widetilde M (\Gamma))_{\rm II} \in {\rm GSC}$. But in order to get this, we need to show first that $(S'_u \, \widetilde M (\Gamma))_{\rm II} \in {\rm GSC}$.

\medskip

D) The main technical tool for proving the GSC theorem~2.3 will be a transformation going, very roughly speaking, from $\Theta^{N+4} (X^2)$ to an object related to $S'_u \, \widetilde M (\Gamma)$, to be a bit more explicit, from $\{ \Theta^{N+4} (X^2)$ with some Holes and DITCHES deleted$\}$ via an infinity of successive steps and taking the following general form,
\setcounter{equation}{26}
\begin{equation}
\label{eq2.27}
\Theta^{N+4} (X^2-H) - {\rm DITCHES} \overset{Z}{=\!\!\!=\!\!\!\Longrightarrow} S'_b (\widetilde M (\Gamma) - H) \, .
\end{equation}
I will explain the notations used in (\ref{eq2.27}). The ``$Z$'' stands for ``zipping'', and the transformation in question  mimicks the zipping process. The left hand side of (\ref{eq2.27}) is ``$\Theta^{N+4} (X^2)$ with Holes and DITCHES deleted''. The Holes and DITCHES are explained first in section~III, then with more details in IV. In IV there will also be a space $S'_u (\widetilde M (\Gamma)-H)$, meaning ``$S'_u \, \widetilde M (\Gamma)$ with Holes''. For the time being, the Holes are still mythical.

\smallskip

But the final object of the construction $Z$ is really another object $S'_b \, \widetilde M (\Gamma)$ defined eventually via a reconstruction formula 
$$
S'_b \, \widetilde M (\Gamma) \equiv S'_b (\widetilde M (\Gamma)-H) + \{\mbox{appropriate 2-handles}\} \, .
$$
The $S'_b (\widetilde M (\Gamma)-H)$ which occurs in the RHS of (\ref{eq2.27}) is a cousin of $S'_u (\widetilde M (\Gamma)-H)$, which will be explicitly constructed in section~IV during the unrolling of the ZIPPING LEMMA, with more details in V, when the zipping lemma will be proved. It is a pivotal object for the whole approach.

\smallskip

Coming back to the $Z$, which mimicks the strategy $X^2 \Longrightarrow fX^2$ for (\ref{eq1.1}), it will be referred to as the ``{\ibf geometric realization of the zipping}'', in high dimensions. But, while the $2^{\rm d}$ zipping strategy is a gigantic quotient-space projection, our $Z$ will be an infinite sequence of inclusion maps, mixed from time to time with another kind of step, which we will call the DITCH-{\ibf jumping}.

\smallskip

Now, for expository purposes, instead of plunging straight into the general construction of the geometric realization of the zipping, we will open a prentice and, in te next section, a simplified TOY MODEL, exhibiting some of the important features of the real life thing, in a simplified context. It should be helpful for understanding the real story.

\smallskip

We will very much refer to it in the more technical section which will follow afterwards.

\medskip

D) We end this section with a brief review of the ARTICULATION of the present paper, and its successive changes of dimensions which are used. We start from the $3^{\rm d}$ REPRESENTATION
$$
Y(\infty) \overset{g(\infty)}{-\!\!\!-\!\!\!-\!\!\!-\!\!\!\longrightarrow} \widetilde M (\Gamma) \sim \Gamma \ \mbox{from (1.6)}. \eqno (*_1)
$$
Here $Y(\infty)$ is a union of $3^{\rm d}$ bicollared handles and what we gain (and this never to be lost again), are local finiteness, $\Gamma$-equivariance and uniformly bounded zipping length. Then we go to an appropriately dense 2-skeleton of $(*_1)$ and change $(*_1)$ into the $2^{\rm d}$ REPRESENTATION
$$
X^2 \overset{f}{-\!\!\!-\!\!\!-\!\!\!-\!\!\!\longrightarrow} \widetilde M (\Gamma) \ \mbox{from (1.1)}. \eqno (*_2) 
$$
All the richness of the structure of double points of maps from dimension two to dimension three is now at our disposal and, more specifically, a clear useful {\ibf zipping strategy} for $(*_2)$ is available too. This means a not everywhere well defined {\ibf zipping flow} on $X^2$, the intersections of which with another similar, {\ibf collapsing flow} (on the same $X^2$), will play an important role.

\smallskip

At this point, partly for reasons of circumventing the lack of local finiteness of $fX^2$ but mostly for opening the door to the high dimensions where our key for $\forall \, \Gamma \in {\rm QSF}$ eventually lies, we go again 3-dimensional, but now in a very different context than in $(*_1)$, and we will construct the locally finite cell-complexes
$$
\Theta^3 (X^2) \quad {\rm and} \quad \Theta^3 (fX^2)^{(')} \, . \eqno (*_3) 
$$
These in turn will be thickened into dimension $N+4$, going to
$$
\Theta^{N+4} (X^2) \quad {\rm and} \quad S'_u \, \widetilde M (\Gamma) \, , \eqno (*_4) 
$$
itself essentially are $(N+4)$ dimensional thickening of the intermediary, 4-dimensional $\Theta^4 (\Theta^3 (fX^2)^{(')} , {\mathcal R})$. All these things have just been done, above. The geometric realization of the zipping, developed in the sections~IV, V below, will make great use of the {\ibf additional} $N$ dimensions, with respect to the three of $(*_3)$ or four of $\Theta^4$ in $(*_4)$. This is how the $N+4$ comes in.

\smallskip

The space $X^2$ is GSC, actually $\Theta^3 (X^2)$ is even collapsible. There is a relatively large contribution
\begin{equation}
\label{eq2.28}
W({\rm COLOUR}) \cap \partial X^2 \ne \emptyset \, , \ \mbox{beyond} \ W({\rm BLACK}) \cap \partial X^2 \ne \emptyset \, .
\end{equation}
This has nothing to do with the Holes $H$ which will occur in (\ref{eq4.7}), but when we will move from $X^2$ to $X^2 - H$, then the contributions to $\partial (X^2 - H)$ coming from (\ref{eq2.28}) and from (\ref{eq4.7}) (i.e. $\partial H ({\rm normal})$), will be physically undistinguishable. By contrast to $X^2$, the $X^2-H$ is highly non-simply connected.

\newpage

\section{A toy model}\label{sec3}
\setcounter{equation}{0}

The situation we will present now is highly simplified from several respects: there are no immortal singula\-rities (\`a la ${\rm Sing} \, \widetilde M(\Gamma)$), there are no triple points and no $p_{\infty\infty}$'s either (hence no difference between $S_u$ and $S'_u$). But then there will not be any group action present now either. We replace (\ref{eq1.1}) by the following map which {\ibf is} our toy model
\begin{equation}
\label{eq3.1}
X^2 \overset{f}{-\!\!\!-\!\!\!-\!\!\!-\!\!\!\longrightarrow} R^3 \, ,
\end{equation}
where

\medskip

\noindent (3.2) \quad $X^2 \equiv \{$the $2^{\rm d}$ region called $R$, homeomorphic to $I \times {\rm int} \, I$ and PROPERLY embedded inside $(z=0) \subset R^3 = \{x,y,z \}$, which is displayed in the figure~3.1.(A)$\} \cup \underset{n=1}{\overset{\infty}{\sum}} D_n^2$, with $R$ and each of the $D_n^2$'s glued along $[\sigma_n , \Sigma_n] + [s_n , S_n]$. Here $D_n^2 = \{$the disk of diameter $[\sigma_n , s_n] = R \cap (x=x_n)$, living in the plane $(x=x_n) \subset R^3$ and displayed in the figure~3.1.(B)$\}$.

\medskip

Figure~3.1 should explain how the spare parts of $X^2$, listed in (3.2) are glued together. Restricted to each of these span parts, the map $f$ is the natural embedding into $R^3$.

$$
\includegraphics[width=165mm]{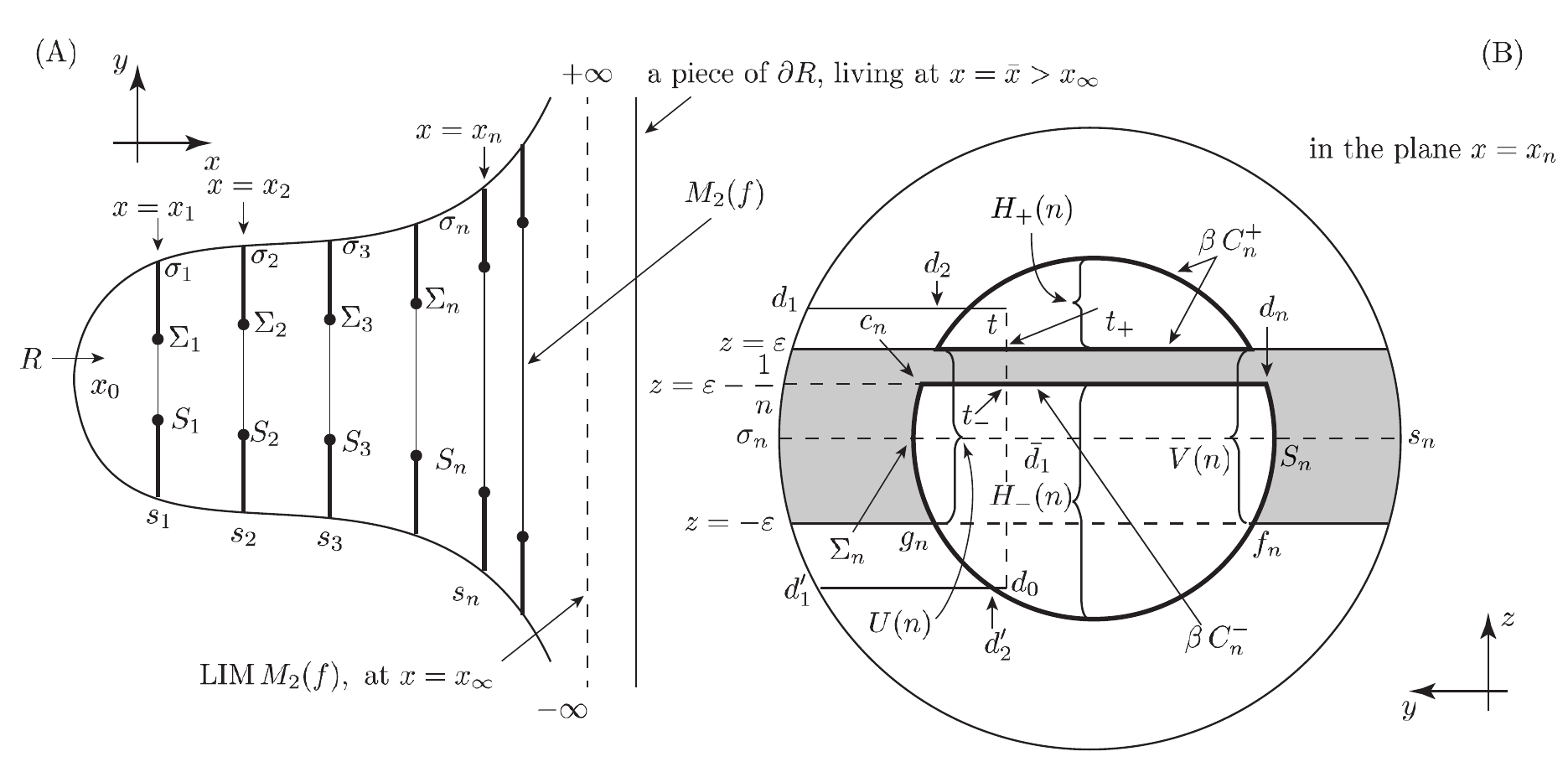}
$$
\label{fig3.1}

\centerline {\bf Figure 3.1.} 

\smallskip

\begin{quote} 
We display here the spare parts $R$ and $D_n^2$ of $X^2$. The $\Sigma_n$, $S_n$ are undrawable singularities of the $f$ (\ref{eq3.1}), and the arcs $U(n) \ni \Sigma_n$, $V(n) \ni S_n$ stand for undrawable singularities of $\Theta^3 (X^2)$. The $H_{\pm} (n)$ are Holes, open sets, drilled out of $D_n^2$.
\end{quote}

\bigskip

We will not need any fins now, so proceeding essentially like in (\ref{eq2.5}) we will define the noncompact smooth 3-manifold, with large boundary
\setcounter{equation}{2}
\begin{equation}
\label{eq3.3}
\Theta^3 (fX^2) \equiv (R \times [-\varepsilon \leq z \leq \varepsilon]) \cup \sum_1^{\infty} D_n^2 \times [-\varepsilon_n \leq x \leq \varepsilon_n] - {\rm LIM} \, M_2 (f) \times \{ z = \pm \, \varepsilon \} \, .
\end{equation}

Next, by analogy to (\ref{eq2.25}) and to (2.1), we introduce the smooth $(N+4)$-manifolds
\begin{equation}
\label{eq3.4}
\Theta^{N+4} (X^2) = \Theta^4 (\Theta^3 (X^2) , {\mathcal R}) \times B^N
\end{equation}
and
$$
\Theta^{N+4} (fX^2) = \{\mbox{the already smooth} \ \Theta^3 (fX^2) \} \times B^{N+1} \, .
$$

In anticipation for what we will do in real life, we will try now to suggest a transformation leading from $\Theta^{N+4} (X^2)$ to $\Theta^{N+4} (fX^2)$ via a succession of steps of the following forms

\medskip

\noindent (3.5.1) \quad Inclusion maps mimicking the infinitely many  elementary zipping operations leading from $X^2$ to $fX^2$. There is no question of proving any GSC property now, but we will still insist that all our steps should be, abstractly speaking, i.e. just by themselves, GSC preserving;

\medskip

\noindent (3.5.2) \quad Deletion of the ${\rm LIM} \, M_2 (f) \times \{ z = \pm \, \varepsilon \} = \partial \sum (\infty)$. This should be realized directly at the level of $\Theta^{N+4} (X^2)$ via an infinite sequence of Whitehead dilatations which sends it to infinity, a GSC-preserving step too.

\medskip

Realizing the more serious step (3.5.1) will need the introduction of some basic tools to be used later on in real life too, namely the

\medskip

\noindent (3.6) \quad  The Holes, the Ditches and the (Partial) Ditch-Filling.

\medskip

But above I had used the word ``try'', the point being that for our project to succeed, we may be forced to throw other, still GSC preserving ingredients, into the definition of $\Theta^{N+4} (fX^2)$ too. In anticipation for that we denote the $\Theta^{N+4} (fX^2)$ from (\ref{eq3.4}) by $\Theta^{N+4} (fX^2)_{\rm I}$ and we will refer to its context as variant~I. The ingredients (3.6) to which we turn next, are there in both variants~I and II (not yet introduced), but more will have to  be thrown in when it will come to II.

\medskip

In figure 3.1.(B) the shaded areas which are outside of the central circle which contains the $H_{\pm} (n)$, stand for identifications to be performed at level $\Theta^3 (X^2)$ leaving us with the singularities $U(n)$, $V(n)$ at $3^{\rm d}$ level. Next, inside each individual $D_n^2 \approx D_n^2 \times [-\varepsilon_n, \varepsilon_n]$ we will drill two Holes $H_n (\pm)$ like in figure~3.1.(B), deleting open sets (unlike the deletion of the closed $p_{\infty\infty} \times (-\varepsilon , \varepsilon)$ in (\ref{eq2.20})).

\smallskip

This leaves us with the boundary curves
\setcounter{equation}{6}
\begin{equation}
\label{eq3.7}
C_n^{\pm} \equiv \partial \, \{{\rm Hole} \ H_{\pm} (n)\} \, ,
\end{equation}
coming with
\begin{equation}
\label{eq3.8}
\left\{\mbox{canonically framed link} \ \sum_n C_n^{\pm} \right\} \overset{\beta}{-\!\!\!-\!\!\!\longrightarrow} \partial \, \Theta^{N+4} (X^2 - H) \, ;
\end{equation}
here $\Theta^{N+4} (X^2-H) \equiv \{ \Theta^{N+4} (X^2)$ with holes deleted$\}$ and $\beta$ is like in figure~3.1.(B). The shaded area in figure~3.1.(B), living inside the central circle and contained in $\varepsilon - \frac1n \leq z \leq \varepsilon$, corresponds to $(D_n^2 - H) \cap f^{-1} \left(R \times \left[ \varepsilon - \frac1n , \varepsilon \right] \right)$, and it stays far from the identifications performed at level $\Theta^{N+4} (X^2-H)$. We also introduce the smooth manifold $\Theta^{N+4} (fX^2-H) \equiv \{\Theta^{N+4} (fX^2)$ from which we delete the holes $H_+ (n)$ and $H_- (n) \mid (z < -\varepsilon)$ (see here the figure~3.1.(B)$\}$.

\smallskip

This comes with its own
\begin{equation}
\label{eq3.9}
\left\{\mbox{canonically framed link} \ \sum_n C_n^{\pm} \right\} \overset{\alpha}{-\!\!\!-\!\!\!\longrightarrow} \partial \, \Theta^{N+4} (fX^2 - H) \, .
\end{equation}
We find here that $\alpha \, C_n^+ = f \, \beta \, C_n^+$, {\ibf BUT} the $\alpha \, C_n^-$ and $f \, \beta \, C_n^-$ are far from each other; the reason for this is that while $\beta \, C_n^- \subset \Theta^{N+4} (X^2-H)$ goes all the way up to $z = \varepsilon - \frac1n$, the $\alpha \, C_n^- \subset \Theta^{N+4} (f X^2-H)$ only reaches up to $z=-\varepsilon$. Both $\alpha \, C_n^-$ and $\beta \, C_n^-$ can be reconstructed from figure~3.2. 

\smallskip

In figure~3.2 we also see the DITCHES, when we restrict ourselves to some fixed value of $y \in [y(V(n))$, $y(U(n))]$. The piece $T(n) \subset \Theta^{N+4} (X^2 - H)$-DITCHES, as defined by figure~3.2 (and then with more detail in the context of (\ref{eq3.12})), lives in the DITCH, without touching the lateral surface $\delta \, \mbox{ditch} \, (n) \mid y$. Here
$$
\mbox{ditch} \, (n) \subset \Theta^{N+4} (X^2-H) \mid R = \Theta^{N+4} (X^2) \mid R \, .
$$

When the gap between $T(n)$ and the lateral surface of the ditch will be filled with material, during the process DIL (\ref{eq3.18}) (which is our present version of the geometrical realization of the zipping), then the curves $\alpha \, C_n^+$ and $\eta \, \beta \, C_n^+$ (as defined by (\ref{eq3.20}) below) will be connected by a short simple-minded isotopy. The filling material mentioned above only occurs inside $\left[\varepsilon - \frac1n \leq z \leq \varepsilon \right] \subset [-\varepsilon \leq z \leq \varepsilon]$, and the same will be true for $T(n)$. We speak hence of a {\ibf PARTIAL DITCH FILLING}.

\bigskip

Here are some EXPLANATIONS CONCERNING THE FIGURE 3.2. The ditch $(n) \mid y$ is a region which is concentrated around the following arc
$$
A(n,y) \equiv (x=x_n , y,-\varepsilon \leq z \leq \varepsilon \, , \ \{ \mbox{some point} \ p \in \partial B^{N+1} , \ \mbox{here} \ t=1\}) \subset R \times [-\varepsilon , \varepsilon] \times B^{N+1} \, ,
$$
and which is contained in $\Theta^{N+4} (fX^2-H) \mid (R \mid y)$. Let us say that
$$
{\rm ditch} \, (n) \mid y = A(n,y) \times b^{N+1} (n) \, ,
$$
with $b^{N+1} (n) \subset B^{N+1}$ concentrated around $p \in \partial B^{N+1}$. With this, for a fixed, generic value $y$, the (A) presents an infinite family of ditches, occurring as {\ibf indentations} inside
$$
\{\mbox{the $R$ in figure 3.1}\} \times [-\varepsilon \leq z \leq \varepsilon] \times B^{N+1} \, .
$$
They make use of the {\ibf additional} dimensions, i.e. of the additional factor $B^{N+1}$ with respect to some already smooth $\Theta^3$, respectively $B^N$ with respect to some smooth $\{\Theta^4$ which has just tamed the undrawable singularities of some $\Theta^3\}$. We insist on this distinction between ``additional dimensions'' as opposed to the mere ``high dimensions''.

\smallskip

In (B) we present
$$
T(n) \mid y \equiv {\rm Tube} \, (n) \mid y \subset {\rm ditch} \, (n) \, ,
$$
and this inclusion will be part of the forthcoming map $j$ in (\ref{eq3.13}). One should imagine here a complete $y$-movie coming with the drawing (B). This movie creates then the
$$
\bigcup_{\overbrace{\mbox{\footnotesize $y ([t_- t_+]$ of $V(n)) \leq y \leq y([t_- , t_+]$ of $U(n))$}}} T(n) \mid y = \{\mbox{The shaded area of $D^2-H$,}
$$
$$
\mbox{living inside the central circle in the figure 3.1.B, viewed now at the level of $\Theta^{N+4} (X^2 - H)$}\} \, .
$$

\smallskip

To be more precise, we actually have
\begin{equation}
\label{eq3.10}
T(n) \mid y = \Theta^{N+4} (X^2-H) \mid \left\{ ((D_n^2 - H) \mid y) \cap \left[ \varepsilon - \frac1n \leq z \leq \varepsilon \right] \right\} \, ,
\end{equation}
for $y ([t_- t_+]$ of $V(n)) \leq y \leq y([t_-  t_+]$ of $U(n))$.

\smallskip

From the (A) one should read an easy diffeomorphism isotopic to the identity $\Theta^{N+4} (fX^2 - H) \mid (R) - \{$DITCHES$\} = \Theta^{N+4} (fX^2 - H) \mid (R)$, and here clearly, we have
\begin{equation}
\label{eq3.11}
\Theta^{N+4} (fX^2-H) \mid R = \Theta^{N+4} (fX^2) \mid R = \Theta^{N+4} (X^2) \mid R = \Theta^{N+4} (X^2-H) \mid R \, ,
\end{equation}
and in these formulae we may as well replace $R$ by $R \mid y$. This ends our explanations for figure 3.2.

$$
\includegraphics[width=175mm]{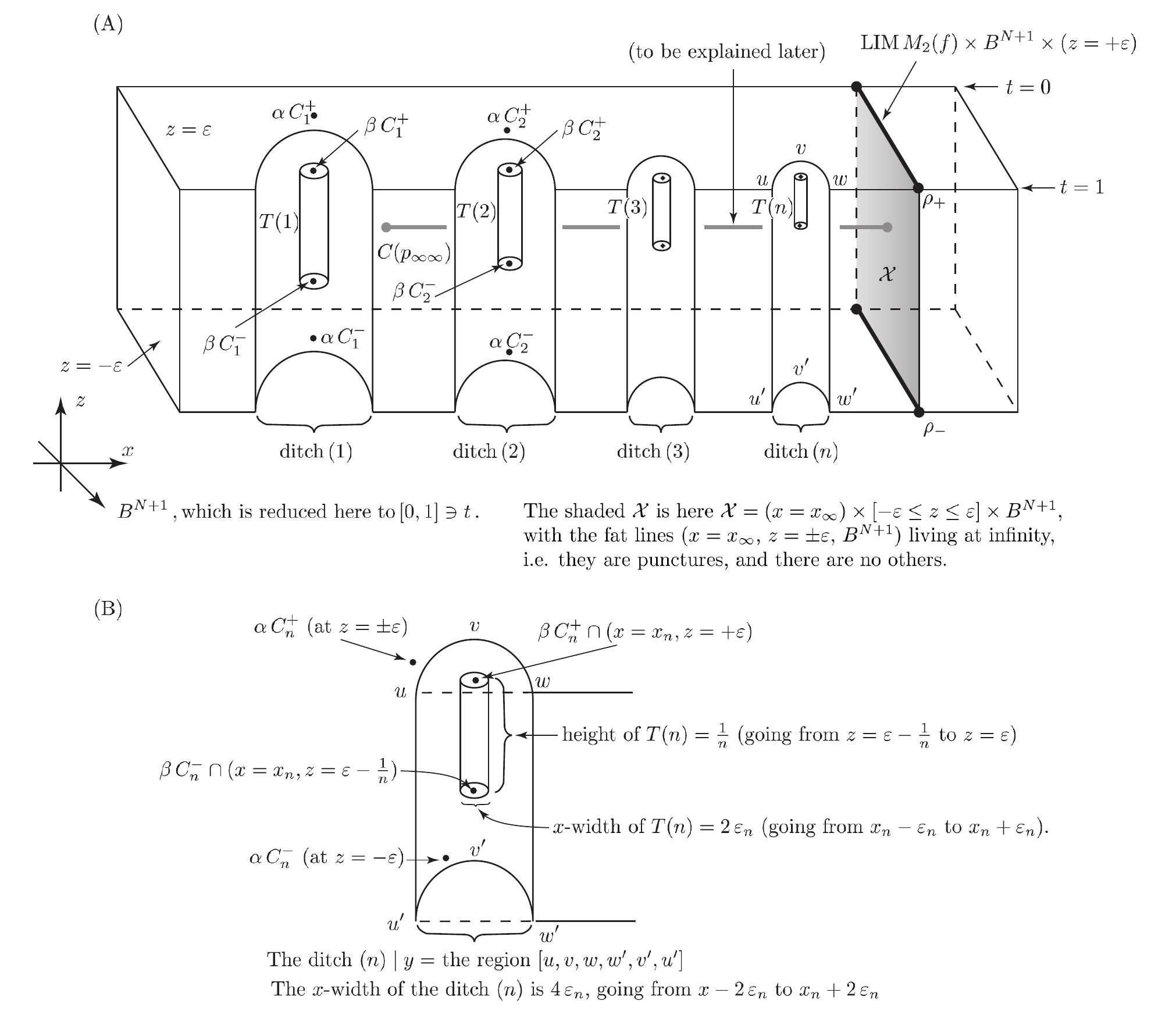}
$$
\label{fig3.2}

\centerline {\bf Figure 3.2.} 

\smallskip

\begin{quote} 
Everything we see both in (A) and in (B), lives at some generic value $y \in [y \, V(n),y \, U(n)]$, with $V(n) , U(n)$ like in figure~3.1.(B). In the present (B) we see the ditch $(n) \mid y$ and, with this ditch $(n) \equiv \bigcup \, {\rm ditch} \ (n) \mid y$, with $y ([t_- , t_+]$ of $V(n)) \leq y \leq y ([t_- , t_+]$ of $U(n))$. Here $[t_- , t_+]$ is like in figure 3.1.(B), i.e. $y([t_- , t_+]$ of $V(n))$ (respectively $y([t_- , t_+]$ of $U(n))$ is slightly larger than $yV(n)$ (respectively slightly smaller than $y \cup (n)$); see the $y$ coordinates in the figure 3.1.(B).
\end{quote}

\bigskip

Both $\Theta^{N+4} (X^2)$ and $\Theta^{N+4} (fX^2)$ can be reconstructed, up to diffeomorphism, starting from the corresponding objects with Holes and then adding 2-handles along (\ref{eq3.8}), respectively along (\ref{eq3.9}).

\smallskip

Also, 
$$
{\rm DITCHES} \equiv \sum_{n=1}^{\infty} {\rm ditch} \, (n) \subset \Theta^{N+4} (X^2 - H) \mid R = \Theta^{N+4} (X^2) \mid R \, .
$$
In the same vein, when it comes to the curves $\alpha \, C_n^{\pm}$, $\beta C_n^{\pm}$, we only see here the $C \mid y$'s.

\smallskip

So, what we  learn from figure~3.2 is that, by making use of the additional dimensions $B^{N+1}$, we define first the restrictions ${\rm ditch} \, (n) \mid y$ for all the
$$
y ([t_- t_+] \ \mbox{of} \ V(n)) \leq y \leq y([t_-  t_+] \ \mbox{of} \ U(n)) \, , \eqno (*)
$$
then the DITCH $(n)$ by summing $(*)$ over $y$, and finally the complete system $\{{\rm DITCHES}\} \subset \Theta^{N+4} (X^2 - H)$, by summing over $n$.

\bigskip

\noindent {\bf Lemma 3.1.} {\it There is a splitting, to be explicitly defined in the proof below}
\begin{equation}
\label{eq3.12}
\Theta^{N+4} (X^2-H) - \{\mbox{DITCHES}\} = \{{\rm MAIN} \ \Theta^{N+4}\} + \sum_n \{4^{\rm d} \, 1\mbox{-handles} \ T(n) \} \, .
\end{equation}
{\it There is also a PROPER embedding}
\begin{equation}
\label{eq3.13}
\xymatrix{
\Theta^{N+4} (X^2 - H) - \{\mbox{DITCHES}\} \ar[rr]^-{j} \ar[d]^-{p} &&\Theta^{N+4} (fX^2-H) \ar[d]^-{q}  \\ 
X^2 - H \ar[rr]^-{f} &&fX^2-H \, . 
}
\end{equation}

\bigskip

\noindent {\bf Proof.} At $2^{\rm d}$ level, we have ${\rm Sing} \, (f) = \{$the $\Sigma_n$'s, $S_n$'s in figure~3.1.A$\} \subset X^2-H$ and its isomorphic image $f \, {\rm Sing} \, (f) \subset fX^2 - H$. Let $x \in {\rm Sing} \, (f)$; since $H_- (n) \subset X^2$ reaches to $z = \varepsilon - \frac1n$, while $H_- (n) \subset fX^2$ only reaches to $z= -\varepsilon$, our $x$ is actually non-singular for $X^2-H$, while it is singular for $fX^2 - H$. We consider canonical neighbourhoods
$$
N = Nbd \, ({\rm Sing} \, (f)) \overset{\rm PROPER \, embedding}{-\!\!\!-\!\!\!-\!\!\!-\!\!\!-\!\!\!-\!\!\!-\!\!\!-\!\!\!-\!\!\!-\!\!\!-\!\!\!-\!\!\!-\!\!\!-\!\!\!-\!\!\!\longrightarrow} X^2 \, .
$$
We use figure~3.1.(B) in order to describe the singular local models $N \mid x \in X^2$, proceeding like in \cite{8}, \cite{18}, \cite{19}. So, let $x = \Sigma_n$, and the same embellishments as those drawn in figure 3.1.B around $\Sigma_n$, should be understood around $S_n$ too. We will take
\begin{equation}
\label{eq3.14}
N \mid x = (R \mid [\sigma_n , \Sigma_n , \bar d_1]) \underset{\overbrace{\mbox{\footnotesize $[\sigma_n , \Sigma_n]$}}}{\cup}
\end{equation}
$$
\{\mbox{The $2$-cell $\Delta^2 = [d_1 \, d_2 \, t \, t_+ \, t_- \, \bar d_1 \, d_0 \, d'_2 \, d'_1]$, cobounding a circle $C$, contained in $D_n^2$}\} \, .
$$
Going, next $3$-dimensional, the $\Theta^3 (X^2) \mid \Sigma_n$ is the obvious $3^{\rm d}$ version of (\ref{eq3.14}) split from the rest of $\Theta^3 (X^2)$ by a surface homeomorphic to $S^1 \times S^1 - {\rm int} \, D^2$, which is a regular neighbourhood $\delta \theta$ of
\begin{equation}
\label{eq3.15}
\partial \, (R \mid [\sigma_n , \Sigma_n , \bar d_1]) \underset{\sigma_2}{\cup} C \, .
\end{equation}
Notice that, at the level of $\Theta^3(X^2)$, the $\Sigma_n , S_n$ become small singular rectangles. Moving to $\Theta^4 (\Theta^3 (X^2),{\mathcal R})$ where ${\mathcal R}$ is a desingularization, we get, like in \cite{8}, \cite{18}, \cite{19}, an ${\mathcal R}$-dependent embedding
$$
\delta\theta \subset \partial \, \Theta^4 (\Theta^3 (X^2) \mid \Sigma_n , {\mathcal R}) = \partial B^4 \, .
$$
When the holes $H$ are now thrown in too, then the local model from (\ref{eq3.14}) undergoes the following change
$$
\Delta^2 \Longrightarrow \Delta^2 -\{\mbox{the two pieces of $\Delta^2$ which live inside $H_{\pm} (n)$, figure 3.1.(B)$\} \subset X^2 - H$} \, .
$$
Accordingly, $\delta\theta$ changes to the now disconnected
\begin{equation}
\label{eq3.16}
\delta (\theta - H) \underset{\rm TOP}{=} S^1 \times I + D^2 ([t_- \, t_+]) \, ,
\end{equation}
and up to isotopy, the canonical embedding
$$
\delta (\theta - H) \subset \partial \, \Theta^4 = \partial \, B^4
$$
is now ${\mathcal R}$-independent. The disks $\underset{\overbrace{\mbox{\footnotesize all $\Sigma_n$, $S_n$}}}{\sum} D^2 ([t_- \, t_+])$, induce now a splitting
$$
\Theta^4 (X^2-H)_{\rm non \, singular} - \{\mbox{DITCHES}\} = \{{\rm MAIN} \, \Theta^4 \} + \sum_n \{\mbox{$4^{\rm d}$ $1$-handles $T(n)$}\} \, .
$$

One gets our desired (\ref{eq3.12}) by multiplying this splitting with $B^N$. With this, here is how the $j$ in (\ref{eq3.13}) is defined, using (\ref{eq3.12}). The $j \mid \{{\rm MAIN} \, \Theta^{N+4}\}$ is the obvious canonical embedding. Next, like in the figure~4.2.(D), we introduce now for each $n$ a small $b^{N+1} (n) \subset B^{N+1}$, corresponding let us say to the indentations in figure~3.2.(A), such that ${\rm diam} \, b^{N+1} (n) \to 0$ when $n \to \infty$. With this, $j \mid (T(n) \mid y)$ sends $T(n) \mid y$ inside
$$
{\rm ditch} \, (n) \mid y - \partial \, {\rm ditch} \, (n) \mid y \subset {\rm DITCHES} \subset \Theta^{N+4} (fX^2 - H) \, ,
$$ 
essentially inside $b^{N+1} (n) \times \left[ \varepsilon - \frac1n \leq z \leq \varepsilon \right]$ as the figure~3.2 suggests to do it. This ends the description of (\ref{eq3.13}). It is the restriction to $\varepsilon - \frac1n \leq z \leq \varepsilon$ occurring above which makes that $j$ is PROPER. To be more precise concerning this issue, in the context of the figure~3.2 we have
$$
\lim_{n=\infty} \, \{\mbox{the $z$-height of} \ T(n) \mid y \} = 0 \, ,
$$
and hence
\begin{equation}
\label{eq3.17}
\lim_{n=\infty} \, j (T(n) \mid y) = \{\mbox{the point} \ (x = x_{\infty} , y , z = \varepsilon , (t=1) \in \partial B^{N+1})\} \, ,
\end{equation}
which lives at infinity. This ends the proof of our lemma. \hfill $\Box$

\bigskip

In the formulae (\ref{eq3.18}), (\ref{eq3.19}) below we work with an $y \in (-\infty , \infty)$ which is fixed. We introduce first the critical rectangle, shaded in figure~3.2.(A)
\begin{equation}
\label{eq3.18}
{\mathcal X} = \{ x = x_{\infty} , - \varepsilon \leq z \leq \varepsilon , t \in B^N \} \, ,
\end{equation}
which has exactly its two fat sides
\begin{equation}
\label{eq3.19}
(f \, {\rm LIM} \, M_2 (f)) \times \{ z = \pm \, \varepsilon \} \times B^{N+1}
\end{equation}
living at infinity. Notice that, in the present very simplified context $f \, {\rm LIM} \, M_2 (f) = {\rm LIM} \, M_2 (f) $, see figure~3.1.(A).

\bigskip

\noindent {\bf Lemma 3.2.} {\it When in the formulae like} (\ref{eq3.4}), {\it the operation of sending the punctures ${\rm LIM} \, M_2 (f) \times \{ z = \pm \, \varepsilon \}$ to infinity is already incorporated into the definition of $\Theta^3$, then all the four manifolds below are {\ibf transversally compact}
$$
\Theta^{N+4} (X^2) \, , \ \Theta^{N+4} (X^2 - H) \, , \ \Theta^{N+4} (fX^2) \, , \ \Theta^{N+4} (fX^2-H) \, .
$$
This means that, when $\Theta^3$ is smooth they are locally of the form $\Theta^3 \times B^{N+1}$ and, when $\Theta^4 (\Theta^3 , {\mathcal R})$ ($=$ smooth $4^{\rm d}$ thickening of $\Theta^3$) is required, then of the form $\Theta^4 \times B^N$.}

\bigskip

In the context of the figure~3.2 we have, inside the LHS of (\ref{eq3.13})
\begin{equation}
\label{eq3.20}
\lim_{n=\infty} \, {\rm ditch} \, (n) \mid y = \{\mbox{the whole line} \ (x=x_{\infty} , y, -\varepsilon \leq z \leq \varepsilon , t=1)\} \, ,
\end{equation}
and this formula should be compound to (\ref{eq3.17}). [Incidentally too, in both formulae the $y$ is fixed.] When we go from $\Theta^{N+4} (X^2-H)$ to $\Theta^{N+4} (X^2)$, respectively from $\Theta^{N+4} (fX^2-H)$ to $\Theta^{N+4} (fX^2)$ then the 2-handles $D^2 (\beta  C_n^{\pm})$, respectively $D^2 (\alpha \, C_n^{\pm})$, are supposed to shoot out of

\begin{equation}
\label{eq3.21}
\beta C_n^{\pm} \subset \partial (\Theta^{N+4} (X^2-H) - {\rm DITCHES}) \cap \partial \, \Theta^{N+4} (X^2-H) \subset 
\end{equation}
$$
\Theta^{N+4} (X^2-H) - {\rm DITCHES} \overset{j}{\longrightarrow} \Theta^{N+4} (fX^2-H) \subset R^3 \times B^{N+1} \, ,
$$
respectively out of
$$
\alpha \, C_n^{\pm} \subset \partial (\Theta^{N+4} (R) - {\rm DITCHES}) - \partial \, {\rm DITCHES} \subset \Theta^{N+4} (R) \subset \Theta^{N+4} (fX^2-H) \, ,
$$
without touching the rest of $\Theta^{N+4} (fX^2-H)$ (see figure~3.2).

\bigskip

\noindent {\bf Lemma 3.3.} 1) {\it Inside $\Theta^{N+4} (fX^2-H)$ there is an infinite system of smooth dilatations and possible additions of handles of index $\lambda > 1$, call it
\begin{equation}
\label{eq3.22}
j (\Theta^{N+4} (X^2-H) - {\rm DITCHES}) \underset{\rm DIL}{=\!\!=\!\!=\!\!\Longrightarrow} \{ j (\Theta^{N+4} (X^2-H) - {\rm DITCHES}) \, ,
\end{equation}
$$
\mbox{with all that part of} \ \bigcup_y \sum_n {\rm ditch} (n) \mid y \ \mbox{not yet occupied by} \ j \bigcup_y \sum_n T (n) \mid y \, ,
$$
$$
\mbox{{\ibf filled in}, but only {\ibf PARTIALLY}, namely only for $\varepsilon - \frac1n \leq z \leq \varepsilon \}$} \, .
$$
Inside $\Theta^{N+4} (fX^2-H) \supset j (\Theta^{N+4} (X^2-H) - {\rm DITCHES})$, the step DIL is PROPER. In anticipation of things to come later, we will introduce the notation
\begin{equation}
\label{eq3.23}
S_b (fX^2-H) \equiv \{\mbox{the long $\{ \ldots\ldots \}$ occurring as the RHS of $(3.22)$.}\}
\end{equation}
This $S_b (fX^2 -H)$ is a smooth $(N+4)$-manifold.

\smallskip

What makes our ditch filling step} (\ref{eq3.13}) $+$ (\ref{eq3.22}) {\it be only {\ibf partial} is that $j \mid T(n)$ is restricted to $\varepsilon - \frac1n \leq z \leq \varepsilon$ and that the same restriction applies, afterwards, to the ditch filling material in} DIL.

\medskip

\noindent 2) {\it The} (\ref{eq3.22}) {\it induces a PROPER embedding
\begin{equation}
\label{eq3.24}
S_b (fX^2-H) \overset{\mathcal J}{\longrightarrow} \Theta^{N+4} (fX^2-H) \, .
\end{equation}

Putting together} (\ref{eq3.13}), (\ref{eq3.22}) {\it and} (\ref{eq3.24}), {\it we get the following commutative diagram of PROPER embeddings}
\begin{equation}
\label{eq3.25}
\xymatrix{
\Theta^{N+4} (X^2 - H)\mbox{-DITCHES} \ar[rr]^-{j} \ar[d]^-{\rm DIL} &&\Theta^{N+4} (fX^2 - H)  \\ 
S_b (fX^2-H) \ar[urr]_-{\mathcal J}
}
\end{equation}

\medskip

\noindent 3) {\it The embedding ${\mathcal J}$ is connected by a not boundary respecting isotopy to a simple-minded diffeomorphism which enters the diagram below
\begin{equation}
\label{eq3.26}
\xymatrix{
&\Theta^{N+3} (\partial X^2) \ar[dl]_b \ar[dr]^{a} \\
S_b (fX^2-H) \ar[rr]_{{\eta, \mbox{\scriptsize simple-minded} \atop \mbox{\footnotesize \scriptsize DIFFEOMORPHISM}}}  &&\Theta^{N+4} (fX^2-H) \, .  \\
&\Sigma \, C_n^{\pm} \ar[ur]_{\alpha} \ar[ul]^{\beta}
}
\end{equation}
Here the $\Theta^{N+3} (\partial X^2) = \{$the regular neighbourhood of $\partial X^2$ (figure~{\rm 3.1.(A)})$\}$ in any of our $\partial \, \Theta^{N+4}$ from lemma~{\rm 3.1}. The embeddings $b,\beta,a,\alpha$ are the obvious ones and ${\rm Im} \, b \cap {\rm Im} \, \beta = \emptyset = {\rm Im} \, a \cap {\rm Im} \, \alpha$. The upper triangle in} (\ref{eq3.26}) {\it commutes. The $\eta$ is suggested in the figure~}3.3.

\medskip

\noindent 4) {\it When restricted to $\underset{n}{\sum} \, C_n^+$, the lower triangle in} (\ref{eq3.26}) {\it also commutes. But when restricted to $\underset{n}{\sum} \, C_n^-$, then the lower triangle in} (\ref{eq3.26}) {\it only} {\ibf commutes up to homotopy}.

\bigskip

\noindent COMMENTS. A) The DIL in (\ref{eq3.22}) is the geometric realization of the zipping in the present toy model context.

\medskip

\noindent B) Figure~3.3 should help understand the failure of the lower triangle in (\ref{eq3.26}) to commute up to PROPER homotopy. Of course, also, with a high enough $N$, PROPER homotopy means isotopy in our context. \hfill $\Box$

\bigskip

We move now to the following two smooth ($N+4)$-manifolds
$$
\Theta^{N+4} (fX^2) = \Theta^{N+4} (fX^2-H) + \{\mbox{the 2-handles} \ D^2 (\alpha \, C_n^{\pm})\} \, ,
$$
$$
S_b (fX^2) \equiv S_b (fX^2-H) + \{\mbox{the 2-handles} \ D^2 (\beta \, C_n^{\pm})\} \, .
$$

\bigskip

\noindent {\bf Lemma 3.4.} 1) {\it There is NO homeomorphism taking the form}
\begin{equation}
\label{eq3.27}
(S_b (fX^2) , \Theta^{N+3} (\partial X^2)) \longrightarrow (\Theta^{N+4} (fX^2) , \Theta^{N+3} (\partial X^2)) \, ,
\end{equation}
{\it and inducing the identity on $\Theta^{N+3}$.}

\smallskip

\noindent 2) {\it Hence there is no such diffeomorphism either, and hence diagram} (\ref{eq3.26}) {\it does NOT commute up to a PROPER homotopy.}

\bigskip

\noindent {\bf Proof.} In the figure~3.1 we can see the following disk, concentric to $D_n^2$ and smaller, namely the
$$
\delta_n^2 = D^2 (\beta \, C_n^+) \cup \left( \bigcup_{\overbrace{\mbox{\scriptsize $y \in[y U (n),y V (n)]$}}} T(n) \mid y \right) \cup D^2 (\beta \, C_n^-) \subset S_b (fX^2) \, .
\eqno (3.27.1)
$$
We denote $\Gamma_n = \partial \delta_n^2$. The space $S_b (fX^2)$ is such that
$$
\delta_n^2 \to \infty \quad {\rm in} \quad S_b (fX^2), \quad \mbox{when} \quad n \to \infty \, .
\eqno (3.27.2)
$$
In the formula, the middle $T(n)$-term lives in the DITCH, disjoined from the contribution of the $R$ part of $X^2$ (3.2) and this is what makes (3.27.2) possible. Of course, we have already $\Gamma_n \subset S_b (fX^2-H)$ and we can look then at ${\mathcal J} \, \Gamma_n \subset \Theta^{N+4} (fX^2-H) \subset \Theta^{N+4} (fX^2)$.

\smallskip

The claim is now that, any conceivable system of singular disks $d_n^2 \subset \Theta^{N+4} (fX^2)$ cobounding ${\mathcal J} \, \Gamma_n$, have to have accumulation points at finite distance. This is forced by that piece of $\partial R$ which lives beyond $x=x_{\infty}$. This proves 1) in our lemma, and the 2) follows. \hfill $\Box$

\bigskip

We did not bother to add subscripts I to our various $\Theta^3 , \Theta^4, \Theta^{N+4}$, but everything discussed so far was in the context of the so-called VARIANT~I (for our toy-model) and we move now to the superior level of VARIANT~II, where the nagging issues of the lemma~3.4 will be superseded.

\smallskip

At any of the $3^{\rm d}$ levels $\Theta^3 (X^2)$, or $\Theta^3 (X^2-H)$, or $\Theta^3 (fX^2)$, or $\Theta^3 (fX^2-H)$, we will introduce the following {\ibf critical rectangle} (which we also call the {\ibf bad rectangle})
\begin{equation}
\label{eq3.28}
{\mathcal R}_{\infty} \equiv (x = x_{\infty} , -\infty < y < + \infty , - \varepsilon \leq z \leq \varepsilon) \subset \Theta^3 \, .
\end{equation}
Notice that, for the critical rectangle ${\mathcal X}$ from (\ref{eq3.18}) and from the figure~3.2 we have at any fixed value of $y$, the equality
$$
\pi {\mathcal X} = {\mathcal R}_{\infty} \mid y \, , \ \mbox{for} \quad \Theta^{N+4} (X^2-H) \overset{\pi}{\longrightarrow} \theta^3 (X^2-H) \, .
$$
We have ${\mathcal R}_{\infty} \cap \{$the undrawable singularities of $\Theta^3 (X^2)\} = \emptyset$. With this,
\begin{equation}
\label{eq3.29}
\Theta^3 (X^2)_{\rm II} \equiv \Theta^3 (X^2) \underset{\overbrace{\mbox{\scriptsize${\mathcal R}_{\infty} = {\mathcal R}_{\infty} \times \{0 \}$}}}{\cup} ({\mathcal R}_{\infty} \times [0,\infty)) \, ,
\end{equation}
and similarly in the other contexts. Starting from here, according to the case we will have
$$
\Theta_{\rm II}^{N+4} \equiv \Theta^3 (X^2)_{\rm II} \times B^{N+1} \quad {\rm or} \quad \Theta_{\rm II}^{N+4} \equiv \Theta^4 (\Theta^3 (X^2)_{\rm II} , {\mathcal R}) \times B^N \, .
$$

The $\Theta_{\rm II}^{N+4}$ is smooth and transversally compact. We can replay now the lemma~3.3 in the context~II, using the larger, still transversally compact objects $\Theta^{N+3} (\partial X^2) \subset \Theta_{\rm II}^{N+4}$, $S_b (\ldots)_{\rm II}$, defined in the obvious way.

\bigskip

\noindent {\bf Lemma 3.5.} {\it When we go to the variant~{\rm II} of lemma~{\rm 3.2}, then {\rm 1)} in the lemma~{\rm 3.3} works just as before, while the higher analogue of the diagram} (\ref{eq3.26}) {\it from the point {\rm 2)} commutes now up to PROPER homotopy, not only for $\underset{n}{\sum} \, C_n^+$ but for $\underset{n}{\sum} \, C_n^-$ too. As a consequence of this, the variant~{\rm II} of} (\ref{eq3.27}) {\it is now a diffeomorphism}
$$
S_b (fX^2)_{\rm II} \underset{\rm DIFF}{=} \Theta^{N+4} (fX^2)_{\rm II} \, .
$$

\bigskip

At the risk of being pedantic, we will give an explicit proof of this lemma. We consider, in the context of the figure~3.1, the following closed curve, at $R$-level rather than at the $D_n^2$ level
\begin{equation}
\label{eq3.30}
\gamma_n \equiv \{ [ c_n , d_n , f_n , g_n] \, , \ \mbox{living at} \ x=x_{\infty}\} \subset R \times [-\varepsilon \leq z \leq \varepsilon] \, .
\end{equation}
This is drawn in the figure~3.1.(B), except that strictly speaking, what we see there is an embedding $\gamma_n \subset D_n^2$. But then the $\gamma_n$ clearly also makes sense at level $R$ (figure~3.1.(A)) as written in (\ref{eq3.30}). The curve $\gamma_n$ clearly makes sense at levels $\Theta^{N+4} (fX^2 - H)_{\rm I}$ and $\Theta^{N+4} (fX^2-H)_{\rm II}$ too, and in both cases we find that $\underset{n=\infty}{\lim} \, \gamma_n = \infty$. In the {\ibf case II}, making use of the ${\mathcal R}_{\infty} \times [0,\infty) \subset \partial \, \Theta^{N+4} (fX^2-H)_{\rm II}$, we can find a family of (singular) disks $\Delta_n^2 \subset \partial \, \Theta^{N+4} (fX^2 - H)_{\rm II}$ such that $\partial \Delta_n^2 = \gamma_n$ and $\underset{n=\infty}{\lim} \, \Delta_n^2 = \infty$. In the context of the variant~II of diagram (\ref{eq3.20}), of which the lower triangle is now written as
\begin{equation}
\label{eq3.31}
\xymatrix{
\partial S_b (fX^2-H)_{\rm II} \ar[rr]_-{\approx}^-{\eta}  &&\partial \, \Theta^{N+4} (fX^2-H)_{\rm II}  \\ 
&\sum C_n^{\pm} \ar[ur]_{\alpha} \ar[ul]^{\beta} \, , 
}
\end{equation}
the $\underset{n}{\sum} \, \Delta_n^2$ provide us with the PROPER homotopy modulo which the $\underset{n}{\sum} \, C_n^-$ part of this triangle commutes (while, for the $\underset{n}{\sum} \, C_n^+$ the commutativity is essentially automatic). To make this little argument work, it was necessary to send the bad rectangle ${\mathcal R}_{\infty}$ to infinity, the way we have just done it. \hfill $\Box$

\bigskip

When we will go next from the toy-model to the real life, we will eventually proceed in the style of something like the variant II. But things are then trykier. In the context of lemma~3.5, the obstruction for homotopy commutativity up to PROPER isotopy was localized inside the following rectangles
\begin{equation}
\label{eq3.32}
\{\mbox{The} \ \Delta_{\infty}^2 (n) \subset {\rm LIM} \, M_2 (f) \times [-\varepsilon , \varepsilon] \ \mbox{which cobound the various} \ \{ \gamma_n \ \mbox{at} \ x=x_{\infty} \}\} \subset {\rm LIM} \, M_2 (f) \times [-\varepsilon , \varepsilon] \, .
\end{equation}
In the real life case, there is no longer anything so simple-minded as the rectangles in (\ref{eq3.32}) to localize our obstruction. So now the key ingredient for getting the desired conclusion, in the likes of the lemma~3.5, will be the group action of $\Gamma$ itself, on the representation space, and its equivariance, i.e. we will rely now upon
$$
\boxed{\mbox{discrete symmetry with compact fundamental domain}}
$$
This ends our discussion of the toy-model and, for purely pedagogical reasons we present now a MOCK-PROOF of our main result, i.e. $\forall \, \Gamma \in {\rm QSF}$. For this purpose, let us start with a smooth closed manifold $M^p$ with $\pi_1 M^p = \Gamma$. For large enough $m$, we have $\pi_1^{\infty} (\widetilde M^p \times R^m)=0$; here our $\Gamma$ which we consider is infinite and hence $\widetilde M^p \times R^m$ is open. It follows from these things that we also have $\widetilde M^p \times R^m \in {\rm GSC}$. An explicit proof for this more or less well known fact, is to be found in \cite{30}, for instance. There is a natural $\Gamma$-action on $\widetilde M^p \times R^m$, trivial on the factor $R^m$. We can find an invariant subset $\widetilde M^p \times B^m$ on which the action of $\Gamma$ is co-compact and this comes with the obvious retraction
$$
\widetilde M^p \times R^m \overset{r}{\longrightarrow} \widetilde M^p \times B^m \, .
$$

Since $\widetilde M^p \times R^m \in {\rm GSC}$ it is certainly WGSC (in the terminology of L.~Funar and D.~Otera \cite{7}, \cite{14}), meaning that it has an exhaustion by compact, simply connected complexes
$$
K_1 \subset K_2 \subset K_3 \ldots \subset \widetilde M^p \times R^m = \bigcup_i K_i \, .
$$
For any compact $k \subset \widetilde M^p \times B^m$ there is a $K_N^i \supset k$. With this in the natural diagram below
$$
\xymatrix{
k \ar@{^{(}->}[rr]^-{i} \ar[dr]_-{j} &&K_N \, , \ar[dl]^-{r} \\ 
&\widetilde M^p \times B^m 
}
$$
we have that $i(k) \cap M_2 (r) \subset k \cap (\widetilde M^p \times \partial B^m)$, since $r M_2 (r) \subset \widetilde M^p \times \partial B^m$. Deceptively, this may look very close to the desired $\widetilde M^p \times B^m \in {\rm QSF}$. But, as long as $k$ touches $\partial (\widetilde M^p \times B^m)$, there is no permissible way to disrectangle it from $M_2 (r)$. This is as far as this naive argument can go. \hfill $\Box$

$$
\includegraphics[width=7cm]{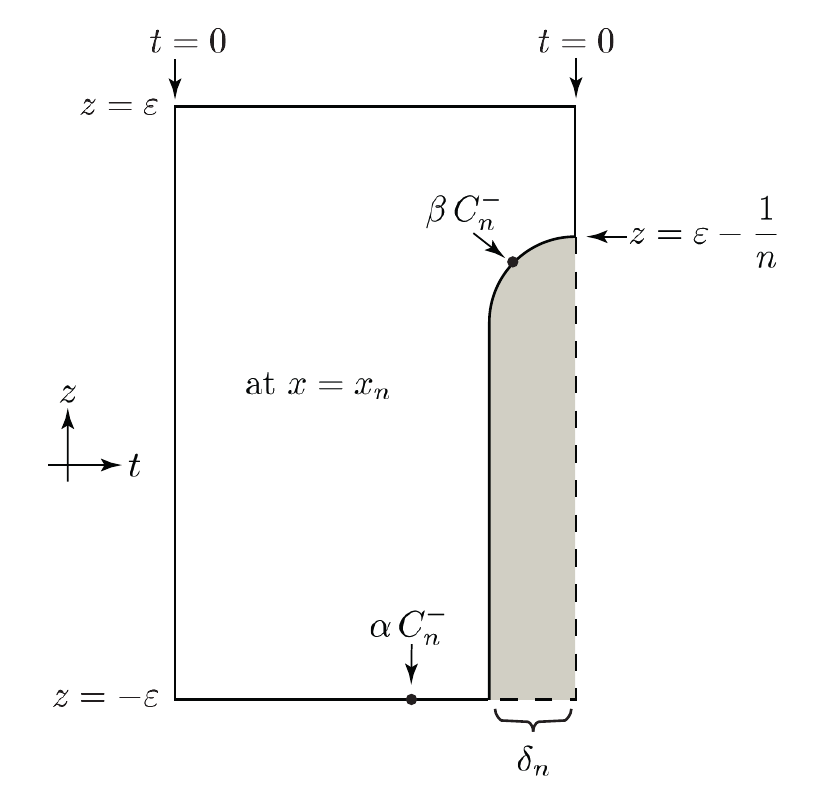}
$$
\label{fig3.3}

\centerline {\bf Figure 3.3.} 

\smallskip

\begin{quote} 
We consider now two $3^{\rm d}$ objects. To begin with the $A_u = (x_1 - k^1 \leq x \leq x_{\infty} + k^2 , z \in [-\varepsilon , \varepsilon], t \in [0,1]) - \{ x = x_{\infty} , z = \pm \, \varepsilon , t \in [0,1])$, for some small positive quantities $k^1 , k^2$. Next, let $A_b = A_u = \underset{n}{\sum} \, \{$the shaded region from the drawing, for $\forall \, x \in [x_n - \varepsilon_n , x_n + \varepsilon_n]$, $\forall \, n\}$. The depth $\delta_n$ is $\delta_n > \varepsilon_n$, with $\underset{n=\infty}{\lim} \, \delta_n = 0$. The $C_n^-$'s are now single points. There is an easy diffeomorphism $A_b \overset{\eta}{\longrightarrow} A_u$. At level $A_u$, we find that for the points marked $\alpha \, C_n^- \in A_u \mid (x=x_n)$, $\beta \, C_n^- \in A_b \mid (x=x_n)$, we have $\underset{n=\infty}{\lim} \, \alpha \, C_n^- = \infty$, $\underset{n=\infty}{\lim} \, \eta \, \beta \, C_n^- = \infty$. Moreover at the level of $\partial A_u$, one can join $\alpha \, C_n^-$ to $\eta \, \beta \, C_n^-$ by a path $\gamma_n$. But there is no way to prevent the $\gamma_n$'s from accumulation at finite distance. This spells the doom for Variant~I, in the context of (\ref{eq3.26}).
\end{quote}

\newpage

\section{The four functors $S_u , S'_u , S_b , S'_b$ and the geometric \\ high-dimensional realization of the zipping}\label{sec4}
\setcounter{equation}{0}

In the last section, in the context of the toy-model, we have introduced the distinction between the variants~I and II. But that distinction will be crucial in real life too. Let us say that everything done in sections~I, II as well as what will be done in the present and the next section will be staying constantly in the context of the VARIANT~I. It is only in the section~VI that the VARIANT~II will be extensively used. In a nutshell, here is how VARIANT~II will come about. In the context of (\ref{eq2.14}) and (\ref{eq2.15}), with $\sum (\infty)^{\wedge} \equiv \{ \sum (\infty)$ from (2.13.1)$\} \ \cup$ fins, with $p_{\infty\infty} (S) \times [-\varepsilon , \varepsilon]$ deleted$\}$, and with $\partial \, \sum (\infty)^{\wedge}$ like in (2.13.1), containing the $p_{\infty\infty} (S) \times [-\varepsilon , \varepsilon]$'s, any $2$-cell
\begin{equation}
\label{eq4.1}
(R_0 , \partial R_0) \subset \left( \sum (\infty)^{\wedge} , \partial \, \sum (\infty)^{\wedge} \right)
\end{equation}
will be called a {\ibf bad rectangle}; bad rectangles are clearly an obstruction for the injectivity of $\pi_1 (\partial \, \sum (\infty)^{\wedge}) \to \pi_1 \, \Theta^3 (fX^2)({\rm I})((2.5)) = \pi_1 \, \sum (\infty)(2.3)$. Without going now into more details, these will come in the section~VI when we will be moving to VARIANT~II, bad rectangles will be treated via the PROPER Whitehead dilatations
\begin{equation}
\label{eq4.2}
\Theta^3 (fX^2)_{\rm I} \, (\equiv {\rm our} \ \Theta^3 (fX^2) \ \mbox{from (\ref{eq2.13}), i.e. from Step III in section II)}
\end{equation}
$$ 
\Longrightarrow \Theta^3 (fX^2)_{\rm II} = \Theta^3 (fX^2)_{\rm I} \underset{\overbrace{\mbox{\scriptsize$\sum R_0 = \sum R_0 \times \{ 0 \}$}}}{\cup} \sum R_0 \times [0,\infty) \, ,
$$
which sends $\sum R_0$ to infinity, not by deletions ($=$ punctures), but by applying to it an infinite, PROPER Whitehead dilatation; this is our constant policy in this paper, for sending things at infinity.

\smallskip

This will be reflected then at all levels $(\Theta^3)'$, $\Theta^{N+4}$ a.s.o. But we will not bother to add subscripts~I (which should be all over the place in the sections~I to V) or II, until we get section VI. So, we ignore the variant~II for the time being, noting only that it should not interfere, when it will come, with what will have achieved before.

\smallskip

So we turn now to the list of ingredients (3.6) and talk about the {\ibf Holes} in real life. As a preparation, we start by listing the individual walls (with $H^{\lambda}$ meaning $\lambda$-bicollared handle)
\begin{equation}
\label{eq4.3}
W({\rm BLUE}) \, , \ W({\rm RED}) = W({\rm RED} \cap H^0) \cup W({\rm RED} - H^0)
\end{equation}
(a splitting to be soon explained below, but the notation should be self-explanatory), $W_{(\infty)} ({\rm BLACK})$.

\smallskip

In a first approximation, it is $\bigcup S_{\infty}^2 \subset \sum_1 (\infty)$ which splits each $W({\rm RED})$ into $W({\rm RED} \cap H^0) \subset \bigcup G_{\infty} H_i^0$ and $W({\rm RED} - H^0) \subset \bigcup G_{\infty} H_i^{\wedge}$. More accurately, if $H^0$ corresponds to $x < x_{\infty}$ then the border separating the $W({\rm RED} \cap H^0)_n$ from $W({\rm RED} - H^0)_n$ should be placed at $x = x_{\infty} + {\mathcal J}_n$ with ${\mathcal J}_n > 0$ converging very fast to zero. Also, each of the $W({\rm BLUE})$, $W({\rm RED} - H^0)$, $W({\rm BLACK} - H^0 - H^1)$ is parallel to a corresponding limit wall in $\sum (\infty)$. We have
\begin{equation}
\label{eq4.4}
f \, {\rm LIM} \, M_2 (f) = f \left[ \sum (W({\rm RED}) + W_{(\infty)} ({\rm BLACK}) \right] \cap \left[ \sum S_{\infty}^2 ({\rm BLUE}) \cup \sum (S^1 \times I)_{\infty} ({\rm RED}) \right] \, .
\end{equation}
In this formula we have $W({\rm RED}) \cap S_{\infty}^2 ({\rm BLUE}) \subset W({\rm RED} \cap H^0)$, in agreement with our separation at $x = x_{\infty} + {\mathcal J}_n$ above.

\smallskip

Turning now to the global $\sum_1 (\infty)$ (\ref{eq1.14}), we have

\medskip

\noindent (4.5.1) \quad $\sum_1 (\infty) \supset \{ \overline{{\rm int} \, \sum (\infty)},$ i.e. the closure of ${\rm int} \, \sum (\infty)$ inside $\sum_1 (\infty)\} = \sum (\infty) \cup \sigma_1 (\infty)$, where the first $\sum (\infty)$ is like in (2.13.1), the second line in (\ref{eq2.3}), and

\medskip

\noindent (4.5.2) \quad $\sum_1 (\infty) - \overline{{\rm int} \, \sum (\infty)} = \{$an infinite disjoined union of open $2$-cells, which will be called the {\ibf ideal Holes}$\}$.

\bigskip

With this, we introduce now an equivariant system of open 2-cells
\setcounter{equation}{5}
\begin{equation}
\label{eq4.6}
\{\mbox{normal Holes}\} \subset X^2 - M_2 (f) \, ,
\end{equation}
to be pushed afterwards into $fX^2$. Here is schematically our complete list of all the Holes of $\Theta^3 (fX^2)$ (2.13.1) and/or of $\Theta^3 (fX^2)'$ (\ref{eq2.17})
\begin{equation}
\label{eq4.7}
\mbox{Holes} = \{\mbox{normal Holes, which are open subsets}\} + \{ H(p_{\infty\infty})
\end{equation}
$$
=p_{\infty\infty} \times (-\varepsilon , \varepsilon)\mbox{'s which are closed subsets (PROPERLY embedded)}\},
$$
with the further subdivisions $\{$normal Holes$\} = \{$completely normal Holes$\} + \{$BLACK Holes$\}$ and also $\{ H(p_{\infty\infty})\} = \{\{ H(p_{\infty\infty} ({\rm proper}))\}$, to be deleted when we make the change $\Theta^3 \Rightarrow (\Theta^3)'$, and not otherwise$\} + \{\{ H(p_{\infty\infty} (S))\}$, to be {\ibf always} deleted, both for $(\Theta^3)'$ and for $\Theta^3 \}$.

\smallskip

The normal Holes are subjected, globally, to the following metric condition
\begin{equation}
\label{eq4.8}
\lim \, \{\mbox{normal Holes}\} = \{\mbox{a subset of the Ideal Holes}\} \, .
\end{equation}
In detail, here is the list of normal Holes.

\medskip

\noindent (4.9) \quad Every maximal smooth $2$-cell in $fW({\rm BLUE})$ or in $fW({\rm RED} - H^0)$ contains its unique {\ibf completely normal} Hole. Next, any $h_i^2 \subset \widetilde M^3 (\Gamma)$ generates the infinite system $\underset{\gamma}{\sum} \, H_i^2 (\gamma) \subset Y(\infty)$ and each $X^2 \mid H_i^2 (\gamma)$ contains a {\ibf unique} $W({\rm BLACK} , {\rm complete})(\gamma)$, see (1.13). For each $i$ and for all $\gamma$'s except for a finite non-void $i$-dependent subset, the $W({\rm BLACK}, {\rm complete})(\gamma)$ carries a central BLACK Hole. These are among the normal Holes in (\ref{eq4.6}). For given $i$ and for $\gamma \to \infty$ these BLACK Holes become larger and larger, so that the metrically correct condition (\ref{eq4.8}) should be satisfied for them too.

\medskip

We move now our discussion to $\Theta^{N+4} (X^2)$, which we think of, from now on, on the lines of (2.17.1) expanded into the splitting
\setcounter{equation}{9}
\begin{equation}
\label{eq4.10}
\Theta^{N+4} (X^2) = \Theta^{N+4} (X^2-H) + \sum \, D^2 (C(H))
\end{equation}
which should mean $\{ \Theta^{N+4} (X^2)$ (2.22), with {\ibf all} the Holes ({\ref{eq4.7}) {\ibf deleted} and then {\ibf compensated} for by adding $2$-handles (of dimension $N+4$), $D^2 (C(H))\}$. For the completely normal holes the $D^2 (C(H))$'s are the obvious smooth $2$-handles corresponding to the framed $\partial ({\rm Holes}) \equiv C(H)$. For the $H(p_{\infty\infty})$'s the $D^2 (C(H))$'s are like in (2.19) and in the figure~2.2, and now we will describe the $D^2 (C(H))$'s for the BLACK Holes. These are displayed in the figure~4.1.

\bigskip

\noindent ADDITIONAL EXPLANATIONS CONCERNING Figure 4.1. All the points where $M_2 (f) \cap W({\rm BLACK})$ is not, locally, a smooth line, are located inside the islands of $p_{\infty\infty}$ and $p_{\infty\infty}(S)$ and in the $S$-regions from figure~1.1. Without loss of generality, the $C(H)$'s are inside the easy region $N_{\infty}^2$. One may also assume the $C(p_{\infty\infty})$'s far from triple points. The immortal singularities in figure~4.1 are the finitely many $c$'s and the infinitely many $c^*$'s, all localized at $x=x_{\infty}$. They never touch the $C(p_{\infty\infty}) \times \left[ - \frac{\varepsilon}4 , \frac{\varepsilon}4 \right]$. These annuli are far from the DITCHES too. The immortal singularities created by the partial zipping of $W(n)$, $W(m)^*$, are not yet there.

\smallskip

For typographical simplicity's sake, the hexagon in figure~1.1 has become here a square, and only three of the $p_{\infty\infty}$'s are drawn. Notice that both $D(H(p_{\infty\infty}))$ and $D^2 (C(H({\rm BLACK})))$ have their attaching curves $C(p_{\infty\infty})$, $C(H)$ surrounding {\ibf from far away} the respective Holes. In the figure~4.1, we are at pre-zipping level, and the dotted (non-fake) lines are here only future punctures. No double points are drawn in figure~4.1. The point marked $c \in \sum (\infty)$ will become one of the finitely many immortal singularities, involving the present $W(n)$ and where this $W(n)$ is subdued, while the infinitely many where it overflows, occur inside the two $p_{\infty\infty} (S)$-islands, with the ``infinite package'' signalized in figure~4.1 cutting the dotted line at points $c^*$. We wanted the $h^1 (H(W(n)))$ ($=$ the $1$-handle dual to the $D^2 (C(H(W(n))))$ in the geometric intersection matrix of $\Theta^{N+4} (X^2)$) to be in the region contained between $\partial W(n)$ and the concentric dotted polygon, region called ``$T_{\infty}$'' in the figure, and we also want that $\sum (\infty)$ should not enter in the ${\rm int} \, T_{\infty}$, reason which explains the location of $C(H(W(n)))$ ($=$ the attaching curve of $D^2 (C(H(W(n))))$, see the drawing (A) and (C). So, now the $D^2 (C(H(W(n))))$ is treated very much like the $D^2 (C(p_{\infty\infty} ))$'s in the figure~2.2. Notice, finally, that in Morse-theory terms, every cocore $(h^1)$ comes with a much larger unstable manifold $W_u (h^1)$, which our figure also tries to suggest. End of explanations for figure~4.1. \hfill $\Box$

$$
\includegraphics[width=165mm]{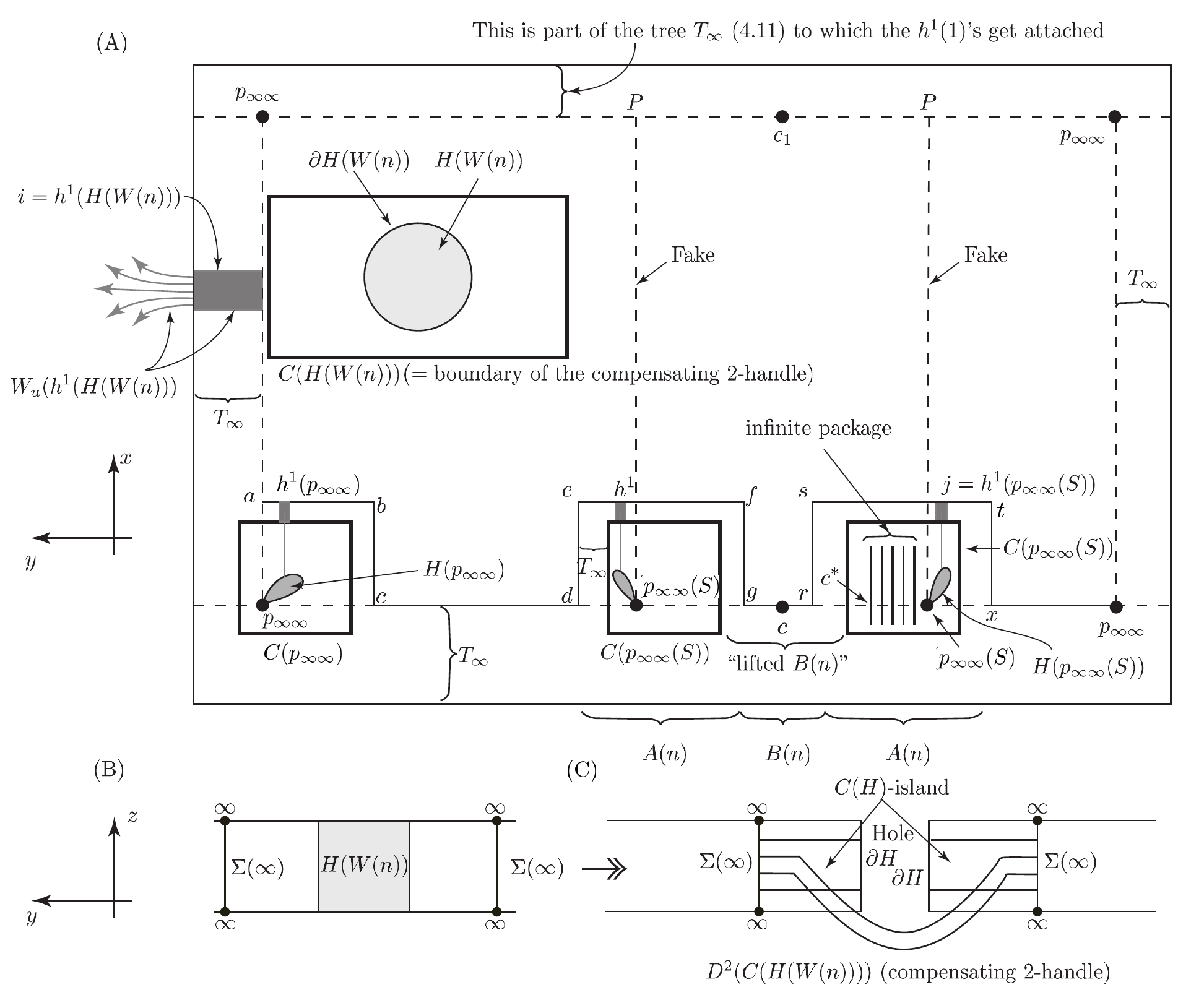}
$$
\label{fig4.1}
\centerline {\bf Figure 4.1.} 

\smallskip

\begin{quote} 
We see here a $W(n) = W({\rm BLACK} {\rm complete})$, and its BLACK Hole $H(W(n))$. Just for better visualization, the holes $p_{\infty\infty} , p_{\infty\infty}(S)$, which in real life are point-like, have been drawn as {\ibf closed} teardrops denoted $H(p_{\infty\infty})$. The various $h^1$'s in the picture are not part of (but rather attached to) $T_{\infty} =$ the tree in (4.11), while all the {\ibf rest} of $\{$the region contained between the $\partial W(n) (=$ the outer rectangle) and the closed curve $[p_{\infty\infty} ,a,b,c,d,e,f,g,c,r,s,t,x,p_{\infty\infty},p_{\infty\infty}]\} - \{$the interiors of the curves $C(p_{\infty\infty})$, $C(p_{\infty\infty} (S))$, $C(p_{\infty\infty} (S))\}$, is part of $T_{\infty}$. LEGEND: $----- =$ future $\sum (\infty)$, real or fake (in the real case the $\partial \sum (\infty)$ is eventually deleted as puncture, but not in the fake case), and we use the word ``future'' because the present figure is drawn at level $X^2$ or (\ref{eq4.10}), prior to any map for zipping thereof, {\Huge{$_-\!\!_-\!\!_-\!\!_- $}} $=$ curves $C(H)$, $C(H(p_{\infty\infty}))$, to which the $2$-handles $D^2 (C(H))$, $D^2 (C(H(p_{\infty\infty})))$ are attached, like in the drawing (C), respectively like in figure 2.2.(A), ``infinite package'' $=$ an infinite package of $\{$germs of $\{$lifted $B(m\geq n)$'s$\}\} \subset C(H(W(m)^*))$'s. When $W(m)^*$ is subdued and the present $W(n)$ is overflowing, their intersection with the dotted lines $(=\sum (\infty))$ are points of type $c^*$, duals of our present point $c$, but belonging now to $W(m)^*$.
\end{quote}

\bigskip

The $\Theta^{N+4} (X^2)$ in (\ref{eq4.10}) will be endowed with a handle body decomposition flexible enough to allow our future manipulations. This includes already all the fins (with their rims). For the time being nothing at all is sent to infinity. The handles will be denoted generically $h_i^{\lambda} (1)$. With this, we will have
\begin{equation}
\label{eq4.11}
\Theta^{N+4} (X^2) = T_{\infty} \ \mbox{(an infinite tree, thickened in dimension} \ N+4) + \underset{i}{\sum} \, h_i^1 (1) \, + 
\end{equation}
$$
\left( \underset{j}{\sum} \, h_j^2 (1)_1 + \underset{k}{\sum} \, h_k^2 (1)_2 \right) + \left( \underset{\ell}{\sum} \, h_{\ell}^3 (1)_1 + \underset{k}{\sum} \, h_k^3 (1)_2 \right) + \ldots \, ,
$$
where

\medskip

\noindent (4.11.1) \quad All the $D^2 (C(H))$ in (\ref{eq4.10}) are among the $h_j^2 (1)_1$.

\medskip

\noindent (4.11.2) \quad In the geometric intersection matrix we have $h_j^2 (1)_1 \cdot h_i^1 (1) = \delta_{ji} +$ nilpotent (easy type), and then, similarly, for all $\lambda \geq 2$ things like $h_{\alpha}^{\lambda +1} (1)_1 \cdot h_{\beta}^{\lambda} (1)_2 =$ easy id $+$ nil, too.

\medskip

\noindent (4.12) \quad As a consequence of all this, $\Theta^{N+4} (X^2)$ certainly is GSC.

\medskip

The more schematic figure~4.1.bis, which accompanies figure~4.1, should help understand the matrix $h_j^2 (1)_1 \cdot h_i^1 (1)$, in the neighbourhood of figure~4.1. 

\smallskip

We will denote now the smooth manifold

\medskip

\noindent (4.13) \quad $S'_u (\widetilde M (\Gamma)-H) = \{ S'_u \, \widetilde M (\Gamma)$ with all the Holes and all the compensating handles $D^2 (C(H))$ {\ibf deleted}$\}$ and similarly we define $\Theta^{N+4} (X^2-H)$.

\smallskip

We will also introduce the following singular object

\medskip

\noindent (4.13.1) \quad $\Theta^3 (fX^2-H)' \equiv \{ \Theta^3 (fX^2)$ with all the holes, including $p_{\infty\infty}$ (all) and also all the compensating disks removed$\}$.

\medskip

Neither the notation ``$\Theta^3 (fX^2-H)$'' nor $S_u (\widetilde M (\Gamma)-H)$ for that matter, will ever be used. With this, we have
$$
S'_u (\widetilde M (\Gamma) - H) = \Theta^4 (\Theta^3 (fX^2 - H)' , {\mathcal R}) \times B^N \, .
\eqno(4.13.2)
$$
The smooth manifold (4.13) comes endowed with an obvious framed link
\setcounter{equation}{13}
\begin{equation}
\label{eq4.14}
\sum_{\overbrace{\mbox{\footnotesize $H \in (4.7)$}}} C(H) \overset{\alpha}{\underset{{\rm PROPER \, embedding}}{-\!\!\!-\!\!\!-\!\!\!-\!\!\!-\!\!\!-\!\!\!-\!\!\!-\!\!\!-\!\!\!-\!\!\!-\!\!\!-\!\!\!-\!\!\!-\!\!\!-\!\!\!\longrightarrow}} S'_u (\widetilde M (\Gamma)-H) \, .
\end{equation}
In this formula, for $H \in \{$completely normal Holes$\}$, we have ${\rm Im} \, \alpha \subset \partial S'_u (\widetilde M (\Gamma) - H)$, while for $H \in \{$BLACK Holes$\} + \{ H(p_{\infty\infty})\}$ we have ${\rm Im} \, \alpha \subset {\rm int} \, S'_u (\widetilde M (\Gamma)-H)$. Notice the obvious reconstruction formula
$$
S'_u \, \widetilde M (\Gamma) = S'_u (\widetilde M (\Gamma)-H) + \sum_{\overbrace{\mbox{\footnotesize $\alpha \, C(H), \, H \in (4.7)$}}} D^2 (C(H)) \, .
\eqno (4.14.1)
$$
We will {\ibf not} give any meaning to ``$S_u (\widetilde M (\Gamma)-H)$'', without an upper $(')$, a symbol we will never use. Enough has been said concerning the Holes, for the time being, and referring to the same list (3.6) we now move next to the DITCHES. These ditches will be {\ibf indentations}, inside some of the thickened walls $W \times [-\varepsilon , \varepsilon] \times B^{N+1}$, making use of additional $N+1$ dimensions. We want our ditches not to interfere with our various mechanisms which concern the $p_{\infty\infty}$'s and so, at least when it comes to the regions $A(p_{\infty\infty}) \subset \Theta^3 (fX^2)'$ (figure~2.4) then our indentations will be contained inside smaller $N+1$ balls $b^{N+1}$, resting on $\partial B^{N+1}$, in such a way that
\setcounter{equation}{14}
\begin{equation}
\label{eq4.15}
b^{N+1} \subset B^{N+1} - \left(\frac12 + \varepsilon \right) B^{N+1} \, .
\end{equation}

$$
\includegraphics[width=15cm]{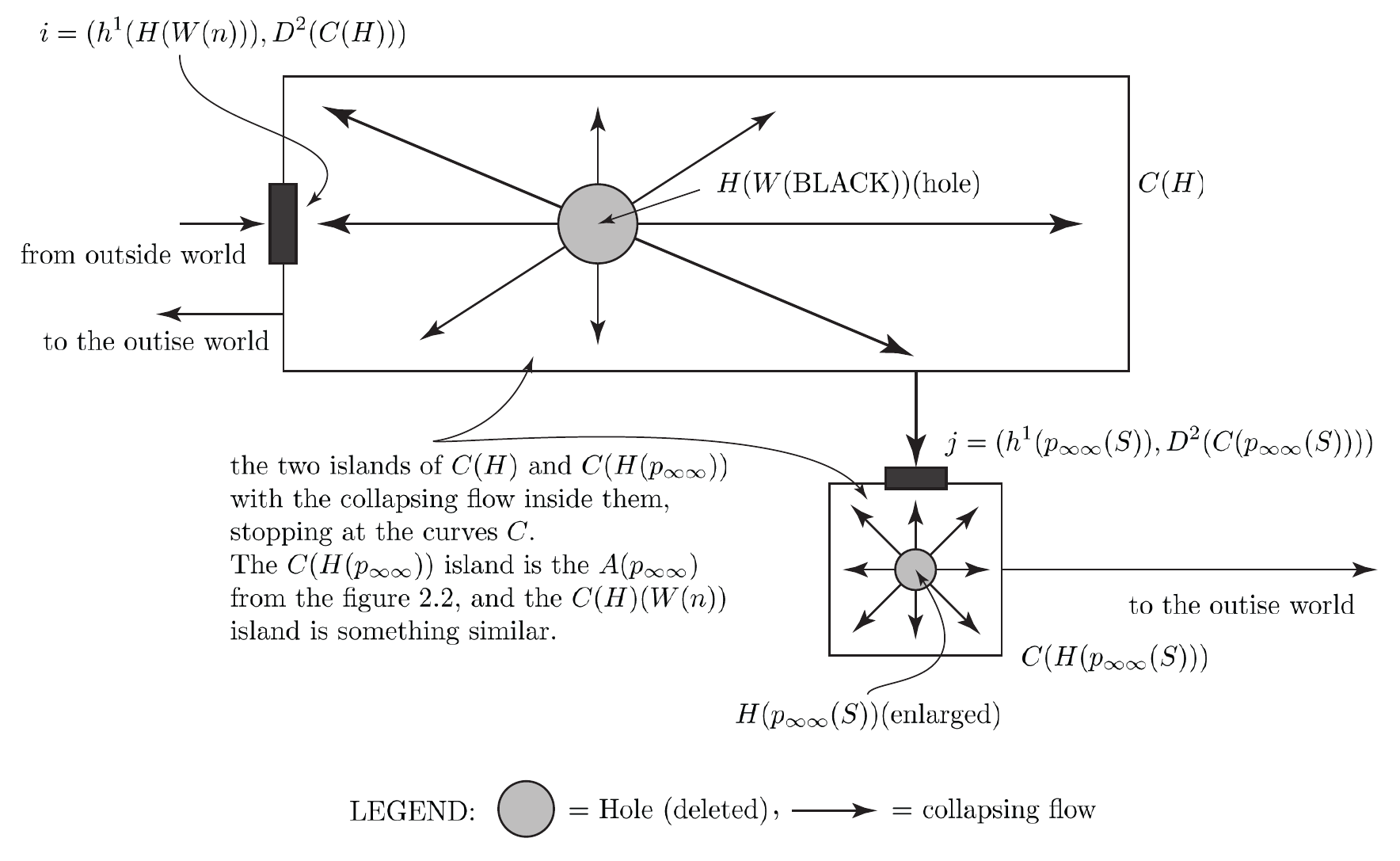}
$$
\label{fig4.1.bis}

\centerline {\bf Figure 4.1.bis} 

\smallskip

\begin{quote} 
This figure should explain the piece of the geometric intersection matrix $h_j^2 (1)_1 \cdot h_j^1 (1) = ({\rm easy}) \, {\rm id} + {\rm nil}$, coming with the figure~4.1. In view of the form of this matrix, it defines an oriented graph, and in a somehow metaphorical language, the reader is invited to think of it as the, not everywhere well-defined ``collapsing flow'' of $X^2$. The vertices of our oriented graph are called {\ibf states}, and they correspond to couples like for instance
$$
i = (h^1 (H(W(n))) , D^2 (C(H (W(n))))), \ {\rm or} \ j = (h^1 (p_{\infty\infty} (S)) , D^2 (H(p_{\infty\infty} (S)))) \, .
$$
So, the off-diagonal part of the geometric intersection matrix of $X^2$ (\ref{eq4.10}) and (4.11.2) defines an abstrat flow, our {\ibf ``collapsing flow''}. This is, of course, not defined on the $2$-cells $h^2 (1)_2$ in (\ref{eq4.11}). Here is a detail of our flow, coming with the present figure
$$
\begin{matrix} {\rm outer} \\ {\rm world} \end{matrix} \left\{ \begin{matrix} -\!\!-\!\!\!\longrightarrow \\ \longleftarrow\!\!\!-\!\!- \end{matrix} \quad i \ \overset{C(H)}{-\!\!\!-\!\!\!-\!\!\!-\!\!\!\longrightarrow} j \longrightarrow \mbox{outer world} \, .
\right.
$$
There is also an internal flow inside the islands, radiating from the corresponding Hole. Exactly one flow line hits the corresponding $i,j$ and this makes it then to the outer world. These smooth lines are the ways for the islands to leak out. We define here the islands as being whatever is inside the curves $C(H)$, $C(H (p_{\infty\infty}))$, with the {\ibf exclusion} of $i,j$ themselves. Very importantly, nothing coming from the outside world ever hits these islands. Finally, when it comes to their leaking out, see also the outer boundary of $A(p_{\infty\infty})$, in the figure~2.2, too.
\end{quote}

\bigskip

The ditches will be carved inside the
$$
\{\mbox{complementary walls}\} \times [-\varepsilon , \varepsilon] \times B^{N+1} \, ,
$$
where we define
\begin{equation}
\label{eq4.16}
\{\mbox{{\ibf complementary} walls}\} \equiv \{ W_{(\infty)} ({\rm BLACK}) \ \mbox{and} \ W({\rm RED} \cap H^0)\} \, .
\end{equation}
These walls never carry completely normal holes, only Black Holes, and $H(p_{\infty\infty})$'s. So all the walls are partitioned now into $\{$walls$\} = \{$complementary walls, which carry DITCHES$\} + \{$the non-complementary walls $W({\rm BLUE})$, $W({\rm RED} - H^0)$, which carry completely normal holes$\}$. The complementary walls are the ones cut by the limit walls $S_{\infty}^2 \cup (S^1 \times I)_{\infty}$, which create ${\rm LIM} \, M_2 (f)$, reason why they need DITCHES.

$$
\includegraphics[width=15cm]{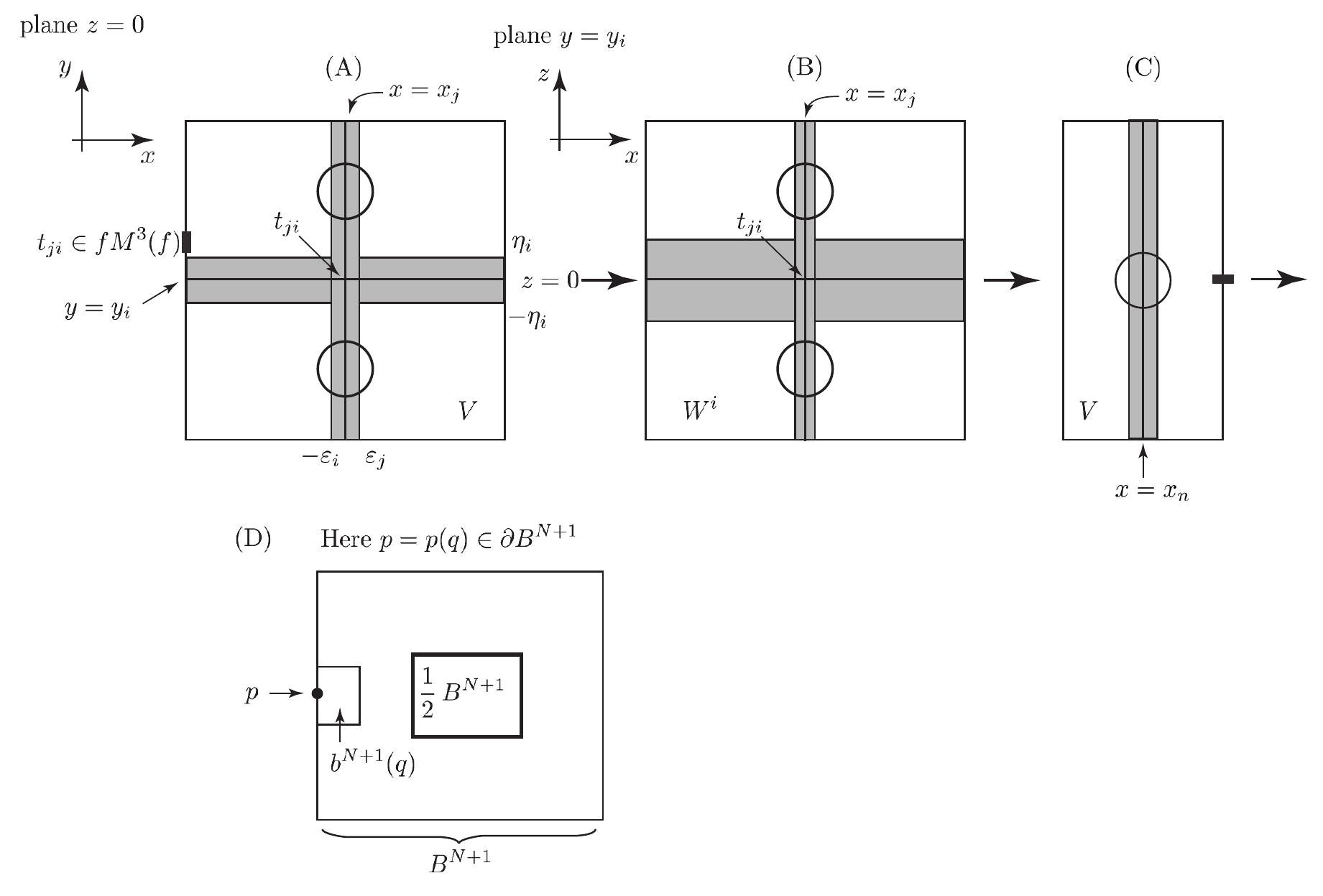}
$$
\label{fig4.2.}
\centerline {\bf Figure 4.2.} 

\smallskip

\begin{quote} 
In (A), (B), (C) we see small representative pieces of complementary walls of $X^2$, in the neighbourhood of pieces of $M_2 (f) \subset X^2$. The (A) and (B) go together, and $x=x_j$ corresponds to a $W_j ({\rm BLUE})$ which cuts transversally through them. The $y=y_i$ in (A) corresponds to the $W^i$ (displayed in (B)),  and, at least in the generic case we have here, $V = W_{(\infty)} ({\rm BLACK})$, $W^i = W({\rm RED} \cap H^0)$. The non-complementary walls, which cut transversally through the displayed complementary walls, carry Holes, signalized by the circles riding on them. But these Holes, two for each circle, live outside of the plane of our drawing. In (A) $+$ (B) we have displayed a triple point $t_{ij}$. The non smooth ramification point is left for figure~4.3.
\end{quote}

\bigskip

In these formulae, for any given $W({\rm RED})$ we have a triple splitting $W({\rm RED}) = W({\rm RED} \cap H^0)$ $\cup$ $W({\rm RED} - H^0) \cup W ({\rm RED} \cap H^0)$ with the splitting lines corresponding to $x = x_{\infty} + \zeta$, with $\zeta > 0$ i.e., with $x=x_{\infty}$, in the ``$\cap \, H^0$'' side, and with
$$
\zeta (W({\rm RED})_n) \to 0 \quad \mbox{when} \quad W({\rm RED})_n \to (S^1 \times I)_{\infty} \, .
$$
With all this, in the figures~4.2 (A, B and C) we have represented small pieces of complementary walls be they $V$ or $W^i$. The arrow drawn on the right side of each drawing, points towards ${\rm LIM} \, M_2 (f)$.

\smallskip

The double lines $fM_2 (f)$ are surrounded by {\ibf shaded zones} of controlled widths, namely the
$$
[x_j - \varepsilon_j , x_j + \varepsilon_j] \, , \ [y_i - \eta_i , y_i + \eta_i] \, , \ [-\varepsilon \leq z \leq \varepsilon] \, .
\eqno (4.17.1)
$$
For reasons of compatibility to be developed later, we will insist on the following kind of metric conditions too
$$
\varepsilon_j \ll \{ \varepsilon \ {\rm and} \ \eta_i \} \, , \quad {\rm and} \quad \lim_{n = \infty} \varepsilon_n = 0 \quad \mbox{when} \quad \lim x_n = x_{\infty} \, .
\eqno (4.17.2)
$$
Every wall will be endowed with a canonical {\ibf transversal orientation}.

\smallskip

The transversal orientations are chosen individually for each $W$. In the case $W({\rm RED})$ they look away from the corresponding handle-core. In the case of $V$ they go from $z = +\varepsilon$ to $z=-\varepsilon$ and this will determine the {\ibf partial} filling of DITCHES, happening close to $z = +\varepsilon$, like in the figure~3.2.

\smallskip

We consider now a context like in the smooth figure~4.2, but with $V+W^i$ replaced by $V \cup W^i$. This may be already so at line $X^2$ OR after some preliminary zipping. As it will be explained later, our zipping strategy will be such that any zipping $V+W^i \Rightarrow V \cup W^i$ will be performed {\ibf before} the zipping of $W_j$. With this, figure~4.2 is replaced by 4.3 and our shaded zones from (4.17.1) become

\medskip

\noindent (4.18) \quad $(V \cup W^i) \cap (x_j - \varepsilon_j \leq x \leq x_j + \varepsilon_j)$ with the transversal arcs $[-\varepsilon \leq z \leq + \varepsilon]$ and $[y_i - \eta_i \leq y \leq y_i + \eta_i]$ melted now into the $U(x_0)$ from figure~4.3.(A).

\medskip

We will denote by $q$ the generic point in a shaded zone along a double line
\setcounter{equation}{18}
\begin{equation}
\label{eq4.19}
q \in \{\mbox{shaded zones}\} \subset \mbox{complementary wall} \, .
\end{equation}
$$
\includegraphics[width=15cm]{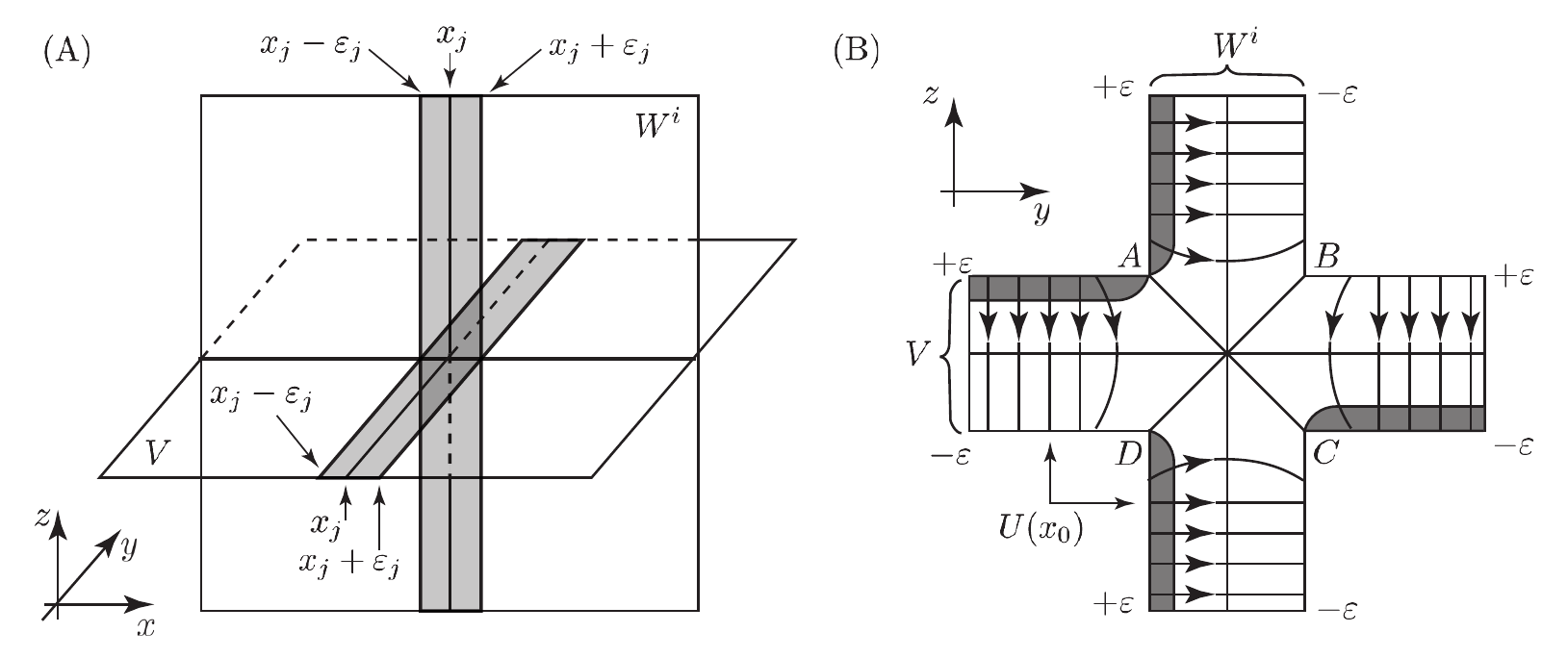}
$$
\label{fig4.3.}
\centerline {\bf Figure 4.3.} 

\smallskip

\begin{quote} 
In (A) we show what becomes of the shaded zones from figures~4.2.(A, B), when $V,W^i$ are glued together. This leave us with an ``$X$'' $\subset \, M_2 (f)$. In (B) we see a section $(x = x_0) \cap (V \cup W^i \, {\rm thickened})$. Here $x_0$ is a generic point in $[x_j - \varepsilon_j , x_j + \varepsilon_j]$. The double shading in (B) corresponds to the PARTIAL ditch filling.
\end{quote}

\bigskip

The ditches, subject of our present discussion, are far from the fins and the singularities and, assuming now also that our generic point $q$ (4.27.4) is a smooth point of $X^2$, then the local contribution of the $W = W$ (complementary) to $\Theta^{N+4} (X^2)$ (and here, as long as we ignore $p_{\infty\infty}$ we might as well write here $\Theta^{N+4} (X^2 - H)$) is
\begin{equation}
\label{eq4.20}
\begin{matrix}
\!\!\!\!\!\!\!\!\!\!\!\!\!\!\!\!\!\!\overbrace{ \qquad \qquad \qquad \quad }^{\Theta^4} \\
W^2 \times [-\varepsilon , \varepsilon] \times I \times B^N \, . \\
\underbrace{ \qquad \qquad \quad }_{\Theta^3} \qquad \!\!\!\!\underbrace{ \qquad \quad }_{B^{N+1}}
\end{matrix}
\end{equation}

For each $q \in (4.19)$ we will attach, like in the figure~4.2.(D) a smaller $q$-dependent $b^{N+1} = b^{N+1} (q)$ verifying (\ref{eq4.15}), such that when $\lim q \in {\rm LIM} \, M_2 (f)$ then $\lim ({\rm diam} \, b^{N+1} (q)) = 0$. The $b^{N+1} (q)$ is concentrated around $p(q) \in \partial B^{N+1} \cap \partial b^{N+1} (q)$, with the $p$ replacing the $t=1$ from the toy-model figure~3.2. When it comes to the DITCHES, do not mix up the ``$t=1$'' $\in \, \partial B^{N+1}$ and the $z = +\varepsilon$, end of $-\varepsilon \leq z \leq \varepsilon$, in terms of which the {\ibf PARTIAL DITCH FILLING} only concerns the portion $\varepsilon - \frac1n \leq z \leq \varepsilon$ of the complete ditch. The precise definition of the ditches will be given a bit later below.

\smallskip

But this is a good place for describing, in very broad terms, {\ibf the main idea of the present paper}, so we will open now a longuish prentice. Completely ideally, what we would like to do, would be to replace the gigantic {\ibf quotient}-space projection
$$
X^2 \underset{\rm ZIPPING}{-\!\!\!-\!\!\!-\!\!\!-\!\!\!-\!\!\!-\!\!\!-\!\!\!\twoheadrightarrow} fX^2 = X^2 / \Psi (f) \, ,
$$
connected to the $2^{\rm d}$ REPRESENTATION of $\Gamma$ from (\ref{eq1.1}) by a process, indicated for the time being only with a question mark
\begin{equation}
\label{eq4.21}
\Theta^{N+4} (X^2) \overset{?}{=\!\!\!=\!\!\!\Longrightarrow} S_u \, \widetilde M (\Gamma)
\end{equation}
consisting of a sequence of {\ibf inclusion} maps of GSC-preserving type. It is shown in \cite{36} how this plan (\ref{eq4.21}) is actually possible, under the very restrictive condition of {\ibf easyness} for the $2^{\rm d}$ representation, i.e. under the following explicit condition
\begin{equation}
\label{eq4.22}
M_2 (f) \subset X^2 \quad \mbox{and} \quad fX^2 \subset \widetilde M(\Gamma) \ \mbox{are {\ibf closed} subsets.}
\end{equation}

But, in the general case which we have to face now the assumption (\ref{eq4.22}) is hopeless. [An ironic side comment: Once we will have proved that any $\Gamma$ is QSF, then we will be able to {\it deduce} that $\Gamma$ always admits easy $2^{\rm d}$ representations too, but that is the object of still another forthcoming paper, a coda to this series.]

\smallskip

So, we are still left to try to do something on the lines of (\ref{eq4.21}), as best as we can. We replace $\Theta^{N+4} (X^2)$ by $\Theta^{N+4} (X^2)$-DITCHES and then try to develop a ditch filling process, imitating the zipping, so as to get to $S_u \, \widetilde M (\Gamma)$, or to something like it. But we want to stay PROPER and, as we have learned from our toy-model, this requires Holes which make {\ibf a partial} ditch-filling possible. The point here is that a ditch filling which is only PARTIAL, is something which can be arranged to be PROPER; see here, for instance, how in the context of the figure 3.2 we have
$$
\lim_{n=\infty} T(n) \in {\rm LIM} \, M_2 (f) \times B^{N+1} \times \{ z = +\varepsilon \} \, ,
\eqno (4.22.0)
$$
where the RHS of the formula above lives at infinity. But then, in order to get the (4.22.0), the hole $H_- (n)$ from the figure~3.1.(B) is necessary.

\smallskip

With things like the DIL in (\ref{eq3.22}) in mind, we will eventually build a process
\begin{equation}
\label{eq4.23}
\Theta^{N+4} (X^2 - H) - {\rm DITCHES} \overset{Z}{=\!\!\!=\!\!\!=\!\!\!\Longrightarrow} S'_b (\widetilde M (\Gamma) - H) \, ,
\end{equation}
when $S'_b (\widetilde M (\Gamma)-H)$ is a smooth $(N+4)$-manifold which is built by the infinite process $Z$ itself, and which is related to $S'_u (\widetilde M (\Gamma) - H)$; and, as it turns out a posteriori, it is actually diffeomorphic to it (this is an easy diffeomorphism, to be compared to the $\eta$ from (\ref{eq3.26}), and see figure~3.3 too.

\smallskip 

But, unlike what happened in the toy-model, in real life we are faced with the following two items which vastly complicate things.

\medskip

\noindent (4.24.1) \quad The $f \, {\rm LIM} \, M_2 (f)$ contains now the obvious complications of the points $p_{\infty\infty}$. In particular the $p_{\infty\infty} (S)$'s make havoc of any normal procedure for getting a locally finite $\Theta^3 (fX^2)$ or $\Theta^n (\ldots)$.

\medskip

\noindent (4.24.2) \quad We also find now a set of triple points $M^3 (f) \subset X^2 \times X^2 \times X^2$. 

\medskip

The presence of the (4.24.1), (4.24.2) make that, when we go to (\ref{eq4.23}), then the relatively simple-minded step of (partial) ditch filling, which are pure GSC-preserving inclusions (like in the case of the TOY MODEL) will have to be supplemented with other trykier steps, like for instance the {\ibf Ditch-jumping}, and others. All this will be developed in great detail in the ZIPPING-LEMMA below. Again, like in the context of the TOY MODEL, there will be now a VARIANT~I (\`a la lemmas~3.2, 3.3) and a variant~II (\`a la Claim (\ref{eq3.23})). But, for the present variant~II, nothing like (\ref{eq3.32}) is available now, and a  whole new machinery will have to be put up, the compactness lemma~4.6 and everything following it to the end of the present section and then section~VI too.

\smallskip

Coming back to (\ref{eq4.23}), it will be shown that the process $Z$ is such that the cell-complex
\setcounter{equation}{24}
\begin{equation}
\label{eq4.25}
S'_b \, \widetilde M(\Gamma) \equiv S'_b (\widetilde M (\Gamma)-H) + \sum_{\overbrace{\mbox{\footnotesize$H \in {\rm (4.7)}$}}} D^2 (C(H)) \, ,
\end{equation}
is GSC. So, we would like to show that $S'_u \, \widetilde M (\Gamma) \underset{\rm DIFF}{=} S'_b \, \widetilde M (\Gamma)$, and here the group action will be necessary. In terms of the diagram from (\ref{eq1.24}) we move to the lower level of the {\ibf compact} $M(\Gamma)$, the fundamental domain of $\Gamma \times \widetilde M (\Gamma) \to \widetilde M (\Gamma)$. We also have to move from Variant~I to Variant~II, make use of the uniformly zipping length (see \cite{29}), and compactify the $S'_u \, M(\Gamma)$. Using these kind of ingredients, it will be shown in section~VI below that
\begin{equation}
\label{eq4.26}
S'_u \, M (\Gamma)_{\rm II} \underset{\rm DIFF}{=} S'_b \, M(\Gamma)_{\rm II} \, .
\end{equation}
The $S'_u, S'_b$ are functorial and so, by taking the universal cover of (\ref{eq4.26}) what we get is
$$
S'_u \, \widetilde M (\Gamma)_{\rm II} \underset{\rm DIFF}{=} S'_b \, \widetilde M(\Gamma)_{\rm II}
$$
proving that the LHS of this formula is GSC too. From here on, there is a short argument leading to the conclusion of our GSC theorem~2.1, where of course ``$S_u \, \widetilde M(\Gamma)$'' should read now $(S_u \, \widetilde M(\Gamma))_{\rm II}$. This $S_u \, \widetilde M (\Gamma)$ comes naturally endowed with a free action
$$
\Gamma \times S_u \, \widetilde M (\Gamma) \to S_u \, \widetilde M (\Gamma)
$$
and, {\ibf if} this action would be co-compact, we could conclude that $\Gamma \in {\rm QSF}$, and this would be the end of our story. But life is not that simple, the action above is not co-compact, the $S_u \, M(\Gamma) = S_u \, \widetilde M(\Gamma)/\Gamma$ is NOT compact, the only compact space which comes in for free in our story is the very initial $M(\Gamma)$, the singular $3$-manifold presenting our $\Gamma$ as its $\pi_1 \, M(\Gamma)$. So, in order to get from $S_u \, \widetilde M (\Gamma)_{\rm II} \in {\rm GSC}$ to $\Gamma \in {\rm QSF}$, a third paper in this series will still be necessary. And then it will be $S_u \, \widetilde M (\Gamma) \in {\rm GSC}$, rather than $S'_u \, \widetilde M (\Gamma) \in {\rm GSC}$ which will have be used. [Remember that the first paper in this trilogy is \cite{29}, while the present one is the second.]

\smallskip

END of the long prentice. \hfill $\Box$

\bigskip

We make now the definition of the DITCHES completely precise, using the toy-model as a guide, when it is suitable.

\medskip

\noindent (4.27.1) \quad In the smooth case of figure~4.2, with triple points absent, and with a generic $q \in \{$shaded areas$\} \subset V + W^i$, we take
$$
{\rm DITCH} = \bigcup_q q \times [-\varepsilon,\varepsilon] \times b^{N+1} (q) \subset (V+W^i) \times [-\varepsilon,\varepsilon] \times B^{N+1} \subset \Theta^{N+4} (X^2 (-H)) \, .
$$

\noindent (4.27.2) \quad After the zipping $V + W^i \Rightarrow V \cup W^i$ which changes figure~4.2, into 4.3, we consider again
$$
q \in \{\mbox{shaded area}\} = [V \cup W^i] \cap [x_j - \varepsilon_j \leq x \leq x_j + \varepsilon_j]
$$
and the retraction $U(x_0) \overset{r}{\longrightarrow} (V \cup W^i) \mid x_0$ of possibly singular fiber, suggested in figure~4.3.(B) and which defines the transversal orientation, from $-\varepsilon$ to $\varepsilon$. With this, in the present situation, we will have
$$
{\rm DITCH} = \bigcup_{\overbrace{\mbox{\footnotesize$x_0 \in [x_j - \varepsilon_j , x_j + \varepsilon_j]$}}} \cup \, (x_0) \times b^{N+1} = \bigcup_{\overbrace{\mbox{\footnotesize$q \in \{{\rm shaded \, area, fig. \, (4.3.(A)}\}$}}} q \times \mbox{``}[-\varepsilon,\varepsilon]\mbox{''} \times b^{N+1} (q) \, .
$$

\noindent [NOTATIONAL COMMENT. The shadings in figure~4.3.(A) are of the same nature as the ones in figure~4.2, and they always mean things like
$$
\{\mbox{generic point $q$}\} \in \{\mbox{shaded zones}\} \subset \{\mbox{complementary walls}\} \, ,
$$
the points $q$ being the ones which occur in the formulae defining ditches. The shadings in figure~4.3.(B), 4.6, 5.7, 5.8, 5.9 concern
$$
\{\mbox{the previous $2^{\rm d}$ shaded zones}\} \times \left[ \varepsilon - \frac1n , \varepsilon \right] \subset {\rm DITCH} \, ,
$$
and are supposed to illustrate the PARTIAL ditch filling.]

\medskip

\noindent (4.27.3) \quad Consider now, in the context of figure~4.2 the zipping of two complementary walls $V,W^i$ and what we will say now should come as an additional refinement with respect to (4.21.1). At the level of the zipped $\Theta^3 (V \cup W^i) \times B^{N+1}$ the two pieces $\Theta^3 (V) \times B^{N+1}$ and $\Theta^3 (W^i) \times B^{N+1}$ are separated by a splitting hypersurface $S$, like in the figure~4.4. Each is supposed to dwell inside the DITCH of the other, in a completely symmetrical manner, topologically speaking, at least; metrically, it may be convenient to break this symmetry.

\medskip

\noindent (4.27.4) \quad The zipping $V + W^i \Rightarrow V \cup W^i$ takes place, as already said, before any BLUE $W_n - {\rm action}$. We can play with the function $q \to b^{N+1} (q)$, so that (4.21.1) ($+$ (4.21.3)) should merge continuously into (4.21.2). Figure~4.4 should serve here as an illustration. More explicitly Figure~4.4 goes with {\ibf the geometric realization of the zipping}, in high dimensions, as explained by lemma~4.1 and its sequels. In our figure ``high dimensions'' means just adding a $t$-coordinate to $(x,y,z)$. In the context of figure~4.4, we have $b_0^{N+1} \subsetneqq b^{N+1} \subset B^{N+1} - \frac12 \, B^{N+1}$.

\medskip

When $W_n ({\rm BLUE})$ is forgotten, figure~4.4 should illustrate, visually, the initial zipping $(V = W({\rm BLACK})) + W^j ({\rm RED}) \Rightarrow V \cup W^i$. This is like in figure~4.2 and the present $b^{N+1}$ {\ibf is} the $b^{N+1} (q)$ from that figure and from (4.27.1). The $S$ mentioned in (4.27.3) is here the RED splitting surface $S = [A,A',B',B,C,C',D,D']$, figure~4.4. Up to some metric rescaling, things should be here like in the (4.20.3).

\smallskip

[To be a bit more explicit concerning this last issue, there is here the following diffeomorphically equivalent and also metrically allowed change of view point, concerning the $W^j ({\rm RED})$: For $\vert z \vert \geq \varepsilon$ expand each of the $[A,A'], [B,B'], [C,C'], [D,D']$ all equal here to $b^{N+1}$, into a full $B^{N+1}$; see the figure~5.5.ter.]

\smallskip

When, after $V + W^j \Rightarrow V \cup W^j$ we move to the next zipping
$$
V \cup W^j + (W_n ({\rm BLUE}) - H) \Longrightarrow V \cup W^j \cup W_n \, ,
$$
then the corresponding $b^{N+1} = b^{N+1} (q)$ is the present {\it smaller} $b^{N+1}_0$. Let us say that the $b^{N+1}$ in the figure is to be used when room for $W^j$ is to be carved out of $V$, while $b_0^{N+1}$ is to be used when room for $W_n - H$ is to be carved out of $V \cup W^j$. The first, initial zipping, heals the
$$
\{\mbox{$S$-wound}\} \equiv \{\mbox{the splitting surface $S$}\} \, ,
$$
leaving, for the purpose of the second zipping the
$$
\{\mbox{$S$-scar}\} \equiv S \cap (W_n ({\rm BLUE}) - H) \times b_0^{N+1} \, , \ \mbox{shaded in the drawing.}
$$

$$
\includegraphics[width=16cm]{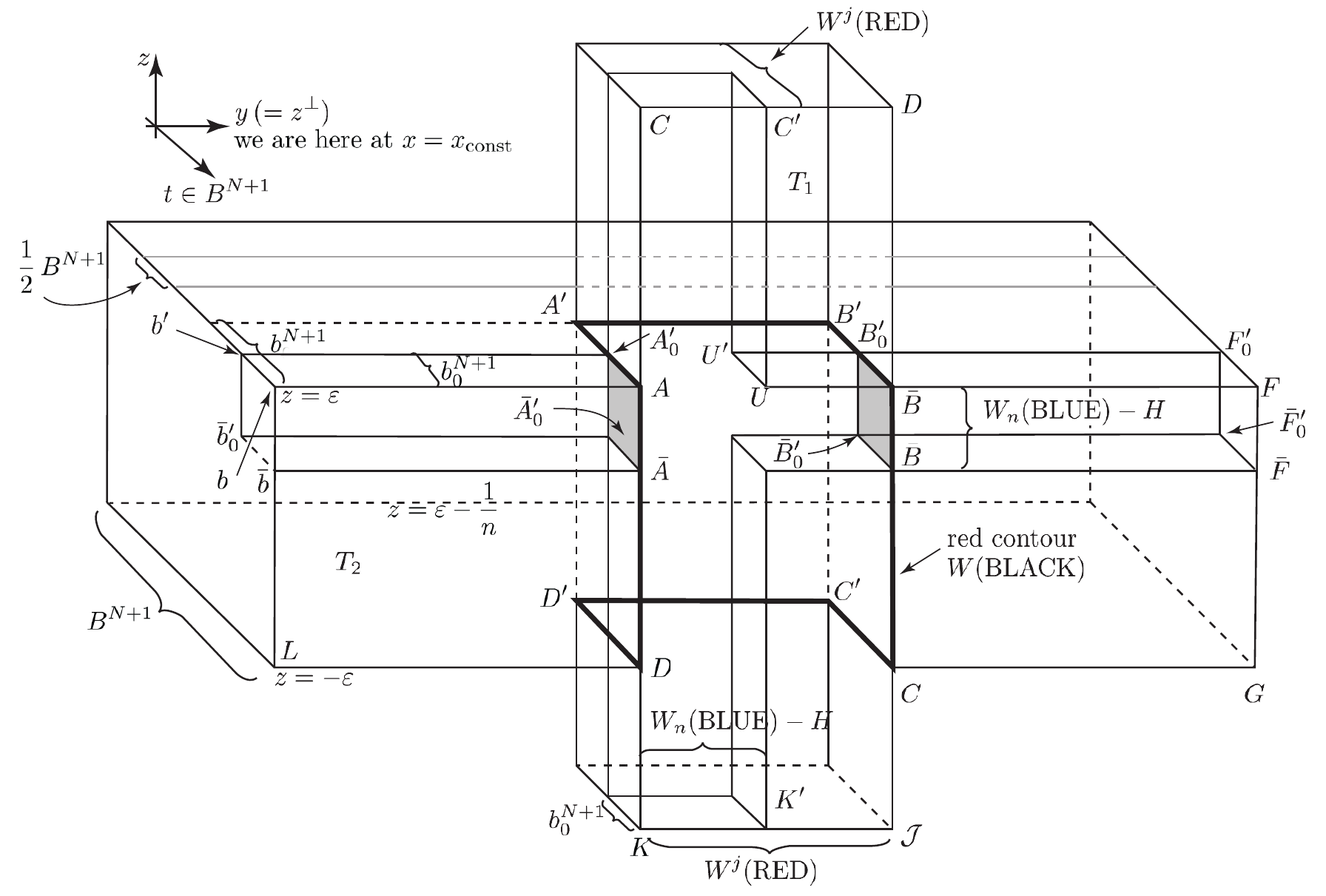}
$$
\label{fig4.4.}

\noindent LEGEND: The RED contour $[A,A',B',B,C,C',D,D']$ is the boundary of a hyper-surface $S$ splitting, in high dimensions, $W({\rm BLACK})$ from $W^j({\rm RED})$.

\medskip

\centerline {\bf Figure 4.4.} 

\smallskip

\begin{quote} 
Illustration for the formulae (4.21.1) to (4.27.4) and for the filling of the ditches. The splitting surface $S$, splits the $\{$thickened $W({\rm BLACK}) \cup W^j ({\rm RED})\}$ into $T_1 \cup T_2$.
\end{quote}

\bigskip

We use the generic notation $h(3)$ for the system of handles corresponding to the ditch-filling material. With this, at the end of the (geometric realization of the) zipping, we will find, inside the grand geometric intersection matrix
\setcounter{equation}{27}
\begin{equation}
\label{eq4.28}
\partial h^2 (3) (\mbox{second zipping}) \cdot \delta h^1 (3)(\mbox{first zipping}) \ne 0 \, ,
\end{equation}
compatibly with GSC.

\smallskip

The figure~4.5 to which we turn next, should help vizualizing the high dimensional geometry of the ditches, in the neighbourhood of the points $p_{\infty\infty}$. Here are EXPLANATIONS for this figure. Both figures 4.5.(A) and (B) refer to the contribution of a $W_{(\infty)} ({\rm BLACK})$ to the $S'_u \, \widetilde M (\Gamma)$, $S'_u \, (\widetilde M (\Gamma) - H)$, $S'_b \, (\widetilde M (\Gamma) - H)$. The coordinate system $(x,y,z)$ is here the same as in the figures~1.1.(B), (C). In 4.5.(A) we see
$$
\left\{ A(p_{\infty\infty}) \times \left( \frac12 + \varepsilon_0 \right) B^{N+1} , \, \mbox{site of the action (2.23)}\right\} \subset \bigl\{ W_{\infty} ({\rm BLACK}) \times B^{N+1} 
$$
$$
- \, \{\mbox{the indentation $b^{N+1}$, site of the zipping action}\}\bigl\} \, ,
$$

\medskip

\noindent making the action $S'_u \Rightarrow S_u$ in (\ref{eq2.23}) and the zipping, independent of each other.

$$
\includegraphics[width=17cm]{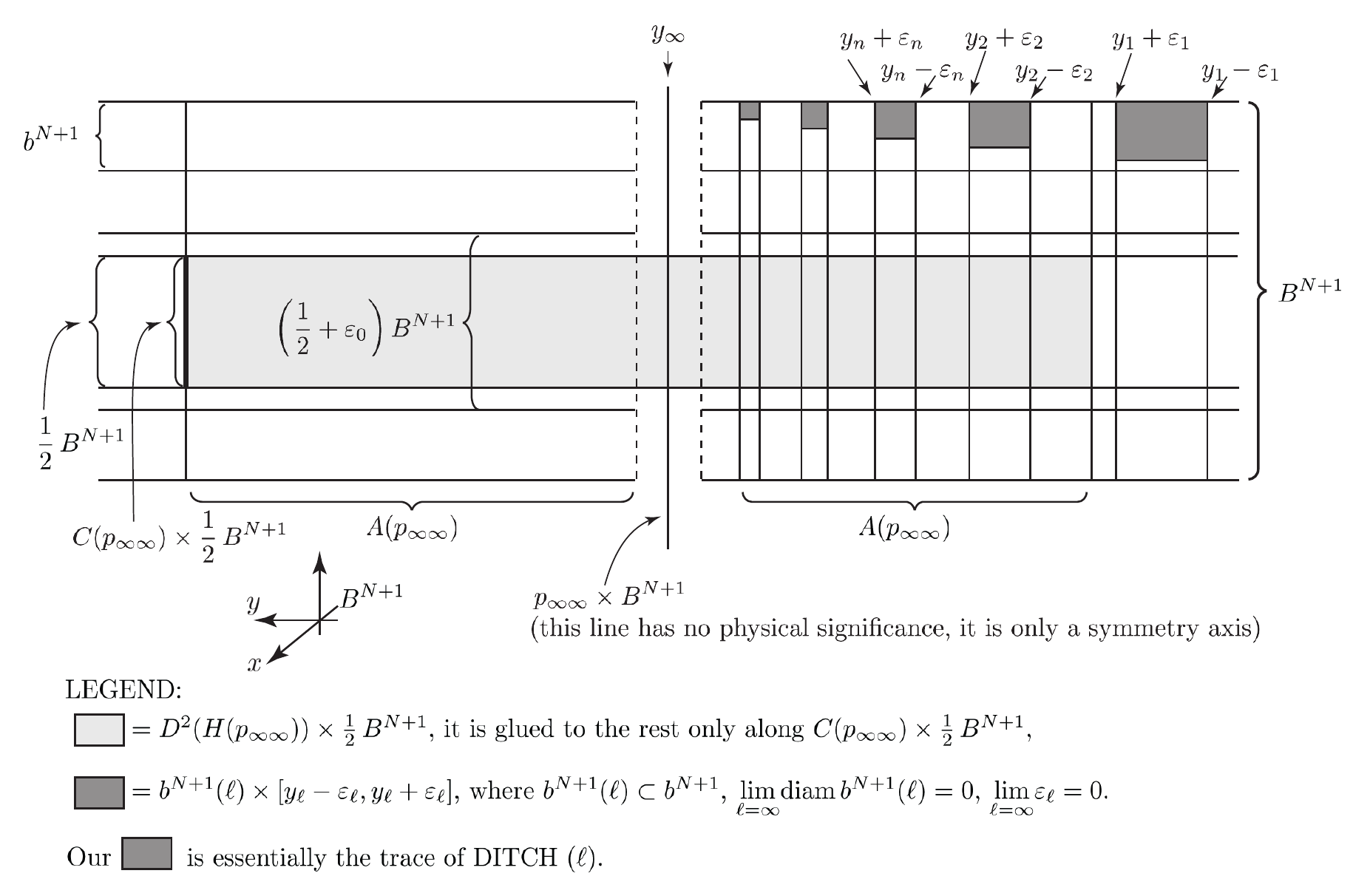}
$$
\label{fig4.5.(A)}
\centerline {\bf Figure 4.5.(A)} 

\smallskip

\begin{quote} 
This figure is a slice through $W_{\infty} ({\rm BLACK}) \times [-\varepsilon \leq z \leq \varepsilon] \times B^{N+1}$ at $x=x_{\infty}$, $z=0$. It may concern a $p_{\infty\infty} ({\rm proper})$ or a $p_{\infty\infty} (S)$. We are here in the same context as in figure~2.2, but now in higher dimensions.
\end{quote}

$$
\includegraphics[width=16cm]{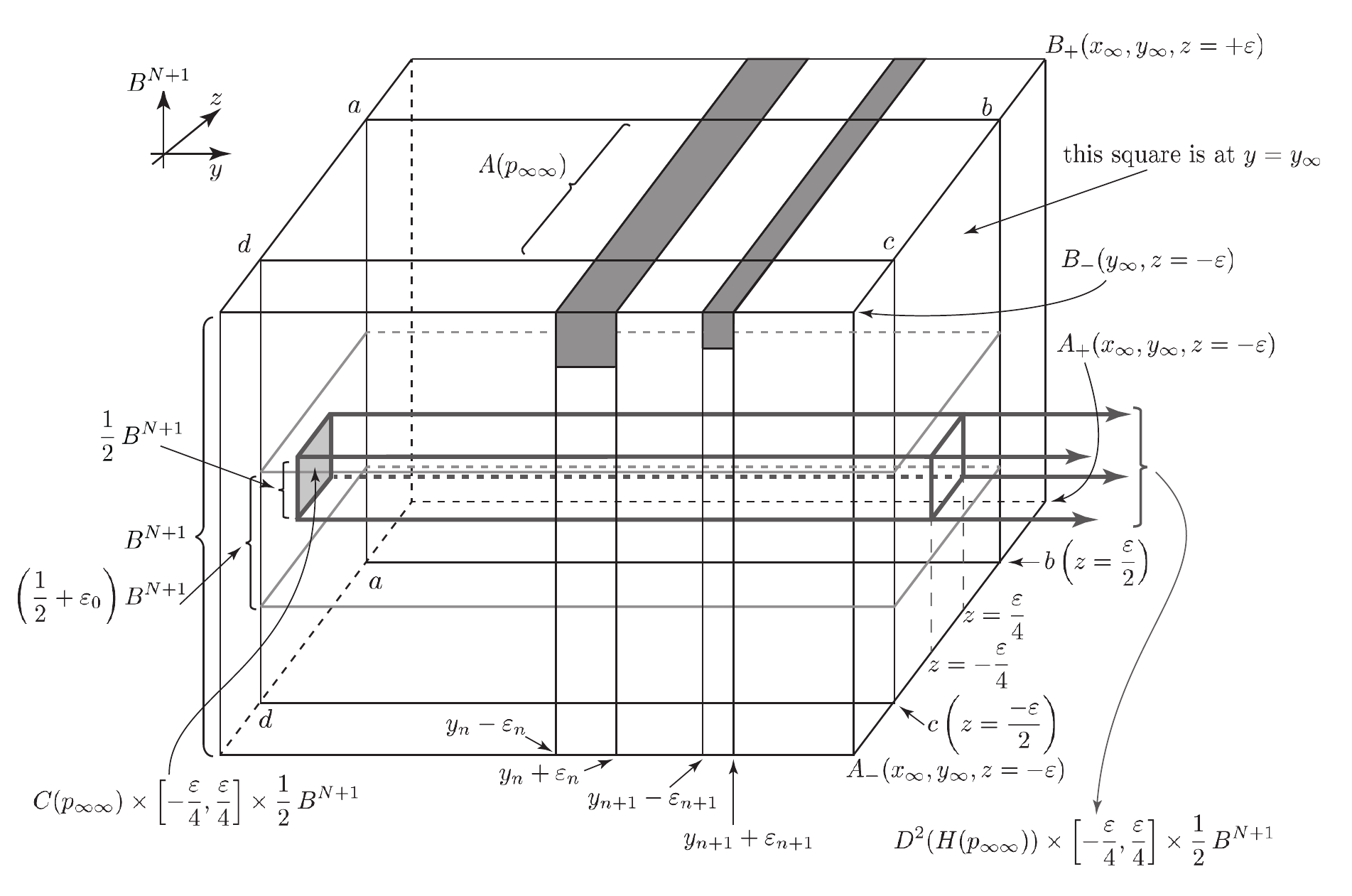}
$$
\label{fig4.5.(B)}
\centerline {\bf Figure 4.5.(B)} 

\smallskip

\begin{quote} 
The $[a,b,c,d]$ is here like in the figure~2.2. The $[A_- , A_+ , B_- , B_+]$ lives at infinity.
\end{quote}

\bigskip

In (B), since we are at the level $S'_{\varepsilon} (\widetilde M (\Gamma)-H)$ the square $[A_- , A_+ , B_- , B_+]$ has actually been deleted, but our (B) is supposed to continue, with less details beyond $y=y_{\infty}$ to $y > y_{\infty}$; so does also the 2-handle $D^2 (H(p_{\infty\infty})) \times \left[ - \frac{\varepsilon}4 , \frac{\varepsilon}4 \right] \times \frac12 \, B^{N+1}$. Notice that, at the end of the collapse $+$ dilatation in (2.22), when the dust has settled then, in the case when $p_{\infty\infty} = p_{\infty\infty} ({\rm proper})$, what we get locally is a diffeomorphic model for $S_u \, \widetilde M (\Gamma)$, where instead of putting back the whole of $[A_- , A_+ , B_- , B_+] = p_{\infty\infty} \times [-\varepsilon \leq z \leq \varepsilon] \times B^{N+1}$ we only put back $p_{\infty\infty} \times [-\varepsilon < z(c) \leq z \leq z(d) < \varepsilon] \times \left( \frac12 + \varepsilon \right) B^{N+1}$, which leaves out the $p_{\infty\infty} \times [-\varepsilon , \varepsilon] \times b^{N+1}$. But this alone is still not of much help when dealing with the sting of $p_{\infty\infty}$. To see how we will actually proceed, watch the zipping lemma~4.1 and its various sequels in this section, up to and including the ``Final arguments''.

\smallskip

For $S'_u \, \widetilde M (\Gamma)$, the contribution of each $W^i ({\rm RED})$ adds to what we can see already in (B), the full contribution
$$
[y_i - \varepsilon_i \leq y \leq y_i + \varepsilon_i] \cap \{ \vert z \vert \geq \varepsilon \} \times B^{N+1} \, . \eqno (*_1)
$$
But then, finally, both (A) and (B) are at $x=x_{\infty}$, hence the contribution of $\partial \, \sum (\infty)$ is to be deleted. Forgetting about the $(*_1)$ above, what this deletion means for (B) is that the whole lateral surface $(z=\pm \, \varepsilon) \times B^{N+1}$ is to go, at $x=x_{\infty}$.

\smallskip

What our figure~4.5 displays is the following fact. At the level of $\Theta^{N+4} (X^2-H) - {\rm DITCHES}$ or $S'_b (\widetilde M (\Gamma) - H)$, when localized at $p_{\infty\infty} ({\rm proper \, or} \ S) \subset W_{(\infty)} ({\rm BLACK}) = \{$let us say the $V$ in model (2.8.i)$\}$, the thickened $V$ contains the region

\bigskip

$$
[A(p_{\infty\infty}) \times B^{N+1} - \{\mbox{indentations}\}] \supset \Bigl\{ C(H(p_{\infty\infty})) \times \frac12 \, B^{N+1} \, , \eqno (4.28.1)
$$
$$
\mbox{attaching zone of the $2$-handle} \ D^2 (H(p_{\infty\infty})) \times \frac12 B^{N+1} \Bigl\} \, ,
$$
where the {\ibf indentations} ($=$ DITCHES) make room for the infinite set of zippings of $V$. In the geometric realization of the zipping, this means ditch-filling operations, involving $W({\rm BLUE})_j \times b^{N+1} (j)$, $W({\rm RED})^i \times b^{N+1} (i)$, all living in the DITCH. Here ${\rm diam} \, b^{N+1} (\ell)$ goes very fast to zero, when $\ell \to \infty$, so that we have a diffeomorphism
$$
A(p_{\infty\infty}) \times B^{N+1} - \{\mbox{indentations}\} \underset{\rm DIFF}{=} A(p_{\infty\infty}) \times B^{N+1} \, .
$$

The indentations protect our $C(H(p_{\infty\infty})) \times \frac12 \, B^{N+1}$ from contacts with the zipping process $Z$.

\smallskip

Now, a similar protection is not necessary for $C(H)$ (figure~4.1), since as we can see in the figure~4.1, $M_2 (f) \cap C(H) = \emptyset$ (to see this, compare actually the figures~4.1 and 1.1.(A)). Details like $[g,c,r]$ in figure~4.1, fit now into the
$$
\bigl\{\mbox{doubly shaded indentations of a figure like 4.5, for $p_{\infty\infty} (S) \in W(m)^*$, $m < n$,}
$$
$$
\mbox{where $W(m)^*$ overflows, not represented in figure 4.1}\bigl\}.
$$

We are ready now to give the {\ibf proof of the CLAIM (2.26)}.

\smallskip

We will choose a smooth function
$$
q \in D^2 \overset{\varepsilon_0}{-\!\!\!-\!\!\!-\!\!\!\longrightarrow} \left[ 0 , \frac\varepsilon2 \right] \, , \eqno (*_2)
$$
with $\varepsilon_0 (q) = \varepsilon_0 (d (q, p_{\infty\infty} = {\rm center \, of} \ D^2))$, increasing from $\varepsilon_0 = 0$ on $\partial D^2$ to $\varepsilon_0 (p_{\infty\infty}) = \frac\varepsilon2$. Our main aim remains to keep the zipping action and the action in (\ref{eq2.23}) far apart from each other. Let $r (B^{N+1}) \equiv$ radius of $B^{N+1}$.

\smallskip

For this we will also ask for the following condition
$$
\Vert b^{N+1} (q) \Vert + \left\Vert \left( \frac12 + \varepsilon_0 (q) \right) B^{N+1} \right\Vert < r(B^{N+1}) \, , \eqno (*_3)
$$
where $b^{N+1} (q)$ is like in the DITCH-defining formulae (4.27.1) and (4.27.2), inside $A \cap D^2 = A \cap (D^2 - p_{\infty\infty})$. This will make that, inside $B^{N+1}$ we will find that
$$
b^{N+1} \cap \left( \frac12 + \varepsilon_0 \right) B^{N+1} = \emptyset \, ; \eqno (*_4)
$$
compare this with lemma~4.4 below. Coming back to lemma 2.2 with the factor $\left( \frac12 + \varepsilon_0 (q)\right) B^{N+1}$ becoming now $q$-dependent, what the collapse $+$ dilatation in (\ref{eq2.23}) do, in the neighbourhood of $p_{\infty\infty} ({\rm proper})$, is to recreate the correct
$$
\Theta^3 (fX^2) \times I \times B^N \subset S_u \, \widetilde M (\Gamma) \, .
$$
\hfill $\Box$

\bigskip

\noindent GEOMETRIC REALIZATION OF THE ZIPPING, in high dimensions. This discussion will occupy all the rest of this section, of which it is the core. There will be actually two levels for the discussion in question, {\ibf upstairs}, at level $\widetilde M (\Gamma)$ and {\ibf downstairs}, at level $M(\Gamma)$. For the time being we will function upstairs and it will be signalized when we will move downstairs.

\smallskip

The central piece of the present discussion will be the ZIPPING LEMMA~4.1, for which we give now some PRELIMINARIES. We will start with the following object, which is a cell-complex
\begin{equation}
\label{eq4.29}
\left[ \Theta^{N+4} (X^2 - H) - {\rm DITCHES} \right] + \sum_{H \in (4.7)} D^2 (C(H)) \, .
\end{equation}
For the time being we will unroll things upstairs.

\smallskip

The $S_0 \equiv [\Theta^{N+4} (X^2-H) - {\rm DITCHES}]$, which occurs in (\ref{eq4.29}) is a smooth $(N+4)$-manifold related to $\Theta^{N+4} (X^2-H)$ by a simple-minded diffeomorphism. Remember that ``$-H$'' means with all Holes, including the $H(p_{\infty\infty})$'s deleted, and with the compensating $D^2 (C(H))$'s deleted too.

\smallskip

The (\ref{eq4.29}) comes with an initial cell-decomposition $h(1)$ like in (\ref{eq4.11}), which is GSC (see (4.12)). Starting from $S_0$ the lemma~4.1 below will develop an infinite process $Z$ of successive steps and here ``$Z$'' stands for ``zipping'', our infinite process {\ibf is} actually the geometric realization of zipping, in high dimensions. During the process $Z$ (see (4.30)) some of the steps involved are inclusions $S_n \subset S_{n+1}$. These steps will add new handles, coming in addition to the system $h(1)$ from (\ref{eq4.11}), and here one should remember that $h(1)$ contains already all the fins, with their rims. Here is the list of the handles which are added:

\medskip

\noindent (4.29.1) \quad There will also be a system of handles $h(2)$, added during a subsequence $S_{n_1} \subset S_{n_1 + 1}$, $S_{n_2} \subset S_{n_2 + 1} , \ldots$ of inclusion steps $S_n \subset S_{n+1}$, and these will send to infinity, via a PROPER Whitehead dilatation, the $\partial \, \sum (\infty)^{\wedge} \supset \{$rims of fins$\}$, see (\ref{eq2.15}).

\medskip

\noindent (4.29.2) \quad Then there will be the main bulk of inclusions $S_n \subset S_{n+1}$ of the process $Z$ (4.30). Each of these inclusions consists of two parts. First, analogously to some of the things which happen during (\ref{eq3.13}) in the context of the toy-model, pieces of our $S_0$ are being send ``in the ditch'', meaning here, of course, in the interiors of the DITCH.

\medskip

This part is, essentially, an isotopy which bends $S_0$, so as to bring it to enter the ditch, via parts of it. Then, comes a red inclusion map which adds more stuff, binding now $S_0$, inside the ditch, to $\partial \, {\rm DITCH} \cap \partial S_0$. This means throwing in a new system of handles (mainly of indices $\lambda \geq 1$) and/or cells. We call this system $h(3)$. We will call this the ``ditch-filling material''.

\smallskip

VERY IMPORTANTLY, like in the context of figure~3.2 from the TOY MODEL, as soon as we are in the dangerous zone $\overline N_{\infty}$, the {\ibf ditch filling is only partial}. This means exactly the following, remember: in terms of formulae (4.27.1), (4.27.2) (which is to be supplemented by the figure~4.3), when a piece of $S_0$ ($+$ the corresponding ditch-filling material $h(3)$) goes into the ditch, this occupies only the part $\left[ \varepsilon - \frac1n \leq z \leq \varepsilon \right]$. Here is one way of conceiving our operations of sending parts of $S_0$ into the DITCH. To begin with think of the standard zipping operation $X^2 \to fX^2$, actually restricted to $X^2 - H \to fX^2 - H$. Locally speaking, there are here two possible cases

\smallskip

i) zipping together two complementary walls $W_{(\infty)} ({\rm BLACK})$ and $W({\rm RED} \cap H^0)$,

\smallskip

ii) zipping together the non-complementary $W_1 = W({\rm BLUE})$ or $W({\rm RED} - H^0)$ (or rather the $W_1 - H$) to $W_2$ (complementary).

\smallskip

Going now from dimension $2$ to $N+4$, {\it ideally} (a word to be explained later on), in case i) each of the two complementary walls involved, goes inside the ditch of the other like in (4.27.3) and/or like in figure~4.4 and, for i), this is true without further qualifications, i.e. the word ``ideally'' may be dropped here. In the case ii), the $W_1 (\mbox{non complementary}) -H$ goes now into the ditch of the $W_2$ (complementary). All this mimicks the zipping in high $d$, in the sense that, in high $d$ the result should look, topologically, like $W_1 \cup W_2$ thickened in dimension $N+4$. Put it also another way, what we just said means an inclusion map
$$
\{\mbox{parts of $S_0$, going in the DITCH}\} \subset \bigl\{\mbox{a shaded region of the complementary $W$}\} \times 
$$
$$
[-\varepsilon \leq z \leq \varepsilon] \times (b^{N+1} \subset B^{N+1}) = \mbox{DITCH}\bigl\} \, ,
$$
and the last piece, between the $\{ \ldots \}$, makes sense both inside $S_0$ and inside the $S'_u (\widetilde M (\Gamma) - H)$.

\smallskip

This should also serve for defining the map $j$ from $S_0$ to $S'_u (\widetilde M (\Gamma)-H)$ occurring in (\ref{eq4.33}) below. We will come back to the word ``ideally'' later on.

\smallskip

Coming back to our successive handle-system $h(\varepsilon)$, with $1 \leq \varepsilon \leq 3$, for each individual $\varepsilon$ we will have that $\partial h^2 (\varepsilon) \cdot \delta h^1 (\varepsilon) =$ easy id $+$ nil. Moreover, we will also find that when $\varepsilon' < \varepsilon''$, then $\partial h^2 (\varepsilon') \cdot \delta h^1 (\varepsilon'') = 0$, while generally speaking we will find that $\partial h^2 (\varepsilon'') \cdot \delta h^1 (\varepsilon') \ne 0$. With these things, there cannot be any violation of GSC coming from our $h(\varepsilon)$'s.

$$
\includegraphics[width=140mm]{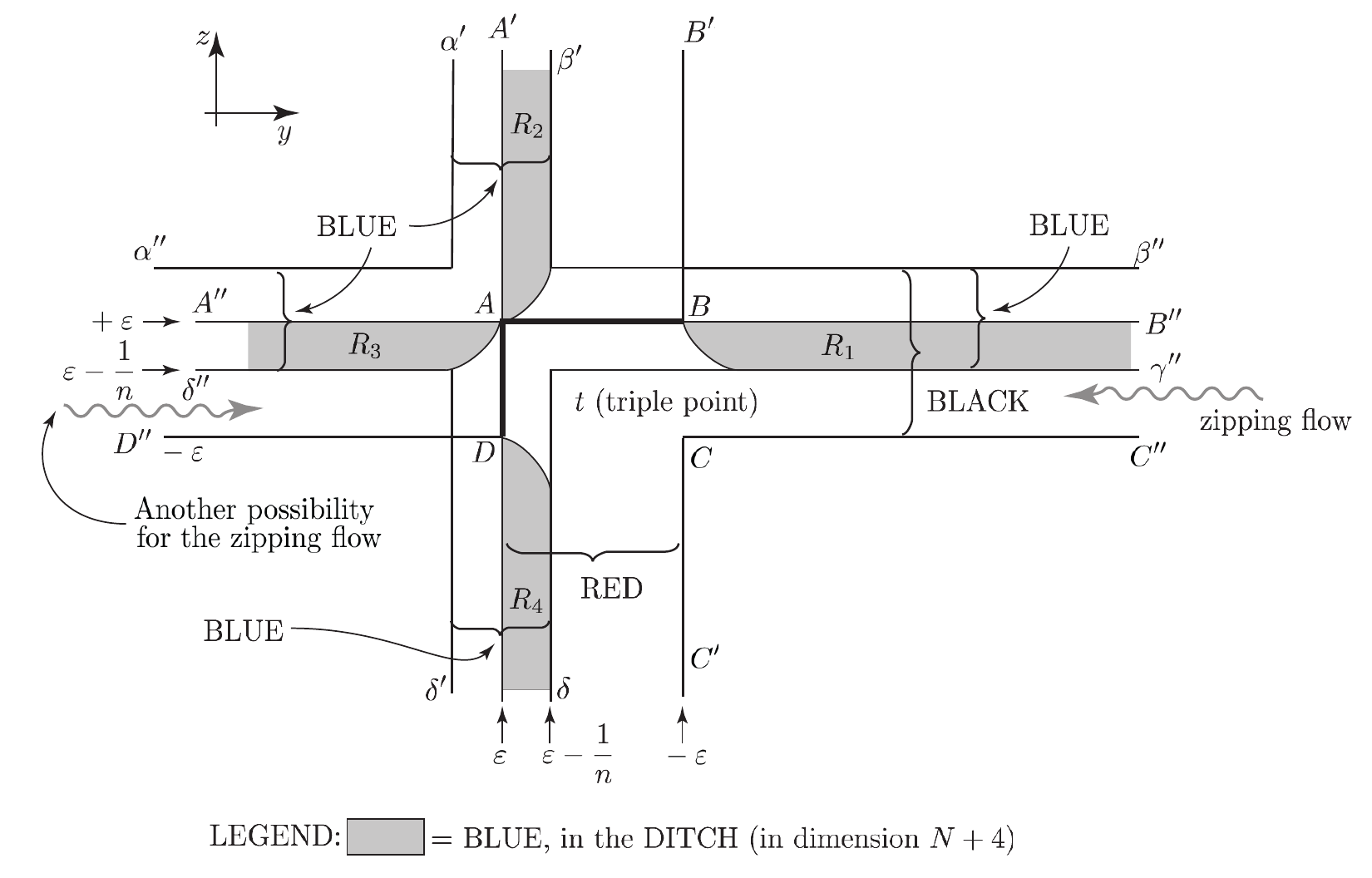}
$$
\label{fig4.6}
\centerline {\bf Figure 4.6.}

\smallskip

\begin{quote} 
The $[AB]$, $[AD]$ are so-called {\ibf special BLUE ${\bf 1}$-handles}. The greek cross stands for $W_n ({\rm BLUE}) - H$ while the latin cross stands for $W({\rm BLACK}) (\mbox{horizontal}) \cup W({\rm RED}) ({\rm vertical})$, already zipped together. The shaded areas of the greek cross live in the $\{{\rm DITCHES}\} \subset (W({\rm BLACK}) \cup W ({\rm RED})) \times [-\varepsilon , \varepsilon] \times b_0^{N+1}$ (see figure~4.4). There are {\ibf no DITCHES} at the level of the square $[ABCD]$, where the triple point $t \in M^3 (f)$ lives. Such DITCHES, occupied by the $[AB] + [AD]$, if they would exist, would come into conflict with the {\ibf metrizability} of the space $S'_b$ introduced below. In our figure we are at $x=x_n$, with a coordinate axis $x$ going transversally through the plane of our drawing, towards $x=x_{\infty}$, where $\sum (\infty)$ lives.

\smallskip

The COLOURS marked in the figure refer to the local models from (2.8). Specifically, the ``RED'' and ``BLACK'' are like in the paradigmatic case (2.8.i). Other combinations are possible for them.

\smallskip

Notice that the whole BLUE cross $[\delta'' , \alpha'' , \alpha' , \beta' , \beta'' , \gamma'' , \delta , \delta']$ (and beyond it, in the plane of $W_n ({\rm BLUE})$, when it lives, are Holes), collapses into $\{$the two special BLUE $1$-handles $[AB]+[CD]\} \cup \{$the four shaded BLUE areas, which during $Z$ will go into the DITCH, namely the $R_{\varepsilon \leq 4}$'s$\}$.

\smallskip

The whole point here is that, as explained above, the special BLUE $1$-handles {\ibf cannot} be sent into the DITCH. Compare this with our term ``ideally'' from above.
\end{quote}

\bigskip

\noindent VERY IMPORTANT REMARK. In the jargon of section~III, in all the present and also the next section, we will be constantly in the context of the VARIANT~I. Later on {\ibf after} all the movements of the ZIPPING LEMMA below will have been performed, then we will move to the context of the VARIANT~II. This will only effectively happen in the section~VI.

\smallskip

Before going on, here are some EXPLANATIONS CONCERNING the NON-METRIZABILITY issue mentioned above, in connection with figure~4.1. Consider a figure like~3.2, but now, with the $T(n)$'s extending all the way from $-\varepsilon$ to $+\varepsilon$, for all $n$'s. They would then accumulate on the segments $[\rho_- , \rho_+]$ from figure 3.2.(A), of which only the extreme points $\rho_{\pm}$ live at infinity, and the open arc $(\rho_- , \rho_+)$ at finite distance. Something similar would happen, if we would use the ditches for housing the special BLUE $1$-handles. Imagine now that the $Z$-process (see (\ref{eq4.30})) would make use of this. If we would still insist that $\underset{n=\infty}{\lim} S_n + \sum \, D^2 (C(H))$ be GSC, as we shall certainly do, this would oblige us to use for $\lim S_n$ a weak topology (consistent with the topology of the individual $S_n$) which would no longer be metrizable, the all-important $S'_u (M(\Gamma) - H)$ (and/or the $S'_u (\widetilde M (\Gamma)  - H)$) from section~VI would no longer be a DIFF, smooth manifold and our technology for proving anything like the GSC theorem~2.3 would be in shambles. The cure for this serious disease consists in: 

\medskip

$\alpha$) Use the {\ibf partial} ditch filling when outside the triple points (inside the difficult regions $\overline N_{\infty}^2$).

\medskip

$\beta$) Use the not yet defined ditch-jumping steps, when one is at triple points (inside the difficult regions $\overline N_{\infty}^2$).

\medskip

$\gamma)$ Last, but not least is also, our treatment of the $p_{\infty\infty}$'s, i.e. delete them at the price of introducing compensating $2$-handles. Concerning the whole issue of the {\ibf non-metrizability barrier}, see also \cite{27}. 

\medskip

\noindent END OF EXPLANATIONS. 

\bigskip

When we go to the variant~II, then the bad rectangles (\ref{eq4.1}) will be sent to infinity, in the manner of (\ref{eq3.29}), using PROPER Whitehead dilatations. This comes with still another system of handles $h(4)$ to be added, {\ibf after} the scenario from the lemma below will be played through. So, we will have $\partial h^2 (\varepsilon) \cdot \delta h^1 (4) = 0$, for $1 \leq \varepsilon \leq 3$. Moreover, internally to $h(4)$, $\partial h^4 (2) \cdot \delta h^4 (1) =$ easy id $+$ nil. So no violation of GSC comes with the transformation
$$
\mbox{Variant  I} \ \Longrightarrow \ \mbox{Variant II, either}.
$$

\bigskip

\noindent THE ZIPPING LEMMA 4.1. 1) {\it Starting from the $S_0$ in} (\ref{eq4.29}), {\it there is an infinite system of transformation, which we will call the Process $Z$ (like zipping)
\begin{equation}
\label{eq4.30}
S_0 \Longrightarrow S_1 \Longrightarrow S_2 \Longrightarrow S_3 \Longrightarrow \ldots \, .
\end{equation}

Each $S_i$ is a smooth $(N+4)$-manifold, coming with an inclusion $S_i \subset S_i + \underset{\overbrace{\mbox{\scriptsize $H \in (4.7)$}}}{\sum} \, D^2 (C(H))$ (which is a cell-complex). This last inclusion is not part of} (\ref{eq4.30}). {\it Each step in} (\ref{eq4.30}) {\it is followed by an appropriate {\ibf admissible subdivision} of the handle-body and/or cell decomposition $h(\varepsilon \leq 3)$, which never subdivide the $D^2 (C(H))$'s. The admissible subdivisions are stellar, barycentric or, more generally, Siebenmann's bisections, never violating the GSC property; remember that more general linear subdivisions may do so. Hence, we will keep $\partial h^2 \cdot \delta h^1 = $ {\rm id} $+$ {\rm nil}, under the admissible subdivisions, and these will also provide us with the necessary flexibility for performing the step $3)$ below.}

\medskip

\noindent 2) {\it Each individual step in} (\ref{eq4.30}) {\it is either a normal inclusion $S_n \subset S_{n+1}$ during which a finite number of new cells $h(2)$ or $h(3)$ is added, like in} (4.29.1), (4.29.2), {\it in a GSC-preserving manner, actually a J.H.C.~Whitehead dilatation, followed possibly by adding handles of index $\lambda > 1$, OR so-called DITCH-jumping steps $S_m \Longrightarrow S_{m+1}$, always attached to triple points $t \in M^3 (f)$. Globally speaking, the whole process $Z$ is a PROPER sequence of Whitehead dilatation, higher $(\lambda \geq 2)$ handle additions and diffeomorphisms.

\smallskip

Coming back to the special step $S_m \Longrightarrow S_{m+1}$, to any $t \in M^3 (f)$ like above, two {\ibf special BLUE} $1$-handles are attached. This is suggested in figure~{\rm 4.6}. [Later on, this figure will be expanded into $5.7 + 5.8$, but we added it also here, to help in understanding our present lemma.] So, to our $t \in M^3 (f)$, the special BLUE $1$-handles $[XY] \equiv \{[BA]$ and/or $[AD]$, figure~{\rm 4.6}$\} \subset W_n ({\rm BLUE}) - H$, are getting attached. We will assume, in the figure~{\rm 4.6}, that the zipping flow comes to the triple point $t$, via the road $B'' \to B$, going from right to left in our figure, at stage $m$ in the process $Z$. We will denote by $[0,m]$ the initial piece of $Z$} (\ref{eq4.30}), {\it up to time $m$ included, i.e. the piece
$$
S_0 \Longrightarrow S_1 \Longrightarrow S_2 \Longrightarrow \ldots \Longrightarrow S_m \, . \eqno (4.30.0)
$$
The next step, $m \Rightarrow m+1$ will be the DITCH jumping step, which we will describe now, actually in the $3)$ below. But still at time $m$ comes the piece $R_1 = [(X \equiv B), (B'',\gamma'')] \subset W_n ({\rm BLUE}) - H)$ see the figure~{\rm 4.6}, which is inductively, already in the DITCH, glued to the rest by the ditch filling, material $h(3) \mid [0,m]$.}

\medskip

\noindent 3) {\it To every special BLUE $1$-handle $[XY]$ there is attached, an embedded shadow arc, $S_j [XY]$, present at the beginning of $m \Rightarrow m+1$, coming with
$$
S_j [XY] \subset \partial \, \bigl[ \{ \Theta^{N+4} (X^2 - H) - \mbox{DITCH} - \{\mbox{special BLUE $1$-handle,}
\eqno (4.30.1)
$$
$$
\mbox{for all times} \ [m,\infty)\} \cup \{(h(2) + h(3)) \mid [0,m] \} \bigl] \, ,
$$
and this arc is liftable to $\partial S'_b (\widetilde M (\Gamma)-H)$, see here} (\ref{eq4.33}) {\it too.

\medskip

\noindent [NOTATIONS: For $[XY] \subset W({\rm BLUE}) - H$, the end-points are $X,Y$. When it comes to $S_j [XY]$, and when we want to stress it, rather than the $[XY]$, then the endpoints may be denoted by $\bar X , \bar Y$. With this we have $\bar X \equiv X \equiv \{$that endpoint of $[XY]$ already reached by the geometric realization of the zipping at time $m\}$. The other end of $S_j [XY]$ is at $\bar Y \equiv \{$the corresponding other endpoint, already reached by the previous zipping BLACK/RED$\}$. In our specific case of $[AB]$ figure~{\rm 4.6}, with the indicated zipping flow, we have $X = \bar X = B$ (with BLACK/BLUE identification), $\bar Y = A$ (with BLACK/RED identification).] So, in our specific example, $\bar X = X = B \in \{$the shaded region $R_1 \equiv (B , (B'' , \gamma''))$ (figure~{\rm 4.6}) consisting entirely of BLUE $h(1)$'s$\}$. One can also see our $A,B$ above in figure~{\rm 4.4}, where the BLUE cross is only a ghost. The four regions $R_i$ will eventually go into the DITCH of $W({\rm BLACK}) \cup W({\rm RED})$, actually in $W \times [-\varepsilon , \varepsilon] \times b_0^{N+1}$, making use of the supplementary dimensions, while, as we shall see, $[AB] + [AD]$ will eventually go very far from the site of figure~{\rm 4.6}, smeared along their corresponding shadow arcs.

\medskip

Imagine now that $X$ is moved slightly into the interior of $R_1$ and, on the shadow arc $S_j [XY]$, which starts at $X$, we pick up a point $X+\varepsilon$ with $[X,X+\varepsilon] \subset R_1$, see here the figure~{\rm 5.10.(A)}. On its road to $A$, the shadow arc will have to cross a lot of handles $\{h(1)$ (of various COLOURS)$\}$, $h(2), h(3)$.

\smallskip

The next mini-step, which is essentially isotopic, will not change the incidence relations $\{ \partial h^{\lambda + 1} (\varepsilon) \cdot \delta h^{\lambda} (\varepsilon), \ \varepsilon \leq 3\}$ at all. Making use of the flexibility which is provided by the admissible subdivision from the end of point {\rm 1)} above, we will {\ibf push} isotopically the BLUE shaded region $R_1$ along the portion $[X+\varepsilon , Y \equiv A] \subset S_j [XY]$, until it engulfs completely the whole $S_j [XY]$, starting from $B$ to the end $A$, and also a bit beyond $A$.

\smallskip

The $h(1)$ (COLOUR) $+ \, h(2) + h(3)$ which are crossed by the shadow arc $S_j [XY]$, get pushed together with $R_1$ so that, as announced, no new incidence relations appear in this process. At this point, before going on, we make the following}

\bigskip

\noindent {\bf Claim (4.31).} {\it Since any individual compact $W({\rm BLUE})$ is contained inside some bicollared $0$-handle of $Y(\infty)$, the BLUE $h(1)$'s involved in the push along the shadow arc above, are all contained inside some bicollared $0$-handle $H_i^0 (\gamma) \subset Y(\infty)$. Once this is so, the following happens. All the trajectories which may start from these $h(1)$'s under discussion now, are always {\rm short}, they are confined inside the same $H_i^0 (\gamma)$. The kind of trajectories we have in mind here, are chains of contacts, in the geometric incidence matrix of $\Theta^{N+4} (X^2)$
$$
\partial h_1^{\lambda + 1} \cdot \delta h_1^{\lambda} = 1  \, , \ \partial h_1^{\lambda + 1} \cdot \delta h_2^{\lambda} = 1 \, , \ \partial h_2^{\lambda + 1} \cdot \delta h_2^{\lambda} = 1 , \ldots , \partial h_k^{\lambda + 1} \cdot \delta h_k^{\lambda} = 1 \, ,
\eqno (4.31.0)
$$
coming with the following items
\begin{enumerate}
\item[\tiny$\bullet$)] In principle at least, we have here $\lambda = 1$.
\item[\tiny$\bullet$$\bullet$)] All the $\partial h_i^{\lambda + 1} \cdot \delta h_i^{\lambda} = 1$ are diagonal contacts in the geometric intersection matrix of $h(1)$ which, remember, has the form {\rm id} $+$ {\rm nil}, while all the $\partial h_i^{\lambda + 1} \cdot \delta h_{i+1}^{\lambda} = 1$ are off-diagonal contacts, of the same matrix.
\item[\tiny$\bullet$$\bullet$$\bullet$)] There are no off-diagonal contacts for $h_k^{\lambda+1}$, and there the trajectory stops for good.

\smallskip

What we have said so far, can be represented graphically as follows, if we use the dictionary $\{$square matrices with entries in $Z_+ \} \Rightarrow\{$oriented graphs$\}$ (i.e. discrete dynamical systems), see figure~{\rm 4.1.bis} too
$$
1 \longrightarrow 2 \longrightarrow 3 \longrightarrow \ldots \longrightarrow k \ \mbox{(no more outgoing arrow here)}
$$
\item[\tiny$\bullet$$\bullet$$\bullet$$\bullet$)] Here of course, the state $1 \equiv (h_1^2 (1) , h_1^1 (1))$ may itself be hit by some other trajectory, from outside of {\rm (4.31.0)}. End of CLAIM. 
\end{enumerate}
\hfill $\Box$

\bigskip

Here is the important consequence of our claim {\rm (4.31)}. We have
$$
[\partial h^2 (1) (\mbox{and here} \ h^2(1) \subset H_i^0 (\gamma) \ \mbox{is getting pushed along the shadow arc})] \cdot \delta h^1 = 0 \, ,
\eqno (4.31.1)
$$
when $h^1$ is outside of $H_i^0 (\gamma)$. As a consequence of this, we find that inside the geometric intersection matrix
$$
\partial h^2 (1) \ (\mbox{from the pushed region} \ R_1) \cdot \delta h^1 (3) \ (\mbox{on which the pushed region rests}) = 0 \, .
$$
We will denote by $S_m \Rightarrow S_{m+\frac12}$ the isotopic push of the BLUE region $R_1$, until it engulfs the whole arc $S_j [BA]$, and a bit beyond $A$. During this pushing operation, the region $R_1$ (or the active part of it) is in contact with $h(3)$, which cushions it like in the figure~{\rm 5.10.(B)}, and it is only indirectly through $h(3)$, that it comes into contact with the outer world. We will introduce the notation $h(1) \mid S_m \equiv \{$That part $R_1 \cup ($the $[XY]$ which is attached to it$) \subset S_m$, consisting entirely of BLUE $h(1)$'s, like in the CLAIM~{\rm (4.31)} and in the formula~{\rm (4.31.1)}, and which is actively involved in covering the shadow arc$\}$.

\smallskip

The next mini-step after this $m \Rightarrow m + \frac12$, we will denote by $m + \frac12 \Rightarrow m + \frac34$, and which we will call the {\ibf smearing}, is a folding map strictly confined inside the $h(1) \mid S_{m+\frac12}$, see the {\rm (4.31.2)} below, and which realizes the identification $\{[XY] = S_j [XY]$ (for $[XY] = [BA])\}$, leaving $B$ fixed. This creates the identification $\{ (Y = A ({\rm BLUE})) = (A({\rm RED} / {\rm BLACK})) = \bar Y \}$, see here the figures~{\rm 4.4, 4.6}. This smearing appears now in the diagram below
$$
\begin{matrix}
S_m &\overset{\rm pushing}{=\!\!=\!\!=\!\!=\!\!=\!\!=\!\!=\!\!\Longrightarrow} &S_{m+\frac12} &\overset{\rm smearing}{=\!\!=\!\!=\!\!=\!\!=\!\!=\!\!=\!\!\Longrightarrow} &S_{m+\frac34} &\Longrightarrow &S_{m+1} \, . \\
\uparrow &&\uparrow &&\uparrow \\
h(1) \mid S_m &=\!\!=\!\!=\!\!=\!\!=\!\!=\!\!=\!\!\Longrightarrow &h(1) \mid S_{m+\frac12} &=\!\!=\!\!=\!\!=\!\!=\!\!=\!\!=\!\!\Longrightarrow &h(1) \mid S_{m+\frac34} 
\end{matrix}
$$
\hfill {\rm (4.31.2)}

\medskip

In this diagram, the double arrows are transformations while the vertical ones are canonical inclusions. The $h(1) \mid S_{m+\frac12} \equiv \{$what becomes of $h(1) \mid S_m$ via the pushing$\}$ and, since the pushing in question is an isotopy, we have a diffeomorphism
$$
h(1) \mid S_m \underset{\rm DIFF}{=\!\!=\!\!=\!\!=\!\!=\!\!=\!\!=\!\!\Longrightarrow} h(1) \mid S_{m+\frac12} \, .
\eqno (4.31.3)
$$
Inductively, the $[0,m]$ from {\rm (4.30.0)} is GSC preserving and $h(1) \mid S_m$ has its handlebody decomposition inherited from the GSC handlebody decomposition of $S_m$. We will use the diffeomorphism from {\rm (4.31.3)} for simple-mindedly {\ibf transporting} the handlebody decomposition of $h(1) \mid S_m$ into a handlebody decomposition of the manifold $h(1) \mid S_{m+\frac12}$. Here the transport only concerns $R_1$.

\smallskip

Outside of $h(1) \mid S_{m+\frac12}$, the $S_{m+\frac12}$ inherits quite naturally a handlebody decomposition from $S_m$. By the (diffeomorphism {\rm (4.31.3)}) $+ \, \{$the transport of handlebody decomposition$\}$, the $h(1) \mid S_{m+\frac12}$ inherits the features from {\rm (4.31.1)} and it follows from all these things that $S_{m+\frac12} + \sum \, D^2 (C(H)) \in {\rm GSC}$. Keep in mind the following fact}

\medskip

\noindent (4.31.3.bis) \quad {\it The transport does not modify the internal incidence relations $\partial h^2 (1) \cdot \delta h^1 (1)$, while it certainly may change the $\partial h^2 (\varepsilon) \cdot \delta h^2 (1)$, for $\varepsilon > 1$. The smearing step which follows next in {\rm (4.31.2)}, bends $[XY]$ modulo $X$ and glues it to $S_j [XY] \subset h(1) \mid S_{n+\frac12}$, blending $[XY]$ into it. In our high dimensions, this can be read as another diffeomorphism, which continues the one in {\rm (4.31.3)}
$$
h(1) \mid S_{m+\frac12} \underset{\rm DIFF}{=\!\!=\!\!=\!\!=\!\!=\!\!=\!\!=\!\!\Longrightarrow} h(1) \mid S_{m+\frac34} \, .
\eqno (4.31.4)
$$

By the same argument of transporting and/or reshuffling handlebody  decompositions as above, we get that $S_{m+\frac34} + \sum \, D^2 (C(H)) \in {\rm GSC}$.

\smallskip

This second transport is a more serious business than the first one. It concerns now 
$$
R_1 \cup \{\mbox{the special BLUE $1$-handle} \ [XY]\} \, ,
$$
and the effect of transporting, which continues to abide to {\rm (4.31.3.bis)}, takes care of the potential nontriviality of the folding map.

\smallskip

Here is now how the ditch-jumping continues from $m+\frac34$ to $m+1$ and it will turn out that this piece is also GSC-preserving.}

\medskip

\noindent i) {\it The two BLUE zones $R_2 , R_3$ which start at $A$ (figure~{\rm 4.6}) are now sent into their respective DITCHES, together with the necessary $h(3)$ ditch-filling material. We have here a {\ibf partial} ditch-filling, again, of course.}

\medskip

\noindent ii) {\it Next, we have to handle the $[XY] = [AD]$ (figure~{\rm 4.6}) too. This will start with a pushing $+$ smearing like in $m \Rightarrow m+\frac34$ above. But the $S_j [AD]$, which is needed for this new pushing $+$ smearing, is now a much trickier thing than before. Notice that we have not said anything, yet, about the explicit location of the shadow arcs. It will turn out that for $S_j [AB]$ this is a relatively simple affair. We just make use of the closest fin, remembering that our concern is to keep the whole $Z$ PROPER. But for $S_j [AB]$ we will have to go quite far from $[AD]$. Out of the zipping flow
$$
\{\mbox{singularities} \ \longrightarrow B \} \, ,
$$
essentially played in reverse, one will be able to extract the continuous arc $S_j [XY]$ for $[XY] \equiv [AD]$. This will be the so-called (zipping)$^{-1}$ process, to be explained with great detail in the next section. Then, like before, $X$ now at $A$, will be pushed slightly, into the interior of one of the BLUE regions $R_2,R_3$ by now already in the ditch, via i) above, then we push $[X,X+\varepsilon] \subset R_{2 \, {\rm or} \, 3}$ along $S_j [XY]$ and then we finally perform the smearing $[XY] = S_j [XY]$.}

\medskip

\noindent iii) {\it We end with a partial ditch filling starting at $D$ (figure~{\rm 4.6}), concerning now $R_4$ and going in the $-z$ direction.

\medskip

For the needs of the long statement of our ZIPPING LEMMA, enough has been said concerning the ditch-jumping steps $S_m \Rightarrow S_{m+1}$.}

\medskip

\noindent 4) {\it For what will come next, it will be necessary to invoke the following list of ingredients: The fins, the fact that $M_2 (f) \subset X^2$ can only accumulate on {\rm LIM }$M_2 (f) \approx \sum (\infty)$, and then finally the process (zipping)$^{-1}$ which will turn out to be PROPER. But it will only be explicitly defined in the next section~{\rm V}. Making use of these ingredients, we can {\ibf choose} the location of the shadow arcs, which are anyway demanded by the steps $m \Rightarrow m+1$ which are DITCH-jumpings, so that they should not accumulate at finite distance, but on $\partial \sum (\infty)^{\wedge}$, which lives at infinity. Symbolically, we may write this as}
$$
\lim_{n=\infty} \, (S_j [XY]_n) = \infty \, .
$$

\medskip

\noindent 5) {\it Each of the individual $h_i^{\lambda} (\varepsilon)$'s, when $\varepsilon \leq 3$, has a first moment $p (h_i^{\lambda} (\varepsilon)) \in Z_+$ {\rm (\ref{eq4.30})} when it occurs, meaning that $p=p(h_i^{\lambda} (\varepsilon))$ is the smallest index such that, inside the process $Z$ from {\rm (\ref{eq4.30})} we should have that
$$
h_i^{\lambda} (\varepsilon) \subset S_p \, .
$$
Then, during the process $Z$ in question, $h_i^{\lambda} (\varepsilon)$ gets finitely many time subdivided, the subdivisions being repositioned during the pushes and smearings of the ditch-jumping steps. During all the ditch-jumping steps, for every given $(\lambda , \varepsilon , i)$ we will find that geometric incidence
$$
\left[ \sum_{{\rm all} \, [XY]} \{\mbox{shadow arcs} \ S_j [XY]\} \right] \cdot h_i^{\lambda} (\varepsilon) < \infty \, .
$$
After the finitely many subdivisions and repositionings, there is a final transformation, defined for every $(\lambda , i, \varepsilon)$
$$
h_i^{\lambda} (\varepsilon)_{\rm initial} \Longrightarrow \sum_{j=1}^{N} h_{ij}^{\lambda} (\varepsilon) \, , \quad N = N(\lambda , \varepsilon , i) \, .
$$
For any given $h_i^{\lambda} (\varepsilon)$ all this happens during a finite span of time, after which the $h_{ij}^{\lambda} (\varepsilon)$'s are no longer subdivided nor touched by the pushes $+$ smearings, during the infinite process {\rm (\ref{eq4.30})}. With this, we have a grand ``geometric incidence matrix'', involving all the stable $h_{ij}^{\lambda} (\varepsilon)$'s. These can be glued now together into a well-defined limit-object for the infinite process $Z$ from {\rm (\ref{eq4.30})}, call it
\setcounter{equation}{31}
\begin{equation}
\label{eq4.32}
S'_b (\widetilde M (\Gamma) - H) \equiv \lim_{n = \infty} S_n \, .
\end{equation}
We claim that this $S'_b (\widetilde M (\Gamma) - H)$ has a natural structure of smooth $(N+4)$-manifold, with large boundary. With all this, one may also interpret now our $Z$ as an embedding map, occurring in the following commutative diagram of PROPER embeddings, to be compared to {\rm (\ref{eq3.25})},}
\begin{equation}
\label{eq4.33}
\xymatrix{
S_0 \equiv \Theta^{N+4} (X^2 - H) - {\rm DITCH} \ar[rr]^-{j} \ar[d]^-Z &&S'_u (\widetilde M (\Gamma) - H) \, .  \\ 
S'_b (\widetilde M (\Gamma) - H) \ar[urr]_-{\mathcal J}
}
\end{equation}

\medskip

\noindent 6) {\it Just like in the toy-model, the embedding ${\mathcal J}$ is actually isotopic, via a simple-minded isotopy of non-boundary respecting injections, to a simple-minded {\ibf diffeomorphism} like in {\rm (\ref{eq3.26})}
\begin{equation}
\label{eq4.34}
S'_b (\widetilde M (\Gamma) - H) \overset{\eta}{\underset{\approx}{-\!\!\!-\!\!\!-\!\!\!-\!\!\!\longrightarrow}} S'_u (\widetilde M (\Gamma) - H) \, .
\end{equation}
[But we prefer to keep ${\mathcal J}$ and $\eta$ distinct, since it is ${\mathcal J}$ which makes the diagram {\rm (\ref{eq4.33})} strictly commutative, and not the $\eta$.]}

\medskip

\noindent 7) {\it Completely analogous to the $\alpha$ from {\rm (\ref{eq4.14})}, there is also second embedded framed link
\begin{equation}
\label{eq4.35}
\sum_{\overbrace{\mbox{\scriptsize$H \in (4.7)$}}} C(H) \overset{\beta}{\longrightarrow} S'_b (\widetilde M (\Gamma) - H) \, ,
\end{equation}
and we have the following diagram, connecting together these various objects
\begin{equation}
\label{eq4.36}
\xymatrix{
S'_b (\widetilde M (\Gamma)-H) \ar[rr]_-{\rm DIFFEOMORPHISM}^-{\eta}  &&S'_u (\widetilde M (\Gamma)-H) \, . \\ 
&\underset{\overbrace{\mbox{\scriptsize$H \in (4.7)$}}}{\sum} C(H) \ar[ur]_{\alpha} \ar[ul]^{\beta}
}
\end{equation}
This diagram has the following features:}

\medskip

\noindent (4.36.1) \quad {\it When it comes to the $H \in \{$BLACK Holes$\} + \{H(p_{\infty\infty})\}$, for which the corresponding curves are such that $C(H) \subset {\rm int} \, S'_{\varepsilon}$, then our diagram commutes strictly (up to isotopy at least).}

\medskip

\noindent (4.36.2) \quad {\it But when it comes to the other $H \in \{$completely normal Holes$\}$, which are such that for the corresponding curves $C(H) \subset \partial S'_{\varepsilon}$, then the {\rm (\ref{eq4.36})} only {\ibf commutes up to homotopy} (at the level of $\partial S'_{\varepsilon}$).}

\medskip

\noindent 8) {\it Making use of {\rm (\ref{eq4.35})} we can define the following cell-complex
\begin{equation}
\label{eq4.37}
S'_b \, \widetilde M (\Gamma) \equiv S'_b (\widetilde M (\Gamma) -H) \cup \sum_{\overbrace{\mbox{\scriptsize$\beta \, C(H), \, {\it for} \, H \in (4.7)$}}} D^2 (C(H)) \, .
\end{equation}
As an essentially formal consequence of the various things which have been already said in this lemma, it follows that $S'_b \, \widetilde M (\Gamma)$ is GSC.}

\bigskip

The proof of the zipping lemma will occupy all of the next section~V.

\bigskip

\noindent COMMENTS. A) Consider, in the neighbourhood of $p_{\infty\infty}$, a
$$
\{W^i ({\rm RED}) \mid (x \leq x_{\infty} + \zeta_i \, , \ {\rm with} \ \zeta_i > 0 \, , \ \zeta_i \to 0 \ {\rm when} \ i \to \infty) \} \approx W^i ({\rm RED} \cap H^0) \, ;
$$
this occupies its corresponding ${\rm DITCH} \subset V = W({\rm BLACK})$, all the way from $z = -\varepsilon$, to $z = +\varepsilon$. The {\ibf partial} ditch filling, inside $V - {\rm DITCH}$ is to be used here only for a piece of $W({\rm RED})$ namely the $\{ W^i ({\rm RED} - H^0) \mid x \geq x_{\infty} + \zeta_i \} -H$, which is a non-complementary wall, see here, also the figure~1.1.(B). Without all these various things, the very important object $S'_b (\widetilde M (\Gamma) - H)$ would have been unmanageably on the wrong side of the non-metrizability barrier, {\ibf if} we would not have drilled out the $p_{\infty\infty} \times [-\varepsilon , \varepsilon]$'s. But, as has been seen, once this is done, a whole can of worms is opened; and the present paper will certainly have to deal with them.

\medskip

\noindent B) Generically, our special BLUE $1$-handles occur essentially as
$$
[XY] \subset C(H(\mbox{completely normal})) \subset W({\rm BLUE}) - H \, ,
$$
and, in terms of (\ref{eq4.33}) $+$ (\ref{eq4.35}) we will find that $j [XY] = \eta \circ \beta [XY] \subset \partial S'_u (\widetilde M (\Gamma) - H)$. Moreover, in terms of what has been said at 2) in the lemma~4.1, $\beta [XY] = S_j [XY] \subset \partial S'_b (\widetilde M (\Gamma)-H)$. In view of these things, we may sometimes use the loose notation ``$[XY] = \beta [XY]$'', a reminder of the smearing step in (4.31.2). In this same vein, the ``$\beta$'' in (\ref{eq4.35}) is here a reminder of the fact that this arrow is ``bizarre'', since it uses the possibly complicated repositionning of $C(H) \subset [XY]$, along the shadow arc $S_j [XY]$. 

\medskip

\noindent C) While the reconstruction formula (4.14.1) is a posteriori, in the sense that it concerns two already defined objects, the superficially analogous formula (\ref{eq4.37}) is a priori, it {\ibf defines} the $S'_b \, \widetilde M (\Gamma)$.

\smallskip

In (4.14.1), the $\alpha \, C(H)$'s (see (\ref{eq4.14})), play exactly the same role as the $\beta \, C(H)$ in (\ref{eq4.37}). \hfill $\Box$

\bigskip

The next lemmas~4.2 to 4.7 are various complements to the ZIPPING LEMMA.

\bigskip

\noindent {\bf Lemma 4.2.} 1) {\it The process $Z$ from {\rm (\ref{eq4.30})} is modelled after a precise zipping strategy for our $X^2 \overset{f}{\longrightarrow} fX^2$ from {\rm (\ref{eq1.1})}, which is {\ibf equivariant}
\begin{equation}
\label{eq4.38}
X^2 \equiv X_0^2 \longrightarrow X_1^2 \longrightarrow X_2^2 \longrightarrow X_3^2 \longrightarrow \ldots \, .
\end{equation}
This strategy will be made explicit in section~{\rm V}.}

\medskip

\noindent 2) {\it So, the whole process $Z$ which mimicks {\rm (\ref{eq4.38})} is equivariant too. It follows that there is a free $\Gamma$-action on $S'_b \, \widetilde M (\Gamma)$ entering into the following commutative diagram, the vertical arrows of which are the natural inclusions
$$
\xymatrix{
&\Gamma &\!\!\!\!\!\!\!\!\!\!\!\!\!\!\!\!\times &\!\!\!\!\!\!\!\!\!\!\!\!\!\!\!\!\!\!\!S'_b (\widetilde M (\Gamma)) \ar[rr]  &&S'_b (\widetilde M(\Gamma))    \\ 
&\Gamma \ar[u]_{\rm id}  &\!\!\!\!\!\!\!\!\!\!\!\!\!\!\!\!\times &\!\!\!\!\!\!\!\!\!\!\!\!\!\!\!\!\!\!\!S'_b (\widetilde M (\Gamma)-H) \ar[u] \ar[rr] &&S'_b (\widetilde M(\Gamma)-H) \ar[u] \, . 
} \eqno (4.38.1)
$$
Whenever it makes sense, the zipping lemma is equivariant.}

\medskip

\noindent 3) {\it There are actually homotopy equivalences
\begin{equation}
\label{eq4.39}
S_i + \sum DC(H) \longrightarrow X_i^2
\end{equation}
which, together with {\rm (\ref{eq4.30})} $+$ {\rm (\ref{eq4.38})} can be put into an infinite, homotopy-commutative ladder.} \hfill $\Box$

\bigskip

The next lemma expresses a very important property of the objects of type $S'_{\varepsilon}$ (where $\varepsilon$ can mean $u$ or $b$) and where the prime $(')$ means that $p_{\infty\infty} ({\rm all})$ are excised, as opposed to the objects of type $S_{\varepsilon}$, where only the $p_{\infty\infty} (S)$ are excised. This propery, expressed here at level $\widetilde M (\Gamma)$, survives at level $M(\Gamma)$ too, and then at the last level it will play a crucial role in the proof of the COMPACTIFICATION LEMMA~4.6, in section~VI.

\bigskip

\noindent {\bf Lemma 4.3.} {\it We consider some arbitrary sequence of points
$$
P_n = (x_n , y_n , z_n , t_n) \in S'_{\varepsilon} (\widetilde M (\Gamma) -H) \mid W_{(\infty)} ({\rm BLACK})
$$
which is such that $\underset{n=\infty}{\lim} \, (x_n , y_n) = (x_{\infty} , y_{\infty}) = p_{\infty\infty}$, $-\varepsilon \leq z_n \leq + \varepsilon$, i.e. which is {\ibf without restrictions} in the values $z_n$, and with $t_n \in b^{N+1} \subset B^{N+1}$ (see {\rm (\ref{eq4.15})}). With this, we have the following items.}

\medskip

\noindent 1) {\it Both at the level $S'_{\varepsilon} (\widetilde M (\Gamma)-H)$ and at the levels
\begin{equation}
\label{eq4.40}
S'_{\varepsilon} (\widetilde M (\Gamma)-H) \underset{\overbrace{\mbox{\scriptsize$\sum C(p_{\infty\infty} (all))$}}}{\cup} \sum \left(D^2 (H(p_{\infty\infty})) \times \frac12 \, B^{N+1}\right) \quad \mbox{or} \quad S'_{\varepsilon} \, \widetilde M (\Gamma) \, ,
\end{equation}
always with an $S'_{\varepsilon}$, we have that
\begin{equation}
\label{eq4.41}
\lim_{n=\infty} P_n = \infty \, .
\end{equation}
}

\medskip

\noindent 2) {\it Our sequence  $\{P_n \}$ stays far both from the $[A(p_{\infty\infty}) \times B^{N+1} - \{$the indentations$\}]$ from {\rm (4.28.1)}, and from the various $2$-handles in} (\ref{eq4.40}).

\bigskip

\noindent {\bf Proof.} This is an easy consequence of the fact that, at level $S'_{\varepsilon}$ the $p_{\infty\infty} \times [-\varepsilon , \varepsilon]$ are completely deleted, and that
$$
b^{N+1} \cap \frac12 \, B^{N+1} = \emptyset \subset B^{N+1} \, .
$$
\hfill $\Box$

\newpage

In the zipping lemma~4.1, at (\ref{eq4.30}) we have defined the $S'_b (\widetilde M (\Gamma) - H)$ and at (\ref{eq4.37}) the $S'_b \, \widetilde M (\Gamma)$. Now, we introduce a slightly smoother cell-complex $S_b \, \widetilde M (\Gamma)$, as follows
\begin{equation}
\label{eq4.42}
S_b \, \widetilde M (\Gamma) \equiv \Bigl\{ S'_b (\widetilde M (\Gamma) - H) \cup \sum_{H\in(4.7)} D^2 (C(H)) \, ,
\end{equation}
so far just like in (\ref{eq4.37}), but then to this the following two operations are, afterwards applied: 

\medskip

\noindent a) first a simple-minded change $\{C(H({\rm BLACK \, hole}))$ (see figure~4.1)$\} \Rightarrow \partial H ({\rm BLACK \, Hole})$ i.e., with the notations from the figure~4.1
$$
\partial \, C(H(W(n))) \Longrightarrow \partial H (W(n)) \, ;
$$
b) Next, the step (\ref{eq2.23}) is applied to the $p_{\infty\infty} ({\rm proper})$, leaving only the $C(H(p_{\infty\infty} (\Gamma))) = \partial D^2 (H(p_{\infty\infty} (\Gamma)))$ singular$\}$. END of (\ref{eq4.42}).

\bigskip

Here the operation a) makes our $S'_b$ be smooth at the BLACK Holes (now filled), and it is already smooth at the usual Holes (filled too). The (2.23) was defined, a priori, in the context $S_u^{(')}$, but it makes perfectly sense here too. It leaves us with an $S_b \, \widetilde M (\Gamma)$ which is a smooth $(N+4)$-manifold, except for mild singularities at the $C(H(p_{\infty\infty} (S)))$'s.

\smallskip

The next lemma, puts together some technicalities which we shall need later on.

\bigskip

\noindent {\bf Lemma 4.4.} 1) {\it Without any loss of generality, the following things happen in the neighbourhood of any $p_{\infty\infty}$ (and see here figures~{\rm 2.2, 4.5} too)}

\medskip

\noindent (4.43.1) \quad {\it $A(p_{\infty\infty}) \mid \{ S'_b$ or $\Theta^{N+4} \} \subset \left[-\frac\varepsilon2 \leq z \leq \frac\varepsilon2 \right] \times B^{N+1} - \{$the indentation in {\rm (4.28.1)}, i.e. the DITCHES in {\rm (4.29.0)}$\}$; here and also in the next items, we ignore those dimensions of the $S'_b$ and/or $\Theta^{N+4}$, when not explicitly mentioned;}

\medskip

\noindent (4.43.2) \quad {\it $\{$Ditch filling and $h(3)\} \mid p_{\infty\infty} \subset b^{N+1} \subset B^{N+1} - \left( \frac12 + \varepsilon_0 \right) B^{N+1}$;}

\medskip

\noindent (4.43.3) \quad {\it $\{$adding fins and sending $\partial \, \sum^{\wedge} (\infty)$ to infinity$\} \subset \{\vert z \vert \geq \varepsilon \}$;}

\medskip

\noindent (4.43.4) \quad {\it $\{$Ditch-jumping and smearing$\} \subset \{ \vert z \vert \geq \frac{3\varepsilon}4 \}$.}

\medskip

\noindent {\it The point, for items} (4.43.3), (4.43.4) {\it is that they are far from $A(p_{\infty\infty} (all))$.}

\medskip

\noindent 2) {\it As a consequence of the {\rm i)} above, the $Z$-action} (\ref{eq4.30}) {\it is disjoined from $A(p_{\infty\infty}) \times \left(\frac12 + \varepsilon_0 \right) B^{N+1}$ which is the site of the action in lemma~{\rm 2.2}, see in particular} (\ref{eq2.23}). {\it This last action concerns the $p_{\infty\infty}$ (proper)'s, and never the $p_{\infty\infty} (S)$'s.}

\medskip

\noindent 3) {\it Starting from any of the two smooth manifolds $S'_{\varepsilon} (\widetilde M (\Gamma) -H)$, where $\varepsilon = u$ or $b$, we have the following chain of transformations}
\setcounter{equation}{43}
\begin{equation}
\label{eq4.44}
S'_{\varepsilon} (\widetilde M (\Gamma) - H) \overset{\textcircled{\mbox{\scriptsize 1}}}{\underset{add \ D^2 (H(normal))}{=\!\!=\!\!=\!\!=\!\!=\!\!=\!\!=\!\!=\!\!=\!\!=\!\!=\!\!=\!\!=\!\!=\!\!=\!\!=\!\!\Longrightarrow}} S'_{\varepsilon} (\widetilde M (\Gamma) - H (p_{\infty\infty} (\mbox{\ibf all}))) \overset{\textcircled{\mbox{\scriptsize 2}}}{\underset{add \ D^2 (H(p_{\infty\infty} (S)))}{=\!\!=\!\!=\!\!=\!\!=\!\!=\!\!=\!\!=\!\!=\!\!=\!\!=\!\!=\!\!=\!\!=\!\!=\!\!=\!\!=\!\!\Longrightarrow}}
\end{equation}
$$
S'_{\varepsilon} (\widetilde M (\Gamma) - H(p_{\infty\infty} (proper)))  \overset{\textcircled{\mbox{\scriptsize 3}}}{\underset{add \ D^2 (H (p_{\infty\infty}  (proper)))}{=\!\!=\!\!=\!\!=\!\!=\!\!=\!\!=\!\!=\!\!=\!\!=\!\!=\!\!=\!\!=\!\!=\!\!=\!\!=\!\!=\!\!=\!\!=\!\!=\!\!=\!\!\Longrightarrow}} S'_{\varepsilon} \, \widetilde M(\Gamma) \overset{\textcircled{\mbox{\scriptsize 4}}}{\underset{(2.23)}{=\!\!=\!\!=\!\!=\!\!=\!\!\Longrightarrow}} S_{\varepsilon} \, \widetilde M (\Gamma) \, .
$$
{\it The notations $(S'_{\varepsilon} (\widetilde M (\Gamma) - H(p_{\infty\infty} (\eta)))$'s mean here that only the $H(p_{\infty\infty} (\eta))$'s and their compensating disks are being removed. Here only the first two objects are smooth $(N+4)$-manifolds, the other are cell-complexes with mild singularities. Moreover, concerning this chain of transformations {\rm (\ref{eq4.44})}, we have the following facts too}

\medskip

\noindent (4.44.1) \quad {\it A lot of pairs of steps in {\rm (\ref{eq4.44})} commute with each other; we have $[\textcircled{\mbox{\scriptsize 1}}, \textcircled{\mbox{\scriptsize 2}}] = [\textcircled{\mbox{\scriptsize 2}}, \textcircled{\mbox{\scriptsize 3}}] = [\textcircled{\mbox{\scriptsize 3}}, \textcircled{\mbox{\scriptsize 1}}] =0$ and then, also $[$adding the $D^2 (H(normal))$'s, $\{$adding $D^2 (H(p_{\infty\infty} (all)))\} + \{\textcircled{\mbox{\scriptsize 4}}$, which only concerns the $p_{\infty\infty} (proper)\}]=0$.}

\medskip

\noindent (4.44.2) \quad {\it The whole {\rm (\ref{eq4.44})}, written above at level $\widetilde M(\Gamma)$, also functions downstairs, at level $M(\Gamma)$.}

\medskip

\noindent (4.44.3) \quad {\it In the context of our {\rm (\ref{eq4.44})}, the two actions, at $\varepsilon = u$ and at $\varepsilon = b$ are {\ibf completely isomorphic}.}

\bigskip

From the fact that $S'_b \, \widetilde M (\Gamma)$ is GSC, see here the 8) in the ZIPPING LEMMA, and making use of the step $\textcircled{\mbox{\scriptsize 4}}$ in (\ref{eq4.44}), for $S_b^{(')}$, we will be able to deduce the following

\bigskip

\noindent {\bf Corollary 4.5.} {\it We have $S_b \, \widetilde M (\Gamma) \in {\rm GSC}$.}

\bigskip

\noindent {\bf Proof.} What one needs to prove is that the transformation
$$
S'_b \, \widetilde M^3 (\Gamma) \underset{(2.23)}{=\!\!=\!\!=\!\!=\!\!=\!\!\Longrightarrow} S_b \, \widetilde M^3 (\Gamma)
\eqno (4.44.4)
$$
does not change the GSC property which $S'_b \, \widetilde M^3 (\Gamma)$ already has, see 8) in the ZIPPING LEMMA. One should remember that the transformation (4.44.4) is localized inside (see figure~4.1.bis too, but now with $p_{\infty\infty} ({\rm proper})$ in focus)
$$
\sum_{p_{\infty\infty} ({\rm proper})} A(p_{\infty\infty}) \times \left( \frac12 + \varepsilon \right) B^{N+1} \cup \{ j = (h^1 (p_{\infty\infty}) , D^2 (H(p_{\infty\infty}))) \} \times \frac12 B^{N+1} \, ,
\eqno (4.44.5)
$$
confined inside $-\frac\varepsilon2 \leq z \leq \frac\varepsilon2$, far from anything interesting, like the zipping action with its (possibly partial) DITCH-filling, the pushing $+$ smearing coming with the special BLUE $1$-handles or the $C(H({\rm normal})) = \partial D^2 (H({\rm normal}))$.

\smallskip

From (\ref{eq4.33}) and (\ref{eq4.44}) we can derive the following diagram which also makes use of the commutativities from (\ref{eq4.44}),
$$
\begin{matrix}
S'_b (\widetilde M^3 (\Gamma) - H) \underset{{\rm add} \ D^2 (H(p_{\infty\infty}))}{=\!\!=\!\!=\!\!=\!\!=\!\!=\!\!=\!\!=\!\!=\!\!=\!\!=\!\!=\!\!=\!\!=\!\!=\!\!=\!\!\Longrightarrow} S'_b (\widetilde M^3 (\Gamma) - H ({\rm normal})) \underset{(2.23)}{=\!\!=\!\!=\!\!=\!\!=\!\!=\!\!\Longrightarrow} S_b (\widetilde M^3 (\Gamma) - H ({\rm normal})) \, . \\
\left\Uparrow Z \qquad \qquad \qquad \qquad \qquad \qquad \qquad \qquad \qquad \qquad \qquad \qquad \qquad \begin{matrix} { \ } \\ {\rm add} \\ D^2 (H({\rm normal})) \\ { \ } \end{matrix} \right\Downarrow \\ \\
S_0 \ ({\rm see} \ (4.33)) \qquad \qquad \qquad \qquad \qquad \qquad \qquad \qquad \qquad \qquad \qquad \qquad \qquad \qquad S_b \, \widetilde M (\Gamma)
\end{matrix}
\eqno (4.44.6)
$$
The ``(\ref{eq2.23})'' is the same in (4.44.4) and (4.44.6). The $Z$ is confined inside $b^{N+1} \times \left\{ \vert z \vert \geq \frac{3\varepsilon}4 \right\}$, far from (4.44.5).

\smallskip

When one considers the transformation of geometric intersection matrix where, in the context (4.44.6) we go from $S_0$ to $S_b \, \widetilde M (\Gamma)$, then one can see that this is without effect on the

\medskip

\noindent $(*)$ $\{\mbox{the internal geometric intersection matrix of (4.44.5)}\} + \{\mbox{the ingoing and outgoing arrows from it}\} ,
$

\medskip

\noindent except that the state $j (p_{\infty\infty} ({\rm proper}))$ gets erased. All this is without impact on the GSC which stays intact under our transformation (4.44.4), as claimed.

\smallskip

Concerning $(*)$, see again the figure~4.1.bis, but for the $p_{\infty\infty} ({\rm proper})$-island, and not the $p_{\infty\infty} (S)$-island, which is not touched by the (4.44.4). \hfill $\Box$

\bigskip

Once the corollary~4.5 has been proved, we turn now our attention back to the $S_b \, \widetilde M (\Gamma)$ from (\ref{eq4.42}), to the actually more inclusive $S_{\varepsilon}^{(')} \widetilde M (\Gamma)$. The following lemma is completely analogous to lemma~2.1.

\bigskip

\noindent {\bf Lemma 4.6.} 1) {\it There is a natural free action}
\begin{equation}
\label{eq4.45}
\Gamma \times S_b^{(')} \widetilde M (\Gamma) \longrightarrow S_b^{(')} \widetilde M (\Gamma) \, .
\end{equation}

\noindent 2) {\it One can, also, define directly $S_b^{(')} M (\Gamma)$ downstairs and this last object has the properties
$$
S_b^{(')} M (\Gamma) = (S_b^{(')} \widetilde M (\Gamma)) / \Gamma \eqno (4.45.1)
$$
and
$$
S_b^{(')} \widetilde M (\Gamma) = (S_b^{(')} M (\Gamma))^{\sim} \, , \eqno (4.45.2)
$$
and, of course, $\pi_1 \, S_b^{(')} M (\Gamma) = \Gamma$.}

\bigskip

So far, our discussion of the geometric realization of the zipping in high dimensions has proceeded UPSTAIRS, at the level of $\widetilde M(\Gamma)$, but now we will move DOWSTAIRS, at the level of the compact $M(\Gamma)$; this crucial {\ibf compactness} reflects, of course, the fact that our $\Gamma$ is a finitely presented group. It will turn out that, once we move dowstairs, then by various COMPACTNESS ARGUMENTS, we can connect $S^{(')}_u M(\Gamma)$ and $S^{(')}_b M(\Gamma)$, actually show that they are diffeomorphic, then via functoriality, this connection can be transported upstairs. With this, we will finally deduce that $S^{(')}_u \widetilde M(\Gamma)$ is GSC, from the known fact that $S^{(')}_b \widetilde M(\Gamma)$ is GSC.

\smallskip

According to our lemma~4.2 above, $Z$ was modelled on an {\ibf equivariant zipping} $X^2 \to fX^2$, but once we have this equivariance, the process $Z$ (\ref{eq4.30}) makes sense downstairs too, leading like in (\ref{eq4.32}) to an object $S'_b (M(\Gamma)-H)$ downstairs when, of course, we have a more obvious $S^{(')}_u (M(\Gamma)-H)$ too. It is important, at this point, that there are here two distinct, but eventually equivalent roads, leading to the same scenario DOWNSTAIRS.

\medskip

A) Everything which we did so far, in this paper, has been $\Gamma$-equivariant and so we can take the viewpoint that the story downstairs is just the quotient by $\Gamma$ of the equivariant story upstairs. This means using the formulae like $S^{(')}_{\varepsilon} M(\Gamma) = (S^{(')}_{\varepsilon} \widetilde M(\Gamma))/\Gamma$ or the same with Holes, as definitions.

\medskip

B) But then, the zipping being equivariant, it makes sense directly downstairs, instead of zipping the upper line in (\ref{eq1.24}) we zip now the lower line and then built directly the analogue of $Z$ (\ref{eq4.30}), downstairs too. Next, one notices that the $S^{(')}_{\varepsilon}$ have good functorial properties like in lemma~2.1, under appropriate maps one can push them forward and pull them backwards. So, one constructs now the $S^{(')}_{\varepsilon} M(\Gamma)$ (possibly with Holes) directly downstairs, notices that $S^{(')}_{\varepsilon} (\widetilde M (\Gamma)(-H)) = S^{(')}_{\varepsilon} (M(\Gamma)(-H))^{\sim}$ and then one deduces the definition from A) as a theorem. In the same vein, the $S_u^{(')}$ certainly has good localization and glueing properties and if one looks at things from the correct viewpoint, so has $S_b^{(')}$ too. So the two viewpoints A) and B) {\ibf are} equivalent, and this is important for the self-consistency of our approach.

\medskip

So, I will start by indicating now, rather briefly, how the zipping lemma~4.1 and its sequels fare downstairs. The starting point is now
$$
[\Theta^{N+4} ((X^2-H)/\Gamma) - {\rm DITCHES}] + \sum D^2 (C(H)) = \Theta^{N+4} (f(X^2-H) - {\rm DITCHES}) + \sum D^2 (C(H)) 
$$
$$
= \{\mbox{the quotient (4.29)} / \Gamma \} \, . \eqno (4.45.3)
$$

From here on, invoking now lemma~4.2 which tells us that the process $Z$ in (\ref{eq4.30}) is equivariant, we can make sense of it DOWNSTAIRS too, and get then the analogue of the ZIPPING LEMMA~4.1, which holds downstairs too, except for the following two items.

\smallskip

The $\Gamma$-action has clearly disappeared, and any reference to GSC (the claim (4.31) included) is no longer with us either. But the following items certainly will survive. To begin with the analogue of the diagram (\ref{eq4.33}), namely the following commutative diagram of PROPER embeddings
\begin{equation}
\label{eq4.46}
\xymatrix{
\Theta^{N+4} ((X^2 - H) / \Gamma) - {\rm DITCHES} \ar[rr]^-{j/\Gamma} \ar[d]^-{Z/\Gamma} &&S'_u (\widetilde M (\Gamma) - H) \, .  \\ 
S'_b (M (\Gamma) - H) \ar[urr]_-{{\mathcal J} /\Gamma} 
}
\end{equation}
Again, like in the context of the map (\ref{eq4.33}) ${\mathcal J} /\Gamma$ is connected, via a simple-minded isotopy to a diffeomorphism, to be denoted again by
$$
S'_b (M(\Gamma)-H) \overset{\eta}{\underset{\approx}{-\!\!\!-\!\!\!-\!\!\!\longrightarrow}} S'_u (M(\Gamma) -H) \, .
$$

The analogue of (\ref{eq4.36}) survives too, but right here we will only retain from it the following $C(H({\rm normal}))$ part, which continues to be {\ibf homotopy-commutative}
\begin{equation}
\label{eq4.47}
\xymatrix{
S'_b (M (\Gamma)-H) \ar[rr]^-{\eta}_-{\approx} &&S'_u (M(\Gamma)-H)  \\ 
\partial S'_b (M (\Gamma)-H) \ar[u] \ar[rr]^-{\eta}_-{\approx} &&\partial S'_u (M(\Gamma)-H) \ar[u] \\
&\underset{H({\rm normal})}{\sum} \{C(H({\rm normal}))+{\rm framing}\} \ar[ur]_{\alpha} \ar[ul]^{\beta}
}
\end{equation}

\noindent (4.47.1) \quad Of course, just by itself the upper square in (\ref{eq4.47}) commutes strictly. Also, the map $\eta$ respects, strictly, the contribution
$$
\sum_{p_{\infty\infty} ({\rm all})} \left[ A(p_{\infty\infty}) \times \left( \frac12 + \varepsilon_0 \right) B^{N+1} \right] \supset \sum_{p_{\infty\infty} ({\rm all})} C(p_{\infty\infty}) \times \frac12 \, B^{N+1} \, . \eqno (*)
$$
This contribution $(*)$ lives buryed deep inside the int $S'_{\varepsilon} (M(\Gamma)-H)$, safely away from the boundary, which is the site of the homotopy which connects $\eta \circ \beta$ and $\alpha$ in (\ref{eq4.47}).

\smallskip

This ends our zipping lemma discussion downstairs, and with this we can move now forward.

\smallskip

We will state below a COMPACTNESS LEMMA, which will be proved later, in section~VI. But I start by listing now the four basic ingredients which are all necessary for this very basic lemma.

\smallskip

\begin{enumerate}
\item[i)] The group $\Gamma$ is finitely presented, hence the $M(\Gamma) = \widetilde M (\Gamma)/\Gamma$ which, loosely speaking is the fundamental domain of $\Gamma$ is COMPACT; the action $\Gamma \times \widetilde M (\Gamma) \to \widetilde M (\Gamma)$ is {\ibf cocompact}.
\item[ii)] We have a uniform upper bound for te zipping length, a fact proved in \cite{29};
\item[iii)] We do not work now with $S_{\varepsilon}$ but with $S'_{\varepsilon}$, which will allow us to take advantage of lemma~4.3, applied DOWNSTAIRS;
\item[iv)] Very essentially, we move from the VARIANT~I which is used throughout the section I to V to VARIANT~II, as this will be implemented in the section VI. With these things, we have the following
\end{enumerate}

\bigskip

\noindent {\bf Compactness Lemma 4.7.} {\it Once we take into account the four ingredients above, the diagram {\rm (\ref{eq4.47})} {\ibf considered now in the context of the} VARIANT~{\rm II} (hence the subscripts {\rm II}) and which we rewrite here as 
$$
\xymatrix{
\partial S'_b (M (\Gamma)-H)_{\rm II} \ar[rr]_-{DIFFEOMORPHISM}^-{\eta}  &&\partial S'_u (M (\Gamma)-H)_{\rm II} \\ 
&\sum \{C(H(normal))+framing\} \ar[ur]_-{\alpha} \ar[ul]^-{\beta}
}
\eqno (4.47.{\rm II})
$$
commutes up to a PROPER homotopy.}

\bigskip

\noindent VERY IMPORTANT REMARK. The $\alpha \, C(H)$, $\eta\beta \, C(H)$ are contained inside $S'_u (\widetilde M (\Gamma)-H)_{\rm I}$, but the homotopy which connects them needs the larger
$$
\partial S'_u (\widetilde M(\Gamma)-H)_{\rm II} \supset S'_u (\widetilde M (\Gamma)-H)_{\rm I} \, .
$$
Let us take now another look at the triple points. Once our $W$'s are endowed with transversal orientations, then at each triple point, we can assign signs $\pm$ to each of the four Holes $H$ adjacent to the triple point, more exactly signs will be associated to the four $H$-corners at the triple point. The signs are corner-dependent, the same $H$ may come with different signs at its various corners. At the triple point $t$ we have a $W({\rm BLUE})$, carrying four holes, and two complementary walls $W_1 , W_2$. When the two orientations, transversal to $W_1 , W_2$ look away from $H$, the corner is ``$+$'', otherwise it is ``$-$'', making that exactly one hole is $H^+$ at $t$ and three are $H^-$. This rule is illustrated in figure~5.7, in the next section. We move now from these corner-dependent signs, to global ones. A hole $H$, or rather its $C(H)$ will be $C^-$ (respectively $C^+$) if at least one of its corners is $-$ (respectively if they are all $+$).

\smallskip

We go back now to the diagrams (\ref{eq4.36}) $+$ (\ref{eq4.47}), in connection with which we make the following

\bigskip

\noindent {\bf Claim (4.48).} Both at the level $S'_{\varepsilon} (\widetilde M (\Gamma)-H)$ (from (\ref{eq4.36})) and at the level $S'_{\varepsilon} (M (\Gamma)-H)$ (from (\ref{eq4.47})) we may assume, without any loss of generality that we have $\alpha = \eta \circ \beta$, when it comes to the curves
$$
\left\{ C(p_{\infty\infty} ({\rm all})) \, , \ C(H(\mbox{BLACK Hole})) \quad \mbox{and} \quad C^+ (H) \right\} \in \sum_H C(H) \, .
$$

So, when it comes to the concern of being PROPER for the homotopy in the diagram (\ref{eq4.37}), it is {\it only} the $\beta \, C^- (H)$ which should worry us at all. The next lemma is valid now, both upstairs and downstairs. But, in its form down in $M(\Gamma)$ it will provide an essential ingredient for the proof of the COMPACTNESS lemma~4.6. So we only state it downstairs here.

\bigskip

\noindent {\bf Lemma 4.8.} (DOWNSTAIRS) {\it We consider now all the various $C^- (H)$'s,
\setcounter{equation}{48}
\begin{equation}
\label{eq4.49}
C^- (H_1) \, , \ C^- (H_2) , \ldots \subset \Theta^3 (X^2-H) / \Gamma = \Theta^3 (fX^2 - H) \, ,
\end{equation}
where remember that ``$-H$'' means not only with all Holes deleted, but with all the $\{ D^2 (C(H))$ $\supset$ $D^2 (C(H$ $(p_{\infty\infty}))) \}$'s removed too.

\smallskip

For each $H_i$ in {\rm (\ref{eq4.49})} there is an embedded arc $\gamma_i \subset \partial S'_u (M(\Gamma)-H)$, which joins $\alpha \, C^- (H_i)$ to $\eta \beta \, C^- (H_i)$ and which is such that, when we go to the closed loop
\begin{equation}
\label{eq4.50}
\Lambda (H_i^-) \equiv \alpha \, C^- (H_i) \underset{\overbrace{\gamma_i}}{\bullet} (\eta \beta \, C^- (H_i))^{-1} \to \partial S'_u (M(\Gamma)-H) \, ,
\end{equation}
then we have the following items}

\medskip

\noindent (4.50.1) \quad {\it As a consequence of the homotopy-commutativity of the lower triangle in {\rm (\ref{eq4.47})}, the loop $\Lambda (H_i^-)$ as defined above, is null-homotopic in $\partial S'_u (M(\Gamma)-H)_{\rm I}$.}

\medskip

\noindent (4.50.2) \quad {\it There is a uniform bound s.t. $\Vert \Lambda (H_i^-) \Vert < N$, for all $i$'s.}

\medskip

\noindent (4.50.3) \quad {\it For every compact $K \subset \partial S'_u (M(\Gamma)-H)$ there are only finitely many $j$'s s.t. $K \cap \Lambda (H_j^-) \ne \emptyset$. In a more schematical language, we have $\underset{n=\infty}{\lim} \, \Lambda (H_n^-) = \infty$, inside $\partial S'_u (M(\Gamma)-H)_{\rm I}$.}

\bigskip

\noindent {\bf Comments and some hints on the proof of the lemma 4.8.} Concerning (4.50.2), once we work downstairs, it is very easy to control the $\Vert \alpha \, C^- (H_i) \Vert$.

\smallskip

The construction of the $\eta\beta \, C^- (H_i)$'s and of the $\gamma_i$'s will be a major objective of the next section. Since one of the main ingredients of the construction is the process ``(zipping)$^{-1}$'', explained in the next section, the control of $\Vert \eta \beta \, C^- (H_i) \Vert$ and $\Vert \gamma_i \Vert$ will need the uniform boundedness of the zipping length, proved in the first paper of this series \cite{29}. It will also need the lemmas~5.1, 5.3 from the next section (actually the various estimates given in those lemmas). So much about (4.50.2). The (4.50.3) will follow from lemma~4.3, in conjunction with the 6) in lemma~5.3, next section. \hfill $\Box$

\bigskip

\noindent FINAL ARGUMENTS FOR THE PROOF OF THEOREM 2.3. ($S_u^{(')} \widetilde M (\Gamma)_{\rm II}$ is GSC). In view of the high dimensions $N+4 \gg 1$, which we are working with, in the context of the diagram (4.47.II), the PROPER homotopy from lemma~4.6 implies smooth ambient isotopy inside $\partial S'_u (M(\Gamma)-H)_{\rm II}$. So, from (4.47.II) $+$ (4.47.1) we can pull now the following commutative diagram, where only smooth objects and smooth maps are involved, with the $d$ below being a DIFFEOMORPHISM, and the two unions being along $C(H({\rm normal}))$
\begin{equation}
\label{eq4.51}
\xymatrix{
{\begin{matrix} S'_b (M(\Gamma)-H)_{\rm II} \ \cup \\ { \ } \\ \sum D^2 (H({\rm normal})) \end{matrix} }\ar[rr]_-{d}  &&{\begin{matrix} S'_u (M(\Gamma)-H) \ \cup \\ { \ } \\ \sum  D^2 (H({\rm normal})). \end{matrix}}  \\
&\underset{p_{\infty\infty} ({\rm all})}{\sum} A(p_{\infty\infty}) \supset \underset{p_{\infty\infty} ({\rm all})}{\sum} C(p_{\infty\infty}) \ar[ur] \ar[ul]
}
\end{equation}
The vertical/oblique arrows are here the canonical inclusions. We have denoted here by $\underset{p_{\infty\infty} ({\rm all})}{\sum} A(p_{\infty\infty})$ the $(*)$ from (4.47.1). In connection with (4.51), both for $\varepsilon = u$ and for $\varepsilon = b$, we have the obvious equality (written here for II, but valid of course for I too)
$$
S'_{\varepsilon} (M(\Gamma)-H)_{\rm II} \cup \sum D^2 (H({\rm normal})) = S'_{\varepsilon} \, M(\Gamma)_{\rm II} - \sum_{\mbox{\footnotesize$\overbrace{p_{\infty\infty} ({\rm all})}$}} H(p_{\infty\infty}) \, .
$$
Also, since (\ref{eq4.51}) commutes, we can put back the $D^2 (H (p_{\infty\infty}))$'s and get another commutative diagram, the upper line being now a diffeomorphism among cell-complexes (which are certainly non-manifolds)
\begin{equation}
\label{eq4.52}
\xymatrix{
S'_b \, M(\Gamma)_{\rm II} \ar[rr]_-{\approx}^-d  &&S'_u \, M(\Gamma)_{\rm II} \\
&\underset{p_{\infty\infty} ({\rm all})}{\sum} (A(p_{\infty\infty}) \cup D(H(p_{\infty\infty}))) \, . \ar[ur] \ar[ul]
}
\end{equation}
We can invoke now the (\ref{eq4.44}) with its features (4.44.i), and then derive from (\ref{eq4.52}) a further diffeomorphism among cell-complexes, namely the following
\begin{equation}
\label{eq4.53}
S_b \, M(\Gamma)_{\rm II} \underset{\rm DIFF}{=} S_u \, M(\Gamma)_{\rm II} \, .
\end{equation}
At this point, we can invoke the functoriality property of the $S_{\varepsilon}$'s (as stated for $S_u^{(')}$ in the lemma~2.1 and for $S_b^{(')}$ in lemma~4.6, and then taking the universal covering spaces of both terms in the formula (\ref{eq4.53}), get the following diffeomorphisms {\ibf UPSTAIRS} again, in the context of the {\ibf VARIANT~II}
\begin{equation}
\label{eq4.54}
S_u \, \widetilde M(\Gamma)_{\rm II} = (S_u (M(\Gamma))_{\rm II})^{\sim} \underset{(4.53)^{\sim}}{=} (S_b (M(\Gamma))_{\rm II})^{\sim} = S_b \, \widetilde M (\Gamma)_{\rm II} \, .
\end{equation}
When one reads the (\ref{eq4.54}) from the end to the other one gets that $S_u \, \widetilde M(\Gamma)_{\rm II} \underset{\rm DIFF}{=} S_b \, \widetilde M(\Gamma)_{\rm II} \in {\rm GSC}$ by corollary~4.5 above, which means our desired GSC theorem~2.3.

\bigskip

\noindent COMMENTS. But $S_u \, \widetilde M(\Gamma)_{\rm II} / \Gamma = S_u \, M (\Gamma)_{\rm II}$ is {\ibf NOT} compact, so even with $S_u \, \widetilde M(\Gamma)_{\rm II} \in {\rm GSC}$, we are not yet done for $\Gamma \in {\rm QSF}$. The action $\Gamma \times S_u \, \widetilde M(\Gamma)_{\rm II} \to S_u \, \widetilde M(\Gamma)_{\rm II}$ fails to be co-compact.

\smallskip

A third paper in this serie will be necessary, for providing us with the ingredient which we are still missing, i.e. the implication
$$
S_u \, \widetilde M(\Gamma)_{\rm II} \in {\rm GSC} \Longrightarrow \Gamma \in {\rm GSF} \, .
$$
In a nutshell, here is the strategy for this implication. We have a little cascade of implications: $\{S_u \, \widetilde M(\Gamma)_{\rm II} \in {\rm GSC}$, the main result of the present paper$\} \Longrightarrow \{$A certain smooth $4$-manifold, on which the $S_u \, \widetilde M(\Gamma)_{\rm II}$ retracts is $4^d$ Dehn exhaustible$\} \Longrightarrow \{$A certain singular $3^{\rm d}$ manifold, result of a further retraction, is $3^{\rm d}$ Dehn exhaustible$\} \Longrightarrow \{$There exists a certain $3$-complex $X^3$ with the following features: $\pi_1 \, X^3 = 0$, there is a free action $\Gamma \times X^3 \to X^3$ with compact quotient $X^3 / \Gamma$ such that $\pi_1 (X^3/\Gamma) = \Gamma$ and then, finally also, $X^3 \in {\rm QSF}\}$.

\smallskip

From the last item, it follows by standard results that $\Gamma \in {\rm QSF}$. The Dehn-exhaustibility, which will be explicitly defined in the next paper of this series, is a forefather of the QSF concept, and it was used by the present author in the very early nineties \cite{23}, \cite{24} and \cite{25}. In the cascade above, the first two arrows need that the $S_u \, \widetilde M(\Gamma)_{\rm II}$ should be transversally compact, i.e. they make use of the fact that we safely stay in the good side of the Stallings barrier. Somehow dually to this, in the present paper we are also striving to stay on the good side of the non metrizability barrier. Concerning these two barriers, which a priori push in opposite directions, making the navigation hazardous, see the introduction to the present paper and also the \cite{27}.

\newpage

\section{The proof of the zipping lemma 4.1}\label{sec5}
\setcounter{equation}{-1}

In the present section we will focus on the closed subset
\begin{equation}
\label{eq5.0}
f \, {\rm LIM} \, M_2 (f) \subset f \, \sum \, W({\rm complementary}) \subset fX^2 \, ,
\end{equation}
with $\sum \, W({\rm complementary}) = \sum \, W_{(\infty)} ({\rm BLACK}) + \sum \, W ({\rm RED} \cap H^0)$.

\smallskip

For each of these $W({\rm complementary})$ we introduce a $\Gamma$-invariant splitting
\begin{equation}
\label{eq5.1}
W({\rm complementary}) = N_{\infty}^2 (W)(\mbox{easy part}) \cup \overline N_{\infty}^2 (W)(\mbox{difficult part}),
\end{equation}
s.t. we have
\begin{equation}
\label{eq5.2}
p_{\infty\infty} = N_{\infty}^2 \cap f \, {\rm LIM} \, M_2 (f) \in \partial N_{\infty}^2 \cap \partial \overline N_{\infty}^2 \, .
\end{equation}

\noindent Moreover, whenever $q_1 , q_2 , \ldots \in N_{\infty}^2$ and $\lim q_n = q_{\infty} \in f \, {\rm LIM} \, M_2 (f)$, then $q_{\infty} = p_{\infty\infty}$.

\medskip

$$
\includegraphics[width=15cm]{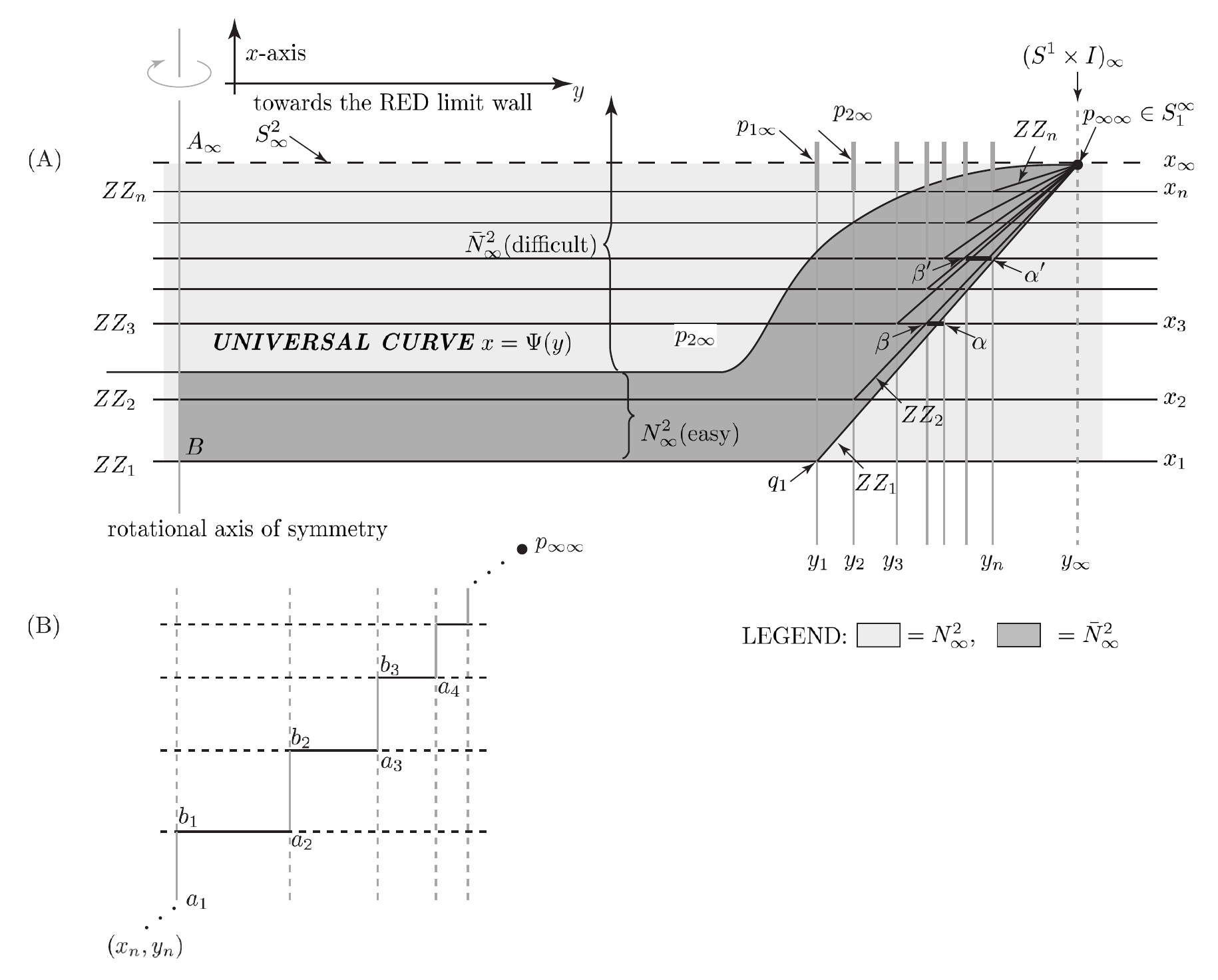}
$$
\label{fig5.1}
\centerline {\bf Figure 5.1.} 

\smallskip

\begin{quote} 
In this figure, the horizontal lines are BLUE, the vertical ones RED.

In (A) we present a piece of $W_{\infty} ({\rm BLACK})_{H_1} = N_{\infty}^2 \cup \partial \overline N_{\infty}^2$, which should be compared to figure~1.3, but be aware that, for typographical commodity's sake, the drawing is not everywhere realistic. The RED lines, dotted or plain do not make it to the $y_1 , y_2 , \ldots , y_{\infty}$ in real life. The dotted piece $[p_{\infty\infty} , y_{\infty}]$ is not realistic, the limit wall $(S^1 \times I)_{\infty}$ stops at $p_{\infty\infty}$, and then goes higher, and similarly the RED walls do not go as far as the $y_1 , y_2 , \ldots$. They actually only start from the various lines $ZZ$, upwards. Similarly, in real life the $W_{\infty} ({\rm BLACK})_{H^1}$ never gets all the way to $[B,A_{\infty}]$ but stops at some $\partial^0 W_{\infty} ({\rm BLACK})_{H^1}$ (figure~1.3), parallel to $[B,A_{\infty}]$, located to the right of it. The lines $ZZ$'s in (A) stand for ``{\ibf zigzags}'', i.e. continuous broken lines of successive red and blue arcs at the level of $X^2$. For typographical reasons, the zigzags are only very schematically represented in (A); they appear more realistically in (B), as well as in the figures~1.2, 1.3. Normally, such a ZZ is a transversal intersection $W_{(\infty)} \cap \partial H^1 \subset \widetilde M (\Gamma)$. The typical $ZZ_n$ in (A) goes from $[B,A_{\infty}]$, along the BLUE level $x=x_n$ to $(x_n , y_n)$, and from there on it really starts zigzagging infinitely many times to the $p_{\infty\infty}$. Inside $H_i^0$, $\underset{n = \infty}{\lim} \, ZZ_n = [A_{\infty} , p_{\infty\infty}] \subset S_{\infty}^2$. In (B) we have $\ldots a_1 \, b_1 \, a_2 \, b_2 \ldots \subset ZZ_n$. The $x$-arrow in the coordinate system corresponds to the increasing $x$'s.
\end{quote}

\bigskip

As drawn, figure~5.1.(A) is part of a $3^{\rm d}$ rotationally symmetric larger object, with symmetry axis $[B,A_{\infty})$. In terms of the figure~1.3, our $W_{\infty} ({\rm BLACK})_{H^1}$ in figure~5.1 has been truncated at $ZZ_1 = \partial H_1^1$ and the $ZZ_i = \partial H_i^1 \cap W_{\infty}$, for $i = 1,2,\ldots$ are drawn too. All the bicollared $H_i^1$'s are attached to a same $H_1^0$, with $ZZ_1$ farthest from $S_{\infty}^2 = \delta H_1^0$. But, on par with figure~5.1, there is also a much denser figure than figure~5.1, {\ibf at the target}, taking into account all the $H_j^0$'s with $g(\infty) \, H_j^0 = g(\infty) \, H_1^0$. In this complete figure we have a doubly infinite family $ZZ_i (j)$, the zigzags from 5.1 being now the $ZZ_i (1)$. There is then a lowest $ZZ_1 (i_0) \subset \partial H_{i_0}^1 ({\rm MAX})$ (figure~1.3), with only finitely many zigzags between $ZZ_i (1) \equiv {\rm our} \ ZZ_1$ and this $ZZ_1 (i_0) \approx \partial H_{i_0}^1 ({\rm MAX})$.

\medskip

\noindent (5.3) \quad For each complete figure~5.1 there is a {\ibf universal curve} $x = \Psi (y)$ inside the $W_{\infty} ({\rm BLACK})$, serving as common boundary $\partial N_{\infty}^2 \cap \partial \overline N_{\infty}^2$.

\medskip

We define ${\mathcal N}_{\infty}^3 \equiv \{$the rotationally invariant $3^{\rm d}$ region between $ZZ_1 (i_0)$ and $x = \Psi (y)\} \subset \widetilde M (\Gamma)$. With this,
$$
N_{\infty}^2 (W_{\infty} ({\rm BLACK})_{H^0}) = W_{\infty} ({\rm BLACK})_{H^1} \cap {\mathcal N}_{\infty}^3 \, , \quad \overline N_{\infty}^2 \equiv \overline{W_{\infty} - N_{\infty}^2} \, .
$$

The universal curve is $C^{\infty}$, infinitely tangent to $S_{\infty}^2$ at $p_{\infty\infty}$. On par with $Y(\infty) \overset{g(\infty)}{-\!\!\!-\!\!\!-\!\!\!\longrightarrow} \widetilde M(\Gamma)$ the whole system of ${\mathcal N}_{\infty}^3$'s and universal curves is $\Gamma$-equivariant.

\bigskip

\noindent {\bf Claim (5.4).} Without any loss of generality, in a complete figure 5.1.(A) we have the following. The part of any $ZZ_i (j)$ living above the universal curve, closer to $x=x_{\infty}$, is always a completely flat BLUE arc; all the actually zigzagging part of $ZZ_i (j)$ lives inside ${\mathcal N}_{\infty}^3$. [Proof. Assume the BLUE levels given. Then we achieve the (5.4) by letting the RED levels converge very fast to their limit position.] This (5.4) will be valid not only for $W_{\infty} ({\rm BLACK})_{H^1} \subset N_{\infty}^2 \cup \overline N_{\infty}^2$, but also for $W_{\infty} ({\rm BLACK})_{H^0}$ and for the $p_{\infty\infty}$-islands of $W({\rm BLACK})$. 

\smallskip

In the case $W_{\infty} ({\rm BLACK})_{H^0}$, which we may assume adjacent to the same $\partial H^1 \cap H^0$ as in the figures~1.3, 5.1, we will take, like above $N_{\infty}^2 (W_{\infty} ({\rm BLACK})_{H^0}) = W_{\infty} ({\rm BLACK})_{H^0} \cap {\mathcal N}_{\infty}^3 =\{$the region contained between $\partial H_{i_0}^1 ({\rm MAX}) \cap W_{\infty} ({\rm BLACK})_{H^0}$ and $x = \Psi (y)\}$, i.e. $\{$the points $(x,y)$ which are such that $\Psi (y) \geq x \geq ZZ_1 (i_0)(y) \}\}$ and also $\overline N_{\infty}^2 (W_{\infty} ({\rm BLACK})_{H^0}) = \overline{W_{\infty} ({\rm BLACK})_{H^0} - N_{\infty}^2}$. The $N_{\infty}^2 (W_{(\infty)} ({\rm BLACK}))$ touches $(S^1 \times I)_{\infty} \cup S_{\infty}^2$ exactly at $p_{\infty\infty} (W_{(\infty)} ({\rm BLACK}))$, while $W_{(\infty)} ({\rm BLACK}) \cap ((S^1 \times I)_{\infty} \cup S_{\infty}^2) \subset \overline N_{\infty}^2 (W_{(\infty)} ({\rm BLACK}))$; in figure~1.2 the $\overline N_{\infty}^2$ is disconnected. We move next to $W({\rm RED} \cap H^0)$, a cylinder bounded by the circles $S^1 = \partial W({\rm RED}) \cap W({\rm BLUE})$ and $S_{\infty}^2 ({\rm BLUE}) \cap W({\rm RED})$. This time there is no point $p_{\infty\infty}$, but the $S^1$ above meets some figure~5.1.(A) and, compatibly with what we have already done before, we shall take now as definition $N_{\infty}^2 (W({\rm RED} \cap H^0)) = W({\rm RED} \cap H^0) \cap {\mathcal N}_{\infty}^3$, with complement the $\overline N_{\infty}^2 (W({\rm RED} \cap H^0))$.

\smallskip

In our figure~5.1, an infinity of pieces $W({\rm RED} \cap H^0)$ (actually their intersection with $W_{\infty} ({\rm BLACK})$), are suggested). They typically go from a bending point of a $ZZ_n$ line (typically like $q_1$ figure~5.1 for $ZZ_1$), up to $S_{\infty}^2$. Remember that pieces like the drawn $[q_1 , y_1]$ are fakes, they are not physically there. This should make it clear that the $N_{\infty}^2 (W({\rm RED} \cap H^0))$ sees only finitely many double lines. But then, our rotationally invariant ${\mathcal N}_{\infty}^3$ (and see here the figures~5.1, 1.2 and 1.3) also  bites into some $W({\rm BLACK})$'s. Generically, these live in $\frac12$-planes $E^2$ containing the same axis of symmetry $[B,A_{\infty})$ from 5.1, but different from the $\frac12$-plane of $W_{\infty} ({\rm BLACK})_{H^1}$, and looking now {\ibf away} from $[B,A_{\infty})$. At the source $X^2$, $W({\rm BLACK})$ is glued to $\{$a bloc $W({\rm RED} \cap H^0) \cup W({\rm BLUE})$, already glued together at the source $X^2\}$, along something like $[a,s_2,b]$, figure~1.3, and to this corner corresponds a $\{p_{\infty\infty}$-island$\} \subset W({\rm BLACK})$. Here, the $p_{\infty\infty}$ of the island in question belongs to the same $S_{\infty}^1 = S_{\infty}^2 \cap (S^1 \times I)_{\infty}$ to which the $p_{\infty\infty}$ in figure~5.1 belongs too. The latter $p_{\infty\infty}$ is like in the figures~1.2, 1.3 and the former like in 1.1. There is an ${\mathcal N}_{\infty}^3$ associated to the $S_{\infty}^1$, both $\Gamma$-invariant and rotationally invariant. Compatibly with all the previous story, we take
\setcounter{equation}{4}
\begin{equation}
\label{eq5.6}
N_{\infty}^2 (\mbox{of the $p_{\infty\infty}$-island}) = {\mathcal N}_{\infty}^3 \cap \{p_{\infty\infty}\mbox{-island}\} \, .
\end{equation}
In this discussion, we talk about a generic $W({\rm BLACK})$ resting on the corner $[a,s_2 , b]$ from figure~1.3 and its $\{p_{\infty\infty}$-island$\}$ from figure~1.1. In 1.1.(A) one sees infinitely many zigzags of equation $x = ZZ(y)$, with $ZZ_{n+1} (y) > ZZ_n (y)$, all infinitely tangent to $x=x_{\infty}$ at $p_{\infty\infty}$, and converging to this $x=x_{\infty}$. The equation (\ref{eq5.6}) above concerns the {\ibf complete} island. Inside the complete island, there are finitely many zigzags (in (\ref{eq5.6})) entering our island through its BLUE side. These are $ZZ_1 = \partial H_{i_0}^1 ({\rm MAX})$ (like in figure~1.3), $ZZ_2 , \ldots , ZZ_{\rho}$; here $ZZ_{\rho}$ is the  highest $ZZ$ entering via BLUE and $ZZ_{\rho+1} = \{ \partial H_i^1 (\gamma_n)$, like in figure~1.3 with respect to the corner $[a,s_2,b] \subset \partial W({\rm BLACK})$, in that figure$\}$. From here on, the $ZZ_{\rho+2} , ZZ_{\rho+3} , \ldots$ enter through the RED side of the island. We also have
\begin{equation}
\label{eq5.7}
\{ ZZ_1 , ZZ_2 , \ldots , ZZ_{\rho+1} \} \subset \{ N_{\infty}^2 (\mbox{of the complete $p_{\infty\infty}$-island}) \subset W({\rm BLACK})\} 
\end{equation}
$$
= \{\mbox{the regions of the points} \ (x \leq \Psi (y) , y) \} \, .
$$
The (\ref{eq5.7}) is no longer valid for the higher $ZZ$'s, while the (5.4) is obeyed by all the zigzags of the complete island.

\smallskip

We move finally to the $\{ \overline S$-region$\} \subset W({\rm BLACK})$, figure~1.1. This corresponds to an immortal singularity which we call $\overline S \subset {\rm Sing} \, \widetilde M (\Gamma)$, reserving the name ``$S$'' for the double infinity of similar, much smaller $S \subset {\rm Sing} \, \Theta^3 (fX^2)$, into which $\bar S$ explodes, like in figure~1.5.(C).

$$
\includegraphics[width=135mm]{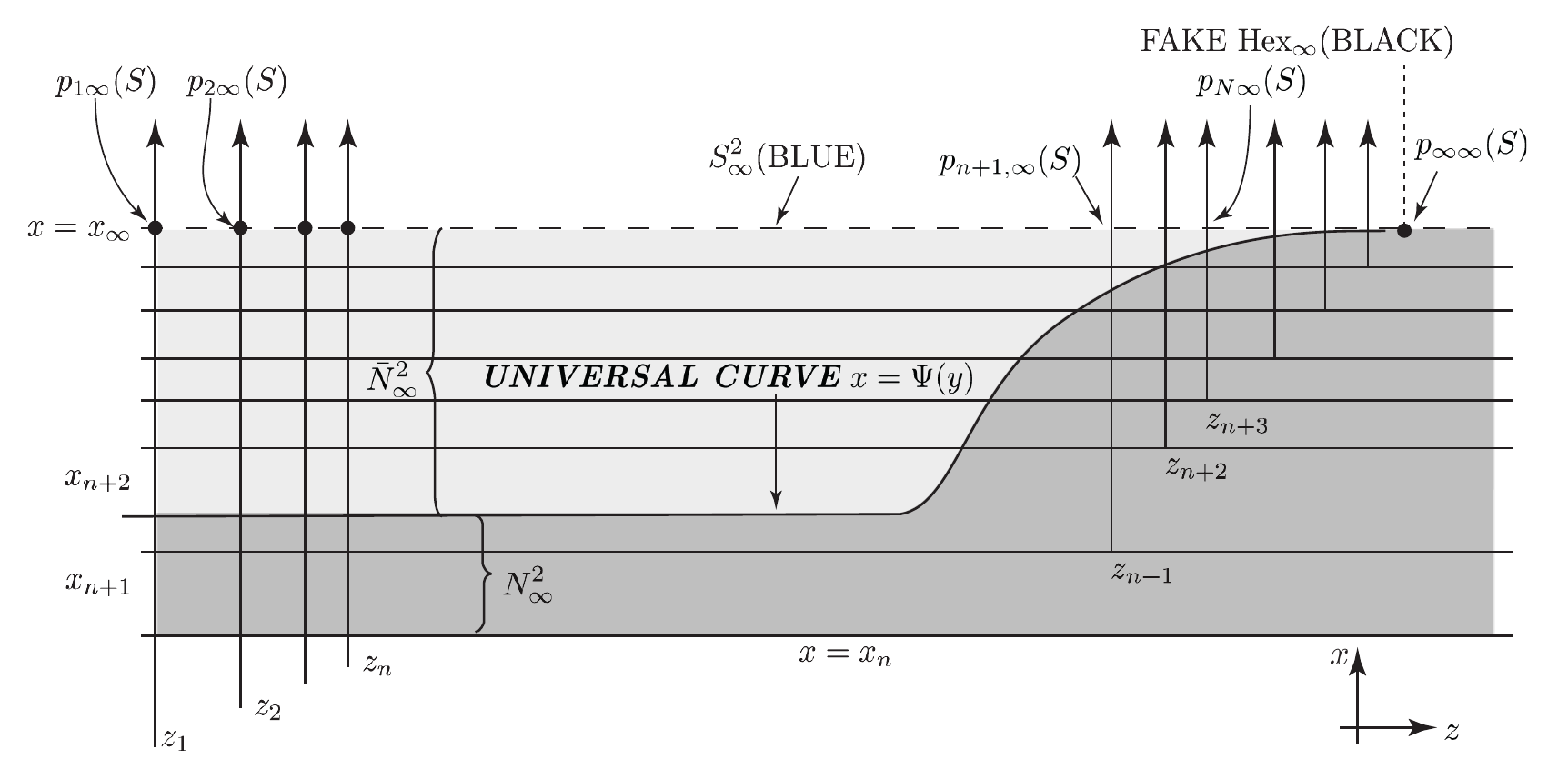}
$$
\label{fig5.1.bis}
\vglue -5mm
\centerline {\bf Figure 5.1.bis} 

\smallskip

\begin{quote} 
We are here in the plane of a $W({\rm BLACK})$ from figure~1.5.(C), corresponding to $y=y_n$ and resting on the BLUE level $x=x_n$. The figure is supposed to be the analogue of a complete figure 5.1.(A), for the case of the immortal $S$'s. We have neither stressed the attaching zones of the various $2$-handles to the present $0$-handle nor have we used the trick of the oblique lines suggesting zigzags, like in figure~5.1.(A). Our $W({\rm BLACK})$ is cut transversally by an infinity of $W({\rm BLACK})^*$'s. Along our $W({\rm BLACK})$, at the intersection sites corresponding to $z_1 , z_2 , \ldots , z_n$ it is $W^*$ which overflows ($\Leftrightarrow W$ subdued), while along the $z_{n+1} , z_{n+2} , \ldots$ it is the other way around. The depths ``$z$'', as drawn, are just a typographical convention.

\smallskip

Legend: (Figure 5.1-bis), we are here in the plane $y=y_n$. The immortal singularities $S$ of $fX^2$ are the various points $p_{n\infty}$, $p_{\infty n}$, and there are $\infty^2$ of them. At $p_{\infty\infty} (S)$, for all $\Theta^3$, $(\Theta^3)'$, $S_{\varepsilon}$, $S'_{\varepsilon}$, there is a Hole $H (p_{\infty\infty} (S))$, inside our $W({\rm BLACK})$.
\end{quote}

\bigskip

Figure~5.1.bis is now the analogue of 5.1. Here there is no longer a previously defined ${\mathcal N}_{\infty}^3$ to hang onto. So, for the whole large $\bar S \subset \widetilde M(\Gamma)$ we will introduce now a {\ibf universal surface}
$$
\{ x = \Psi (y,z) \} \subset \{ x \leq x_{\infty} , (y,z) \in \bar S \ (\mbox{figure 1.5.(C))}\} \, .
$$
All the $3^{\rm d}$ box occurring in the RHS of the formula above, will be referred to, again, as $\bar S$. With this, we will take now
$$
{\mathcal N}_{\infty}^3 = \{ x \leq \Psi (y,z) , y , z \} \subset \bar S \, , \ N_{\infty}^2 (W({\rm BLACK}) \mid \bar S) = W({\rm BLACK}) \cap {\mathcal N}_{\infty}^3 \, .
$$
This ${\mathcal N}_{\infty}^3$, i.e. the universal surface from which it stems, is supposed to satisfy (5.4) as well as the (5.7), (5.8) below. The $\bar S$ is generated by two bicollared handles $H_{\rm I}^2 , H_j^2$ and, for the double infinity of $W({\rm BLACK})$'s stemming from them, we will have

\medskip

\noindent (5.7) \quad $\partial W({\rm BLACK}) \cap \bar S \subset {\mathcal N}_{\infty}^3 \supset \{$the $\{$attaching zones $\partial H_i^2 , \partial H_j^2 \} \cap \bar S \}$.

\medskip

\noindent (5.8) \quad Consider some level $(x=x_{\ell}) = W({\rm BLUE})_{\ell}$ cutting through $\bar S$ and, for expository purposes we will pretend that all the doubly infinite family of $W({\rm BLACK})$'s makes it at least until $x=x_{\ell}$. Inside the large square $\bar S \cap W_{\ell} ({\rm BLUE})$, there is a concentric smaller square $W({\rm BLUE})_{\ell} \cap \overline{\mathcal N}_{\infty}^3$ (figure~5.2), inheriting a finite checker-board. Inside it, we make the choices overflowing/subdued which the figure~5.2 suggest and, afterwards, these choices are continued for $\ell \to \infty$. The ``$O(3)$-lines'' in figure~5.2 refer to our zipping strategy (\ref{eq4.38}).

\smallskip

In the lemma which follows next, we consider a complementary wall $V = W_{(\infty)} ({\rm BLACK})$. Then, going to the context (\ref{eq2.7}), (2.8), for any arc $L_n = V \cap W_n ({\rm BLUE})$, we introduce the quantities
$$
\alpha (n) = {\rm dist} \, (L_n , S_{\infty}^2) \, , \ \beta_n = \# \, \{{\rm points} \ t_{in} \in \overline N^2 \cap L_n \} < \infty \, .
$$

\bigskip

\noindent {\bf Lemma 5.1.} {\it Without loss of generality, there is a uniform bound $P_0$ s.t. for any pair $(V , n \in Z_+)$ we should have}
\setcounter{equation}{8}
\begin{equation}
\label{eq5.10}
\alpha (n) \cdot \beta (n) < P_0 \, .
\end{equation}

\bigskip

\noindent {\bf Proof.} In the context of the figure~5.1, set $\alpha (n) = \vert x_{\infty} - x_n \vert$ and $t_{in} = (x_n , y_i)$. Here the $x_1 , x_2 , \ldots$ are given in the beginning, once and for all.

\smallskip

More precisely, we have infinitely many charts $\{ V_i \}$ where the BLUE walls $V_1 , V_2 , \ldots$ live at levels $x_1 , x_2 , \ldots$, and we will take the sequence $x_1 , x_2 , \ldots$, given once and for all, independent of $\{ V_i \}$. Hence, $\alpha (n) \equiv \vert x_{\infty} - x_n \vert$ will have a universal meaning. This comes with the following universal rule
$$
\{ (x_N , y_i) \in \overline N_{\infty}^2 \} \Longrightarrow \{(x_{N+1} , y_i) \in \overline N_{\infty}^2 \} \, .
$$
For any individual $\{ V_i \}$, we are free to push the vertical lines $y = y_1$, $y=y_2$, $y=y_3,\ldots$ (see figure~5.1) closer and closer to the God-given $y=y_{\infty}$. With these things, it may be assumed that there is a $\{ V_i \}$ independent level $n-1$ such that, for $j < n-1$ we have $\beta (j) = 0$ and that the $\{ V_i \}$-independent quantity $\beta (n-1)$ is positive, $\beta (n-1) > 0$. We define then
$$
P_0 \equiv \beta (n-1) \cdot \vert x_{\infty} - x_{n-1} \vert + 1 \, .
$$

Next, consider the various triple points which are accounted for by $\beta (n-1)$
$$
(x_{n-1} , y_1) , (x_{n-1} , y_2) , \ldots , (x_{n-1} , y_q) , \quad \mbox{when} \ q \equiv \beta (n-1) \, .
$$
From now on, we will keep the values $y_1 , y_2 , \ldots , y_q$ fixed, but we will allow ourselves to bring the next $y$'s closer to $y_{\infty}$. We move now from $\vert x_{\infty} - x_{n-1} \vert$ to the next, smaller $\vert x_{\infty} - x_n \vert$ and let here ${\mathcal J} = {\mathcal J} (n)$ be the largest integer such that ${\mathcal J} (n) < P_0 \cdot \vert x_{\infty} - x_n \vert^{-1}$; we clearly have here ${\mathcal J} (n) \geq q \equiv \beta (n-1)$. With this, we will move all the vertical lines $y_{q+1} , y_{q+2} , \ldots$ closer to $y_{\infty}$, until we have achieved that
$$
(x_n , y_{\mathcal J}) , (x_n , y_{{\mathcal J} + 1}) , \ldots \in N_{\infty}^2 \, .
$$
It follows that $\beta (n) < P_0 \cdot \vert x_{\infty} - x_n \vert^{-1}$ too, and also, by now, an obvious inductive process has been initiated. 

\hfill $\Box$
$$
\includegraphics[width=135mm]{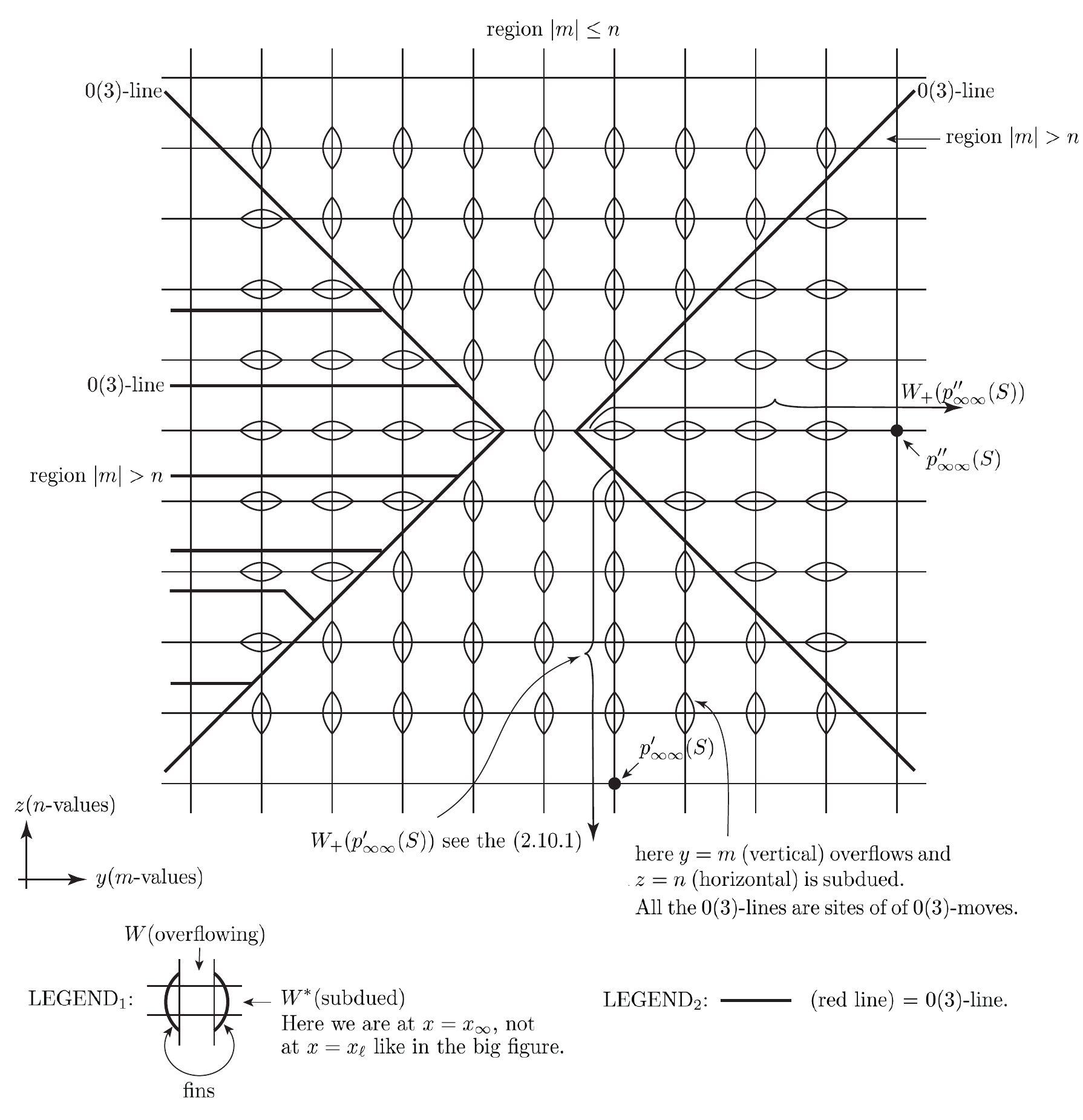}
$$
\label{fig5.2.}
\vglue -9mm
\centerline {\bf Figure 5.2.} 

\smallskip

\begin{quote} 
We see here how the choice overflowing/subdued is made inside the Square $W({\rm BLUE})_{\ell} \cap \bar {\mathcal N}_{\infty}^3$, which our present figure represents.

\smallskip

There are no intersection points $W(m) \cap W (n)^*$ on the boundary of the square nor on the $0(3)$-lines. Around each intersection point, the location of the fins (at $x = x_{\infty}$) is suggested by the round lines. 

\smallskip

The pattern suggested here is to be continuated consistently, for $\ell \to \infty$.

\smallskip

In the LEGEND above, we have thickened the walls, as it is the case in real life.

\smallskip

We have here two main angles, all the lines of which are at $45^{\rm o}$, all of them $0(3)$-lines, but then much more branches, each of them horizontal or vertical, resting on the main angles and going to infinity. They are also $0(3)$-lines, and only few are represented here, for illustration. The main oblique $0(3)$ angles define the distribution of choices subdued/overflowing. Each overflowing site is surrounded by two curved lines, symbolizing the fins (which actually live much higher, at $x=x_{\infty}$, not at the level $x=x_{\ell}$ of the figure). The additional $0(3)$ lines, which are either vertical or horizontal always cut through segments both the endpoints of which are overflowing. After all the $\{ W (n) , W^* (m)\}$ have been zipped together, until they create an immortal singularity $S$ at $x=x_{\infty}$, we start performing the zipping of $W_{\ell} ({\rm BLUE})$ with these $W({\rm BLACK})$'s, stopping momentarily the zipping in question when we meet the boundary of the square. Inside the $\bar S \cap W_{\ell} ({\rm BLUE})$ of our figure, starting from the singularities which have been momentarily created at the boundary, the zipping continues among the zipping flow lines
$$
\sum \{ W(n) , W^* (n)\} \cap \bar S \cap W_{\ell} ({\rm BLUE}) - \{\mbox{the $0(3)$ lines}\} \, ,
$$
with $0(3)$ moves at the lines in question.
\end{quote}

\bigskip

We consider now a $T_i$ like in (\ref{eq2.7}). The following lemma is pretty trivial and it will be largely improved later on.

\bigskip

\noindent {\bf Lemma 5.2.} {\it For each $T_i$ there is a quantity $N = N(T_i)$ such that for all $j > N$ we have $t_{ij} \in M_3 (f)$. Moreover, we can (and will) always drive the zipping flow of $f$ so that it should glue together $V \mid t_{ij}$ and $W^i \mid t_{ij}$ {\ibf before} any action of $W_j ({\rm BLUE}) \mid t_{ij}$; (we refer here to the models {\rm (2.8)}).}

\bigskip

In the context of this lemma and of figure~5.1, for $V = N_{\infty}^2 \cup \overline N_{\infty}^2$, $\partial N_{\infty}^2 \cap \partial \overline N_{\infty}^2 = \{$universal curve$\}$, let
\begin{equation}
\label{eq5.11}
T_i \cap \overline N_{\infty}^2 = [I(T_i) \in \partial N_{\infty}^2 \cap \partial \overline N_{\infty}^2 , p_{i\infty}] \, .
\end{equation}
In agreement with (5.4), any $t_{ij} \in [I(T_i) , p_{i\infty}]$ is actually a triple point, and not a ramification point.

\bigskip

Our next immediate aim is now to render explicitly precise the zipping strategy (\ref{eq4.38}). As a preparation for that, we start with some reminders. To begin with, for the convenience of the reader, we have redrawn a generic zipping path
$$
(x_t , y_t) \in \widetilde M^2 (f) \, , \quad 0 \geq t \geq -\infty \, ,
$$
going from the singularities (let us say at $t=-\infty$) to $(x,y) \in M^2 (f)$ (at, let us say $t=0$), in the figure~5.3; this is, essentially, the same as figure~1.1 in \cite{29}.

$$
\includegraphics[width=12cm]{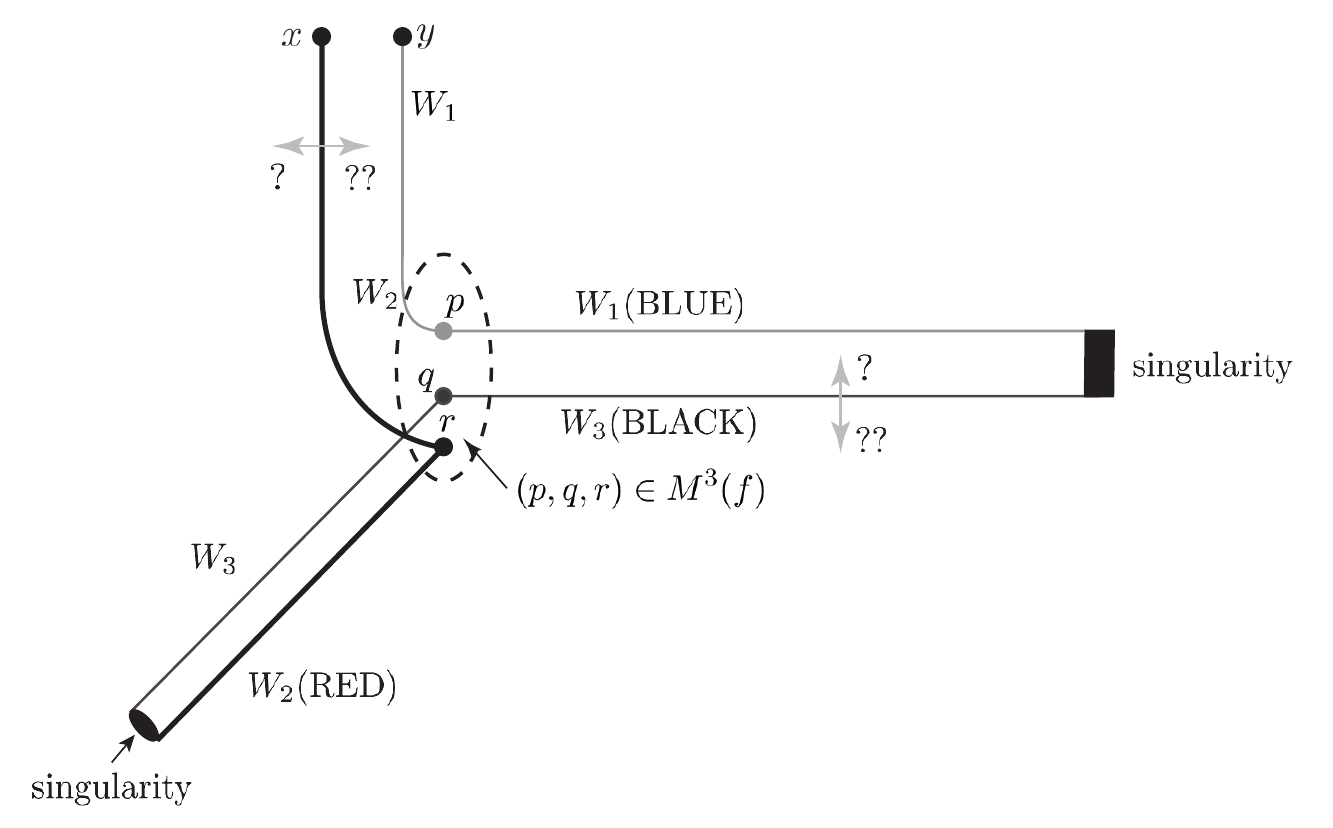}
$$
\label{fig5.3.}
\centerline {\bf Figure 5.3.} 

\smallskip

\begin{quote} 
A schematical view of a zipping path for $(x,y) \in M^2 (f)$, inside $\widehat M^2 (f) = M^2 (f) \cup {\rm Sing} (f)$. Here, for the sake of the argument, imagine that the wall $W_1$ carries Holes, while $W_2$, $W_3$ carry ditches. The arrow stands for the transversal orientations to the complementary walls, going from $+\varepsilon$ to $-\varepsilon$. Since we have not specified the $\pm \, \varepsilon$, we added the question marks. The point is that, with an $n=n(W_j)$, where $W_j$ is non complementary and with $n \to \infty$, (i.e. with ${\rm dist} \, (W_j , \mbox{limit wall}) \to 0$), ditches only get filled between the levels $\varepsilon - \frac1n$ and $\varepsilon$. At the triple point the $W_1 , W_2 , W_3$ meet transversally inside the smooth part of $\widetilde M (\Gamma)$. There are no ditches at the triple point and the special DITCH jumping step is required there (see lemma~4.1).
\end{quote}

\bigskip

As a preliminary for what is coming next, we go back now to the three cases from (2.8), which we will rewrite as follows:
\begin{equation}
\label{eq5.12}
{\rm i)} \ (V({\rm BLACK}) , W^i ({\rm RED}) , W_j ({\rm BLUE})) \, , \quad  {\rm ii)} \ (V({\rm BLACK}) , W^i ({\rm BLACK}), W_j ({\rm BLUE})) \, ,
\end{equation}
$$
{\rm iii)} \ (V({\rm RED}) , W^i ({\rm BLACK}), W_j ({\rm BLUE})) \, .
$$
As a side remark, notice that 
$$
\partial X^2 = \sum \{ \partial X^2 \mid W (\mbox{BLACK, reduced}) \} + \sum \{\partial X^2 \mid W_{\infty} (\mbox{BLACK})_{H^{\varepsilon \leq 1}} \} \, .
$$
It is only the last piece in the formula above which creates $\partial X^2 \cap M_2 (f)$; this involves $W_{\infty} \cap (W({\rm BLUE}) + W({\rm RED}))$. See figures~1.2, 1.3, here.

\smallskip

In the figure~5.4 we have displayed a small area of $\sum (\infty)$, corresponding, in principle, to some $p_{i\infty}$ attached to (\ref{eq5.12}). Our {\ibf policy for locating the fins} is explicitly shown in figure~5.4. Except in (iii)$_2$, (iii)$_3$ when the two branches are already glued in $X^2$, in all the other drawings we are {\ibf after} the preliminary zipping $V+W^i \Rightarrow V \cup W^i$ and the wiggly line, living at some $x(t_{ij}) < x_{\infty}$, should suggest the flow line for the next zipping $W_j ({\rm BLUE}) + V \cup W^i \Rightarrow W_j \cup V \cup W^i$. Figure~5.4 is supposed to be compatible with the lemma~5.3 below. Also, when conflicting arrows have been smeared along some wiggly line, that should mean that in real life one is to be chosen, but we do not yet know which.

\smallskip

The next lemma~5.3, which together with its proof makes (\ref{eq4.38}) explicit, should be considered as an addition to the $2^{\rm d}$ representation theorem~1.1.

\bigskip

\noindent {\bf Lemma 5.3.} {\it For the representation $X^2 \overset{f}{\longrightarrow} \widetilde M (\Gamma)$ one can construct a zipping strategy {\rm (\ref{eq4.38})}, which is equivariant, and such that}

\medskip

\noindent 1) {\it For every $(x,y) \in M^2 (f)$ there is a zipping path $\lambda (x,y) \subset \widehat M^2 (f)$ which has the property that for all $g \in \Gamma$, we have
$$
\lambda (g \cdot (x,y)) = g \, \lambda (x,y) \quad \mbox{(this expresses the equivariance),}
$$
and, moreover, s.t. there is a uniform bound $M$ coming with}
$$
\Vert \lambda (x,y) \Vert < M \, , \quad \forall \, (x,y) \, .
$$

\medskip

\noindent 2) {\it For every $W({\rm BLACK})$, complete or not, there is one arc $[\alpha (\infty),\beta) \subset W({\rm BLACK})$ transversal to the double lines, like in the figure~{\rm 1.1}, with $\alpha (\infty) \in f \, {\rm LIM} \, M_2 (f)$. All the points $[\alpha (\infty),\beta) \, \cap \, M_2 (f)$ will correspond to $0(3)$ moves for our zipping strategy, so that we will have
\begin{equation}
\label{eq5.13}
\lambda (x,y) \cap [\alpha (\infty) , \beta) = \emptyset \, .
\end{equation}

Consider now, for any $W({\rm BLACK})$ the subsets $fM_3 (f) \cap W \subset fM_2 (f) \cap W \subset W$. By resolving every triple point into two disjoined smooth lines, we can perceive $fM_2 (f) \cap W$ as the image of an immersion, which we will denote by $M^2 (f) \cap W \overset{f}{\longrightarrow} W$, where the abstract $M^2 (f) \cap W$ resolves the $fM_2 (f) \cap W = fM^2 (f) \cap W$. With this, we will have that (see figure~{\rm 1.1})
$$
\pi_1 (M^2 (f) \cap W({\rm BLACK}) - [\alpha (\infty) , \beta)) = \pi_1 (f \, {\rm LIM} \, M_2 (f) \cap W({\rm BLACK}) - [\alpha (\infty) , \beta)) = 0 \, .
$$
The FAKE ${\rm LIM} \, M_2 (f)$ is NOT to be taken into account here, and it will be mute throughout this paper.}

\medskip

\noindent 3) {\it At any triple point $t_{ij} \in (V \cap W^i \cap W_j) \cap \overline N_{\infty}^2$, the following things will happen}

\smallskip
\begin{enumerate}
\item[3.1.] {\it The zipping flow will perform the step $V + W^i \Rightarrow V \cup W^i$, before any BLUE $W_j$-action, i.e. before $V \cup W^i + W_j \Rightarrow V \cup W^i \cup W_j$.}

\smallskip

\item[3.2.] {\it For this second zipping move $V \cup W^i + W_j \Rightarrow V \cup W^i \cup W_j$, one always gets to the site $t_{ij}$ via the road ${\rm BLUE} \cap {\rm BLACK}$ (i.e. $W_j \cap V$), and {\rm never} via the road ${\rm BLUE} \cap {\rm RED}$.}
\end{enumerate}

\medskip

\noindent 4) {\it Independently of the $M$ from $1)$ and of the $P_0$ from {\rm (\ref{eq5.10})}, there is also a third uniform bound $P$, s.t.}
$$
\# \, \{ \lambda (x,y) \cap [\{\{ p_{\infty\infty} \, \mbox{islands} \} \subset W_{(\infty)} ({\rm BLACK}) \} + \{ \bar S \, \mbox{regions} \subset W({\rm BLACK})\}] < P \, .
$$

\medskip

\noindent 5) {\it We consider now the $1$-skeleton of the $3$-dimensional $Y(\infty)$ {\rm (\ref{eq1.6})}, i.e. the
$$
Y(\infty)^{(1)} = \bigcup_{\overbrace{\mbox{\footnotesize$i,\gamma;\lambda \leq 1$}}} H_i^{\lambda} (\gamma) \subset Y(\infty) \, .
$$
Next, for the bicollared handles $H_i^{\lambda} (\gamma)$ we consider the completions $\widehat H_i^{\lambda} (\gamma) = H_i^{\lambda} (\gamma) \cup \delta H_i^{\lambda} (\gamma)$ and, finally, we introduce the {\ibf ``ideal surface''}}
\begin{equation}
\label{eq5.14}
\delta_{\infty} \, Y(\infty)^{(1)} \equiv \bigcup_{\overbrace{\mbox{\footnotesize$i,\gamma ; \lambda \leq 1$}}} \delta H_i^{\lambda} (\gamma) \subset \widehat Y (\infty)^{(1)} \equiv \bigcup_{\overbrace{\mbox{\footnotesize$i,\gamma ; \lambda \leq 1$}}} \widehat H_i^{\lambda} (\gamma) \overset{f}{\longrightarrow} \widetilde M (\Gamma) \, ,
\end{equation}
{\it and here we have (see {\rm (\ref{eq1.10})} and {\rm (\ref{eq1.14})})
\begin{equation}
\label{eq5.15}
\Sigma_1 (\infty) = f (\delta_{\infty} \, Y(\infty)^{(1)}) \supsetneqq \partial (\widetilde M (\Gamma)^{(1)}) \equiv \partial \{1\mbox{-skeleton of} \ \widetilde M (\Gamma) \} \, .
\end{equation}
In the formula {\rm (\ref{eq5.14})}, we also have $f \delta H_i^{\lambda} (\gamma) = \delta h_i^{\lambda}$, $\delta h_i^0 \subset f \delta_{\infty} \, Y(\infty)^{(1)}$, but then we find that $\delta h_i^0 - \underset{j}{\sum} \, h_j^1 = \delta h_i^0 \cap \partial (\widetilde M (\Gamma)^{(1)})$ which should explain the $\supsetneqq$ occurring in {\rm (\ref{eq5.15})}; the $\partial (\widetilde M^{(1)})$ misses the $\delta h_i^0 \cap h^1$. The $\delta_{\infty} \, Y(\infty)^{(1)} \cap X^2$ splits $X^2$ into a main piece $X_0^2 \supset {\rm int} \, X_0^2 \supset M_2 (f)$ and the rest, $X^2 - X_0^2$ which consists of the various pieces of $W_{(\infty)} ({\rm BLACK})$ living on the other side of $\sum (\infty)$, at the level of the figures~{\rm 1.2, 1.3}. There is also an intersection $f \delta_{\infty} \, Y(\infty)^{(1)} \cap f ({\rm int} \, X_0^2) \ne \emptyset$, to be compared to the ``$\supsetneqq$'' in  {\rm (\ref{eq5.15})}. But this {\rm 5)} is only a preliminary for what comes next.}

\medskip

\noindent 6) {\it Let $(x,y) \in M^2 (f)$ and let $\lambda (x,y)$ be its zipping path from $1)$. This comes with}
$$
\lambda (x,y) \subset {\rm int} \, X_0^2 \times {\rm int} \, X_0^2 \ \begin{matrix} \nearrow \ ^{\mbox{$X_0^2 \times X_0^2$}} \hfill \\ \searrow \ _{\mbox{$\widehat Y (\infty)^{(1)} \times \widehat Y (\infty)^{(1)} \supset \delta_{\infty} \, Y(\infty)^{(1)} \times \delta_{\infty} \, Y(\infty)^{(1)} \, .$}} \end{matrix}
$$

\medskip

\noindent (5.15) \quad {\it Fix now a compact $K \subset \widetilde M (\Gamma)$; then, for any $\eta > 0$ there is an $\varepsilon > 0$ such that if $f(x) = f(y) \in K$ and ${\rm dist} \, (f(x) , f\delta_{\infty} \, Y(\infty)^{(1)}) < \varepsilon$, then we also have, for all $(u_t , v_t) \in \lambda (x,y)$}
$$
\min \left\{{\rm dist} (u_t , \delta \, Y(\infty)^{(1)} \cap {\rm int} \, X_0^2) \, , \ {\rm dist} (v_t , \delta(\infty)^{(1)} \cap {\rm int} \, X_0^2) \right\} < \eta \, .
\eqno (*)
$$
{\it All our metric structures here are descending from the equivariant metric of $\widetilde M (\Gamma)$ and, for any $(u,v) \in M^2 (f)$ we have that
$$
\left[ \min \left\{{\rm dist} (u , \delta \, Y(\infty)^{(1)} \cap {\rm int} \, X_0^2) \, , \ {\rm dist} (v , \delta \, (\infty)^{(1)} \cap {\rm int} \, X_0^2) \right\} \mbox{(like in $(*)$ above)} \right] =
$$
$$
\left\{ {\rm dist} \left(f(u) , \underset{\overbrace{\mbox{\footnotesize $i;\lambda \leq 1$}}}{\bigcup} \delta \, h_i^{\lambda}\right) \, , \ \mbox{in} \ \widetilde M (\Gamma) \right\} .
$$
This ends our item} (5.15).

\smallskip

{\it But then, here is still another way of reading {\rm (5.15)}, which should be useful. If, at level $\widetilde M(\Gamma)$, a double point $(x,y)$ belongs to the given compact $K$ and is at distance $< \varepsilon$ from $\underset{i;\lambda \leq 1}{\bigcup} \delta \, h_i^{\lambda}$, then considered now inside $\Theta^3 (fX^2)$, the image of the zipping path $\lambda (x,y)$ is at distance $< \eta$ from $\sum (\infty) \subset \Theta^3 (fX^2)$. Or, again, now in a more impressionistic language, in appropriate topologies, we find that
$$
\lim_{n = \infty} (x_n , y_n) = \infty \ \mbox{implies that} \ \lim_{n=\infty} \lambda (x_n , y_n) = \infty \, , \ \mbox{too.}
$$
}

\bigskip

Before going into the proof of this lemma, I will offer some comments. The regions concerned by 4) are the same as the ones in lemma~5.1, but $P,P_0$ are independent of each other. With this, the $M$ in 1) is the one from the uniformly bounded zipping length in \cite{29}, and the uniform bound $\Vert \Lambda (H_i^-) \Vert < N$ in the lemma 4.7 is, essentially, the
$$
N \cong kM +2\pi \max \{ \mbox{diameter of Hole}\} + \{\mbox{the product} \ P \cdot P_0 \} \, ,
\eqno (5.15.1)
$$
where $k$ is some small controlled quantity, where $P$ accounts for the number of special regions of type 4) in lemma~5.3, which the zipping path $\lambda (x,y)$ may have to cross, and $P_0$ (lemma~5.1) accounts for the additional complications which each individual crossing may being.

\smallskip

We have written ``$\cong$'' in the formula (5.15.1) above, since some additional ingredients may still have to be thrown in, before we can write ``$=$''. These come from the fact that in the formula (\ref{eq4.50}), which defines the curve $\Lambda (H_i^-)$, the hole $H_i$ occurs twice and, moreover, both $\beta C^- (H_i)$ {\ibf and} $\gamma_i$ may involve the not yet defined process (zipping)$^{-1} [XY]$. It is, anyway, the (zipping)$^{-1} [XY]$ which brings the contribution $M$ (or some fixed multiple of it) into (5.15.1), and it is its passage through the special regions which counts, rather than the one of $\lambda (x,y)$ which is only a first approximation of (zipping)$^{-1} [XY]$. But the general idea should be clear: as long as we control the quantities $M,P$ and $P_0$, we have a good uniform bound for the lengths or the curves $\Lambda (H_i^-)$ from (\ref{eq4.50}).

\smallskip

We come back now to figure~5.4. The drawings (i), (ii), (iii)$_1$ refer to triple points $t_{ij}$, and they should illustrate the point 3.2) in lemma~5.3. The (iii)$_2$ does not concern triple points at all. Then, there should be a (iii)$_3$, which is not explicitly drawn, corresponding to a situation $V({\rm RED}) \cup W_{(\infty)}^i ({\rm BLACK})$ glued already at the source, let us say for instance a RED portion of $\partial H_i^1 (\gamma_n)$ in figure~1.3 glued already at level $X^2 \mid H_i^1$. This does not correspond then to a triple point either. This drawing (iii)$_3$ should be like $\{({\rm iii})_1$, without ``$t_{ij}$''$\}$. Notice (and see here the figure~1.3), that the (iii)$_3$ should always live in the ${\mathcal N}_{\infty}^3$.

$$
\includegraphics[width=140mm]{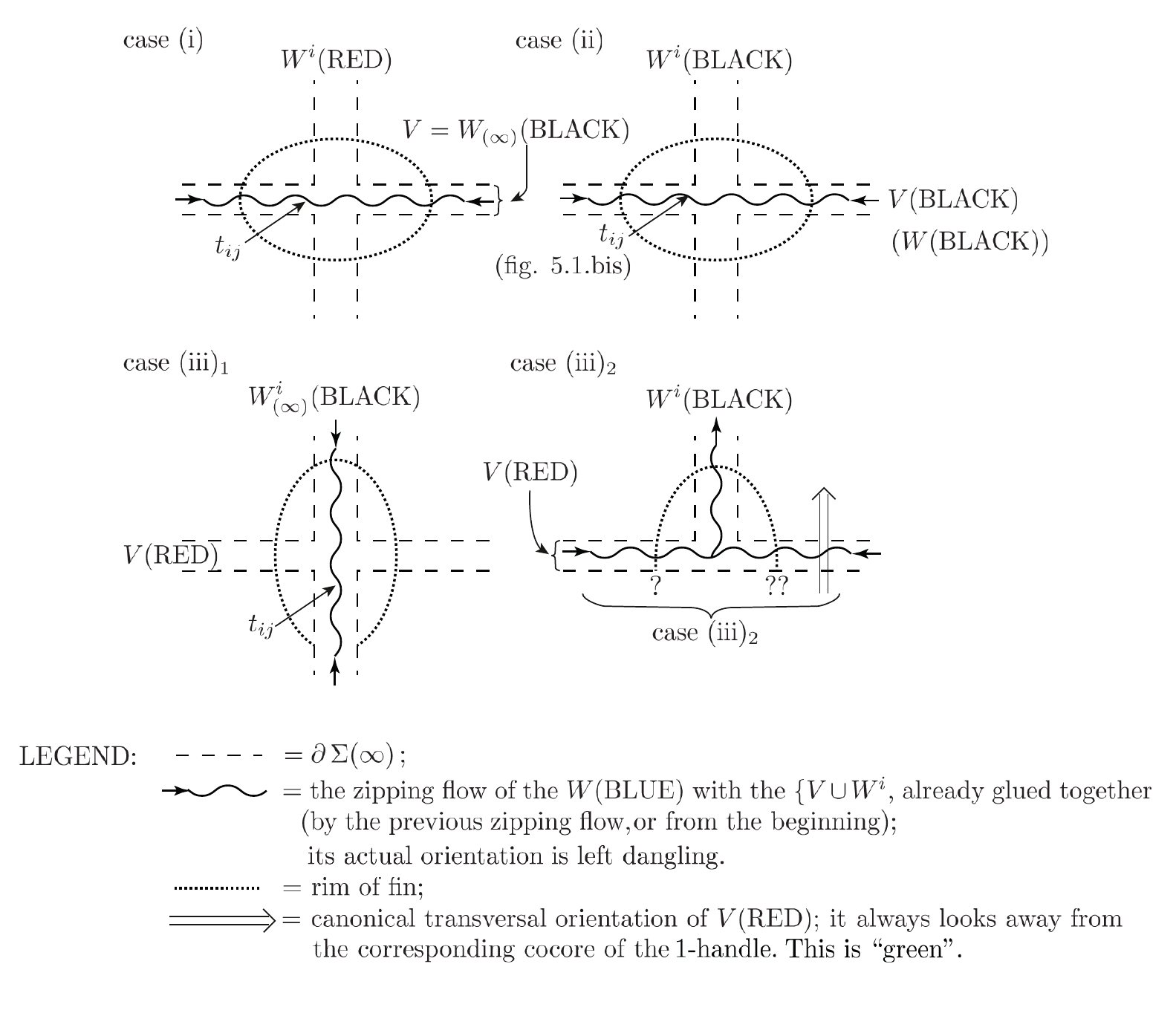}
$$
\label{fig5.4.}
\centerline {\bf Figure 5.4.} 

\smallskip

\begin{quote} 
We display here our policy for locating the fins, corresponding to the three cases (i), (ii), (iii) from (\ref{eq5.12}).
\end{quote}

\bigskip

\noindent {\bf An addendum to lemma 5.3.} In the context of 3) the $V , W^i$ are both complementary walls without Holes but with DITCHES; the preliminary zipping from 3.1) can be performed simple-mindedly with these ditches completely filled. Also inside $N_{\infty}^2 , {\mathcal N}_{\infty}^3$ the zipping is always simple-minded, with completely filled ditches and, a priori, at least, without restrictions for the order of operations. This is not so for $\overline N_{\infty}^2$. Next, remember that $W({\rm RED} \cap H^0) \cap \{ p_{\infty\infty} \} = \emptyset$; but for the other complementary walls, i.e. for the $W_{(\infty)} ({\rm BLACK})$'s we find that $(W_{(\infty)} \cap {\rm LIM} \, M_2 (f)) \cap {\mathcal N}_{\infty}^3 = N_{\infty}^2 (W_{(\infty)}) \cap {\rm LIM} \, M_2 (f) = \{$the points $p_{\infty\infty}$, which are all isolated$\}$.

\smallskip

When we are in the neighbourhood of a $p_{\infty\infty}$, we work with completely filled ditches both for $W_{(\infty)}$ $({\rm BLACK}) \cap W({\rm RED} \cap H^0)$ and for $N_{\infty}^2 (W_{(\infty)} ({\rm BLACK}))$. This opens potential dangers, taken care of by the fact that we make use of $S'_{\varepsilon} (\widetilde M (\Gamma) - H)$ and/or $S'_{\varepsilon} (M(\Gamma) -H)$. This also allows us to appeal to the lemma~4.3.

$$
\includegraphics[width=13cm]{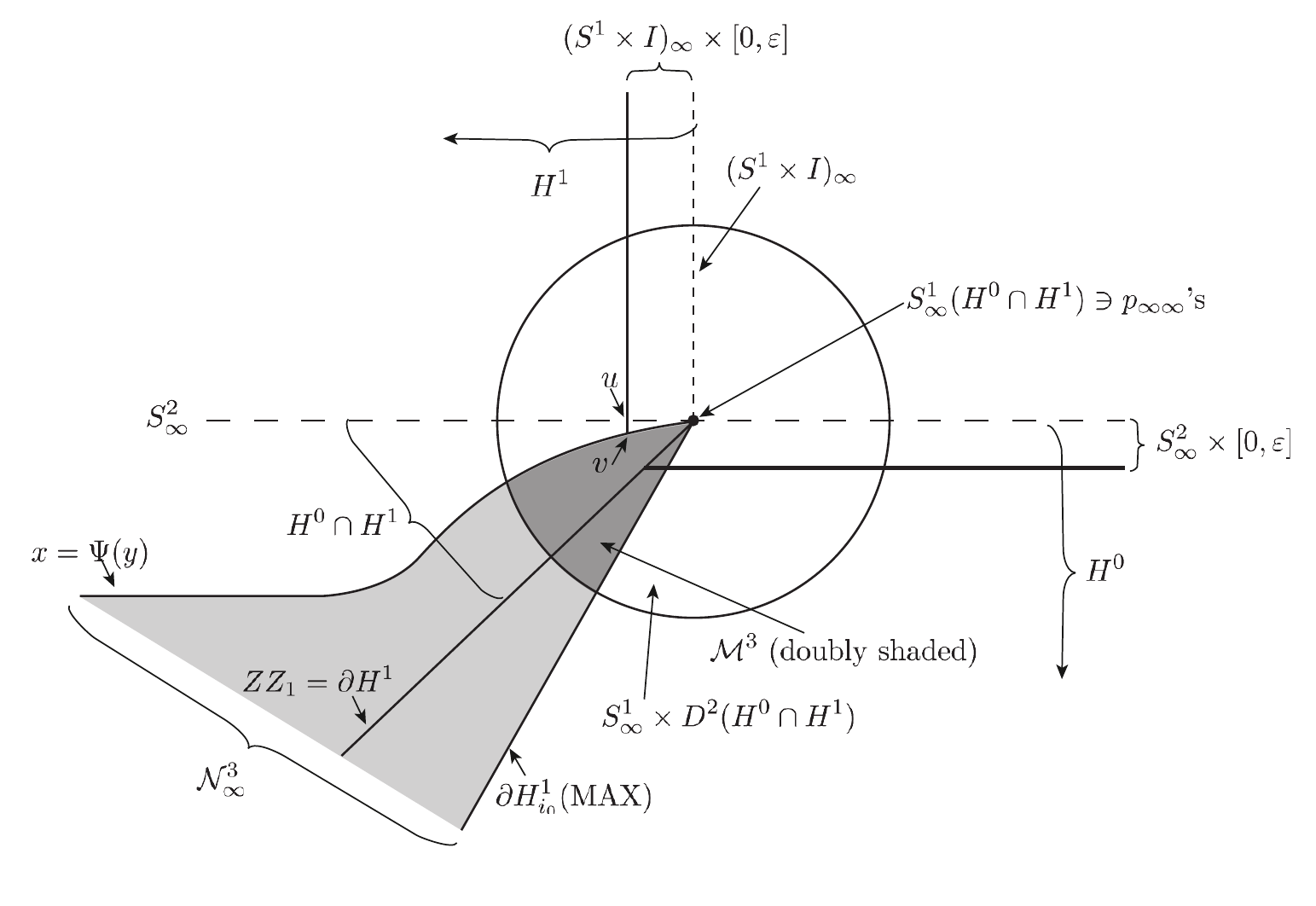}
$$
\label{fig5.5.}
\centerline {\bf Figure 5.5.} 

\smallskip

\begin{quote} 
This figure refers to the figures 1.1, 1.2, 1.3, 5.1 and there is only one such figure, both $\Gamma$-invariant and rotationally invariant (when appropriately extended), for each circle $S_{\infty}^1$. The $S_{\infty}^1 \times D^2 (H^0 \cap H^1)$ is a very thin tubular neighbourhood of $S_{\infty}^1$, biting out of ${\mathcal N}_{\infty}^3$ the doubly shaded region ${\mathcal M}^3$.
\end{quote}

\bigskip

\noindent {\bf The proof of lemma 5.3.} For every pair of adjacent handles $\{ H^0 , H^1 \} \subset Y(\infty)$ we consider the corresponding circle at infinity $S_{\infty}^1 = S_{\infty}^1 (H^0 \cap H^1)$ and then, at the level of $\widetilde M (\Gamma)$ and/or of $Y(\infty)$, a $\Gamma$-invariant tubular neighbourhood of it $S_{\infty}^1 \times D^2 (H^0 \cap H^1)$, like in figure~5.5 and, using it, the following $\Gamma$-invariant region of the $\{$1-skeleton of $Y(\infty)\} \mid X^2$, (the notations are here like in figure~5.5)
\setcounter{equation}{15}
\begin{equation}
\label{eq5.17}
{\mathcal R}^3 = \sum_{H^0 \cap H^1} {\mathcal M}^3 \cup \sum_{H^0} S_{\infty}^2 \times [0,\varepsilon] \cup \sum_{H^1} (S^1 \times I)_{\infty} \times [0,\varepsilon] \, .
\end{equation}

At this point, I will also refer to the section~4 of \cite{29}, the first paper of the present series, where in the $3^{\rm d}$ context of $Y(\infty)$, a zipping strategy of uniformly bounded zipping length has been described; see here the lemma~4.1 in the paper \cite{29}. The essential part of the strategy in question was concentrated inside $Y(\infty)^{(1)}$. With this, our present proof will proceed in several steps.

\bigskip

\noindent {\bf Step I.} We start with the following

\bigskip

\noindent {\bf Claim (5.17)} Let $(x,y) \in M^2 \left( f \mid \sum {\mathcal M}^3 \cap X^2 \right)$ where $\sum$ is over all the blocs $\{ H^0 , H^1 \}$ and which is such that $(x,y)$ does not involve the $H^2$'s; these will be dealt with later in our proof. For this $(x,y) \in {\mathcal M}^3$, there is then a $\Gamma$-invariant, uniformly bounded zipping path ($=$ strategy),  [which stays confined inside ${\mathcal R}^3 \cap X^2$ (see here (\ref{eq5.17}))] constructed like in \cite{29}, and this construction is such that for our $(x,y) \in {\mathcal M}^3$ the controlled zipping path can be chosen inside the ${\mathcal R}^3$ from (\ref{eq5.17}).

\smallskip

Inside the doubly shaded ${\mathcal M}^3$, displayed in figure~5.5, there are no restrictions for the strategy in question: Ditches can be filled completely, no order is imposed for the zipping steps, a.s.o. Now, when we are inside the regions $S_{\infty}^2\times [0,\varepsilon]$, $(S^1 \times I)_{\infty} \times [0,\varepsilon]$, which should be understood as resting on ${\mathcal M}^3$ but with interiors disjoined from it, then the zipping takes always the following form, for our presently discussed strategy in Step~I
$$
W (\mbox{of NATURAL COLOUR}) \cap W_{\infty} ({\rm BLACK}) \, ,
$$
and ditches should now be filled only partially (for the $W_{\infty} ({\rm BLACK})$'s). Also, inside $((S^1 \times I)_{\infty} \times [0,\varepsilon])$ $\cap$ $(x < x_{\infty})$ when along each double line infinitely many triple points will be met, on the way, these will be ignored, for the time being. This {\ibf ENDS} our claim~(5.17), and the {\ibf proof} follows. \hfill $\Box$

\bigskip

\noindent {\bf Proof.} The idea here is to adapt the strategy from lemma~4.1 in \cite{29}, to our present $2$-dimensional situation. This is made possible by the fact that in our present $X^2$, the security walls $W_{\infty} ({\rm BLACK})$ {\ibf overflow}, in particular the $W_{\infty} ({\rm BLACK})_{H^0}$ in figure~1.2 extends on the $H^1$ side of $\partial H_i^1 (\gamma_n)$, thus catching a whole infinite RED/BLUE checkerboard, and then also something similar for $W_{\infty} ({\rm BLACK})_{H^1}$ in the figure~1.3. Remember that this overflowing of the security walls $W_{\infty} ({\rm BLACK})$ is one important instance where the present paper deviates from \cite{31}, which nevertheless stays, largely, very adaptable for our present purposes.

\smallskip

Here is an illustration. Consider a bloc $h_A^1 \cup h_L^0 \cup h_B^1 \subset \widetilde M (\Gamma)$, covered by two blocs $H_a^1 \cup H_{\ell}^0 \cup H_b^1 + H_{\alpha}^1 \cup H_{\lambda}^0 \cup H_{\beta}^1 \subset X^2$, and which are such that, at the other end of $(H_a^1 , H_{\alpha}^1)$, a zipping of type $W_{\infty} (H_a^1) \cap \{$very high level $W({\rm RED})_{\alpha}\}$ has already started. Our aim now is to exhibit a $2$-dimensional zipping flow line, abiding to (5.17), confined inside ${\mathcal R}^3$ and going like
$$
A \longrightarrow L \longrightarrow B \, .
$$
The initial zipping above will reach a spot $(W({\rm RED})_{\alpha} \underset{\overbrace{\mbox{\footnotesize$\partial W ({\rm RED})_{\alpha}$}}}{\cup} W({\rm BLUE})_{\lambda}) \cap W_{\infty} (H_a^1)$, where from the obvious newly created singularity, we continue with a zipping $W({\rm BLUE})_{\lambda} \cap W_{\infty} (H_a^1)$. This will reach a spot $W({\rm BLUE})_{\lambda} \cap (W_{\infty} (H_a^1) \cup W({\rm RED})_a)$, with the ``$\cup$'' created by some short, ${\mathcal N}_{\infty}^3$-confined vertical zipping from a singularity $s'$ in a figure of type~1.3. We can continue now with a zipping $W({\rm BLUE})_{\lambda} \cap W({\rm RED})_a$, until we reach a spot $W({\rm BLUE})_{\lambda} \cap (W({\rm RED})_a \cup W_{\infty} (H_{\lambda}^0))$, with the ``$\cup$'' created now by an ${\mathcal N}_{\infty}^3$-confined short vertical zipping from a singularity $s'$ in a figure~1.2. This zipping concerns the {\ibf overflow} of $W_{\infty} (H_{\ell}^0)$ mentioned above. Next, we move with a zipping $W({\rm BLUE})_{\lambda} \cap W_{\infty} (H_{\ell}^0)$ from the side $A$ of $h_L^0$ to the side $B$. Notice that the last described move has happened close to $p_{\infty\infty} (W_{\infty} (H_{\ell}^0))$ and also that we move now from ${\mathcal N}_{\infty}^3$ to $\overline N_{\infty}^2$. Several such moves, back and forth between ${\mathcal N}_{\infty}^3$ and $\overline N_{\infty}^2$, in the neighbourhood of the $p_{\infty\infty}$'s, will be part of the zipping flow story. For the convenience of the reader, in figure~5.5.bis we have presented in the form of a chart, that part of the proof of the claim~(5.17), developed so far. Anyway, when continued, our zipping $W({\rm BLUE})_{\lambda} \cap W_{\infty} (H_{\ell}^0)$ will encounter (see figure~5.5.bis) a spot $W({\rm BLUE})_{\lambda} \cap (W_{\infty} (H_{\ell}^0) \cup W ({\rm RED})_b)$, where the ``$\cup$'' has been created, again by a vertical zipping from some $s'$, in a figure~1.2. This will unleash a new zipping $W({\rm BLUE})_{\lambda} \cap W ({\rm RED})_b$, getting to a spot $(W({\rm BLUE})_{\lambda} \cap W_{\infty} (H_{\beta}^1)) \cap W({\rm RED})_b$, with the ``$\cup$'' coming from a short, ${\mathcal N}_{\infty}^3$-confined horizontal zipping from some $s'$, figure~1.3. So, by now we have managed to create a contact $W_{\infty} (H_{\beta}^1) \cap W({\rm RED})_b$ and hence gotten from $A$ to $B$. 

$$
\includegraphics[width=14cm]{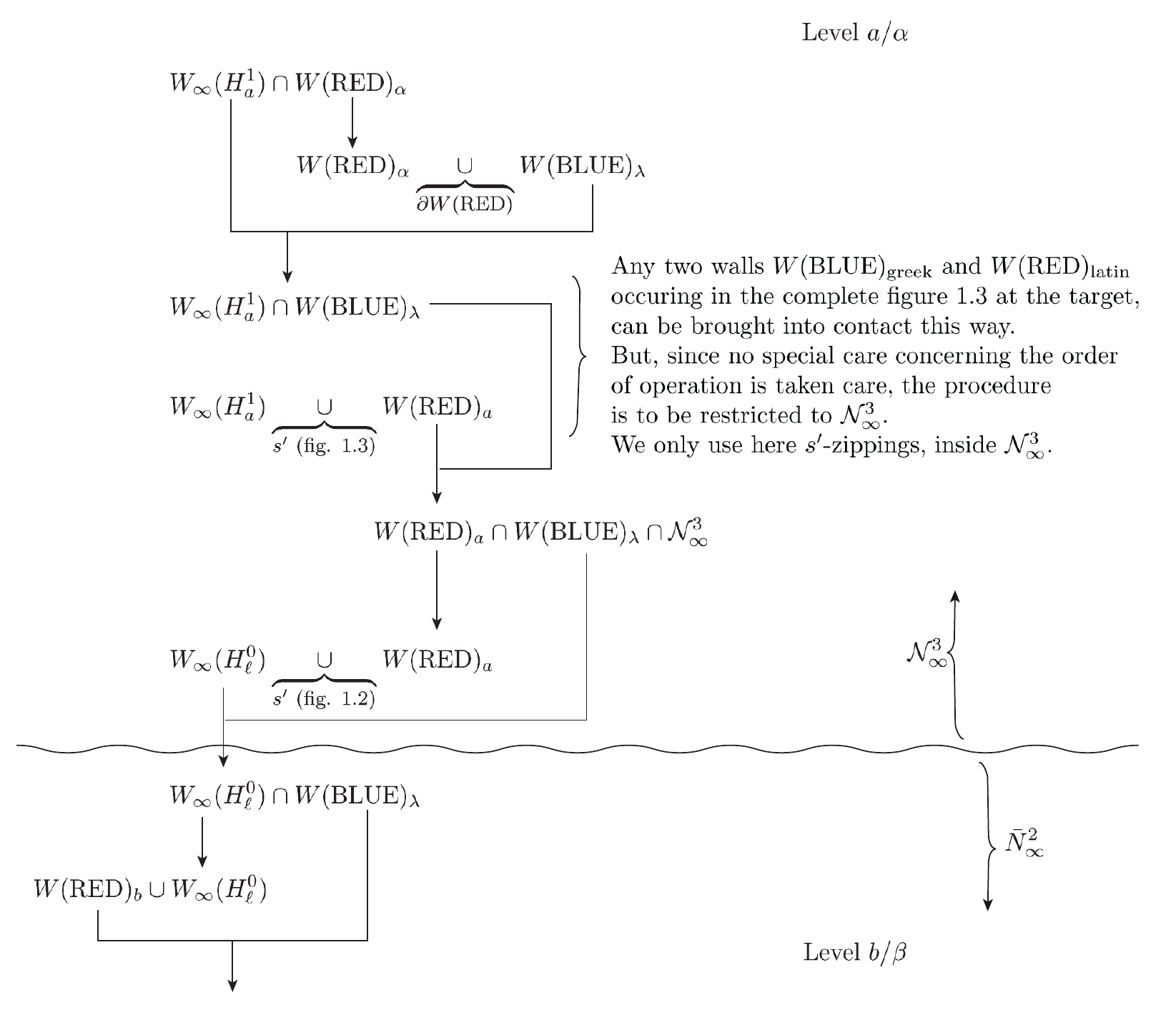}
$$
\label{fig5.5.bis}
\centerline {\bf Figure 5.5.bis} 

\bigskip

This ends our discussion of (5.17). From now on, using (5.17), the rest of the zipping in lemma~5.3 will aways take place in a neighbourhood of controled diameter of ${\mathcal R}^3$, in a way which should be compatible with the zipping strategy for $Y(\infty) \overset{g(\infty)}{-\!\!\!-\!\!\!-\!\!\!-\!\!\!\longrightarrow} \widetilde M (\Gamma)$ in \cite{29}. All the $\lambda (x,y)$'s in 1) will be created this way.

\bigskip

\noindent {\bf Step II.} We move now from ${\mathcal R}^3$ to the zigzags occurring below, and from which one should see the figures~1.2, 1.3, 5.1,
$$
{\mathcal R}^3 \supset \sum_{H^0 \cap H^1} {\mathcal M}^3 \subset \sum_{H^0 \cap H^1} {\mathcal N}^3_{\infty} \supset \{\mbox{all the zigzags} \ ZZ_n\} \, .
$$
The $\sum {\mathcal M}^3$ which, with the exception of the $W({\rm BLACK})$'s, is by now completely zipped contains all the singularities, which at the present level of the construction control the double lines
$$
ZZ_n (\mbox{union of BLUE and RED pieces glued in} \ X^2) \cap W_{\infty} ({\rm BLACK})_{H^{\varepsilon}} \, , \eqno (*)
$$ 
and hence all these $(*)$'s can now be zipped. This brings to life the singularities which correspond to the following double lines
$$
W({\rm RED}) \cap W_{\infty} ({\rm BLACK}) \, , \quad W({\rm BLUE}) \cap W_{\infty} ({\rm BLACK}) \mid {\mathcal N}_{\infty}^3 \, . \eqno (*)_1
$$ 
So, staying completely inside ${\mathcal N}_{\infty}^3$, at least when the $W({\rm BLUE})$'s are concerned, we can zip the $(*)_1$.

\bigskip

\noindent {\bf Step III.} Starting from the zipping-flow chart in figure~5.5.bis, one can see that all the double lines $W({\rm RED})$ $\cap$ $W({\rm BLUE}) \cap {\mathcal N}_{\infty}^3$ can by now be zipped. Moreover, all the singularities which command the double lines $W ({\rm BLACK}) \cap ZZ_n$ are all confined inside ${\mathcal N}_{\infty}^3$ (see here, for instance, the figure~1.1.(B) where the lines in question are on the left side of the universal curve). So, these singularities can be brought now to life and we can zip now all the double lines $W ({\rm BLACK}) \cap ZZ_n$. It is only now that $W ({\rm BLACK})$ enters our game. [Typically, the $[a,s_2,b]$ in the figure~1.3, but in a plane different from the one of $W_{\infty} ({\rm BLACK})_{H^1}$ (of the figure in question), is a corner of the attaching zone of $W ({\rm BLACK})$ to the rest of $X^2$.] So, by now we have access to all the singularities commanding lines $W ({\rm BLACK}) \cap W ({\rm RED})$, all of them in ${\mathcal N}_{\infty}^3$. The $W ({\rm BLACK}) \cap W ({\rm RED})$ can be zipped now, the 3.1) is established and 4), 5), 6) are left to the reader.

\bigskip

\noindent {\bf Step IV.} We take care now of 2) and of 3.2), with the exclusion of the $\{\bar S\mbox{-regions}\} \subset W ({\rm BLACK})$ which are left for the next step. The zipping $W({\rm RED} \cap H^0) \cap (W_{\infty} ({\rm BLACK}) +W ({\rm BLACK}))$ has created by now all the figures~5.4.(iii)$_1$. But our discussion has, by now also certainly taken us beyond ${\mathcal N}_{\infty}^3$. It may be assumed that on each circle $W({\rm RED}) \cap W({\rm BLUE})$ the points of intersection with $W_{\infty} ({\rm BLACK})$ are a dense enough system, such that each figure~5.4.(iii)$_2$ should be sandwiched between two 5.4.(iii)$_1$'s, with $W_{\infty} ({\rm BLACK})$ in the position marked ``$W_{(\infty)}^i ({\rm BLACK})$'', in figure~5.4.(iii)$_1$.

\smallskip

On each $W_{\infty} ({\rm BLACK})$ the zippings $(*)_1$ from Step~II have already been done and we can zip up now completely the $W_{\infty} ({\rm BLACK})$'s in agreement with 3.2). So, for each circle $W({\rm RED}) \cap W({\rm BLUE})$, the figures~5.4.(iii)$_1$ for the various $W_{\infty} ({\rm BLACK})$'s, where $W_{\infty} ({\rm BLACK}) \cup W({\rm BLUE})$ is already zipped, create a complete system of singularities involving $W({\rm RED})$ and $W({\rm BLUE})$. Now we want to zip completely the $W ({\rm BLACK})$'s too. So, let us consider a generic line $\partial W ({\rm BLACK}) \cap W({\rm RED}) \cap \overline N_{\infty}^2$, like for instance the $[s_2,a] \cap \overline N_{\infty}^2$, in the figure~1.3. Such a line, glued already at the source, will navigate through an infinity of $W({\rm BLUE})$'s. Then, for each given $W({\rm BLACK})$ and each $W({\rm BLUE})_n$, with $n \to \infty$, one can navigate along $W({\rm BLUE})_n \cap W({\rm RED})$, where $\partial W ({\rm BLACK})$ rests on {\it this} $W({\rm RED})$, from a BLUE/RED singularity created on a site (iii)$_1$ to the site (iii)$_2$ concerned by {\ibf our} $W({\rm BLACK})$ and then complete the remaining zipping $W({\rm BLUE})_n \cap W({\rm BLACK})$ through the corresponding $p_{\infty\infty}$-island, in a manner compatible with 3.2) in our lemma. This goes along the arrows wiggly and/or double green (they any way point in the same direction), from figure~5.4.(iii)$_2$.

\medskip

\noindent (5.18) \quad Continuing now our navigation along the $W({\rm BLUE})_n \cap W({\rm BLACK})$, we go beyond the $p_{\infty\infty}$-island, until we get to the right of the already zipped $ZZ_1 = \partial H_{i_0}^1 ({\rm MAX})$ in figure~1.3. This kind of thing goes on for all the $\{p_{\infty\infty}\mbox{-islands}\}$ and we can happily decide that
$$
[\beta , \alpha (\infty)] \cap fM_2 (f) \subset \{0(3)\mbox{-moves for the zipping flow}\} \, ,
$$
thereby taking care of 2) in our lemma. In order to finish the zipping of $W({\rm BLACK})$ and the proof of our lemma too, we still have to take care of the $\{\bar S\mbox{-regions}\} \subset W({\rm BLACK})$. This is our STEP~V to be developed next. We may assume that, at this stage in the game, for every $\bar S$-region like in the figure~1.1 or in 1.5.(C), the $W(n) , W(m)^*$ have all been zipped already up to exactly $x=x_{\infty}$, and then an immortal singularity $S(n,m)$ has been created. For each fixed $W_{\ell} ({\rm BLUE})$ the zipping has temporarily stopped at the boundary of the square $W_{\ell} ({\rm BLUE}) \cap \overline{\mathcal N}_{\infty}^3 (S)$, figure~5.2.

\medskip

\noindent (5.19) \quad Any point $\partial \{\mbox{square above}\} \cap (W(n) + W(m)^*)$ is, at this level, a singularity, created by the zipping of $W({\rm BLUE})_{\ell} - \overline{\mathcal N}_{\infty}^3$, and from these singularities, the zipping inside $W({\rm BLUE})_{\ell} \cap \overline{\mathcal N}_{\infty}^3$ will start, with a STOP at the $0(3)$-lines system from figure~5.2.

\medskip

In full argument with 3.2) in our lemma, these things also establish the equation below, which tells us where to locate the fins inside the $S$-region
$$
V = \mbox{overflowing}, \quad W^i = \mbox{subdued}. \eqno (*)
$$
This equation $(*)$ may need some explanations. According to the (ii) in figure~5.2, which concerns us now, the location of the fins $F_{\pm}$ on $W$ OR $W^*$, in the case of an encounter
$$
t_{\infty} = \{ W({\rm BLACK}) \cap W^* ({\rm BLACK}) \cap S_{\infty}^2 \} \, ,
$$
obeys to the following rule. The $W,W^*$ are already zipped together when we come to $t_{\infty}$, when the action of $W({\rm BLUE})_{\ell}$ occurs. Corresponding to $t_{\infty}$ there is a triple point, displayed in our figure~5.2,
$$
t_{\ell} = W \cap W^* \cap W({\rm BLUE})_{\ell} \, .
$$
In order to get to $t_{\ell}$, the $W_{\ell} ({\rm BLUE})$ action uses either $W$ OR $W^*$ as branch $V$ in (5.12) (ii). And it is {\ibf this} $V$ which gets the fins $F_{\pm}$ {\ibf at} $t_{\infty}$. This makes that $(*)$ completely agrees with what the figure~5.2 tells us.

\smallskip

This ends the proof of our lemma~5.3 and we continue with the proof of the zipping lemma, dealing now with the process $Z$ from (\ref{eq4.30}) and with the (partial) DITCH filling. For expository purposes, for a short while the triple points will be kept out of focus in our discussions.

\smallskip

Anyway, each complementary wall comes with a canonical transversal orientation, going conventionally from $z = +\varepsilon$ to $z=-\varepsilon$. When the ditches were explicitly defined, in (4.27.1) to (4.27.4), they had carried, corresponding to these transversal orientations, the canonical factors $[-\varepsilon \leq z \leq \varepsilon]$. There are two possibilities for the DITCH FILLING, namely

\medskip

\noindent (5.19.1) \quad Either the complete filling, for all the $[-\varepsilon \leq z \leq \varepsilon]$, OR the {\ibf partial} ditch-filling restricted to $\varepsilon - \frac1n \leq z \leq \varepsilon$, with an $n \to \infty$ to be made explicit later.

\medskip

The partial ditch filling is to be used in the situation of a zipping $W({\rm complementary}) \cap W$(non-comple\-mentary) $\cap$ $\overline N_{\infty}^2$. The full ditch filling means the full $q \times [-\varepsilon \leq z \leq \varepsilon] \times b^{N+1} (q)$, like in (4.27.1), or (4.27.2). When two complementary walls, let's say $W({\rm BLACK})$ and $W^j ({\rm RED} \cap H^0)$ are to be zipped together, then in $(N+4)$-dimensions their respective contributions 
$$
\{W ({\rm BLACK}) (\subset (x,y)) \} \times [-\varepsilon \leq z \leq \varepsilon] \times B^{N+1} \ - 
$$
$$
\left\{\mbox{the DITCH}  = \{\mbox{the shaded area} \, (-\varepsilon_j \leq y-y_j \leq \varepsilon_j , x)\} \times [-\varepsilon \leq z \leq \varepsilon] \times \frac12 \, B_{\rm black}^{N+1} \ \mbox{(see figure~5.5.ter)} \right\}
$$ 
and 
$$
\{ W^j ({\rm RED})(\subset (x,z)\} \times [-\varepsilon_j \leq y-y_j \leq \varepsilon_j] \times B^{N+1} \ - 
$$
$$
\left\{\mbox{the DITCH} = [\mbox{the shaded area} \, (-\varepsilon \leq z \leq \varepsilon , x)] \times [-\varepsilon_j \leq y-y_j \leq \varepsilon_j] \times \frac12 \, B_{\rm red}^{N+1} \right\}
$$ 

\smallskip

\noindent can be clamped together into a perfectly symmetrical fit, like in figure~5.5.ter. This figure is a more symmetrical version of the figure~4.4, with the $W_n ({\rm BLUE})$ ignored now. Notice that richer symmetry acquired now by the red splitting surface $S$. What we see in figure~5.5.ter corresponds, of course, to (4.27.3), and it is valid both in $N_{\infty}^2$ and $\overline N_{\infty}^2$.

\smallskip

Inside $N_{\infty}^2$ there is no prescribed a priori order for the elementary zipping steps, let us say that they all ``commute'', in particular, when we are at a triple point we are free to proceed arbitrarily. Inside $\overline N_{\infty}^2$ they do not commute any longer, and their order is fixed by the lemma~5.3. Moreover, when we are inside $\overline N_{\infty}^2$, then for double lines which involve (in $(N+4)$ dimension)
$$
W({\rm complementary}) \times [-\varepsilon \leq z \leq \varepsilon] \times B^{N+1} - \{{\rm DITCH}\}
$$
and
$$
W_n = W(\mbox{BLUE OR (RED} - H^0)) - \{{\rm Hole} \ H\} \, ,
$$
the ditch gets filled only up to $q \times \left[ \varepsilon - \frac1n \leq z \leq \frac1n \right] \times b^{n+1} (q)$ and, in the formula above, one should not mix up the ``$H$'' for Holes and ``$H^{\lambda}$'' for $\lambda$-handles.

$$
\includegraphics[width=140mm]{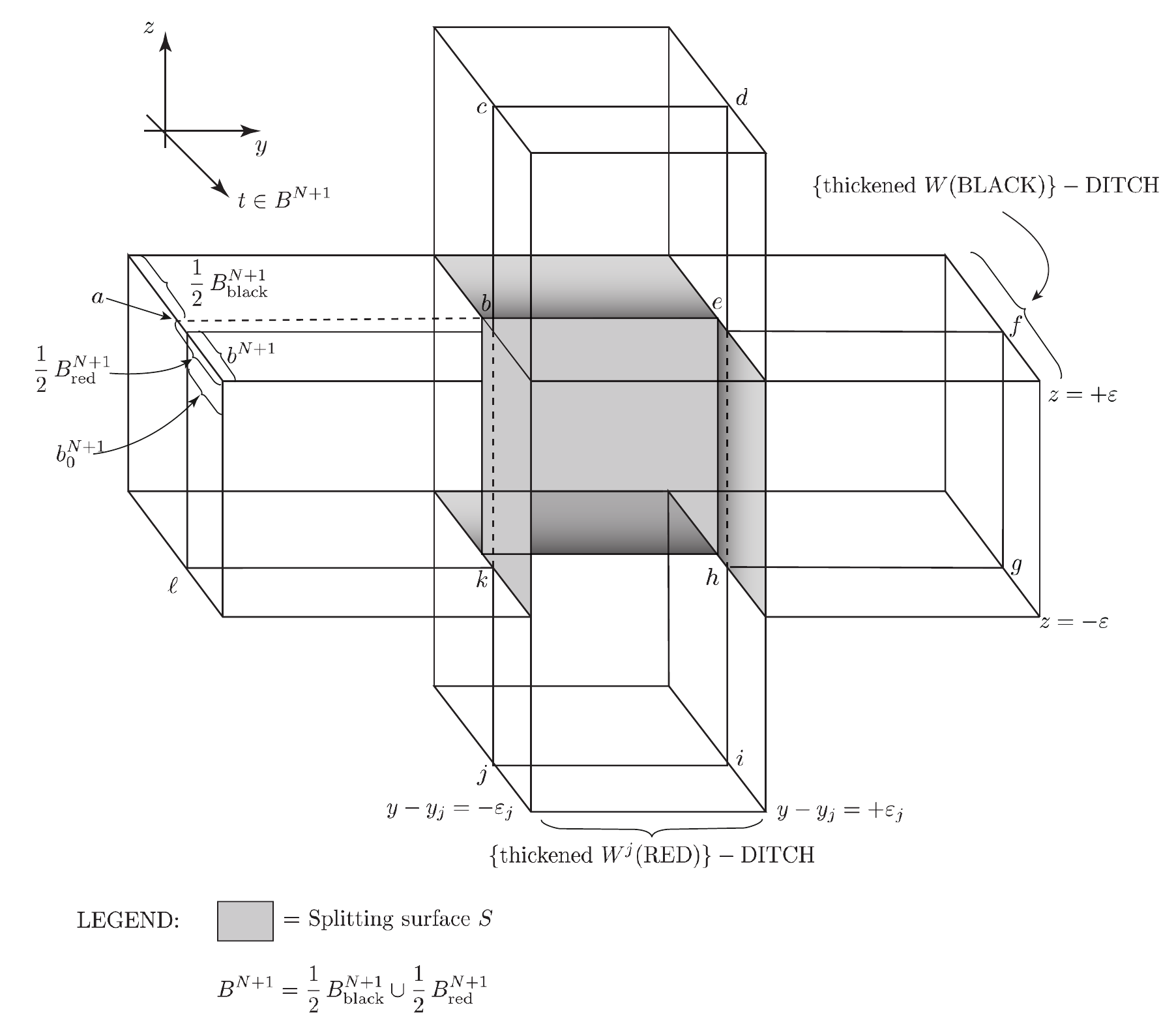}
$$
\label{fig5.5.ter}
\centerline {\bf Figure 5.5.ter} 

\smallskip

\begin{quote} 
We illustrate here the geometric realization, in dimension $N+4$ of the zipping together of two complementary walls. This figure should be compared of course with 4.4. The big cross $[a,b,c,\ldots ,j,k,\ell]$ is a piece of $\sum (\infty)$, reconstructed now at the level of $S_b^{(')} (\widetilde M (\Gamma)(-H))$. It is made out of pieces $h(1)$ and $h(3)$ (along the dotted lines).
\end{quote}

\bigskip

It is understood here that

\medskip

\noindent (5.20) \quad When of the same colour then $W_{n+1}$ is closer to the limit wall than $W_n$, and

\medskip

\noindent (5.21) \quad The two procedures for the geometrical realization of the sipping, for $N_{\infty}^2$ and for $\overline N_{\infty}^2$ are separated by $\{$the universal curve$\} = \partial N_{\infty}^2 \cap \partial \overline N_{\infty}^2$. This may be crossed at will by the zipping flow lines going along the $L_n = M_2 (f) \cap W_n$; figure~5.6 illustrates these things.

$$
\includegraphics[width=145mm]{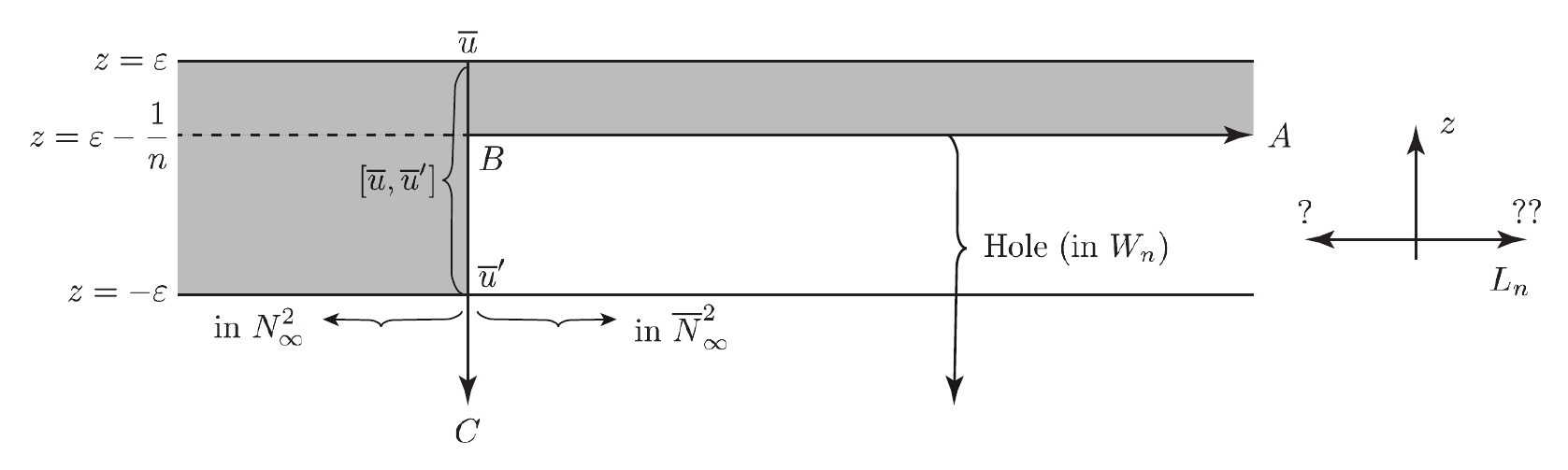}
$$
\label{fig5.6.}
\centerline {\bf Figure 5.6.} 

\smallskip

\begin{quote} 
We are here in the plane $x = x_n$ of a non-complementary wall $W_n$, superposed at the target $\widetilde M (\Gamma)$ of $f$, to a $W({\rm complementary}) \times [-\varepsilon \leq z \leq \varepsilon]$. The axis $L_n$ in the coordinate system should suggest the zipping flow which, in terms of $X^2$, goes along the intersection $W_n \cap W({\rm complementary}) = (x=x_n , z=0)$. Inside the $W({\rm complementary})$, the $[\bar u , \bar u']$ is the trace of the (universal curve) $= \, \partial N_{\infty}^2 \, \cap \, \partial \overline N_{\infty}^2$. When the width $[-\varepsilon , \varepsilon]$ is collapsed to zero, then the whole present arc $[\bar u , \bar u']$ could become the point $\bar u$ in figure~1.1.(B). Inside $W_n$, the line $[A,B,C]$ is in $\partial ({\rm Hole})$. But, a priori, the $[BC] \subset \partial ({\rm Hole})$ could also be located just very slightly to the right of the line $[\bar u , \bar u']$. The width of $W_n$ is $[x_n - \varepsilon_n \leq x \leq x_n + \varepsilon_n]$, and with all this, our ditch is $\underset{q}{\bigcup} \, q \times [-\varepsilon \leq z \leq \varepsilon] \times b^{N+1} (q) \subset W({\rm complementary}) \times [-\varepsilon , \varepsilon] \times B^{N+1}$, the union being over $q \in (x \in [x_n - \varepsilon_n , x_n + \varepsilon_n] , y , z=0) \subset W({\rm complementary})$.

\smallskip

The shaded area in the figure suggests the part of the DITCH which is filled, by the $\{ [(W_n - H) \subset (y,z)] \times [x_n - \varepsilon_n \leq x \leq x_n + \varepsilon_n] \times [$a reduced version of $b^{N+1}]\} \cup \{$The DITCH filling material $h(3)\}$, like in the process $Z$ from (\ref{eq4.30}). Here the reduced version of $b^{N+1}$ fits inside ${\rm int} \, b^{N+1} (q)$, glued to the outside by $h(3)$. When the Hole $H$ is forgotten, our whole figure is in the plane of $W_n$, or at least its projection is.
\end{quote}

\bigskip

We introduce now the triple points in the discussion too, and then some additional modulations with respect to the (5.20), (5.21) above become necessary. The proof of 2) in the ZIPPING LEMMA starts now and so we move back to the figure~4.6 which, with more details and embellishments, occurs below as figures~5.7 $+$ 5.8.

\smallskip

The triple points which are of interest for us now, are the
$$
fM^3 (f) \cap \overline N_{\infty}^2 = \{ t_{ij} \in [I(T_i)], p_{i\infty}\} \, .
$$
So, we consider here the situation of a $t_{in} \in \overline N_{\infty}^2$, after the initial zipping $V + W^i \Longrightarrow V \cup W^i$ has already been performed and when $W_n ({\rm BLUE}) - H$ encounters the line $V \cap W^i$. All the three cases (i), (ii), (iii) from (\ref{eq5.12}) are concerned here, but it is essentially only the paradigmatic case (i) which will be dealt with in some detail. In terms of $\Theta^3 (X^2) \overset{f}{\longrightarrow} \widetilde M (\Gamma)$ and of the figure~4.3, what we have displayed in the figure~5.7, inside the reference plane $x=x_0$, is the superposition of the two independent contributions of $U(x_0) \subset V({\rm BLACK}) \cup W^i ({\rm RED})$ and $W({\rm BLUE}) - H$, via the map $f$ above.

\smallskip

[Our present discussion does not involve the $p_{\infty\infty}$'s and so, the special refinements involving the $H(p_{\infty\infty})$, $D^2 (H (p_{\infty\infty}))$'s can be ignored here.] Like in the formula (4.27.2), the DITCH is here something like
$$
U(x_0) \times b^{N+1} \ \mbox{where} \ b^{N+1} (V+W^i \Longrightarrow V \cup W^i) = \{{\rm the} \ b^{N+1} \ \mbox{in figure 4.4}\} \, ,
$$
while $b^{N+1} (W_n + V \cup W^i \Longrightarrow W_n \cup V \cup W^i) = \{$the $b_0^{N+1} \subsetneqq b^{N+1}$, in the same figure~4.4$\}$. The $U(x_0) \subset W({\rm BLACK}) \cup W({\rm RED})$ can be seen in figure 4.3.(B).

\smallskip

Actually, the figure 5.7 below is supposed to be coherent with 4.4.

$$
\includegraphics[width=155mm]{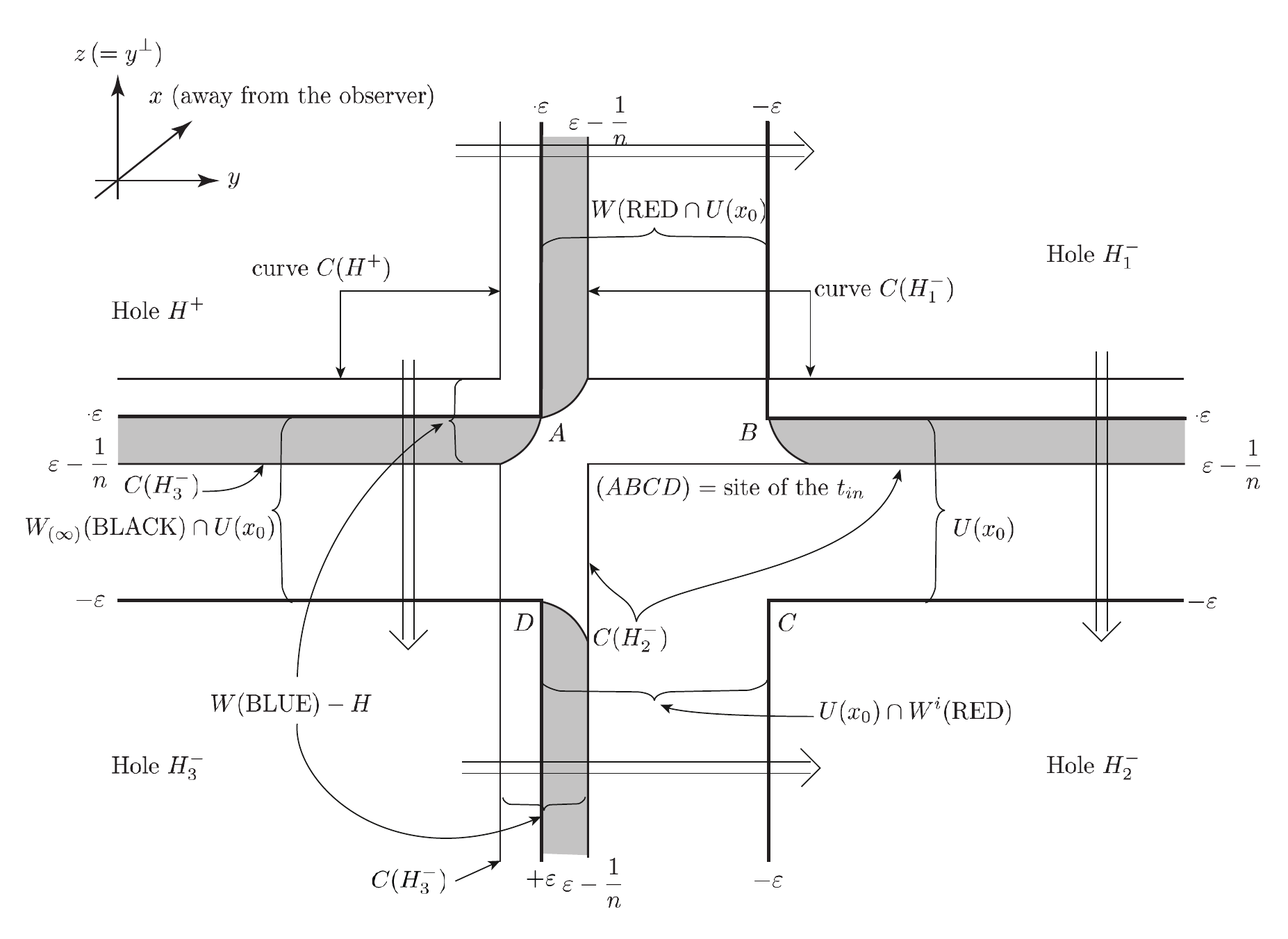}
$$
\label{fig5.7.}

\noindent LEGEND: $\includegraphics[width=5mm]{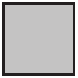}$ ditch filling material, contribution of $W_n ({\rm BLUE}) - H$ seen here as projected by the $\{$extension of the $f$ from (\ref{eq1.1}) to $\Theta^3 (X^2) \overset{f}{\longrightarrow} \widetilde M (\Gamma) \}$ to the plane $x=x_0$ of this figure. For the fate of the unshaded part of $U(x_0) \cap (W({\rm BLUE})-H)$, see what will be coming next; $\Longrightarrow \ =$ canonical transversal orientation for $U(x_0)$.

\smallskip

\centerline {\bf Figure 5.7.} 

\smallskip

\begin{quote} 
We are here at some level $x=x_0$, inside $\widetilde M (\Gamma)$, and we can see the relative situations of $U(x_0) \subset V \cup W^i$ and $W_n ({\rm BLUE}) - H$, in the neighbourhood of the triple point $t_{in}$, before any modulations concerning the special BLUE $1$-handles are in effect. The present $2^{\rm d}$ image is a projection by the map $f$, which superposes the two contributions involved. The $U \approx V({\rm BLACK}) \cup W^i ({\rm RED})$ and the $W_n ({\rm BLUE})$ have not yet been zipped together and, as a preparation for the modulations which are to come, we will start by collapsing the $W_n ({\rm BLUE}) - H$ into the thinner version presented in figure~5.8.(A). But the important modulation is the one described schematically by the transformation Fig. 5.8.(A) $\Rightarrow$ Fig. 5.8.(B). This is a transformation of the special BLUE handles $[AB]+[AD]$, which also drags along the four curves $C = \partial ({\rm Hole})$.
\end{quote} 

$$
\includegraphics[width=165mm]{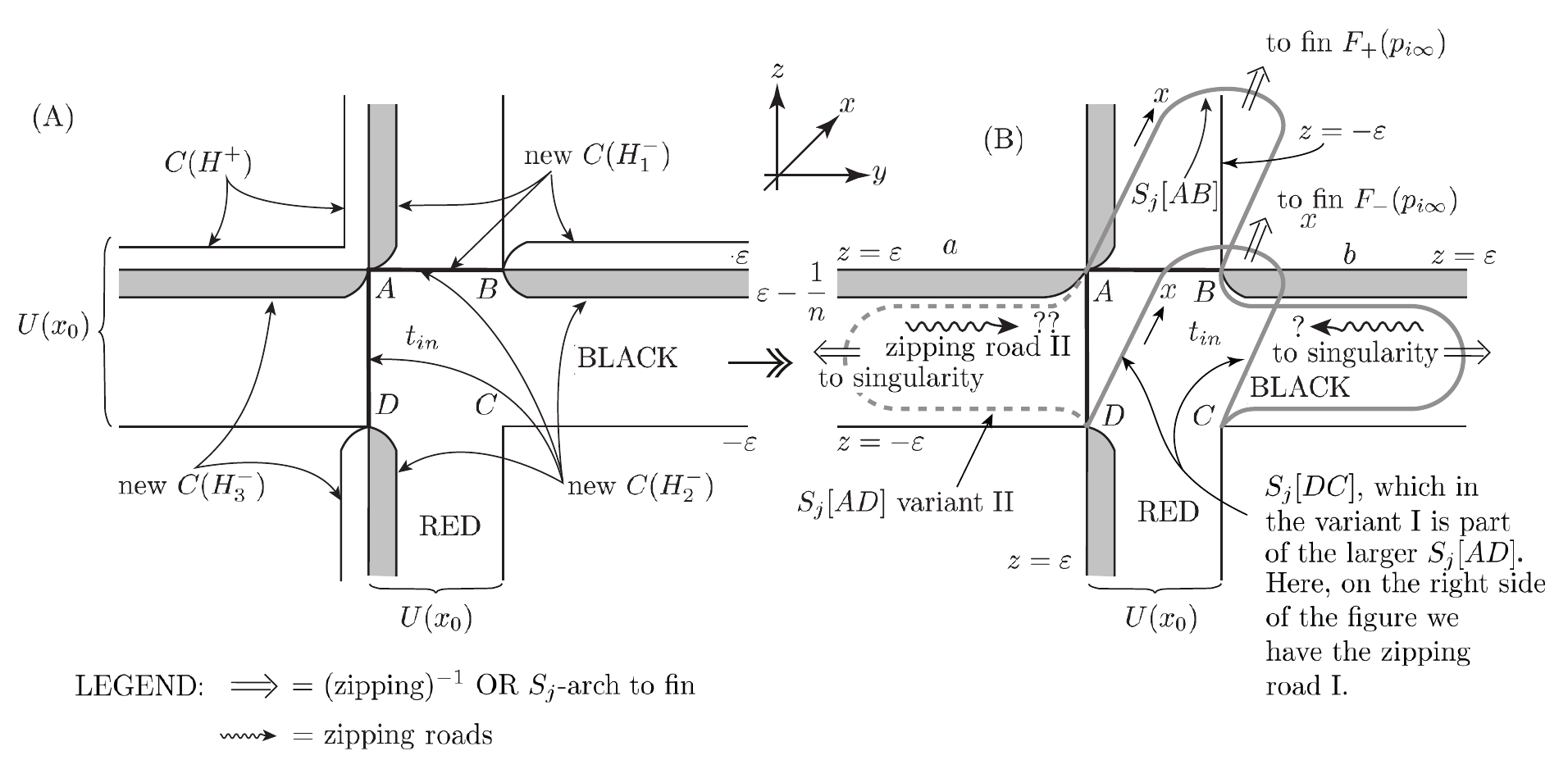}
$$
\label{fig5.8.}
\centerline {\bf Figure 5.8.} 

\smallskip

\begin{quote} 
In (B) we see, in wiggly lines, two possible zipping roads which lead to $t_{in}$, the variants~I (?) and II (??). The double arrows point to the fins $F_{\pm} (p_{i\infty})$, respectively to the singularities which command the zipping roads. When they point to the fins, then they go in the positive direction $x \to x_{\infty}$. Like in the figure~4.3, the BLUE cross $W_n ({\rm BLUE}) - H$, visible in the figure~5.7, is to be thickened in the $x$-direction, with width $x_n - \varepsilon_n \leq x \leq x_n + \varepsilon_n$. The shadow arcs $S_j [AB]$, $S_j [DC]$ (here we have a virtual $[XY] = [DC]$), which is necessary for the more complex $S_j [AD]$ (in the case of variant~I), go to the fins. For $S_j [AD]$, the (zipping)$^{-1}$, to be described later, is necessary.
\end{quote}

\bigskip

In both figures 5.7 and 5.8, the positive $x$-direction looks away from the observer. In figure~5.7, outside the $[ABCD]$, the shaded pieces $P_2(A)$ (which is double), $P(B)$, $P(D)$, are pieces of $W_n ({\rm BLUE})-H$, confined inside the thin bands $\left[ \varepsilon - \frac1n , \varepsilon \right]$.  [Notational remark. Comparing figures~5.7 $+$ 5.8 and 4.6, here is how the respective notations are supposed to fit together $P_2 (A) = R_2 \cup R_3$, $P(B) = R_1$, $P(D) = R_4$.]

\smallskip

In our figure, these three pieces appear superposed to $U(x_0)$. In real life, everything is smoothly thickened in dimension $N+4$ and our three pieces above live in the ditch (see (4.27.2) and the figure~4.4) and, with a $b_0^{N+1} \subset b^{N+1} \subset B^{N+1} -\frac12 \, B^{N+1}$ (figure~4.4) and an $\frac12 \, b_0^{N+1} \subset b_0^{N+1}$, concentric to $b_0^{N+1}$, we get, precisely inside the ``official''
$$
\mbox{``DITCH''} = \underset{\overbrace{\mbox{\scriptsize$x_n - \varepsilon_n \leq x_0 \leq x_n + \varepsilon_n$}}}{\bigcup} U(x_0) \times b_0^{n+1} \, ,
$$
the following smaller pieces, the only ones to be used in the ZIPPING LEMMA, when we are in the neighbourhood of $M_3 (f) \cap \overline N_{\infty}^2$. So, the actual {\ibf effective} DITCH is here

\newpage

\setcounter{equation}{21}
\begin{equation}
\label{eq5.23}
{\rm DITCH} \equiv \bigcup \underset{\overbrace{\mbox{\scriptsize$x_n - \varepsilon_n \leq x_0 \leq x_n + \varepsilon_n$}}}{\bigcup} \overline{\{ U (x_0)-\mbox{the square} \, [ABCD] \}} \times b_0^{N+1} \supset
\end{equation}
$$
\left[\underset{\overbrace{\mbox{\scriptsize$x_n - \varepsilon_n \leq x_0 \leq x_n + \varepsilon_n$}}}{\bigcup} (P_2 (A) + P(B) + P(D))\right] \times \frac12 \, b_0^{N+1} 
\left\{\begin{matrix}
\mbox{This is the part of $W({\rm BLUE})_n-H$} \\
\mbox{which, in the neighbourhood of the} \\
\mbox{triple point, lives inside the DITCH.}
\end{matrix}
\right\}
$$

\bigskip

The reason for excluding the square $[ABCD]$ in (\ref{eq5.23}) is, of course, that the metrizability requirement for $S'_b$ does not allow to have there a ditch housing the special BLUE $1$-handles $\beta [XY] = \{[AB],[AD]\}$.
$$
\includegraphics[width=155mm]{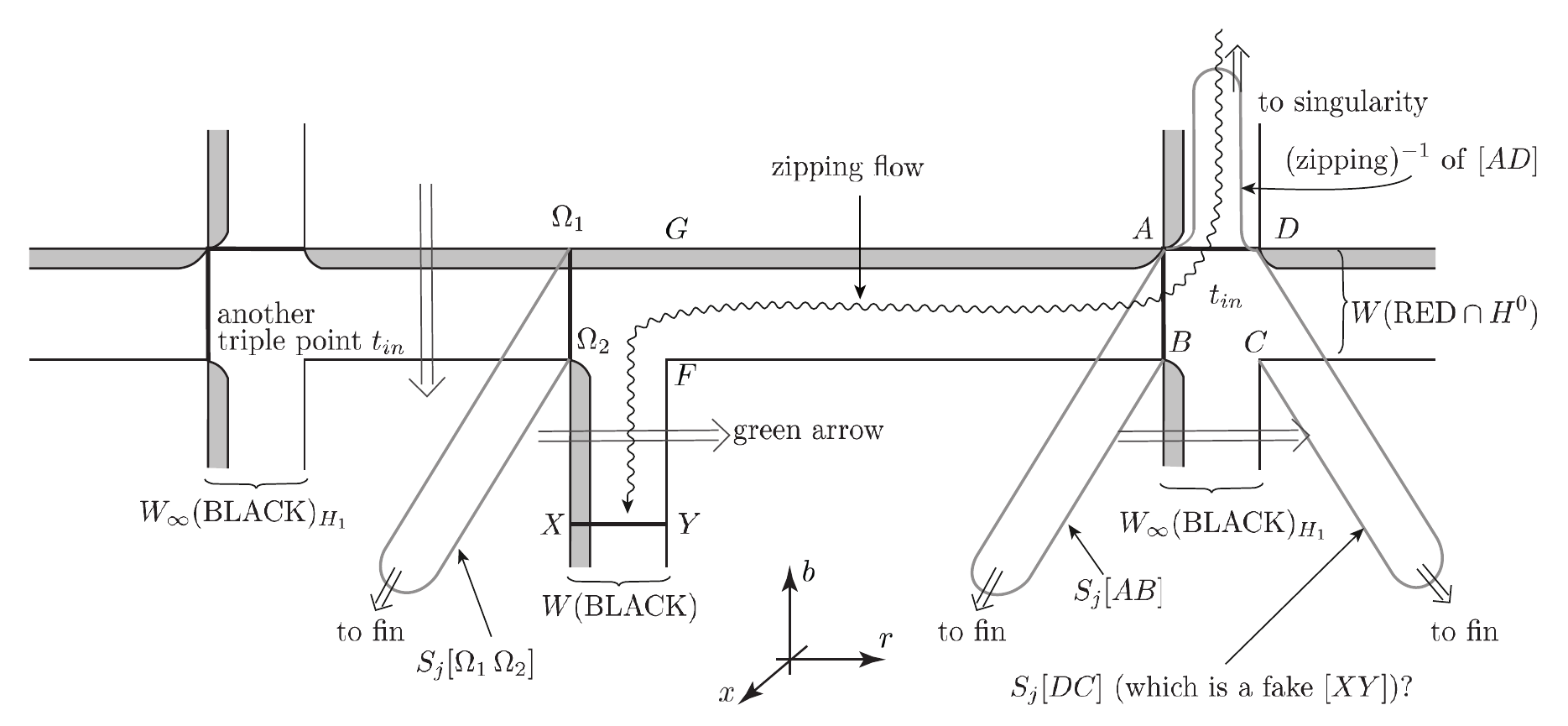}
$$
\label{fig5.9.}
\centerline {\bf Figure 5.9.} 

\smallskip

\begin{quote} 
In terms of our figure 5.4, the figures 5.7 $+$ 5.8 had referred to the case (i), while in the present figure the $[\Omega_1 \, \Omega_2 \, F \, G]$, $[ABCD]$ refer to the cases (iii)$_2$, (iii)$_1$, respectively. The direction $x$ goes to those fins of figure~5.4, which are attached to the black sheet at $x=x_{\infty}$. The green double arrow $(\Longrightarrow)$ stand for the transversal orientations, like in figure~5.7; they go away from the core, in the case $W({\rm RED})$. The transversal orientation determines how the (ditch filling part of the geometric realization of the) zipping reaches  the BLACK walls. The figure suggests only one of the variants which could occur.
\end{quote}

\bigskip

There are two kinds of special BLUE $1$-handles $[XY]$ namely real, like $[AB]$, $[AD]$, or {\ibf fake}, which will actually also need shadow arcs, like the $[BC]$ or $[DC]$ in the same figure 5.8.(A). There is also a second partition of the set of BLUE $1$-handles, transversal to the partition real/fake, and independent of it. Namely, we can have the RED, respectively the BLACK $[XY]$, which in a figure like 5.8.(A) go across $W({\rm RED} \cap H^0)$, respectively across $W_{(\infty)} ({\rm BLACK})$. We discuss first the easier case of the shadow arcs for the RED $[XY]$'s. At a first schematical level, in terms of the figure~5.8.(B), here is the definition of the shadow $S_j$, in the RED case
\begin{eqnarray}
\label{eq5.24}
S_j [AB] &=& [A, \mbox{rim of the fin} \ F_+ (p_{i\infty}),B] \, , \nonumber \\
S_j [DC] &=& [D, \mbox{rim of} \ F_- (p_{i\infty}),C] \, .
\end{eqnarray}

But, before going on, let us notice the following basic potential danger coming with the shadow arcs. At this point we will use the things said in the comment B after lemma~4.1. So, we have $[XY] \subset C(H(\mbox{completely normal}))$ and this $C$ is a $\partial h^2 (1)$. Now, normally, as a simple glance at figure~5.8.(B) may immediately suggest, the shadow arcs $S_j[XY]$ should cross the ditch filling material $h(3)$ which, in terms of the geometric intersection matrix should come now with contacts
$$
\beta [XY] \cdot \delta h^1 (3) = \partial h^2 (1) \cdot \delta h^1 (3) \ne \emptyset \, . \eqno (*_1)
$$
Remember, at this point, that $X^2$ supports two not everywhere well-defined flows, namely the collapsing flow with arrows $\partial h_i^2 (1) \cdot \delta h_j^1 (1)$, and the zipping flow. By itself, each of these two flows may not be terribly complicated, but then we certainly are bound to have transversal intersections

\bigskip

\noindent (5.24) \quad $\{$collapsing flow lines$\} \pitchfork \{$zipping flow lines$\} \ne \emptyset$ inside $X^2$, and what we are discussing now are bad closed loops potentially created by the union of the two flows, once (5.24) is taken into account.

\bigskip

What the (5.24) might bring for us, are the following kind of contacts:
$$
\partial h^2 (3) \cdot \delta h_i^1 ({\rm RED}) \ne 0 \, , \ \mbox{when there is also collapsing flow trajectory,} \eqno (*_2)
$$
$$
\delta h_i^1 ({\rm RED}) \overset{\partial h^2 (1) \cdot \delta h^1(1)}{-\!\!\!-\!\!\!-\!\!\!-\!\!\!-\!\!\!-\!\!\!-\!\!\!-\!\!\!-\!\!\!-\!\!\!\longrightarrow} [XY] ({\rm BLUE}) \, ,
$$
with  an $[XY]$ occurring in $(*_1)$.

\smallskip

Clearly, the combination of $(*_1)$ and $(*_2)$ would mean the doom, for our crucial property $S'_b \, \widetilde M (\Gamma) \in {\rm GSC}$. To prevent this disaster from happening, we will make use of the {\ibf isotopic push} from 3)  in the ZIPPING LEMMA~4.1. This will avoid the $(*_1)$, altogether, and hence it will prevent the bad loop to exist.

\smallskip

We will show now, in more detail, how this works for the $[AB]$ in (\ref{eq5.24}), in the case of the variant~I for the zipping, in the figure~5.8.(B). Remembering here that the figures~5.7 $+$ 5.8 go with 4.4 we redraw in figure~5.10 several details of interest for us now, from the figure~4.4. Here are the detailed

\bigskip

\noindent EXPLANATIONS FOR THE FIGURE 5.10. In (A) one sees the $\left[B , B -\frac1n, B'_0 , B'_0 - \frac1n; F,F-\frac1n , F'_0 , F'_0 - \frac1n \right]$, which is our partial DITCH filling from $z = \varepsilon$ (at $B$), to $z = \varepsilon - \frac1n$ (at $B-\frac1n$), of the effective DITCH (\ref{eq5.23}), by material in $h(1) \mid {\rm BLUE}$, which in figure~4.4 hits the shaded scar $[B , \overline B , B'_0 , \overline B'_0]$. We are here at time $m$ in the process (\ref{eq4.38}). Similarly, in (C), the $D'_2 \equiv [A,\bar A , A'_0 , \bar A'_0 ; b , \overline b , b' , \overline b'_0]$ and $D''_2 = [A,U,U',A'_0;a,u,u',a'_0]$ are two other effective ditches which are {\ibf empty} at this same time $m$, to the point that the dotted lines $[Ab]$, $[Aa]$ (this points to C, figure~4.4) are ghostly, i.e. they are not physically present. The physical $[U,u]$ points to $C'$ (see figure~4.4). In (A), notice that at the finer resolution level of our present figure, the $B$ from figure~5.8.(B) becomes now $B$ AND $X$. The $[X,X+\varepsilon , X_1 , X_2]$ in (A), is a piece of $S_j [XY]$ for $[XY] = [BA]$.

$$
\includegraphics[width=145mm]{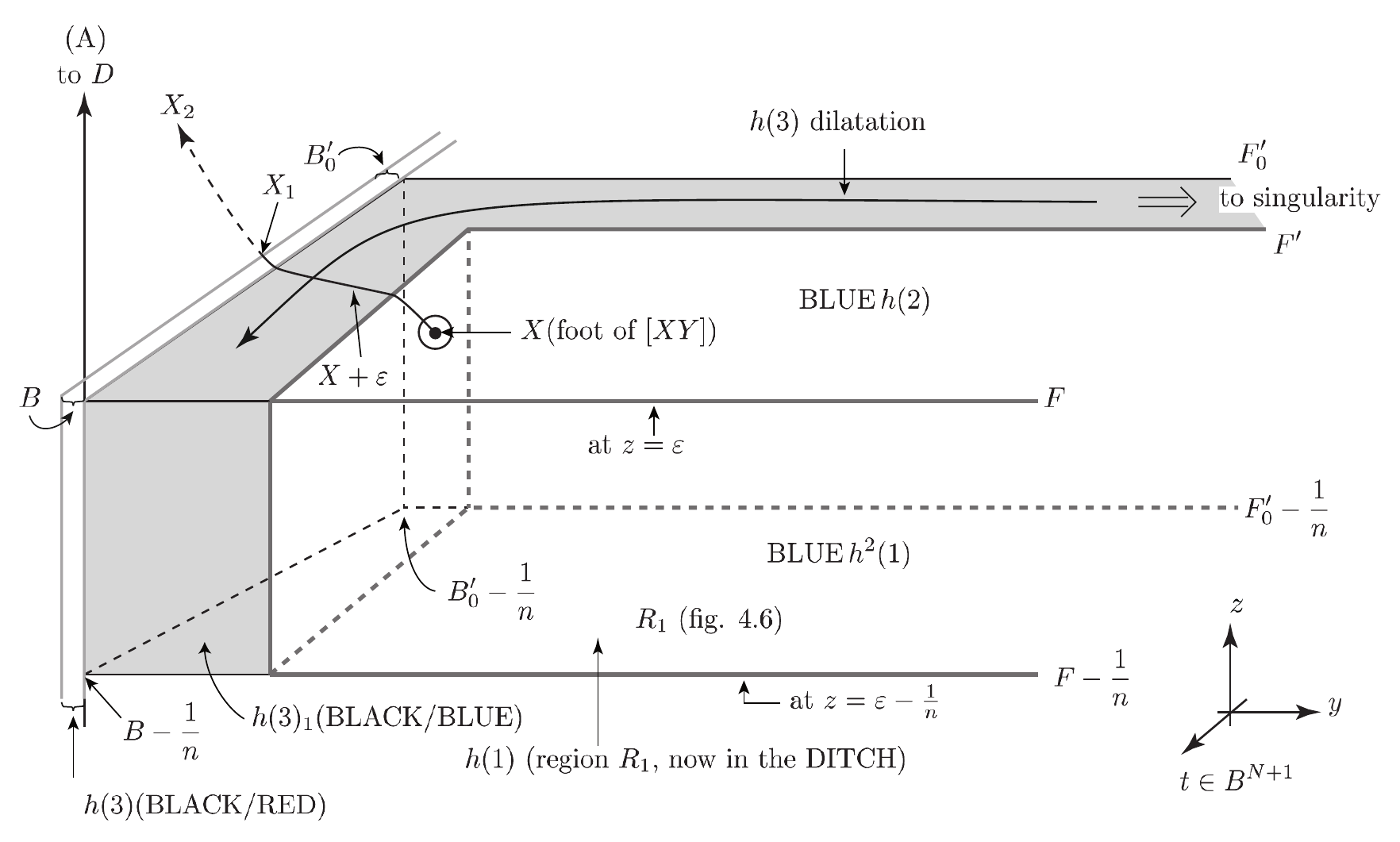}
$$

\bigskip

$$
\includegraphics[width=145mm]{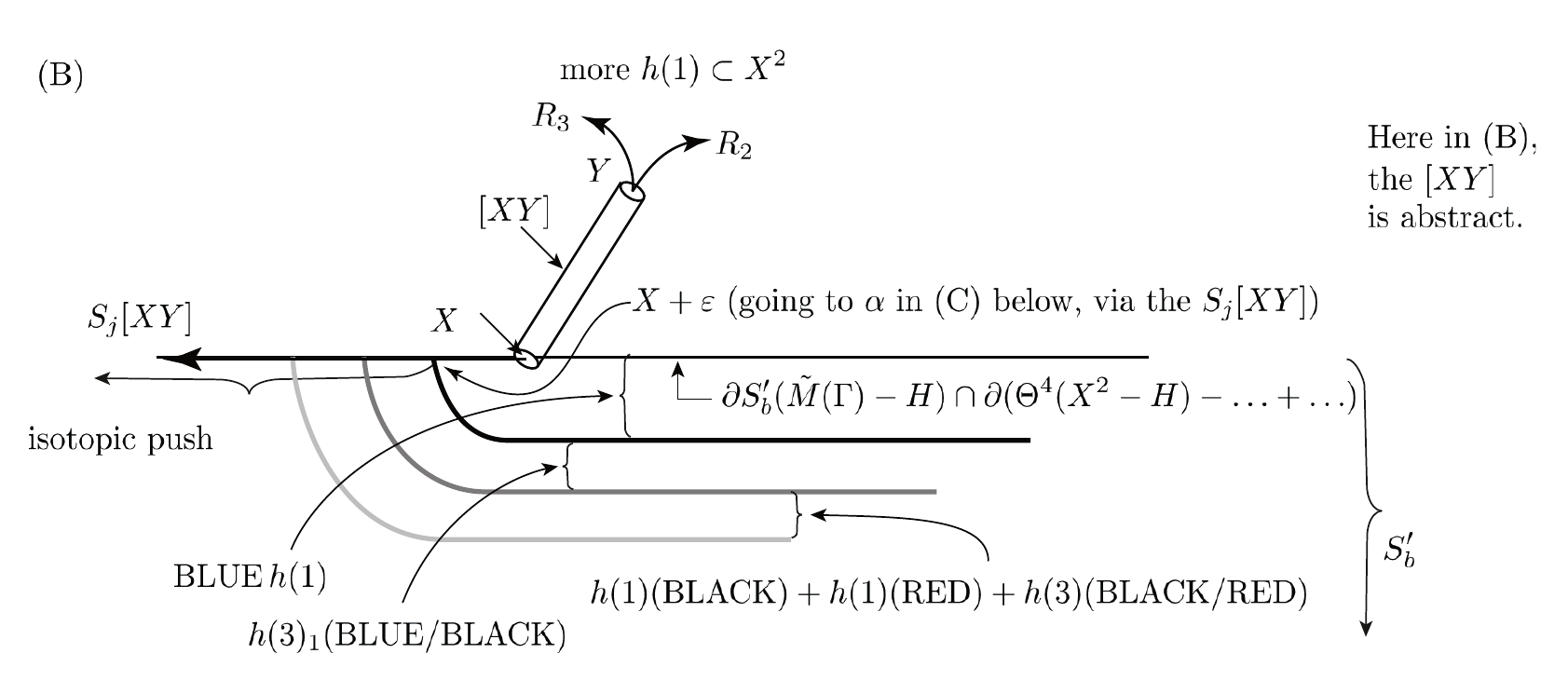}
$$

\bigskip

$$
\includegraphics[width=145mm]{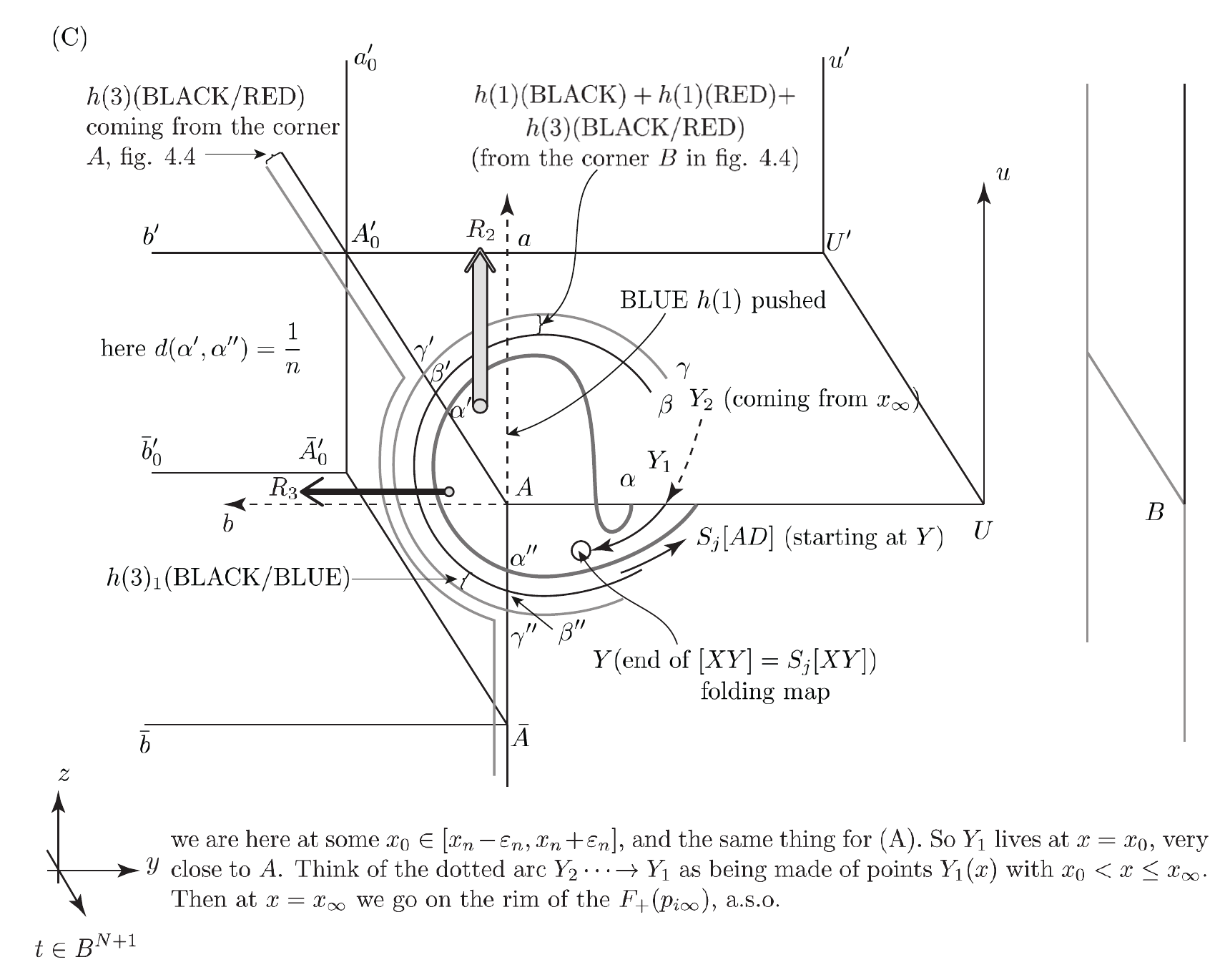}
$$

$$
\includegraphics[width=8cm]{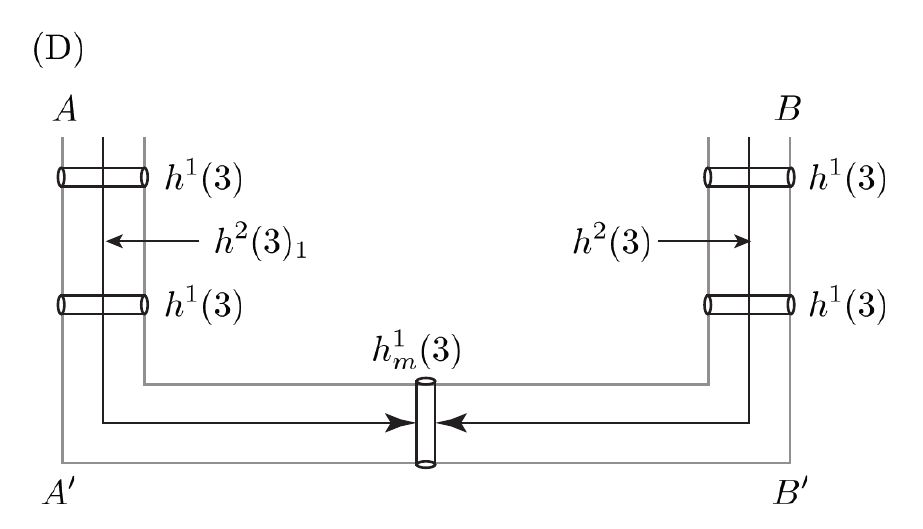}
$$
\label{fig5.10.}

\centerline {\bf Figures 5.10.} 

\smallskip

\begin{quote} 
The notations refer here mostly to figure~4.4 (but see also 5.8). The notations $R_1 , R_2 , R_3$ are like in the figure~4.6.

\smallskip

In (D) we see, schematically, the $S$-wound from the figure~4.4, and its cell-decomposition. This lives at $x=x_n$. The two arrows suggest the collapsing flow, at $x=x_n$, in the arch $S$. This continues then from $h_m^1 (3)$ on towards $x < x_n$. Remember that the collapsing flow means here $\overset{\partial h^i (3) \cdot \delta h^1 (3)}{-\!\!\!-\!\!\!-\!\!\!-\!\!\!-\!\!\!-\!\!\!-\!\!\!-\!\!\!-\!\!\!-\!\!\!\longrightarrow}$. When going into the action
$$
W({\rm BLACK}) + W^j ({\rm RED}) \Longrightarrow W({\rm BLACK}) \cup W^j ({\rm RED}) \, ,
$$
then the two $W$ (complementary) which are involved, come already equipped with disks $D_1 \Bigl( \frac12 \, B_{\rm RED}^{N+1}$, $\frac12 \, B_{\rm BLACK}^{N+1} \Bigl)$, see figure~5.5.ter, and also with the disks $D_2 (b_0^{N+1})$. These are necessary for the next zipping, with $W_n ({\rm RED}) - H$, see figure~4.4. It is understood that the ditch-filling $h(3)({\rm BLACK} / {\rm RED})$, which heals the wound $S$, in figure~5.5.ter, is added in such a way that the $D_2$'s rest on it, as our present figure (C) suggests they should.
\end{quote}

\smallskip

Concerning the dotted lines $[A,a], [A,b]$, in the drawing (C): Before our isotopic pushing action reaches the corner $A$, the volumes $D'_2 , D''_2$ above are just empty ditches, to be filled with material now. Concerning the drawing (A), the dotted part $[X_1 , X_2]$ lives outside the space of our figure (A), which is at $x=x_n$, and it is going along $x_n \to x_{\infty}$. The (B) is also at time $m$, before any pushing action of $h^2 (1)({\rm BLUE})$ has started. This (B) is a reminder, on the one hand of the fact that our $S_j [XY]$ really lives in high dimensions, inside $\partial S'_b (\widetilde M (\Gamma) - H)$ and that, on the other hand, what one sees in (C) when the push along $S_j [BA] = S_j [XY]$ has been completed on the surface, has depth too. This push sends the $X+\varepsilon$ from (B) into the area $(\alpha , \alpha' , \alpha'')$ from (C). The BLUE regions $R_2 , R_3$ in (C) grow out of the
$$
\{ \alpha \, \alpha' \, \alpha'' \} \subset \{\mbox{pushed BLUE} \ h(1) \} \, ,
$$
after the smeering ($=$ the folding map) has been performed.

\smallskip

From the viewpoint of their increasing time of first occurrence in the process $Z$ (\ref{eq4.38}), we have

\medskip

\noindent (5.25) \quad $h(1) < h(3) ({\rm BLACK}/{\rm RED}) < h(3)_1 ({\rm BLACK}/{\rm BLUE}) < \{h(3)_2 ({\rm BLACK}/{\rm BLUE}) \cup h(3)_2 ({\rm RED} /$ ${\rm BLUE})$, not yet represented graphically, but contributing eventually as binding material of $R_2 + R_3$ with the rest, in the partial ditch filling corresponding to (C)$\}$.

\medskip

In (C), once the push is finished (see 3) in lemma~4.1), then the arches $[\alpha , \alpha']$, $[\beta , \beta']$, $[\gamma , \gamma']$ are contained in $[A \, U \, U' A'_0]$ and the $[\alpha' , \alpha'']$, $[\beta' , \beta'']$, $[\gamma' , \gamma'']$ are contained in $[A \, \bar A \, \bar A'_0 , A'_0]$.

\smallskip

Once the pushing has extended the $h(1)$ from $X+\varepsilon$ in (B) to the $(\alpha , \alpha' , \alpha'')$ in (C), proceeding along $S_j [XY]$, what the smearing does is to lean $[XY]$ on $S_j [XY]$ and blend it into it, which realizes the diffeomorphisms (4.31.3) $+$ (4.31.4), final result of the present mini-step pushing $+$ smearing. The fact that now, with the smearing $[XY] = S_j [XY]$ (all happening inside $h(1)$, so as to preserve GSC), the $Y$ has reached the position from (C) is having as consequence the occurrence of the shaded $R_1 , R_3$ in (C), which are sent at this point into the interiors of their respective ditches. The impact of pushing $+$ smearing on $h(1)$ (which internally stays GSC and which does not acquire arrows to the outside), has already been discussed. I claim that the same pushing $+$ smearing does not induce any changes on the matrix $\partial h^2 (3) \cdot \delta h^1 (3)$. This follows from the following geometrical fact. Inside any of the bands $h(3) ({\rm BLACK} / {\rm RED})$ occurring once in (A) and twice in (C), the respective $1$-handles $h^1(3)$ are transversal to the band in question, like it is suggested in the figure~5.10.(D). What this also means, is that when like in (C) the $h(3)$ is pushed along by BLUE $h(1)$, then the $h(1)$ in question never has a chance of coming into contact with the $\delta h^1 (3)$'s. Notice, also, that our $h(3)({\rm BLACK}/{\rm RED})$ corresponds to the $S$-wound from the figure~4.4.

\smallskip

When we complete now the partial DITCH-filling, at the level of (C) then the $h(3)_2$ from (5.25) finally comes in too. In terms of the collapsing flow, we will have now new arrows in the geometric incidence matrix, of the form
$$
\{\mbox{higher terms in (5.25)}\} \Longrightarrow \{\mbox{lower terms in (5.25)}\} \, ,
$$
which is OK for GSC.

\smallskip

This ends our EXPLANATIONS and we move now to the BLACK $[XY]$'s, for which we have to construct the shadow arcs $S_j [XY]$. For this, some preliminaries will be necessary.

\smallskip

Lemma 5.3 has provided us, for any $(x_0 , y_0) \in M^2 (f) \subset (X^2 - H) \times (X^2 - H)$ with a zipping path $\lambda (x_0 , y_0) \subset \widehat M^2 (f) \equiv M^2 (f) \cup {\rm Diag} ({\rm Sing} (f))$. So far, it has not mattered much whether $X^2 \times X^2$ was meaning the usual product, with $(x,y) \ne (y,x)$ when $x \ne y$, OR the symmetric product, with $(x,y) = (y,x)$. For what comes next, it will be convenient to work rather with the {\ibf symmetric product}, which allows us to perceive a map
$$
M_2 (f) \cup {\rm Sing} (f) \overset{2}{-\!\!\!-\!\!\!-\!\!\!-\!\!\!\longrightarrow} M^2 (f) \cup {\rm Diag} ({\rm Sing} (f)) \, ,
$$
such that, for $(x,y) \in M^2 (f)$, we have $2^{-1} (x,y) = \{ x,y \}$, and which, outside of the singularities, covers exactly twice its image. We will also think now of the zipping strategy $\lambda (x_0 , y_0)$ as being a {\ibf subset} of $\widehat M^2 (f) \equiv M^2 (f) \cup {\rm Diag} ({\rm Sing} (f))$, forgetting the parametrization by $(x_t , y_t) \in \lambda (x_0 , y_0)$, for a short while.

\bigskip

\noindent {\bf Claim (5.25.1)} Given the $\lambda (x,y) \subset \widehat M^2 (f)$, there is a commutative diagram of {\ibf continuous} maps
$$
\xymatrix{
I \ar[rr]^-{\widetilde\lambda^{-1}} \ar[drr]^{\lambda^{-1}} &&\lambda \subset \widehat M^2 (f)\subset (X^2 - H) \times (X^2 - H)/(Z/2Z) \\ 
&&\widehat M^2 (f) \equiv M_2 (f) \cup {\rm Sing} (f) \subset X^2 - H \, , \ar[u]_-{2}
}
$$
with the following feature.

\smallskip

\noindent a) We start from $\lambda^{-1}$, which is a homeomorphism on its image, having the feature that the composite {\ibf map} $2 \circ \lambda^{-1}$ has as its image our $\lambda \subset \widehat M^2 (f)$. This {\ibf defines} then our map $\widetilde\lambda^{-1}$. The $\lambda^{-1}$  consists of finitely many successive continuous pieces on each of which, individually, either the natural time flows on I and on $\lambda$ go in the same direction OR in opposite direction (a time reversal).

\smallskip

\noindent b) $\lambda^{-1} (0) = x_0$, $\lambda^{-1} (1) = y_0$. We will sometimes write $\lambda^{-1} = \lambda^{-1} (x_0 , y_0)$. \hfill $\Box$

\bigskip

The claim (5.25.1) is easily visualizable in the figure~5.3, which can unambiguously be read as the image of a {\ibf continuous} path going from $x$ to $y$, or from $y$ to $x$ (a free choice). Also, in a very unprecise and metaphorical manner, if one thinks of the zipping path $\lambda (x_0 , y_0)$ as a continuous path in $\widehat M^2 (f)$, parametrized by $0 \geq t \geq -\infty$, with $\lambda (x_0 , y_0)$ ($t=-\infty$) $=$ singularity, and $\lambda (x_0 , y_0)$ ($t=0$) $=$ $(x_0 ,y_0)$, then $\widetilde\lambda^{-1}$ can be thought of as being the same path as $\lambda$, but with a time reversal, going now back from the singularities to $(x_0,y_0)$. Point 6) in lemma~5.2 may be reformulated now as follows

\medskip

\noindent (5.26) \quad Let $\underset{n=\infty}{\lim} (x_0 (n) , y_0(n)) = \infty$, in $X^2 \times X^2$. Then, with appropriate uniform estimates, and measured in the metric $\min \{ \ldots ; \ldots \}$ from the 6) in lemma~5.3, we also have $\underset{n=\infty}{\lim} \lambda (x_0 (n) , y_0(n)) = \infty$ ($\Leftrightarrow$ $\underset{n=\infty}{\lim} (\lambda^{-1} (x_0 (n) , y_0(n)) = \infty)$.

\medskip

For given $(x_0 , y_0) \in M^2 (f)$ the zipping path $\lambda (x,y)$ involves only a finite initial portion of the process $Z$ in (\ref{eq4.30}), call that $Z \mid [0,m]$, with $m=m(x_0,y_0)$. 

\smallskip

We go now high-dimensional and we introduce the notations

\medskip

\noindent (5.27) \quad $\Theta^{N+4} \equiv \Theta^{N+4} (X^2 - H) - \{{\rm DITCH}\} = S_0$ like in (\ref{eq4.33}), and $\Theta_1^{N+4} \equiv \Theta^{N+4} \cup \{h(2) + h(3)\} \mid [0,m]$, coming with
$$
\Theta_1^{N+4} \underset{{\rm the \, embedding} \, Z \, {\rm in \, (4.33)}}{-\!\!\!-\!\!\!-\!\!\!-\!\!\!-\!\!\!-\!\!\!-\!\!\!-\!\!\!-\!\!\!-\!\!\!-\!\!\!-\!\!\!-\!\!\!-\!\!\!-\!\!\!-\!\!\!-\!\!\!-\!\!\!-\!\!\!\longrightarrow} S'_b (\widetilde M (\Gamma)-H) \, .
$$

The reason for denoting now the $S_0$ from (\ref{eq4.33}) by $\Theta^{N+4}$ is to stress that it is a regular neighbourhood of $X^2-H$; we will thus consider for it the canonical retraction $\Theta^{N+4} \overset{p}{\longrightarrow} X^2-H$.

\smallskip

In a similar vein as in (5.27) we also introduce now the $(N+3)$-manifold $\Theta_2^{N+3} \equiv \partial \, \Theta^{N+4} - \bigl\{$that piece of $\partial \, \Theta^{N+4}$ which stays hidden behind the ditch-filling material $h(3)$, glued to it, when we go to $\Theta^{N+4} \overset{j}{\longrightarrow} S'_u (\widetilde M (\Gamma)-H)$ (\ref{eq4.33}). The hidden piece above can be visualized, in the context of the toy model, in figure~3.2.(B) where it is the curved surface $[u,v,w,u',v',w'] \mid \left[ \varepsilon - \frac1n , \varepsilon \right] \bigl\}$.

\bigskip

\noindent {\bf Lemma 5.4.} {\it With $N$ sufficiently high, for each $u \in X^2-H$, the fiber $\Phi (u) \equiv p^{-1} (u) \cap \Theta_2^{N+3}$ is connected.}

\bigskip

\noindent {\bf Proof.} Ignoring the discrete set ${\rm Sing} (f) \subset X^2-H$, the fiber of
\setcounter{equation}{27}
\begin{equation}
\label{eq5.29}
\Theta^3 (X^2-H) \to X^2-H
\end{equation}
is a compact tree with a single branching point, and at this stage exactly the ends of the tree live on the boundary. The $p$ above factorizes through (\ref{eq5.29}) and with $N$ high enough these ends get connected at the level of $p^{-1} (u) \cap \partial \, \Theta^{N+4}$; the contribution of that $\{ \ldots \}$ occurring in the definition of $\Theta_2^{N+3}$ when it gets deleted (and which should be visualizable in figure~3.2) does not make here any difference. [The codimension of $\Phi (u)$ in the very high-dimensional $\Theta_2^{N+3}$ is two, while $\{ \ldots \}$ is formally $1$-dimensional, in the sense that it is the very thin regular neighbourhood of a $1$-dimensional set.] \hfill $\Box$

\bigskip

We consider next the following commutative diagram, where the arrow $j$ is like in (\ref{eq4.33}), when the $j_2$ factorizes through $S'_b (\widetilde M (\Gamma) - H)$, with a $j_1$ which is not everywhere well-defined (it is certainly not quite well-defined in the BLUE 1-handles $[XY] \mid [0,m]$)
\begin{equation}
\label{eq5.30}
\xymatrix{
X^2-H \ar[d]^-{f} &&\Theta^{N+4} \ar[ll]_-{p} \ar[rr]^-{j_1} \ar[dr]_-{j} &&\Theta^{N+4}_1 \ar[dl]^{\qquad j_2 \equiv {\mathcal J} \mid \Theta_1^{N+4} \ ({\rm see} \, (4.33))}  \\ 
fX^2 - H &&&S'_u (\widetilde M (\Gamma) - H). \ar[lll]^-{q}  
}
\end{equation}
This diagram should be completed with the following
\begin{equation}
\label{eq5.31}
\xymatrix{
\Theta^{N+4} \ar[d]^{Z} \ar[rr]^-{j_1} &&\Theta_1^{N+4} \ar[d]^-{j_2} \ar[dll]_{[Z]} \\
S'_b(\widetilde M (\Gamma)-H) \ar[rr]_-{\mathcal J} &&S'_u (\widetilde M (\Gamma)-H) \, ,
}
\end{equation}
where $[Z]$ is the obvious extension of the $Z$ (\ref{eq4.30}). The combination of (\ref{eq5.30}) and (\ref{eq5.31}) makes that we will be able to talk about embeddings into $S'_{\varepsilon} (\widetilde M (\Gamma)-H)$, in the context of the lemma~5.5 below.

\smallskip

With all this, we look now at the map $[0,1] \overset{\lambda^{-1} (x_0 , y_0)}{-\!\!\!-\!\!\!-\!\!\!-\!\!\!-\!\!\!-\!\!\!-\!\!\!-\!\!\!\longrightarrow} X^2-H$ from the CLAIM (5.25.1) and, from now on, whenever there is no danger of confusion, the image $\lambda^{-1} (x_0 , y_0)$ (I) of this map, will be denoted just by $\lambda^{-1} (x_0 , y_0)$. This can be lifted then from $X^2-H$ to an embedding $\lambda^{-1} (x_0 , y_0) \subset \Theta^{N+4}$. Because of lemma~5.4, we can take this embeding as being of the form
$$
\lambda^{-1} (x_0 , y_0) \subset \Theta^{N+3}_2 \, .
$$

\bigskip

\noindent {\bf Lemma 5.5.} {\it For each BLACK special BLUE arc $[XY]$ (like the real $[AD]$ or the fake $[BC]$, in the figure~}5.8.(A){\it), there is a {\ibf continuous} path
\begin{equation}
\label{eq5.32}
({\rm zipping})^{-1} [XY] \subset \partial \, \Theta_1^{N+4} \, , \ \mbox{connecting $X$ and $Y$},
\end{equation}
with the following features.}

\medskip

\noindent (5.31.1) \quad {\it Notice, to begin with, that for the $[XY]$ under scrutiny, one of the endpoints, call it the $X$, corresponds to some double point, typically for $X=B$, $(X({\rm BLUE}) , X({\rm BLACK})) \in M^2 (f)$ (and then the same thing is not necessarily so for $Y$ too). By lemma~{\rm 5.3} to this double point corresponds a zipping path $\lambda$, coming with $2 \Vert \lambda \Vert = \Vert \lambda^{-1} \Vert$, and with this we will have the estimate below, with some universal constants $C_1 , C_2$ (themselves determined, essentially, by our former $P_0$} (\ref{eq5.12}) {\it and by lemma~{\rm 5.3})}
$$
\Vert ({\rm zipping})^{-1} [XY] \Vert \leq C_1 \cdot \Vert \lambda^{-1} (X({\rm BLUE}) , X({\rm BLACK})) \Vert + C_2 \, .
$$

\medskip

\noindent (5.31.2) \quad {\it We can make use of the diagram} (\ref{eq5.30}) {\it and, without loosing the feature} (5.31.1), {\it transport the $({\rm zipping})^{-1} [XY]$ from} (\ref{eq5.32}) {\it to $({\rm zipping})^{-1} [XY] \subset \partial S'_{\varepsilon} (\widetilde M (\Gamma)-H)$. 

\smallskip

In terms of} (\ref{eq4.15}) {\it we may ask that, whenever it makes sense, the $({\rm zipping})^{-1} [XY]$ should be localized inside $b^{N+1} \subset B^{N+1}$. This condition will be particularly important when $\lambda^{-1}$ enters a neighbourhood of $p_{\infty\infty} \in W_{(\infty)} ({\rm BLACK})$, and it will allow us then to invoke lemma~{\rm 4.3}, when that will be needed, inside the proof of the present lemma.}

\medskip

\noindent (5.31.3) \quad {\it In the context of} (5.31.2) {\it we also have the following. If $\underset{n=\infty}{\rm lim} X(n) = \infty = \underset{n=\infty}{\rm lim} Y(n)$ then we also have $\underset{n=\infty}{\rm lim} (({\rm zipping})^{-1} [X(n) , Y(n)]) = \infty$.

\medskip

Here, starting with the endpoints, various uniform estimates can be also imposed. The convergence which is meant above is the following one, more mundane than in {\rm (5.26)}, namely the following: For every compact $K \subset S'_{\varepsilon} (\widetilde M (\Gamma)-H)$, there is an $n_0 \in Z_+$ s.t., when $n > n_0$ then we also have $({\rm zipping})^{-1} [X(n) , Y(n)] \subset \partial S'_0 - K$.}

\bigskip

Before going into the proof of this lemma we give a comment.

\bigskip

\noindent {\bf Comment (5.31.3).} Ideally, the $({\rm zipping})^{-1} [XY]$ should cover the $\lambda^{-1} (x_0 , y_0)$, where (in our specific case) $(x_0 , y_0) = X({\rm BLUE}) , X({\rm BLACK}))$. This means the following.

\smallskip

Both $({\rm zipping})^{-1} [XY]$ and $\lambda^{-1} (x_0 , y_0)$ are arcs parametrized by $[0 \leq t \leq 1]$. With this, in terms of the diagram
$$
\xymatrix{
\lambda^{-1} (x_0 , y_0) \ar[r] &X^2 - H \ar[d] \\
&fX^2 - H &S'_{\varepsilon}(\widetilde M (\Gamma)-H) \ar[l] &({\rm zipping})^{-1} [XY] \ar[l]
}
$$
points parametrized by the same $t \in I$, should live in the same spot in $fX^2-H$. But this ideal plan will have to be subjected to quite some modifications, as it will be seen in the proof below.

\bigskip

\noindent {\bf Proof.} In order to make the exposition easier, we will forget about the map ``2'' from (5.25.1), and revert to a more impressionistic and/or heuristic viewpoint (already mentioned before), where the $\lambda^{-1} (x_0,y_0)$ is conceived now as a path in $(X^2-H) \times (X^2-H)$ defined by $\lambda^{-1} (x_0,y_0) (t) = \lambda (x_0,y_0)(-t)$, the $t$ going here from $t=0$ to $t = -\infty$. Our $({\rm zipping})^{-1} [XY]$ will cover, essentially, {\ibf this} $\lambda^{-1} (x_0 , y_0)$, actually not quite so, as we shall see. But anyway, with this change of viewpoint, the $({\rm zipping})^{-1} [XY]$ should live now in $\partial \, \Theta_1^{N+4} \, \times \, \partial \, \Theta_1^{N+4}$ and/or in $\partial S'_{\varepsilon} (\widetilde M (\Gamma)-H) \, \times \, \partial S'_{\varepsilon} (\widetilde M (\Gamma)-H)$. Keep in mind that, in the discussion which follows $[XY] = [BC]$ (figure~5.8), the zipping is like in variant I of the figure in question and $\lambda^{-1} = \lambda^{-1} (B({\rm BLUE}) , B({\rm BLACK}))$. [As a final remark, once our arguments will have bee completely unrolled in the framework of the impressionistic/heuristic viewpoint of time reversal along the zipping path, the reader should have no trouble of coaching them in the more accurate frame work of the CLAIM~(5.25.1).]

\smallskip

For our $\lambda^{-1} (t)$ with $\lambda^{-1} (0) = B$, there is a first time $t_0 < 0$, such that, while we have $\lambda^{-1} \mid [0,t_0] \subset \overline{\mathcal N}_{\infty}^3$, at $t_0$, the $\lambda^{-1}$ enters ${\mathcal N}_{\infty}^3$. For the same $0 \geq t \geq t_0$, the $(X_t , Y_t) \in ({\rm zipping})^{-1} [BC]$ is a continuous path in $\partial S'_{\varepsilon} (\widetilde M (\Gamma) - H) \times \partial S'_{\varepsilon} (\widetilde M (\Gamma)-H)$ (the level we consider now), such that $X_0 = B$, $Y_0 = C$, covering $\lambda^{-1} \mid [0,t_0]$ modulo the prescription which will be given in (5.33) below, and such that, outside of the encountered triple points we should have
\begin{equation}
\label{eq5.33}
z (X_t ({\rm BLUE})) \in \left[ \varepsilon , \varepsilon - \frac1n \right] \, , \quad z (Y_t({\rm BLACK})) = - \varepsilon \, ,
\end{equation}
where the ``$z$'' makes sense by factorizing the map $p$ in $X^2 - H$ through $\Theta^3 (X^2-H) \supset W_{(\infty)} ({\rm BLACK}) \times [-\varepsilon \leq z \leq \varepsilon]$, and also it should satisfy the following modulating prescriptions.

\medskip

\noindent (5.33) \quad Every time our $({\rm zipping})^{-1} [XY]$ encounters a RED arc $[X'Y']$ (like the real $[AB]$ or the fake $[DC]$, in figure~5.8.(A)), then it makes use of the already defined $S_j [X'Y']$ (which goes to the fins). With these things, at least for the part $[0,t_0]$ of $({\rm zipping})^{-1}$, the (5.31.1) to (5.31.3) in our lemma are well satisfied.

\medskip

When we get inside ${\mathcal N}_{\infty}^3$, at $t=t_0$, then the $W({\rm BLUE})$ does no longer carry Holes, we can move freely inside $-\varepsilon \leq z \leq \varepsilon$, and we can realize, immediately beyond $t_0$, the following condition
\setcounter{equation}{33}
\begin{equation}
\label{eq5.35}
z (X_{t_0-0}) = z (Y_{t_0-0}) \in \{ \pm \, \varepsilon \} \, .
\end{equation}
One should notice here that, because of lemma~4.3, which puts no restriction on $z \in [-\varepsilon , \varepsilon]$, realizing (\ref{eq5.35}) stays compatible with the (5.31.3) in our lemma~5.5.

\bigskip

\noindent [Comment. Several times in our argument we will need to move back and forth between conditions (\ref{eq5.33}) and (\ref{eq5.35}). Every time this happens, we have to be within the  juridiction of lemma~4.3 and the appropriate version of lemma~5.4.]

\bigskip

By now we have

\medskip

\noindent (5.35) \quad $(X_{t_0-0} , Y_{t_0-0}) \in M^2 (f)$ with a $\lambda^{-1} (X_{t_0-0} , Y_{t_0-\varepsilon})$ which is essentially the $\lambda^{-1} \mid [t_0 , -\infty]$ (in our heuristic view).

\medskip

It is, from now on, this $\lambda^{-1} (X_{t_0-0} , Y_{t_0-0})$ which has to be lifted to $\partial \, \Theta_1^{N+4} \times \partial \, \Theta_1^{N+4}$, or rather to $\partial S'_{\varepsilon} (\widetilde M (\Gamma)-H) \times \partial S'_{\varepsilon} (\widetilde M (\Gamma)-H)$, respecting (5.31.1) $+$ (5.31.3). The $\lambda^{-1}$ consists of successive pieces in $\lambda^{-1} \cap {\mathcal N}_{\infty}^3$, $\lambda^{-1} \cap \overline{\mathcal N}_{\infty}^3$. The pieces $\lambda^{-1} \cap \overline{\mathcal N}_{\infty}^3$ can belong to one of the following three types:

\medskip

I) $W_{\infty} ({\rm BLACK})_{H^0} \cap W({\rm BLUE})$, II) $W_{\infty} ({\rm BLACK})_{H^1} \cap W({\rm RED} - H^0)$ and III) $W_{\infty} ({\rm BLACK})_{H^1} \cap W({\rm RED} \cap H_0) \subset H^0$, which is the very short piece $[u,v]$, barely visible in the figure~5.5. There, it stretches from ${\mathcal M}^3$ to $S_{\infty}^2$. Just before we get to the pieces I or II above, we are still inside ${\mathcal N}_{\infty}^3$, and then (\ref{eq5.35}) can be rechanged into (\ref{eq5.33}). 

\smallskip

Then, the I) and II) can be treated just like we did with the piece $[0,t_0]$ above, and nothing more will be said concerning them. When we get to III), and this is supposed to be followed by a journey of type II), then we change (\ref{eq5.35}) into (\ref{eq5.33}), and worry about the branches $z = +\varepsilon$ and $z=-\varepsilon$, independently.

\smallskip

The III) goes through infinitely many triple points, like the one in figure~5.8, and each time our $\lambda^{-1}$ crosses another $(W_n ({\rm BLUE}) - H)$, coming with $n \to \infty$. This then also means that $\lambda^{-1}$ goes by infinitely many ${\rm BLACK} \, [XY]$'s and ${\rm RED} \, [X'Y']$'s, each of them being a special BLUE $1$-handle of our $W_n ({\rm BLUE}) - H$; see here the figure~5.8. In view of our  high dimensions, these $[XY], [X'Y']$'s can be ignored. More importantly still, going now through the available
$$
({\rm Hole})_n \subset W_n ({\rm BLUE}) \, ,
$$
one can lift the corresponding arc $\lambda^{-1} \cap \overline{\mathcal N}_{\infty}^3$ to an arc of $({\rm zipping})^{-1} [XY]$, staying compatible with (5.31.1) to (5.31.3). Our $N$ is high, and so this lift can be made disjoined of the $C(H({\rm BLUE}))$'s.

\smallskip

We move now to the $\lambda^{-1} \cap {\mathcal N}_{\infty}^3$'s, which may be either of type $W_{(\infty)} ({\rm BLACK}) \cap (W({\rm BLUE}) + W({\rm RED} \cap H^0))$ or $W({\rm RED} \cap H^0) \cap W({\rm BLUE})$. We may assume the arcs $\lambda^{-1} \cap {\mathcal N}_{\infty}^3$ confined inside the ${\mathcal M}^3$'s, see here figure~5.5 and formula (\ref{eq5.17}), and as soon as we enter such a piece $\lambda^{-1} \cap {\mathcal N}_{\infty}^3$ we change again (\ref{eq5.33}) to (\ref{eq5.35}).

\smallskip

We start by lifting our $\lambda^{-1} \cap {\mathcal N}_{\infty}^3$ to $\Theta^3 (X^2-H) \times \Theta^3 (X^2-H)$, thus disregarding the fact that our arc, even if of very short length may go to uncontrollable many points where $X^2$ is not a smooth $2$-manifold. Since with our high $N$ comes a natural embedding $\Theta^3 (X^2-H) \subset \partial S'_{\varepsilon} (\widetilde M (\Gamma)-H)$, this automatically takes case of (5.31.2). Once we are inside ${\mathcal N}_{\infty}^3$ we do not have to worry about DITCHES, Holes and partial DITCH fillings, nor about the $[X'Y']$ which may be in the way. We can certainly impose something like (\ref{eq5.33}), at the moment we leave ${\mathcal N}_{\infty}^3$ for the next $\bar {\mathcal N}_{\infty}^3$. There are clearly no problems with (5.31.1), and because of lemma~4.3, not with (5.31.3) either. This ends our proof. \hfill $\Box$

\bigskip

We finally can go to the $S_j [XY]$ for the ${\rm BLACK} \, [XY]$'s. Here are two typical cases, connected to figure~5.8.
\setcounter{equation}{35}
\begin{equation}
\label{eq5.37}
\{ S_j [AD] \, , \ \mbox{in the case of variant I in 5.8.(B)}\} = S_j [AB]  \underset{\mbox{\footnotesize$\overbrace B$}}{\cup} ({\rm zipping})^{-1} [BC] \underset{\mbox{\footnotesize$\overbrace C$}}{\cup} S_j [CD] \, ,
\end{equation}
and then also
$$
\{ S_j [AD] \, , \ \mbox{in the case of variant II}\} = ({\rm zipping})^{-1} [AD] \, .
$$
With this, the pushes and smearing are now just like in the RED case $[XY] = [BA]$ which was explained at length in the context of figure~5.10.

\smallskip

As another typical case, along with (\ref{eq5.37}), look now at figure~5.9 and at the ${\rm BLACK} \, [XY]$ which we see there. We can take now $S_j [XY] = \{({\rm zipping})^{-1} [XY]$ which makes use of the $S_j [\Omega_1 \, \Omega_2]({\rm RED})$, and which we stop at $[AB]\} \cup S_j [AB]({\rm RED})$. Other choices are possible too.

\smallskip

To stay with the same kind of argument, we will stop temporarily the discussion of lemma~4.1 and move to the

\bigskip

\noindent CONSTRUCTIONS OF THE ARCS $\gamma_i$, FROM LEMMA 4.8. Notice that, if it would only be a matter of having (\ref{eq4.50}) $+$ (4.50.1), then we could proceed just like in the standard textbooks of algebraic topology, and take some simple-minded arc joining $\alpha \, C^- (H_i)$ to $\eta \, \beta \, C^- (H_i)$ inside $\partial S'_u (M(\Gamma)-H)$. This could be just some arc of type $[XY]$. But then, we also insist of having (4.50.3) and this without violation of (4.50.2). All this can be happily satisfied if, forgetting about the standard textbooks, we use lemma~5.5 and take
\begin{equation}
\label{eq5.38}
\gamma_i = ({\rm zipping})^{-1} [XY] \, , \quad \mbox{with $[XY]$ like just above.}
\end{equation}

\bigskip

\noindent {\bf Remark.} When we had an $[XY]$ which was an actual special BLUE $1$-handle, like it was the case for the $[AB],[AD]$ in the figure~5.8 then the $j[XY]=\beta[XY]$ contributed to the final geometric intersection matrix (of the still to be defined $S'_b \widetilde M (S)$).  But this is so even ever for the contribution of the $S_j$ or $({\rm zipping})^{-1}$ of the fake 1-handles, like $[BC],[DC]$ in figure~5.8. [Of course, as such, these arcs have nothing to do with the geometric intersection matrix, but their $({\rm zipping})^{-1}$ contributes to the real $\beta [XY]$'s.]

\smallskip

Now, when it comes to our arcs $\gamma_i$ from (\ref{eq5.38}), there is absolutely no contribution to the geometric intersection matrices. \hfill $\Box$

\bigskip

We go back to our lemma~4.1 and, in particular, we will prove now the CLAIM~(4.31). Remember here that not only is $X^2$ GSC, but every individual $X^2 \mid H_i^{\lambda}$ is also GSC (see (1.11.1)). Next, here is a general recipee for collapsing away most of $X^2$, after some redundant 2-cells might have been deleted.

\medskip

i) Going to $Y(\infty)$, consider infinite sequences of states $j_0 = (H_{j_0}^2 , H_{j_0}^1)$, $j_1 = (H_{j_1}^2 , H_{j_1}^1) , \ldots$, where inside $j_n$ there is the canonical diagonal contact $\partial H_{j_n}^2 \cdot \delta H_{j_n}^1 = 1$, and which also come with trajectories $\partial H_{j_{n+1}}^2 \cdot \delta H_{j_n}^1 = 1$.

\medskip

ii) Next, going 2-dimensional, each $H_j^2$ comes with its unique complete $W_j ({\rm BLACK})$ and, inside $X^2$, $W_{j_{n+1}} ({\rm BLACK})$ is incident to
$$
W_{i_n} ({\rm RED} - H^0) \subset \left[ W_{i_n} ({\rm RED}) \cup R(\infty)_{n,n-1} \cup \ldots \cup W_{i_n - N(n)} ({\rm RED}) \right] \subset X^2 \mid H_{j_n}^2 \, ,
$$
where the next $W_{j_n} ({\rm BLACK})$ is incident to $W_{i_n -N(n)} ({\rm RED} - H^0)$, with a pattern which by now should be clear. Suppressing a lot of $R(\infty)$'s and $W_{(\infty)} ({\rm BLACK})$'s, let us denote by $W_{i_n} ({\rm RED})$ the whole $[\ldots]$ above.

\medskip

iii) With this, collapse away all the pieces
$$
W_{i_0} ({\rm BLACK}) \cup W_{i_0} ({\rm RED} - H^0) \cup W_{i_1} ({\rm BLACK}) \cup \ldots \subset X^2 \, .
$$
This leaves us with a collection of $X^2 \mid H_i^0$'s connected by arcs in a tree-like manner. Our claim~(4.31) should be now transparent. \hfill $\Box$

\newpage

\section{Compactification}
\setcounter{equation}{0}

This section contains the proof of the COMPACTNESS LEMMA~4.7. Our action will happen mostly downstairs, at level $S'_u (M(\Gamma)-H)$. Of course, lemma~4.7 (see the (4.47.II)) mentions both $S'_b (M(\Gamma) - H)$ and $S'_u (M(\Gamma)-H)$, but then, via $\eta$ everything will be happily transferred on $S'_u (M(\Gamma)-H)$. Many of the things we will say make sense on $S'_u (\widetilde M (\Gamma)-H)$ too, and at least for reason of smoothness of exposition, the $S'_u (\widetilde M (\Gamma)-H)$ will be mentioned too. But, very importantly, we work here throughout with $S'_{\varepsilon}$ and never with $S_{\varepsilon}$, because we will want to be able to invoke the lemma~4.3.

\smallskip

Inside the $\sigma_1 (\infty)$ from (\ref{eq1.19}) we consider the more mundane $\sigma_2 (\infty)$ defined below
\begin{equation}
\label{eq6.1}
\sigma_1 (\infty) \subset \sigma_2 (\infty) \equiv \bigcup \left( \bigcap \, \{\mbox{limit walls}\}\right) \subset \widetilde M (\Gamma) \, .
\end{equation}

The $\sigma_2 (\infty)$ is a connected graph
\begin{equation}
\label{eq6.2}
\sigma_2 (\infty) = \sum {\rm circles} \, (S_{\infty}^2 \cap (S^1 \times I)_{\infty}) \cup \sum {\rm arcs} \, (\partial \, {\rm Hex}_{\infty} \cap \, (S_{\infty}^2 \cup (S^1 \times I)_{\infty})) \, ,
\end{equation}
where remember that $S_{\infty}^2 ({\rm BLUE})$, $(S^1 \times I)_{\infty} ({\rm RED})$, ${\rm Hex}_{\infty} ({\rm BLACK})$ are the limit walls of the respective colours.

\smallskip

The vertices of $\sigma_2 (\infty)$ are the $\{p_{\infty\infty}(\infty) ({\rm proper})\}$ (see figures~1.4 and 6.1), which are the intersections of the three kinds of limit walls and the $\{p_{\infty\infty}(\infty) (S)\}$ from figure~1.5.(C); they are produced by the two pairs of BLACK limit walls associated to our $\bar S \subset {\rm Sing} \, (\widetilde M(\Gamma))$. The $p_{\infty\infty}(\infty)$'s are accumulation points of $p_{\infty\infty}$'s, which occur inside the interiors of the edges of $\sigma_2 (\infty)$. When one moves from $\sigma_2 (\infty)$ to $\sigma_1 (\infty)$, then one has to add more vertices, namely the endpoints of the $\lim$ LIM living inside $\sigma_2 (\infty)$, figure~1.4.

\smallskip

We also have

\medskip

\noindent (6.3) \quad $\sigma_1 (\infty) \cap \Theta^3 (fX^2)({\rm II})$ (see (\ref{eq2.12})) $= \sigma_2 (\infty) \cap \Theta^3 (fX^2)({\rm II}) \subset \Bigl\{ \overset{\circ}{\sum} \, (\infty) \cup \underset{S}{\sum} \, p_{\infty\infty} (S) \times (-\varepsilon , \varepsilon)$, i.e. the ${\rm int} \, \sum (\infty)$, as defined in (2.13.1), with the contribution of the $p_{\infty\infty} (S)$'s restored (the one of the $p_{\infty\infty} ({\rm proper})$ was never deleted)$\Bigl\}$.

\medskip

The set defined in (6.3) is a collection of arcs, like in the figure~6.1. In this same connection, we introduce the closed subset
\setcounter{equation}{3}
\begin{equation}
\label{eq6.4}
\sigma (\infty) \equiv \sigma_1 (\infty) - \left\{{\rm int} \, \sum (\infty) \ \mbox{with the $p_{\infty\infty}(S)$ restored} \right\} = \overline{\sigma_1 (\infty) - \Theta^3 (fX^2)} ({\rm II}) \subset \widetilde M (\Gamma) \, .
\end{equation}
The little arcs $[x_{\infty} , y_{\infty}]$ in figure~6.1 are typical connected components of this $\sigma (\infty)$. More complicated such are gotten by assembling, at the level of figure~1.4, dotted lines in $\sigma_2 (\infty) - \Theta^3 (fX^2)({\rm II})$ and $\lim {\rm LIM}$'s.

\smallskip

At this point, we introduce an IMPORTANT EXPOSITORY TWIST. For expository purposes, we have decided to tackle the bad rectangles (\ref{eq4.1}) as late as possible in this story.

\smallskip

It is only in the present proof of the compactness lemma~4.7, that the bad rectangles $R_0$ have to be included into our soup of ingredients which are to be sent to infinity, when we construct the $S^{(')}_{\varepsilon}$'s. Our present $R_0$'s are somehow analoguous (although not quite), to the critical rectangle from (\ref{eq3.22}). Like we have already done it, in a toy-model context, in the section~III, we will distinguish now again, between a variant (or context) I where $R_0$'s are ignored and the real life variant (or context) II when the $R_0$'s are also sent to infinity. As long as one ignores our present compactness lemma~4.7 everything which was done or said,  so far, in this paper is valid in both contexts, reason for not distinguishing them until now.

$$
\includegraphics[width=160mm]{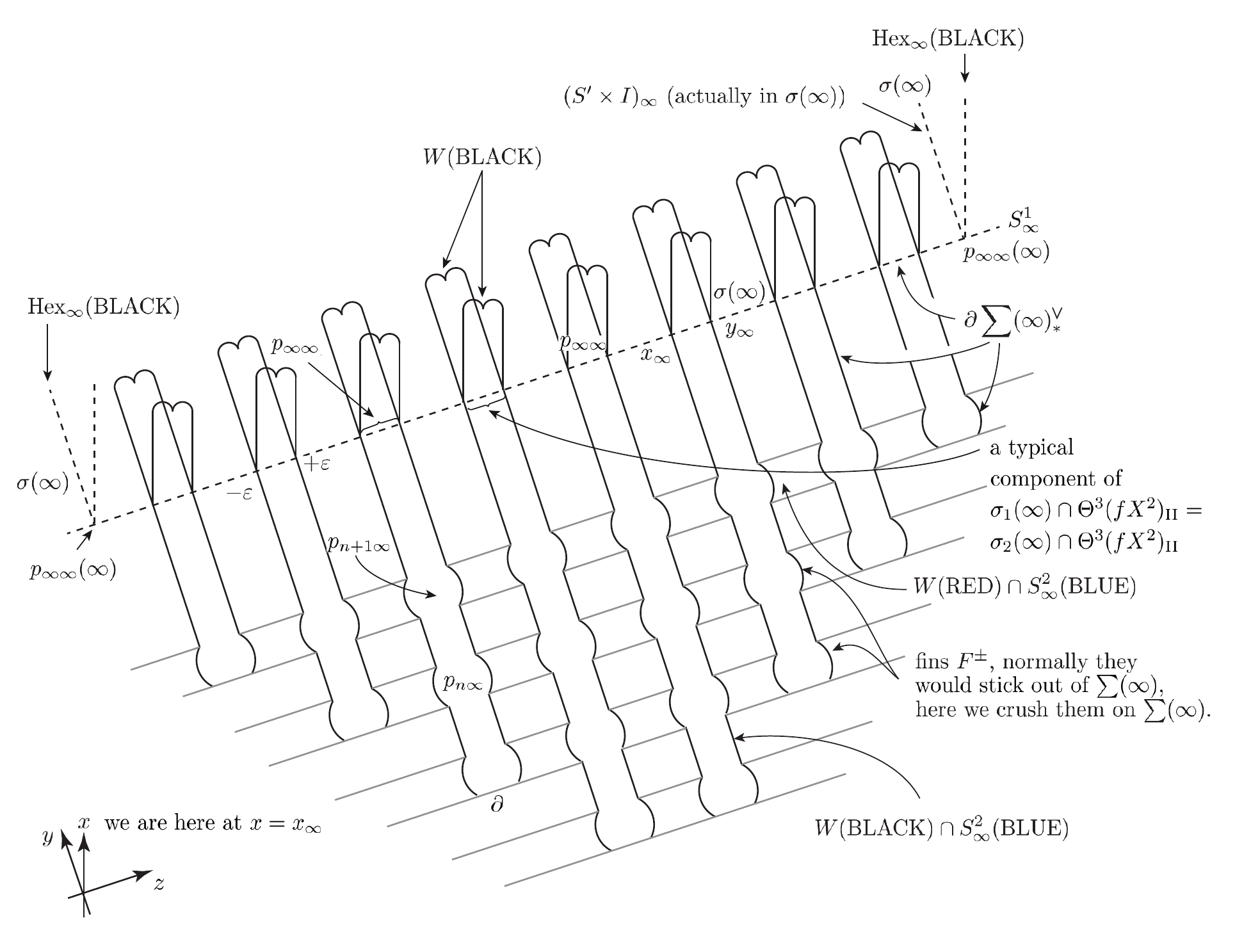}
$$
\label{fig6.1.}
\centerline {\bf Figure 6.1.} 

\smallskip

\begin{quote} 
We see here a detail of $\sum (\infty)$, at its most complicated. We see, mostly, a piece of $S_{\infty}^2 ({\rm BLUE})$, with a piece of $\sigma_2 (\infty)$ (or actually of $\sigma_1(\infty)$, since no $\lim {\rm LIM}$ is present here) occurring as dotted lines. The figure is, actually, an enlarged thickened detail of figure~1.4. The coordinate system is the same as in figure~1.1.(B). The details $W_{\infty} ({\rm BLACK})$ from figure~1.4 give rise to figures similar to this one, except that there, instead of a pair of points $p_{\infty\infty} (\infty)$, we have just one, endpoint of an arc $\lim {\rm LIM}$; see the figure~1.4. The present line $[\partial , p_{n\infty} , p_{n+1\infty} , p_{\infty\infty}]$ occurs in figure~1.1.(A) too. Notice how the $p_{\infty\infty}$'s accumulates on $p_{\infty\infty} (\infty)$. See here also figure~3.1 in \cite{31}. A small piece of the ``surface'' of infinite genus ${\rm int} \, \sum (\infty)$ is visible here. This is smooth except for ramification points along the $p_{\infty\infty} ({\rm proper}) \times (-\varepsilon , \varepsilon)$. These ramifications are of type $\{{\rm figure} \, Y \} \times R$, embeddable in $R^3$. When fins get included, and we move from $\sum(\infty)$ to $\sum(\infty)^{\wedge}$, then additional ramification points of the same type are to be included. Often, we will find it convenient to embed the $\sum(\infty)^{\wedge}$ {\ibf inside} $\sum (\infty)$, via the trick already used in this drawing of crushing each fin on its diameter $c$ $\sum(\infty)$. The $\includegraphics[width=5mm]{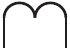}$ means truncation for typographical reasons.
\end{quote}

\bigskip

\noindent {\bf A remark concerning Figure 6.1.} Our figure above is an enlarged, detailed version of that part of figure~1.4, which concerns $H^2$. Then there should also be a figure~6.1.bis (not drawn explicitly), concerning the $W_{\infty} ({\rm BLACK})_{H^0}$'s rather than $W({\rm BLACK})$'s and 6.1.ter, concerning the $W_{\infty} ({\rm BLACK})_{H^1}$'s. The 6.1.bis is again on an $S_{\infty}^2 ({\rm BLUE})$ background, there is no longer the ${\rm Hex}_{\infty} ({\rm BLACK})$, but a $\lim {\rm LIM}$ line. The 6.1.ter, drawn on a background $(S^1 \times I)_{\infty}$ similarly has a $\lim {\rm LIM}$ line, and then continues without one on the $S_{\infty}^2 ({\rm BLACK})$ background. \hfill $\Box$

\bigskip

At the level of section~II we had a cell-complex $\Theta^3 (fX^2)^{(')}$ (see (\ref{eq2.13}) and (\ref{eq2.17})) at the level of which all of $\partial \sum (\infty)^{\wedge}$ was already sent to infinity, and we will call, from now on, this object $\Theta^3 (fX^2)^{(')}_{\rm I}$.

\smallskip

By sending the bad rectangles $R_0$ to infinity too (and see below for the details of this) we go from $\Theta^3 (fX^2)^{(')}_{\rm I}$ to the next level $\Theta^3 (fX^2)^{(')}_{\rm II}$, with which will come a
\begin{equation}
\label{eq6.5}
S'_{\varepsilon} (M(\Gamma)-H)_{\rm II} \supsetneqq S'_{\varepsilon} (M(\Gamma)-H)_{\rm I} \quad (\equiv \mbox{our previous} \ S'_{\varepsilon} (M(\Gamma)-H)) \, ,
\end{equation}
but it is only with $S'_u (M(\Gamma) - H)_{\rm II}$ that lemma~4.7 is valid (we will prove it in this section) and the diagram (\ref{eq4.47}) is PROPERLY homotopy commutative, but otherwise, everything said so far is valid in both contexts.

\smallskip

In (\ref{eq2.3}) we had already introduced a $\sum (\infty)$, but at the time, this was just a provisional starting point, soon to be superseded by the $\bigl\{ \sum (\infty)$ from (2.13.1), where the contribution of $p_{\infty\infty} (S)$ was removed$\bigl\}$. This was the object present in the ZIPPING LEMMA~4.1 and in its various complements, from 4.2 to 4.6. For $S_{\varepsilon} (\widetilde M (\Gamma)-H)$ the bad locus was this last $\sum (\infty)$, while for $S'_{\varepsilon} (\widetilde M (\Gamma)-H)$ it was $\bigl\{\sum (\infty)$ with the contribution $p_{\infty\infty}$ (all) removed$\bigl\}$. All these $\sum (\infty)$'s, came with $\partial \, \sum (\infty) \ne \emptyset$. The ballet of successive $\sum (\infty)$'s will continue now.
$$
\includegraphics[width=115mm]{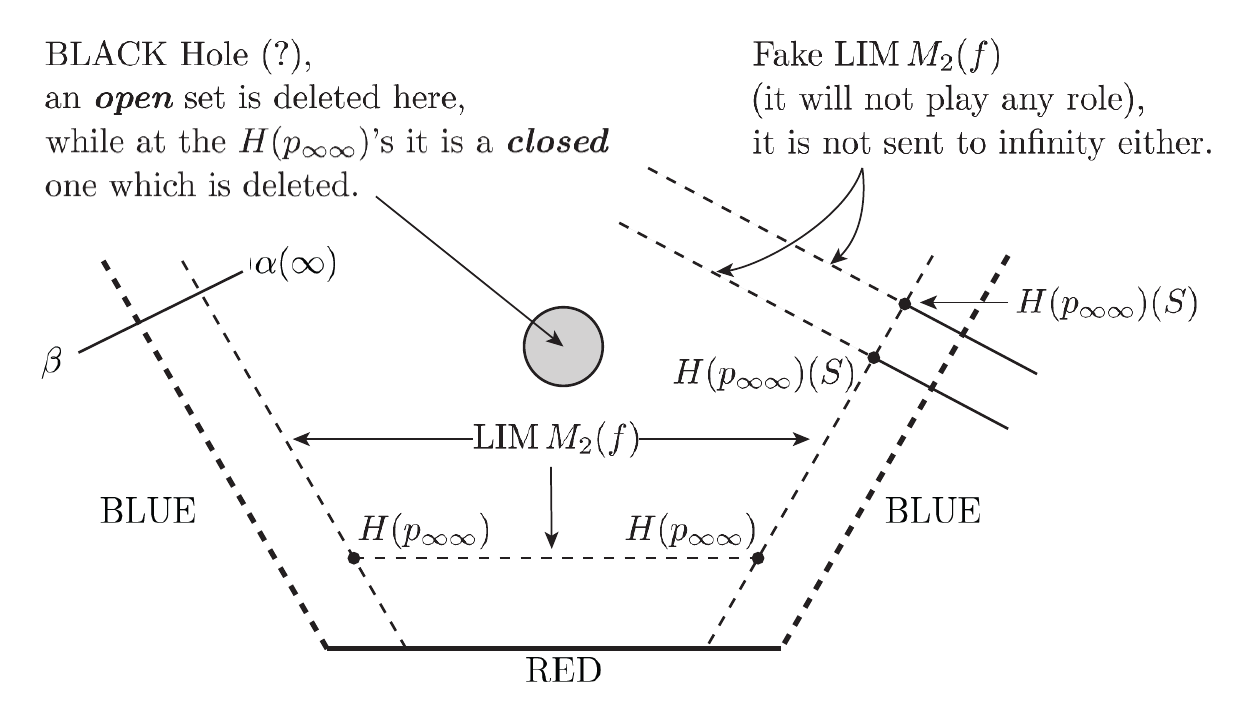}
$$
\label{fig6.2.}
\centerline {\bf Figure 6.2.} 

\smallskip

\begin{quote} 
A detail of figure 1.1.(A), without any $M_2 (f)$ drawn in. Like in (\ref{eq5.13}), for each $W({\rm BLACK})$ there is a unique {\ibf cut arc} $[\beta , \alpha (\infty))$ like in (\ref{eq5.13}), drawn here too.
\end{quote}

\bigskip

The next formula (6.6) presents the appropriate version(s) of $\sum (\infty)$ when one starts adding things at infinity so as to compactify (after quotienting by $\Gamma$). Keep in mind that the (6.6) coming below is a change of perspective with respect to our previous triplet
$$
\left( \Sigma (\infty) \ \mbox{(2.13.1)} \, ; \quad {\rm int} \, \Sigma (\infty) \ \mbox{(2.13.1)} \, , \quad \partial \, \Sigma (\infty) \ \mbox{(2.3)} \right) .
$$

\noindent (6.6) \quad We introduce now the 4-tuple
$$
\left( \Sigma (\infty)_* \, ; \quad {\rm int} \, \Sigma (\infty)_* \supsetneqq \overset{\circ}{\sum} (\infty)_* \, , \quad \partial \, \Sigma (\infty)_* \right) ,
$$
where 

\medskip

\noindent $\sum (\infty)_* \equiv \{$the $\sum (\infty)$ from (2.13.1), with $\underset{S}{\sum} \, p_{\infty\infty}(S) \times [-\varepsilon , \varepsilon]$ put back$\} = \{$the original $\sum (\infty)$ from (\ref{eq2.3})$\}$,

\medskip

\noindent ${\rm int} \, \sum (\infty)_* \equiv \{$the $\sum (\infty)$ (\ref{eq2.3})$\} - \{$its $\partial \, \sum (\infty)\}$,

\medskip

\noindent $\overset{\circ}{\sum} \, (\infty)_* \equiv \Bigl\{ \overset{\circ}{\sum} \, (\infty) \equiv {\rm int} \, \sum (\infty)$ (2.13.1)$\Bigl\} - \underset{p_{\infty\infty} ({\rm all})}{\sum} p_{\infty\infty} \times [-\varepsilon , \varepsilon]$,

\medskip

\noindent (where one should remember that at level $\sum (\infty)$ (2.13.1) the $\underset{p_{\infty\infty} (S)}{\sum} p_{\infty\infty} (S) \times [-\varepsilon , \varepsilon]$ had already been deleted)

\medskip

\noindent $\partial \, \sum (\infty)_* \equiv [\partial \, \sum (\infty)$ (from (\ref{eq2.3})$] \cup \underset{p_{\infty\infty} ({\rm all})}{\sum} p_{\infty\infty} \times [-\varepsilon , \varepsilon]$, with the two terms being glued together along $\underset{p_{\infty\infty} ({\rm all})}{\sum} p_{\infty\infty} \times \{ \pm \, \varepsilon \}$.

\medskip

\noindent This ends our formula (6.6).

\bigskip

At this point, formula (4.5.1) is to be expanded now into the following disjoined decomposition (see here (6.3) $+$ (\ref{eq6.4}) too)
\setcounter{equation}{6}
\begin{equation}
\label{eq6.7}
\overline{\sum (\infty)_*} = \sum^0 (\infty)_* + \partial \, \sum (\infty)_* + \sigma (\infty) \, ,
\end{equation}
where the closure is inside $\sum (\infty)_1$, and to which ``topological relations'' of the form $\lim x_n = x_{\infty}$, should also be added. When fins are thrown into the game too, then (6.6) extends to things like

\medskip

\noindent (6.8.1) \quad $\sum (\infty)_*^{\wedge} \equiv \sum (\infty)_* \cup ({\rm fins}) \subset \Theta^3 (fX^2-H)'_{\rm I}$ coming with $\partial \, \sum (\infty)_*^{\wedge} = \partial \, \sum (\infty)_* \cup \sum (\mbox{rims of fins}) = \partial \, \sum (\infty)^{\wedge} \cup \sum \, p_{\infty\infty} ({\rm all}) \times [-\varepsilon , \varepsilon]$, with $\sum (\infty)$ any of the objects in (\ref{eq2.3}) or in (2.13.1), and with ${\rm int} \, \bigl( \overset{\circ}{\sum}  (\infty)_* \cup {\rm fins} \bigl) \equiv \overset{\circ}{\sum} (\infty)_* \cup (\mbox{fins-rims}) = \sum (\infty)_*^{\wedge} - \partial \, \sum (\infty)_*^{\wedge}$.

\medskip

Also, the (\ref{eq6.7}) should now be further expanded into
$$
\overline{\sum (\infty)_*^{\wedge}} = {\rm int} \, \left( \sum (\infty)_* \cup {\rm fins} \right) + \partial \, \sum (\infty)_*^{\wedge} + \sigma (\infty) \, . \leqno (6.8.2)
$$
Finally cut arcs are defined like in (\ref{eq5.13}) and in figure~6.2. They are far from the fins and a further variation on (6.6) are the following useful objects
\setcounter{equation}{8}
\begin{equation}
\label{eq6.9}
\sum (\infty)_* ({\rm cut}) \equiv \sum (\infty)_* - \{\mbox{cut arcs}\} \subset \sum (\infty)_*^{\wedge} ({\rm cut}) \, .
\end{equation}

\noindent {\bf Remark.} Notice that, according to (6.6) our present $\sum (\infty)_*^{\wedge}$ is actually $\{$the old simple-minded $\sum (\infty)$ from (\ref{eq2.3})$\} \cup \{{\rm fins}\}$. \hfill $\Box$

\bigskip

The system of cut arcs is $\Gamma$-equivariant and everything said functions downstairs too. It should be stressed that the cut arcs do not correspond to any physical deletion or puncture. They serve for the following two related purposes: Firstly, they will be essential for the $\pi_1$-injectivity issues to be developed below and connected with this, all the part of the zipping flow which is necessary for the $\beta$ (\ref{eq4.36}), with its $\beta [XY]$, and for the $\Lambda (H_i^-)$'s (\ref{eq4.50}), will be driven so as to accumulate only on $\overline{\sum (\infty)_*^{\wedge}} \, ({\rm cut})$, avoiding the cut arcs (which only see $0(3)$'s).

\smallskip

We consider now the following commutative diagram of inclusion maps
\begin{equation}
\label{eq6.10}
\xymatrix{
\overset{\circ}{\sum} (\infty)({\rm cut}) \ar[rr]^-{\alpha} &&\Theta^3 (fX^2 - H ({\rm normal})) \\ 
\overset{\circ}{\sum} \, (\infty)_* ({\rm cut}) \ar[rr]^-{\alpha} \ar[u] &&\Theta^3 (fX^2-H)' \, .\ar[u]
}
\end{equation}
[In the first line, all the $p_{\infty\infty} (S)$'s are deleted here, but NOT the $p_{\infty\infty} ({\rm proper})$'s. NO $D^2 (H(p_{\infty\infty}))$ is deleted, either.]

\smallskip

\noindent [In the second line, ALL the $H$'s, including the $H(p_{\infty\infty})$ and the corresponding $D^2 (H)$'s are deleted now.]

\smallskip

In the context of this diagram, we have the

\bigskip

\noindent {\bf Lemma 6.1.} {\it The following map inject
$$
\pi_1  \left( \sum^0 (\infty) ({\rm cut}) \right) \overset{\alpha_*}{-\!\!\!-\!\!\!-\!\!\!-\!\!\!-\!\!\!-\!\!\!-\!\!\!\longrightarrow} \pi_1 \, \Theta^3 (fX^2 - H ({\rm normal}))
\eqno (6.11.1)
$$
$$
\pi_1  \left( \sum^0 (\infty) ({\rm cut}) / \Gamma \right) \overset{(f\alpha)_*}{-\!\!\!-\!\!\!-\!\!\!-\!\!\!-\!\!\!-\!\!\!-\!\!\!\longrightarrow} \pi_1 \, (\pi \, \Theta^3 (fX^2 - H ({\rm normal}))) \, ,
\eqno (6.11.2)
$$
where, of course we also have $\pi_1 \, \Theta^3 = \Theta^3 / \Gamma$. It follows from} (6.11.1), (6.11.2) {\it that the same injectivities as above hold for the lower $\alpha$ in} (\ref{eq6.10}).

\bigskip

\noindent {\bf Proof.} Both upstairs and downstairs, the proofs are easy applications of Van~Kampen, independent of each other. One proves first the $\pi_1$-injectivity at the local level of the individual $0$-handles and $1$-handles. For the 2-handles the cut arcs are necessary for taking care of those $W (\mbox{BLACK complete})$ not carrying Black Holes. Finally, one glues together all the local data, and we invoke Van~Kampen again. \hfill $\Box$

\bigskip

With our high $N \gg 1$, we get a more or less canonical embedding
\setcounter{equation}{11}
\begin{equation}
\label{eq6.12}
{\rm int} \, \left(\sum (\infty)_* \cup {\rm fins} \right) \subset \Theta^3 (fX^2 - H)'_{\rm I} \ \overset{\mbox{\footnotesize smooth embedding}}{\underset{i}{-\!\!\!-\!\!\!-\!\!\!-\!\!\!-\!\!\!-\!\!\!-\!\!\!-\!\!\!-\!\!\!-\!\!\!-\!\!\!-\!\!\!-\!\!\!-\!\!\!-\!\!\!\longrightarrow}} \ \partial S'_u (\widetilde M (S)-H)_{\rm I} \, .
\end{equation}
Starting from (\ref{eq6.12}), one has the following fact

\medskip

\noindent (6.13) \quad  One can glue $\left( \partial \, \sum (\infty)_*^{\wedge} \cup \sigma (\infty)\right) / \Gamma$ (see here (6.8)) to the infinity of $S'_u (M(\Gamma) - H)_{\rm I}$. This way one gets the following space $S'_u (M(\Gamma) - H)_{\rm I}^{\wedge} \equiv S'_u (M(\Gamma)-H)_{\rm I} \cup \left[\left( \partial \, \sum (\infty)_*^{\wedge} \cup \sigma (\infty)\right) / \Gamma \right]$. We will introduce here the notations
$$
\partial_{\infty} \equiv \left[ \partial \, \sum (\infty)_*^{\wedge} \cup \sigma (\infty)\right] / \Gamma \supset \partial_{\infty} ({\rm cut}) \equiv \left[ \partial \, \sum (\infty)_*^{\wedge} \cup \sigma (\infty)\right]({\rm cut}) / \Gamma \, .
$$

\bigskip

\noindent {\bf Lemma 6.2.} 1) {\it Provided we take the metrically correct definition for the $H ({\rm normal})$, leading to} (\ref{eq4.8}), {\it the $S'_u (M(\Gamma)-H)_{\rm I}^{\wedge}$ above is exactly the {\ibf closure} of}
$$
S'_u (M(\Gamma)-H)_{\rm I} \subset \{\mbox{the compact metric space} \ M(\Gamma) \times B^L , \quad L \ {\rm large}\} \, .
$$

\noindent 2) {\it The space $S'_u (M(\Gamma)-H)_{\rm I}^{\wedge}$ is hence {\ibf compact}, and so is also the $\partial S'_u (M(\Gamma)-H)^{\wedge}_{\rm I} \equiv \partial S'_u (M(\Gamma)-H)_{\rm I} \, \cup \, \partial_{\infty}$.}

\medskip

\noindent 3) {\it Let  $\Lambda_n = \Lambda (H_n^-)$, be like in the lemma~{\rm 4.8.} For every subsequence $\Lambda_{i_1} , \Lambda_{i_2} , \ldots$ of $\Lambda_1 , \Lambda_2, \ldots$, there is a sub-sub-sequence $\Lambda_{j_1} , \Lambda_{j_2} , \ldots$ and also a closed continuous curve
\setcounter{equation}{13}
\begin{equation}
\label{eq6.14}
\Lambda_{\infty} \subset \partial_{\infty} ({\rm cut}) \subset \partial_{\infty} = \partial S'_u (M(\Gamma) - H)_{\rm I}^{\wedge} - \partial S'_u (M(\Gamma)-H)_{\rm I} \, ,
\end{equation}
which is such that we have
\begin{equation}
\label{eq6.15}
\lim_{n=\infty} \Lambda_{j_n} = \Lambda_{\infty} \, , \ \mbox{uniform convergence in} \ S'_u (M(\Gamma) - H)_{\rm I}^{\wedge} \, .
\end{equation}
}

\bigskip

\noindent {\bf Proof.} Inside the already compact $M(\Gamma)$, the $X^2 / \Gamma \subset M(\Gamma)$ naturally compactifies to its closure, which is
$$
(X^2 / \Gamma)^{\wedge} = (X^2/\Gamma) \cup \left\{ \Sigma_1 (\infty) = \bigcup \ \{\mbox{limit walls}\} \ \mbox{from (1.14)} \right\} \, . \eqno (*)
$$
Here the $\{ H({\rm normal})\}$ accumulate on the $\{$ideals Holes$\} \subset \sum_1 (\infty)$, when the correct metric choices have been made, while $\{ H(p_{\infty\infty})\}$ accumulates on $\{ p_{\infty\infty} (\infty)\} \subset \sigma_2 (\infty)$ (\ref{eq6.1}). The 1) $+$ 2) follow from the things. Concerning 3), it follows from (4.50.3) that the $\Lambda_n$'s can only accumulates on $\partial_{\infty}$, more precisely we will have
\begin{equation}
\label{eq6.16}
\lim_{n = \infty} {\rm dist} (\Lambda_n , \partial_{\infty}) = 0 \, .
\end{equation}
We will prove now the 3) and for the sake of rigour of our exposition, we will be now more pedantic than we usually are. From (4.50.2), which is essentially a consequence of the uniformly bounded zipping length, established in \cite{29}, there is a uniform bound $N$ such that $\Vert \Lambda_n \Vert \leq N$. Here, the $\Lambda_n$'s are closed curves living inside the smooth $(N+3)$-manifold $\partial S'_u (M(\Gamma)-H)$, which may be assumed, without loss of generality to be $C^{\infty}$, and $\Vert \ldots \Vert$ is the length thereby calculated; so the $\partial S'_u (M(\Gamma) - H)$ has a metric structure $d$, and the length is measured with respect to it. Remember how this comes about. To begin with, $M(\Gamma)$ has a metric $d$ which is riemannian on each of its smooth 3-cells, with compatibilities on the common 2-faces. This $d$ lifts as a $\Gamma$-equivariant metric on $\widetilde M (\Gamma)$ and from there on it lifts on our various objects of interest, like $Y(\infty)$, $X^2$ and $S_{\varepsilon}^{(')} (\widetilde M (\Gamma))$ and $S_{\varepsilon}^{(')} (\widetilde M (\Gamma)-H)$. Only quasi-isometry classes count here, and we may assume that, on $S'_u (\widetilde M (\Gamma)-H)$ our $d$ is a $\Gamma$-invariant riemannian metric. This descends, afterwards, on $S'_u (M(\Gamma)-H)$, the object with which we work now.

\smallskip

Now, by adding well-chosen internal zigzags for each individual $\Lambda_n$ we may assume that for all $n$'s, we have strict equality $\Vert \Lambda_n \Vert = L$. These zigzags are internal reparametrizations which do not change the image and we may assume that there are uniform upper and lower bounds $M_1$, $\delta_1$ such that
$$
\# \, (\mbox{zigzags of} \ \Lambda_n) \leq M_1 \, , \ {\rm dist} \, (\mbox{two consecutive zigzags}) \geq \delta_1 \, .
$$

From now on, our $\Lambda_n$'s are parametrized by arc-length, and so with a universal $L > 0$ they are given by maps
$$
R_+ \supset [0,L] \overset{\Lambda_n \in \, C^{\infty}}{-\!\!\!-\!\!\!-\!\!\!-\!\!\!-\!\!\!-\!\!\!-\!\!\!-\!\!\!\longrightarrow} \partial S'_u (M(\Gamma) - H) \, , \ \mbox{such that for } \ x,y \in [0,L]
$$

\begin{equation}
\label{eq6.17}
\mbox{length of} \ \Lambda_n \mid [x,y] = \vert x-y \vert \, .
\end{equation}

\bigskip

\noindent {\bf CLAIM (6.17.1).} Our $\Lambda_n$'s may be chosen such that there are two universal constants $C_1 , C_2$, independent of $n$, such that
$$
C_2 \, d (\Lambda_n (x) , \Lambda_n (y)) \geq \vert x-y \vert \geq C_1 \, d (\Lambda _n (x) , \Lambda_n (y)) \, , \quad \forall \, n \, .
$$

\bigskip

\noindent {\bf Proof.} One may happily take $C_1 = 1$. The meaning of the first inequality is that the amount of contorsions that $\Lambda_n$ can have, is uniformly bounded; one may express this also in terms of control of ${\rm grad} \, \Lambda_n$. Now, our $\Lambda_n$ consists of three pieces, namely the $\alpha \, C^- (H_n) = \partial H_n^-$, $\eta \beta \, C^- (H_n)$ and the $\gamma_n$ (see (\ref{eq4.50})).

\smallskip

The infinitely many $H_n^-$'s come in infinite families, each corresponding to an ideal hole (see (\ref{eq4.8})), with all the $\alpha \, C^- (H_n)$'s in one family looking more or less like the boundary of the corresponding ideal hole; one can happily assume that here everything is fairly round, without uncontrolled contorsions.

\smallskip

Next, we consider the zipping paths $\lambda (x_0 , y_0)$.

\smallskip

The zipping paths $\lambda$ have actually been constructed in the last section of \cite{29}, where we also forced for them a uniformly bounded length. Each $\lambda$ comes with accidents, seable in figure~5.3, namely singularities (in ${\rm Sing} (f)$) and triple points (in $M^3 (f)$). The construction in \cite{29} can be easily seen to come with the following features:

\medskip

i) There is a uniform upper bound for the number of accidents which each $\lambda$ sees. There is also a uniform lower bound for the distance between them.

\medskip

ii) Between two accidents we can take our $\lambda$'s to be essentially straight, i.e. with a uniformly bounded ${\rm grad} \, \lambda$. All this means that, just like for $\alpha \, C^- (H_n)$, there is no uncontrollably much contorsion coming with $\lambda (x_0 , y_0)$. The $\beta \, C^- (H_n)$ consists of a piece which is essentially a shifted copy of $\alpha \, C^- (H_n)$, and see here the figure~5.7, plus a piece like $\lambda$. The $\gamma_n$ is also essentially like $\lambda$. Incidentally, we have talked here about ${\rm Sing} (f)$ and not about the immortal singularities; these are branching points of spaces and not contorsions of $\lambda$. \hfill $\Box$

\bigskip

At this point, we choose an infinite, increasing sequence of finite subsets $F_1 \subset F_2 \subset F_3 \ldots \subset [0,L]$ which is such that $\underset{i}{\overset{\infty}{\bigcup}} \, F_i$ is dense in $[0,L]$, i.e. $\overline{\underset{i}{\bigcup} \, F_i} = [0,L]$.

\smallskip

Via the diagonal argument we get the following item:

\medskip

\noindent (6.17.2) \quad There exists a subsequence $\Lambda_{j_1} , \Lambda_{j_2} , \ldots$ of $\Lambda_1 , \Lambda_2 , \ldots$ which is such that the 
$\Lambda_{j_1} , \Lambda_{j_2},\ldots$, when restricted to $\underset{1}{\overset{\infty}{\bigcup}} \, F_i$ converges {\ibf uniformly}.

\medskip

If $x \in \underset{1}{\overset{\infty}{\bigcup}} \, F_i$, we let $x_n \equiv \Lambda_{j_n} (x)$, $x_{\infty} \equiv \underset{x=\infty}{\lim} x_n$ (given by the (6.17.2)).

\bigskip

\noindent {\bf Claim (6.17.3).} With the $C_1 , C_2$ from (6.17.1), for each pair $x,y \in \underset{i}{\bigcup} \, F_i$, we have the inequalities
$$
C_2 \, d (x_{\infty} , y_{\infty}) \geq \vert x-y \vert \geq C_1 \, d (x_{\infty} , y_{\infty}) \, , \ \mbox{for} \ x,y \in \bigcup_{\rm I} F_i \, .
$$

\bigskip

\noindent {\bf Proof.} For any fixed $n$, the (6.17.1) tells us that we do have
$$
C_2 \, d (x_n , y_n) \geq \vert x-y \vert \geq C_1 \, d (x_n , y_n) \, ,
$$
and since $d (x_{\infty} , y_{\infty}) = \underset{x=\infty}{\lim} \, d (x_n , y_n)$, our claim (6.17.3) is proved by taking the limit. \hfill $\Box$

\bigskip

Let now $\alpha \in [0,L]$ be some arbitrary point, for which we choose a sequence $x(n) \in \underset{i}{\bigcup} \, F_i$ s.t. $\underset{m=\infty}{\lim} \, x(m) = \alpha$, in $[0,L]$. So, the sequence $x(1) , x(2) , \ldots$ is Cauchy in $[0,L] \subset R_+$. Via (6.17.1), the sequence $x(1)_n , x(2)_n \ldots$ $\in \Lambda_n$ is then Cauchy in $\partial S'_u (M(\Gamma) - H)$. We also have, via (6.17.2) that, for all $m$'s
$$
\lim_{n=\infty} x(m)_n = x(m)_{\infty} \, , \quad \mbox{{\ibf uniformly} in $m$}.
$$
It follows that $x(1)_{\infty} , x(2)_{\infty} , \ldots \in \partial_{\infty}$ is itself Cauchy, inside the compact space $\partial_{\infty}$, where it has a limit $\underset{m=\infty}{\lim} \, x(m)_{\infty} \equiv \alpha_{\infty} \in \partial_{\infty}$.

\smallskip

Let now $x'(m)$ be another sequence in $\bigcup F_i \subset [0,L]$ with the same limit $\underset{m=\infty}{\lim} \, x'(m) = \alpha$. We claim that, we have then
$$
\lim_{m=\infty} x(m)_{\infty} = \lim_{m=\infty} x'(m)_{\infty} \, , \eqno (6.17.4)
$$
i.e. the $\alpha_{\infty}$ does not depend on the particular sequence in $\underset{i}{\bigcup} \, F_i$ which approximates the $\alpha \in R_+$.

\bigskip

\noindent {\bf Proof.} With $m,n \in Z_+$ we have
$$
d(x(m)_{\infty} , x'(m)_{\infty}) \leq d(x(m)_{\infty} , x(m)_n) + d(x(m)_n , x'(m)_n) + d(x'(m)_n , x'(m)_{\infty}) \, .
$$
Because of the uniform convergence in (6.17.2), if $n$ is large enough, then independently of $m$, the two extreme terms in the RHS of the inequality above are $< \varepsilon$. We fix, from now on, such a large $n$. The middle term is also controlled now, since $x(m)_n \equiv \Lambda_{j_n} (x(m))$, $x'(m)_n \equiv \Lambda_{j_n} (x'(m))$ and hence, by (6.17.1), the $d(x(m)_n , x'(m)_n)$ is controlled by $\vert x(m) - x'(m) \vert$, which goes to zero when $m \to \infty$. \hfill $\Box$

\bigskip

So, by now we have defined a new map
$$
[0,L] \underset{\Lambda_{\infty}}{-\!\!\!-\!\!\!-\!\!\!\longrightarrow} \partial_{\infty} \subset \partial S'_u (M(\Gamma) - H)^{\wedge} \, ,
$$
by $\Lambda_{\infty} (x) = x_{\infty}$ when $x \in \underset{1}{\overset{\infty}{\bigcup}} \, F_i$ and $\Lambda_{\infty} (\alpha) = \alpha_{\infty}$ when $\alpha \in [0,L] - \,\underset{1}{\overset{\infty}{\bigcup}} \, F_i$. Each individual $\Lambda_{j_n}$ is a {\ibf closed} curve, i.e. it comes with $\Lambda_{j_n} (0) = \Lambda_{j_n} (L)$. There is no harm in choosing, at the very beginning, the $\bigcup F_i \subset [0,L]$ s.t. it contains both $0$ and $L$. Then
$$
\Lambda_{\infty} (0) \equiv 0_{\infty} = \lim_{n=\infty} \, 0_n = \lim_{n=\infty} \, \Lambda_{j_n} (0) =  \lim_{n=\infty} \, \Lambda_{j_n} (1) =  \lim_{n=\infty} \, 
1_n = 1_{\infty} \equiv \Lambda_{\infty} (L) \, .
$$

\bigskip

\noindent {\bf Claim (6.17.5).} The map $\Lambda_{\infty}$ is continuous i.e. it defines a continuous closed curve in $\partial_{\infty} ({\rm cut})$.

\bigskip

\noindent {\bf Proof.} Our claim follows, once we show that, for every $\alpha , \beta \in R_+$, for our $\alpha_{\infty} = \Lambda_{\infty} (\alpha)$ we have, by analogy with (6.17.1), the double inequality
$$
C_2 \, d (\Lambda_{\infty} (\alpha) , \Lambda_{\infty} (\beta)) \geq \vert \alpha - \beta \vert \geq C_1 \, d (\Lambda_{\infty} (\infty) , \Lambda_{\infty} (\beta)) \, . \eqno (6.17.6)
$$
For $\alpha , \beta \in \bigcup F_i$ this (6.17.6) is just our (6.17.3) above. For the general case, we take $x(m) , y(m) \in \bigcup F_i$ with
$$
\lim_{m=\infty} x(m) = \alpha \, , \quad \lim_{m=\infty} y(m) = \beta \, .
$$
The (6.17.3) tells us now that, at $\underset{i}{\bigcup} \, F_i$ level we have
$$
C_2 \, d (x(m)_{\infty} , y(m)_{\infty}) \geq \vert x(m) - y(m) \vert \geq C_1 \, d (x(m)_{\infty} , y(m)_{\infty}) \, . \eqno (*)
$$
Here, we plug in the definition of $\Lambda_{\infty}$, i.e. we take
$$
\lim_{m=\infty} x(m)_{\infty} = \alpha_{\infty} = \Lambda_{\infty} (\alpha) \, , \quad \lim_{m=\infty} y(m)_{\infty} = \beta_{\infty} = \Lambda_{\infty} (\beta) \, .
$$
If we let now $m \to \infty$ in the context of $(*)$, we get our desired (6.17.6). \hfill $\Box$

\bigskip

In order to complete the proof of 3) in our lemma~6.3, it remains to show that $\Lambda_{j_n}$ converges uniformly to $\Lambda_{\infty}$ (i.e. to prove the (\ref{eq6.15})).

\smallskip

By our diagonal argument used in (6.17.2), we know that for all $\varepsilon > 0$, there exists an $N \in Z_+$ such that, if $x \in \underset{1}{\overset{\infty}{\bigcup}} \, F_i$ and $m,n \geq N$, then $d(\Lambda_{j_m} (x) , \Lambda_{j_n} (x)) < \varepsilon$. We denote $x_n \equiv \Lambda_{j_n} (x)$ and we know now that there is an $x_{\infty}$ such that, with a possible minor modification  of $\varepsilon$, we have $d(x_n , x_{\infty}) < \varepsilon$. 

\smallskip

Let now $\alpha \in [0,L]$ and chose $x \in \underset{i}{\bigcup} \, F_i$ such that $\vert \alpha - x \vert < \frac\varepsilon2$. With $\alpha_n \equiv \Lambda_{j_n} (\alpha)$, $\alpha_{\infty} \equiv \Lambda_{\infty} (\alpha)$, $x_n$ like above and $x_{\infty} \equiv \Lambda_{\infty} (x)$, we have then
$$
d(\alpha_n , \alpha_{\infty}) \leq d (\alpha_n , x_n) + d(x_n , x_{\infty}) + d(x_{\infty} , \alpha_{\infty}) \leq \left( 1+\frac1{C_1} \right) \varepsilon \, .
$$
The proof of lemma~6.3 is now finished. \hfill $\Box$

\bigskip

One should notice that $S'_u (M(\Gamma) - H)_{\rm I}$ is a smooth manifold and so will also be $S'_u (M(\Gamma)-H)_{\rm II}$, but with one dimension higher, as we shall soon see.

\smallskip

We had originally defined bad rectangles in (\ref{eq4.1}) and actually never mentioned them again since, until now. But then since $\sum (\infty)^{\wedge}$ has been superseded by $\partial \, \sum (\infty)_*^{\wedge} \supset \underset{p_{\infty\infty} ({\rm all})}{\sum} \, p_{\infty\infty} \times [-\varepsilon , \varepsilon]$, this will bring about a complete revision, concerning the bad rectangles $R_0$ too. We come now with the following commutative diagram, which will supersede the (\ref{eq4.1})
\setcounter{equation}{17}
\begin{equation}
\label{eq6.18}
\xymatrix{
(R_0 , \partial R_0) \ar[rr] &&\left( \sum (\infty)_* \cup \underset{F}{\sum} \frac12 D^2 (F) \equiv \sum (\infty)_*^{\wedge} , \partial \, \sum (\infty)_*^{\wedge} \right) \\ 
{\rm int} \, R_0 \ar[rr] \ar[u] &&{\rm int} \left( \overset{\circ}{\sum} \, (\infty)_* \cup {\rm fins}\right) \ \mbox{(see (6.8.1)).} \ar[u]
}
\end{equation}
In this formula, we refer to (6.6). We should also notice that, while $\sum (\infty) \varsubsetneqq \sum (\infty)_*$, we have $\overset{\circ}{\sum} \, (\infty)_* \varsubsetneqq {\rm int} \left( \sum (\infty) (2.3) \right) \equiv {\rm int} \, \sum (\infty)_*$, hence also ${\rm int} \left( \sum (\infty)_* \cup {\rm fins} \right) \varsubsetneqq {\rm int} \left( \sum (\infty) \cup {\rm fins} \right)$. Of course, we also have $\partial \, \sum (\infty)$ (like in (\ref{eq2.3})) $\varsubsetneqq \partial \, \sum (\infty)_*$, and the same kind of thing, when hats are being added.

\smallskip

At this point, we may as well assume that there are enough arcs $p_{\infty\infty} \times [-\varepsilon , \varepsilon]$ and fins too, so that we should have
$$
\Sigma (\infty)_* \cup \sum_F \frac12 \, D^2 (F) = \bigcup \, \{\mbox{bad rectangles $R_0$}\} \, , \eqno (6.18.0)
$$
an equality which should be understood ``with multiplicities'', by which we mean that there is a natural {\ibf surjection}, restricting to an injection (\ref{eq6.17}) for each individual $R_0$
$$
\sum_{R_0} R_0 -\!\!\!-\!\!\!\twoheadrightarrow \Sigma (\infty)_* \cup \sum_F \frac12 \, D^2 (F) = \Sigma (\infty)_*^{\wedge} \, , \eqno (6.18.1)
$$
and then in the same spirit, a second equality with multiplicities
$$
\partial \, \sum (\infty)_*^{\wedge} = \bigcup_{R_0} \, \partial R_0 \, , \eqno (6.18.2)
$$
and then also
$$
{\rm int} \, \left(\sum (\infty)_* \cup {\rm fins}\right) = \bigcup_{R_0} {\rm int} \, R_0 \, , \ \mbox{see here (6.8.1)} \, . \eqno (6.18.3)
$$
One may think of the RHS of (6.18.1) as living inside
$$
\Theta^3 (\pi f X^2 - H)'_{\rm I} \cup \partial \, \sum (\infty)_*^{\wedge} \subset \partial S'_u (M(\gamma) - H)_{\rm I}^{\wedge} \, .
$$

The various bad rectangles $R_0$ belong to the types below

\medskip

\noindent (6.19) \quad $\left[ [S_{\infty}^2 \ {\rm or} \ (S^1 \times I)_{\infty}] \cap W_{(\infty)} ({\rm BLACK}) \right] \times [-\varepsilon, \varepsilon]$ and in the case $W({\rm BLACK})$, these rectangles may be ``Cut'', like in (\ref{eq6.9}); then, similarly to this, we also have $[S_{\infty}^2 \cap W({\rm RED})] \times [-\varepsilon, \varepsilon]$, producing $R_0$'s too.

\medskip

In these formulae, the $[-\varepsilon, \varepsilon]$ is the thickness of $\sum (\infty)$, i.e. the thickness of the whole line $[ \ldots ]$ inside the corresponding limit wall. Of course, also, (6.19) needs to be amplified by the contribution from the fins. Each individual rectangle $R_0$ has two long sides, corresponding to the $\pm \, \varepsilon$ above and living in $\partial \sum (\infty)$ and then also two short sides which may be rims of fins, or $p_{\infty\infty} \times [-\varepsilon, \varepsilon]$'s, living anyway at infinity from the viewpoint of $\Theta^3 (fX^2 - H)'_{\rm I}$. When rims of fins are pushed inside $\overset{\circ}{\sum}  (\infty)$, a short side of $R'_0$ may rest on a long one of $R''_0$, creating a smooth figure $Y$. At $p_{\infty\infty} \times [-\varepsilon, \varepsilon]$, when $\overset{\circ}{\sum}  (\infty)$ is not smooth we may similarly get a non smooth $Y$. See here the figures~5.4 too.

\bigskip

\noindent THE TRANSFORMATION FROM $\Theta^3 (fX^2-H)'_{\rm I}$ TO $\Theta^3 (fX^2-H)'_{\rm II}$. To begin with, for each individual
$$
{\rm int} \, R_0 \subset {\rm int} \left(\sum (\infty)_* \cup {\rm fins}\right) \subset \Theta^3 (fX^2-H)'_{\rm I}
$$
we add a whole copy of $({\rm int} \, R_0) \times [0,\infty)$, starting from ${\rm int} \, R_0 \times \{ 0 \} = R_0$. It should be kept in mind here that, later on, at the completed (compactified) level, we will add
$$
\{ R_0 \times [0,\infty] \ \mbox{with each fiber $x \times [0,\infty]$, where $x \in \partial R_0$, is crushed into $x\}$},
$$
and this will allow us to write, along with (6.18.2)
\setcounter{equation}{19}
\begin{equation}
\label{eq6.20}
\partial \, \sum (\infty)_*^{\wedge} = \bigcup \, \partial R_0 \times \{ 0 \} \, .
\end{equation}
We move on now to the higher level of the transformation
$$
S'_{\varepsilon} (\widetilde M (\Gamma) - H)_{\rm I} \Longrightarrow S'_{\varepsilon} (\widetilde M (\Gamma) - H)_{\rm II} , \mbox{where $\widetilde M (\Gamma)$ may be happily changed to $M(\Gamma)$}.
$$
This will happen via the following successive steps.

\medskip

\noindent (6.21.1) \quad Consider, to begin with, the following inclusions
$$
\partial S'_b (\widetilde M (\Gamma)-H)_{\rm I} \supset {\rm int} \, \left(\sum (\infty)_* \cup {\rm fins}\right) \subset \Theta^3 (fX^2 - H)'_{\rm I} \subset \partial S'_u (\widetilde M (\Gamma)-H)_{\rm I}
$$
which induce the following inclusions, now at infinity
$$
\{\mbox{infinity of $S'_b$}\} \supset \partial \, \sum (\infty)_*^{\wedge} \subset \{\mbox{infinity of $S'_u$}\} \, .
$$

\medskip

For the first inclusion in the formula (6.21.1), see the figure~5.5.ter. 

\smallskip

Starting from these things, and totally disregarding the $\Theta^3 (fX^2 - H)'_{\rm I}$ part when dealing with $S'_b$ we add, in both cases, i.e. to $S'_{\varepsilon}$, the $\underset{R_0}{\sum} \, ({\rm int} \, R_0) \times [0 \leq u < \infty)$ along
$$
\sum_{R_0} {\rm int} \, R_0 \overset{\underset{R_0}{\sum} \ {\rm (individual \ inclusions)}}{-\!\!\!-\!\!\!-\!\!\!-\!\!\!-\!\!\!-\!\!\!-\!\!\!-\!\!\!-\!\!\!-\!\!\!-\!\!\!-\!\!\!-\!\!\!-\!\!\!-\!\!\!-\!\!\!-\!\!\!-\!\!\!-\!\!\!\longrightarrow} {\rm int} \, \left(\Sigma (\infty)_* \cup {\rm fins}\right) \, .
$$
Our first step towards the $S'_u (\widetilde M (\Gamma) - H)_{\rm II}$ will be to go from the $S'_u (\widetilde M (\Gamma) - H)_{\rm I}$ to
$$
S'_u (\widetilde M (\Gamma) - H)_{\rm I} \ \underset{\overbrace{\mbox{\footnotesize${\rm int} \left(\sum (\infty)_* \cup {\rm fins}\right) \subset \partial S'_u (\widetilde M (\Gamma) -H)_{\rm I}$}}}{\cup} \ \sum_{R_0} ({\rm int} \, R_0 \times [0,\infty)) \, . \eqno (*)
$$

The $S'_u (\widetilde M (\Gamma) -H)_{\rm I}$ is a smooth manifold, while $(*)$ is not, hence the next step. [Everything said above also makes sense for $S'_b$.]

\medskip

\noindent (6.21.2) \quad By taking a transversally compact regular neighbourhood around the added pieces one gets a smooth $(N+4)$-dimensional space $S'_{\varepsilon} (\widetilde M (\Gamma) -H)_{\rm II}$ (provisional). Here $S'_u (\widetilde M (\Gamma) -H)_{\rm II}$ (provisional) is a transversally compact thickening of $\Theta^3 (fX^2 - H)'_{\rm II} \equiv \{$the smooth $\Theta^3 (fX^2)$ with all Holes $(\supset H(p_{\infty\infty} ({\rm all})))$ deleted and with the compensating $D^2$'s too, no prime $(')$ is necessary right here$\} \cup \Bigl( \underset{R_0}{\sum} \, ({\rm int} \, R_0 \times [0,\infty)$, added along ${\rm int} \, (\sum (\infty)_* \cup {\rm fins})$, which it houses$\Bigl)$. But at the level of $S'_{\varepsilon} (\widetilde M (\Gamma) -H)_{\rm II}$ (provisional) ($\varepsilon = u$ OR $b$) the ${\rm int} \, (\sum (\infty)_* \cup {\rm fins})$ has been pushed towards the interior of our $S'_{\varepsilon} (\widetilde M (\Gamma) -H)_{\rm II}$ (provisional), and we want to bring it back to the boundary, where it should belong.

\medskip

\noindent (6.21.3) \quad So, we introduce the next, real-life object
$$
S'_{\varepsilon} (\widetilde M (\Gamma) -H)_{\rm II} \equiv [S'_{\varepsilon} (\widetilde M (\Gamma) -H)_{\rm II} \, ({\rm provisional})] \times [0,1]
$$
a definition which will get some modulations, afterwards. Anyway, we have now a canonical embedding, which extends the one in (\ref{eq6.12})
$$
\sum_{R_0} {\rm int} \, R_0 \times [0,\infty) \subset \Theta^3 (fX^2-H)'_{\rm II} \subset [S'_u (\widetilde M (\Gamma)-H)_{\rm II} ({\rm provisional})] \times \{1\} \subset
$$
$$
\partial S'_u (\widetilde M (\Gamma) -H)_{\rm II} \, , \ \mbox{and something similar for $\varepsilon = b$, with the $\Theta^3$ left out.}
$$

So we have now, superseding (\ref{eq6.5}), the following embedding which certainly is {\ibf not} codimension zero, when read from the first to the fourth term
$$
S'_{\varepsilon} (\widetilde M (\Gamma) -H)_{\rm I} \subset S'_{\varepsilon} (\widetilde M (\Gamma) -H)_{\rm II} ({\rm provisional}) \subset \partial S'_{\varepsilon} (\widetilde M (\Gamma) -H)_{\rm II} \subset S'_{\varepsilon} (\widetilde M (\Gamma) -H)_{\rm II} \, .
$$

At face value, what we should find now is a copy of $\partial_{\infty} \times [0,1]$ living at the infinity of $S'_u (\widetilde M (\Gamma) -H)_{\rm II}$. But then, via a very simple passage to the quotient at infinity, we can crush every $x \times [0,1]$, when $x \in \partial_{\infty}$, into $x$, and hence replace the $\partial_{\infty} \times [0,1]$ (where the factor $[0,1]$ is here the same as in the beginning of (6.21.3)) by the old $\partial_{\infty}$. This little passage to the quotient will be understood to be always there, from now on.

\smallskip

All these things were developed upstairs, at level $\widetilde M (\Gamma)$ and they are $\Gamma$-equivariant. So, one way or another, they make sense downstairs for $M (\Gamma)$ too.

\smallskip

So far we have been constantly working in the context $(\ldots - H)'$, but now we will drop the prime (and the ``$-H$''). We introduce first
\begin{eqnarray}
\sum^0 (\infty)^{\wedge} &\equiv &\left\{ {\rm int} \left( \sum (\infty)_*^{\wedge} \right) - \sum_{p_{\infty\infty}(S)} p_{\infty\infty} \times [-\varepsilon,\varepsilon] \right\} \nonumber \\
&= &\left\{ {\rm int} \, \sum (\infty)_* {\rm (6.6)} - \sum_{p_{\infty\infty}(S)} p_{\infty\infty} \times [-\varepsilon,\varepsilon] \right\} \cup ({\rm fins} - {\rm rims}) \varsupsetneqq {\rm int} \left( \sum^0 (\infty)_* \cup {\rm fins} \right) \, ,  \nonumber
\end{eqnarray}
where, while in $\overset{\circ}{\sum} \, (\infty)_*$ the whole contribution $p_{\infty\infty} ({\rm all})$ is deleted, in $\overset{\circ}{\sum} \, (\infty)^{\wedge}$ it is only the one of $p_{\infty\infty} (S)$ which is now deleted.

\smallskip

With this, we change the diagram (\ref{eq6.18}) into
$$
\xymatrix{
\left( \underset{R_0}{\sum} \ R_0 , \underset{R_0}{\sum} \ \partial R_0 \right) \ar[rr] &&\left( \sum (\infty)_*^{\wedge}, \partial \, \sum (\infty)_*^{\wedge} \right) - \underset{p_{\infty\infty}(S)}{\sum} \, p_{\infty\infty} \times (-\varepsilon,\varepsilon) \, . \\ 
\underset{R_0}{\sum} \ {\rm int} \, R_0 \ar[rr] \ar[u] &&\overset{0}{\sum} (\infty)^{\wedge} \ar[u]
}
$$

\smallskip

I claim, also, that things like the ZIPPING LEMMA and its various complements which were developed in the context I, remain perfectly valid in the context II, too.

\medskip

\noindent (6.21.4) \quad What will follow now is a useful, alternative description of that basic piece which was added going from I to II. Start with the inclusion

\medskip

\noindent $(*_1)$ \quad $A \equiv {\rm int} \Bigl( \overset{\circ}{\sum} \, (\infty)_* \cup {\rm fins}\Bigl) \, \cup \, \underset{R_0}{\sum} \, {\rm int} \, R_0 \times [0,\infty) \subset B \equiv \Bigl\{ \sum (\infty)_* \cup \underset{F}{\sum} \, \frac12 \, D^2(F) \cup \underset{R_0}{\bigcup} \, R_0 \times [0,\infty]$ where  each $\partial R_0 \times [0,\infty]$ is crushed into the corresponding $\partial R_0 \times \{0\} = \partial R_0 \subset \partial \, \sum (\infty)_*^{\wedge}$ (\ref{eq6.20})$\Bigl\}$.

\medskip

Notice that it is the $A$ in $(*_1)$ above which has been added during the step (6.21.1), and which lives now inside $\partial S'_u (\widetilde M (\Gamma)-H)_{\rm II}$, and/or (when ``$/\Gamma$'' included) inside $\partial S'_u (\widetilde M (\Gamma)-H)_{\rm II}$.

\smallskip

Similarly (with ``$/\Gamma$'' is included) the $B$ will soon be housed inside the compactification $\partial S'_u (\widetilde M (\Gamma)-H)_{\rm II}^{\wedge}$. But the point we want to make is that we have an equality

\medskip

\noindent $(*_2)$ \quad $B \underset{\rm TOP}{=} \Bigl\{ \sum (\infty) \times [0,\infty]$, with each $x \times [0,\infty]$, for $x \in \partial \, \overset{\wedge}{\sum} (\infty)_*$ crushed into $x \Bigl\}$. Here the inclusion $\partial \, \overset{\wedge}{\sum} (\infty)_* \subset \sum (\infty)_* \times [0,\infty]$ comes from (6.6).

\medskip

This equality $(*_2)$ is a simple local matter. In the neighbouhood of the $p_{\infty\infty} \times [-\varepsilon , \varepsilon]$, the situation is suggested in the figure~6.3, while in the neighbourhood of each $\frac12 \, D^2 (F)$ one can read it from the figure~6.5.

\smallskip

There is an obvious inclusion $\underset{R_0}{\bigcup} \, R_0 \times \{ \infty \} \subset B$, which combined with the $(*_2)$ yields the homeomorphism
$$
\bigcup_{R_0} \, R_0 \times \{\infty\} = \sum (\infty)_* \times \{\infty\} = \sum (\infty)_* \, . \eqno (*_3)
$$
$$
\includegraphics[width=11cm]{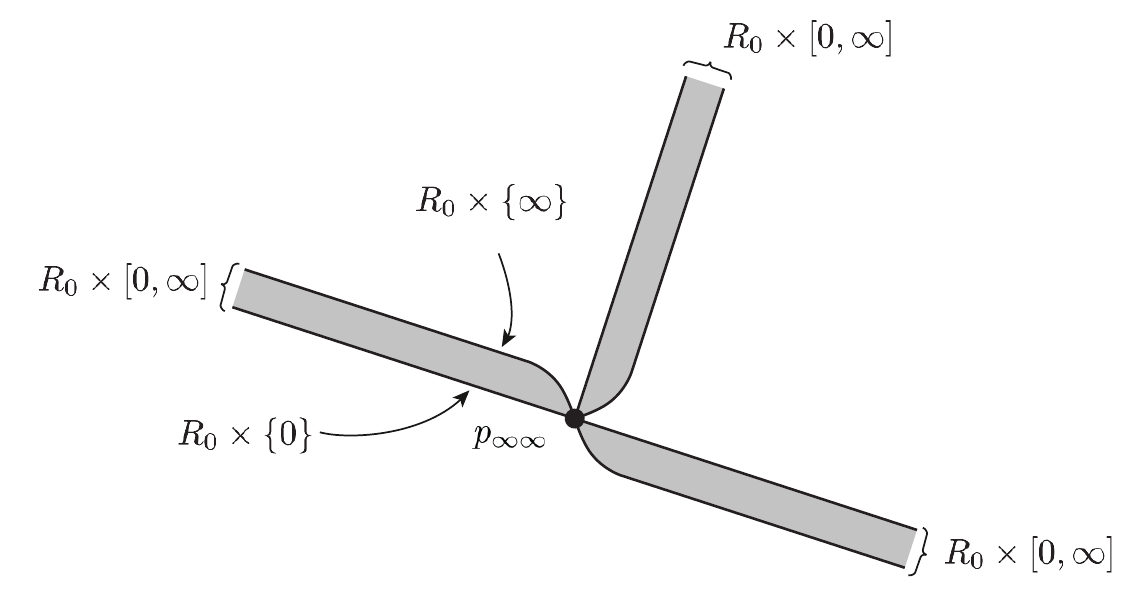}
$$
\label{fig6.3}
\centerline {\bf Figure 6.3.} 

\smallskip

\begin{quote} 
A schematical representation of $R_0 \times [0,\infty]$, in the neighbourhood of $p_{\infty\infty} \times [-\varepsilon , \varepsilon]$. The factor $[-\varepsilon , \varepsilon]$ is projected here down to a point.
\end{quote}

\bigskip

For further purposes (i.e. from the next and third paper in this series), we still have to phrase these things yet another way.

\bigskip

\noindent {\bf Complement (6.21.5)}, concerning the transformation
$$
\mbox{Variant I} \Longrightarrow \mbox{Variant II} \, .
$$
In the context of the variant II, we start by introducing the $\Theta^3 (fX^2)_{\rm II}$, which should not be mixed up with the $\Theta^3 (fX^2)({\rm II})$ from (\ref{eq2.12}). This is $\Theta^3 (fX^2)_{\rm II} \equiv \{ \Theta^3 (fX^2)$ (\ref{eq2.13}), where the contribution $p_{\infty\infty} (S)$ but not one of $p_{\infty\infty} ({\rm proper})$, is deleted$\} \cup \underset{R_0}{\sum} \, {\rm int} \, R_0 \times [0,\infty)$, the union being, like before, along ${\rm int} \left( \sum \, (\infty)_* \cup {\rm fins} \right)$. We have the decomposition
\begin{eqnarray}
\Theta^3 (fX^2)_{\rm II} &= &\left[\left\{ \Theta^3 (fX^2)(2.12) - \sum_{p_{\infty\infty} (S)} p_{\infty\infty} (S) \times (-\varepsilon , \varepsilon) \right\} \cup \sum_{R_0} {\rm int} \, R_0 \times [0,\infty) \right] \nonumber \\
&+ &\sum_{p_{\infty\infty} (S)} D^2 (H(p_{\infty\infty} (S))) \times \left[ - \frac\varepsilon4 , \frac\varepsilon4 \right] \, , \nonumber
\end{eqnarray} our ``$+$'' standing for addition of $2$-handles. Also, I will denote by $[\Theta^3]_{\rm II}$ the piece between brackets ($[\ldots]$) in the formula above, just before the addition of the $2$-handles.

\smallskip

With this comes now the variant II of the definition (2.19), namely the following
$$
\Theta^4 (\Theta^3 (fX^2) , R)_{\rm II} \equiv \Theta^4 ([\Theta^3]_{\rm II} , R) + \sum_{p_{\infty\infty} (S)} D^2 (H(p_{\infty\infty} (S))) \times \left[ - \frac\varepsilon4 , \frac\varepsilon4 \right] \times I \underset{{\rm projection} \, \pi}{-\!\!\!-\!\!\!-\!\!\!-\!\!\!-\!\!\!-\!\!\!-\!\!\!\longrightarrow} \Theta^3 (fX^2)_{\rm II} \, , \eqno (*_4)
$$
where $\pi \mid \Theta^4 ([\Theta^3]_{\rm II} , R) = \{$the natural retraction of the {\ibf smooth} $4^{\rm d}$ regular neighbourhood of $[\Theta^3]_{\rm II}$, on the singular space $[\Theta^3]_{\rm II}\}$, while $\pi \mid \{$the 2-handle part$\}$ is the obvious
$$
D^2 \times \left[ - \frac\varepsilon4 , \frac\varepsilon4 \right]  \times I \longrightarrow D^2 \times \left[ - \frac\varepsilon4 , \frac\varepsilon4 \right] \, .
$$
In $(*_4)$, the desingularization $R$ concerns the $\sum {\rm int} \, R_0 \times [0,\infty)$ too.

\smallskip

The main point of the present complement (6.21.5) is the following variant II version of the definition (2.19.bis)
$$
(S_u \, \widetilde M (\Gamma))_{\rm II} = \Theta^4 ([\Theta^3]_{\rm II} , R) \times B^N + \sum_{p_{\infty\infty} (S)} D^2 (H(p_{\infty\infty} (S))) \times \left[ - \frac\varepsilon4 , \frac\varepsilon4 \right] \times \frac12 \, B^{N+1} \, ,
$$
for very high $N$. This ends our complement (6.21.5). \hfill $\Box$

\bigskip

This also ends our description of the passage from I to II. Notice, also, that, with our trick of multiplying by $[0,1]$ and the quotienting at infinity $\partial_{\infty} \times [0,1]$ by the same $[0,1]$, the $\partial_{\infty}$ lives now at the infinity of $S'_u (M(\Gamma) - H)_{\rm II}$ and/or $\partial \, S'_u (M(\Gamma) - H)_{\rm II}$, compatibly with lemma~6.2. That lemma itself, melts now into the next

\bigskip

\noindent {\bf Lemma 6.3.} 1) {\it One can glue
\setcounter{equation}{21}
\begin{equation}
\label{eq6.22}
\left\{ \left[ \sum_{R_0} R_0 \times \{ \infty \} \right] \cup \sigma (\infty) \right\} / \, \Gamma \supset \partial_{\infty}
\end{equation}
at the infinity of $S'_u (M(\Gamma) - H)_{\rm II}$, thereby producing a compactification $S'_u (M(\Gamma) - H)_{\rm II}^{\wedge}$ of boundary $S'_u (M(\Gamma) - H)_{\rm II}^{\wedge} = \partial S'_u (M(\Gamma) - H)_{\rm II}^{\wedge} \cup \{$the contribution {\rm (6.22)}$\}$.}

\medskip

2) {\it With this, the analogue of {\rm 3)} from lemma~{\rm 6.2} remains valid, i.e. for every subsequence $\Lambda_{i_1} , \Lambda_{i_2} , \ldots$ of the $\Lambda_1 = \Lambda (H_1^-)$, $\Lambda_2 = \Lambda (H_2^-), \ldots$ there is a sub-sub-sequence $\Lambda_{j_1} , \Lambda_{j_2} , \ldots$ and a closed continuous curve $\Lambda_{\infty} \subset \partial_{\infty}$, which is such that we have
\begin{equation}
\label{eq6.23}
\lim_{n=\infty} \Lambda_{j_n} = \Lambda_{\infty} \, , \ \mbox{uniform convergence inside} \ S'_u (M(\Gamma)-H)_{\rm II}^{\wedge} \, .
\end{equation}
The $\Lambda_{j_n} , \Lambda_{\infty}$ are here {\ibf the same} as in lemma~{\rm 6.2}, coming with embeddings into $(\partial S'_u)_{\rm II}^{\wedge}$ which factorize via}
$$
\partial S'_u (M(\Gamma)-H)_{\rm I}^{\wedge} \subset \partial S'_u (M(\Gamma)-H)_{\rm II}^{\wedge} \, .
$$

\bigskip

For a while now we revert to the variant I and we will move to context II only when we will be forced to do so. Consider next the
\begin{eqnarray}
\label{eq6.24}
fX^2 - H(\mbox{completely normal}) - H(p_{\infty\infty}) &= &fX^2 - (H - \{\mbox{BLACK Holes}\})  \\
&= &\left[ \sum_n (W_n ({\rm BLUE}) - H) + \sum_n (W_n ({\rm RED} - H^0)-H) \right] \nonumber \\
&\cup &\Biggl[\sum_n \left(W_{(\infty)} ({\rm BLACK})_n - H(p_{\infty\infty} ({\rm all}))\right) \nonumber \\ 
&&\quad \cup \, \sum_n W ({\rm RED} - H^0)_n \Biggl] \nonumber
\end{eqnarray}
and here we draw the separation line between $W ({\rm RED} \cap H^0)$ and the rest of $W ({\rm RED})_n$ at $x = x_{\infty} + y_n$ with $y_n > 0$ and $\underset{n=\infty}{\lim} \, y_n = 0$. With this
$$
\overset{\circ}{\sum} (\infty)_* \subset \sum_n (W_{(\infty)} ({\rm BLACK})_n - H(p_{\infty\infty})) \cup \sum_n W({\rm RED} \cap H^0)_n \, .
$$
For each $x \in \overset{\circ}{\sum} \, (\infty)$ (where $\sum p_{\infty\infty} ({\rm all}) \times (- \varepsilon , \varepsilon)$ is deleted) there is a locally defined {\ibf transversal parameter} $(- \, \bar\varepsilon , \bar\varepsilon)$ to $\overset{\circ}{\sum} \, (\infty)$ (see (6.6)) inside the corresponding $W \times (-\varepsilon , \varepsilon)$ and here $\bar\varepsilon = \bar\varepsilon (x)$ is a $C^{\infty}$ function of $x$. To begin with, we will work here with the $\sum (\infty)$ from (2.13.1) which, remember, contains the contribution of the $p_{\infty\infty} ({\rm proper})$ but NOT of the $p_{\infty\infty} (S)$'s.

\smallskip

But then, it will be assumed that our $\bar\varepsilon$ is ``extended'' to another positive function
$$
\sum (\infty)_*^{\wedge} \ (\mbox{see (6.8.1)}) \overset{\varepsilon_0}{-\!\!\!-\!\!\!\longrightarrow} R_+ \eqno (6.24.1)
$$
such that $\varepsilon_0 \mid \overset{\circ}{\sum} \, (\infty)^{\wedge} \in C^{\infty}$ and $\varepsilon_0 (x) > \bar\varepsilon (x)$ in $\overset{\circ}{\sum} \, (\infty) - \sum p_{\infty\infty} \times (-\varepsilon , \varepsilon)$, when $\bar\varepsilon$ is defined. Here, we use the notation $ \overset{\circ}{\sum} \, (\infty)^{\wedge}$ from (6.21.5). Making use of $\varepsilon_0$ we can define the following neighbourhoods of $\sum (\infty)_*^{\wedge}$
\begin{equation}
\label{eq6.25}
{\mathcal N}^3 \equiv \sum (\infty)_*^{\wedge} \times (-\varepsilon_0 , \varepsilon_0) \, ,
\end{equation}
and in the same vein
$$
({\mathcal N}^3)' \equiv {\mathcal N}^3 - \partial \, \sum (\infty)_*^{\wedge} \times (\varepsilon_0 = 0) \subset \Theta^3 (fX^2-H)'_{\rm I} =
$$
$$
\{\Theta^3 (fX^2)' \ \mbox{(2.17), with all the holes and with all the $D^2$'s removed}\} \eqno (6.25.1)
$$

It should be understood that (6.24.1) to (6.25.1) are $\Gamma$-equivariant and the notations ${\mathcal N}^3$, $({\mathcal N}^3)'$ will be happily used when everything is quotiented by $\Gamma$ and $\sum (\infty)_*^{\wedge}$ is replaced by $\sum (\infty)_*^{\wedge} / \Gamma$. This kind of thing will be current practice in this section.

\smallskip

These objects enter the following commutative diagram, where the middle $r$ is the obvious {\ibf retraction}, the other $r$ is its restriction to $({\mathcal N}^3)'$, all the other maps are inclusions and, from now on we will normally work downstairs. So, the notation ``$\Theta^3 (fX^2-H)_{\rm I}$'' will be often used instead of $\Theta^3 (\pi fX^2-H)_{\rm I} = \Theta^3 (fX^2-H)_{\rm I}/\Gamma$, with the $\pi$ like in (\ref{eq1.24}) 
\begin{equation}
\label{eq6.26}
\xymatrix{
\partial S'_u (M(\Gamma) - H)_{\rm I}^{\wedge} - \overset{\circ}{\sigma} (\infty) \ar[rr]_-{{\rm via} \, S'_u (M(\Gamma) - H)_{\rm I}^{\wedge}} &&&\!\!\!\!\!\!\!\!\!\!\!\!\!\!\!\!\!\!\!\!\!\!\!\!\!\!\!\!\!\!\!\!\!\!\!\!\!\!\!\!\partial S'_u (M(\Gamma) - H)_{\rm II}^{\wedge}  \\ 
{\mathcal N}^3 \ar[u]^i \ar[rr]^-r &&\sum (\infty)_*^{\wedge} \ \mbox{(6.8.1)} \\
({\mathcal N}^3)' \ar[u] \ar[urr]^-{r} \ar[rr] &&\Theta^3 (fX^2 - H)'_{\rm I} / \Gamma \, .
}
\end{equation}
Here $\sigma(\infty)$ is like in (\ref{eq6.4}) and $\overset{\circ}{\sigma} (\infty) \equiv \sigma (\infty) - \partial \, \sigma (\infty)$. Also,
$$
\partial S'_u (M(\Gamma) - H)_{\rm I}^{\wedge} - \overset{\circ}{\sigma} (\infty) =\partial S'_u (M(\Gamma) - H)_{\rm I} \cup \partial \, \sum (\infty)_*^{\wedge} \, .
$$
In connection with the lower $r$, we also have the following commutative diagram, both of the lower arrows of which are inclusions
\begin{equation}
\xymatrix{
({\mathcal N}^3)' \ar[rr]^-r  &&\sum (\infty)_*^{\wedge} \\
&\overset{\circ}{\sum} (\infty)_* \cup ({\rm fins} - {\rm rims}) = {\rm int} \left(\overset{\circ}{\sum} (\infty)_* \cup {\rm fins}\right)  \, . \ar[ur] \ar[ul]
} \nonumber
\end{equation}

Here, for $\overset{\circ}{\sum} (\infty)_* \cup ({\rm fins} - {\rm rims})$ all the $p_{\infty\infty}$'s are deleted. Then, in connection with the formula (\ref{eq6.26}), the $\sum (\infty)_*^{\wedge}$ lifts to an embedding into both $(\partial S'_u)_{\rm I}^{\wedge}$ and $(\partial S'_u)_{\rm II}^{\wedge}$. In this context, we have the inclusions
$$
\Theta^3 (fX^2 - H)'_{\rm I} / \Gamma \subset \partial S'_u (M(\Gamma)-H)_{\rm I} \subset \partial S'_u (M(\Gamma) - H)_{\rm II} \, ,
$$
compatible with the rest of (\ref{eq6.26}).

\bigskip

\noindent {\bf Lemma 6.4.} {\it In the context of the diagram {\rm (\ref{eq6.26})}, the following map injects}
\begin{equation}
\label{eq6.27}
\pi_1 \, {\mathcal N}^3 ({\rm cut}) \overset{i_*}{-\!\!\!-\!\!\!\longrightarrow} \pi_1 [\partial S'_u (M(\Gamma) - H)_{\rm I}^{\wedge} - \overset{\circ}{\sigma} (\infty)] \, .
\end{equation}

\bigskip

\noindent {\bf Proof.} We will start with an extension of lemma~6.1. We will make use now of the ${\rm int} \, \sum (\infty)_*$ from (6.6), where all the contribution $p_{\infty\infty} ({\rm all})$ is present. Here, also
$$
\overset{\circ}{\sum} (\infty)_* \, (p_{\infty\infty} ({\rm all}) \, {\rm deleted}) \varsubsetneqq \overset{\circ}{\sum} (\infty) (6.6) = {\rm int} \left\{ \sum (\infty) (2.13.1) , \ p_{\infty\infty} (S) \, {\rm deleted} \right\} \, .
$$

Also, a ``$/\Gamma$'' should be everywhere understood, even if not written down explicitly. We also introduce the following non-locally finite space
$$
[\Theta^3 (fX^2 - H({\rm normal}))] \equiv \{{\rm the} \ \Theta^3 (fX^2-H ({\rm normal})), \, \mbox{{\ibf with} the contribution $p_{\infty\infty} (S)$ restored}\}.
\eqno (6.27.1)
$$

Remember that the $p_{\infty\infty} ({\rm proper})$'s had never been removed from $\Theta^3 (fX^2)$. So, all the $p_{\infty\infty}$'s are now back, in (6.27.1). I claim next, that using again Van~Kampen like in lemma~6.1, one can prove that
\begin{equation}
\label{eq6.28}
\pi_1 \, ({\rm int} \, \sum (\infty)_*) ({\rm cut})) \overset{[f\alpha]_*}{-\!\!\!-\!\!\!-\!\!\!\longrightarrow} \pi_1 [\Theta^3 (fX^2 - H ({\rm normal})] , \ {\rm injects} \, .
\end{equation}
We then go to the following commutative diagram of inclusions
\begin{equation}
\label{eq6.29}
\xymatrix{
\left( {\rm int} \, \sum (\infty)_* \right)({\rm cut}) \ar[rr]^-{\gamma} \ar[d]^-{[f\alpha]}  &&{\mathcal N}^3 ({\rm cut}) \ar[d]^-i \\
[\Theta^3 (fX^2 - H({\rm normal}))] \ar[rr]_-{\beta} &&\partial S'_u (M(\Gamma)-H)_{\rm I}^{\wedge} - \overset{\circ}{\sigma} (\infty) \, .
} 
\end{equation}
Here $\gamma$ is a homotopy equivalence, according to (\ref{eq6.28}) the $[f\alpha]$ $\pi_1$-injects and so, it suffices to show that $\beta_*$ is injective in $\pi_1$ too. So, let us start with the following homotopy equivalence
$$
\Theta^3 (fX^2 - H)'_{\rm I} = \Theta^3 (fX^2-H)' \underset{{\rm inclusion} \ i_0}{-\!\!\!-\!\!\!-\!\!\!-\!\!\!-\!\!\!-\!\!\!-\!\!\!\longrightarrow} S'_u (M(\Gamma)-H)_{\rm I} = \Theta^4 (\Theta^3 (fX^2-H)',R) \times B^{N+1} \, .
\eqno (6.29.1)
$$

The $\Theta^4$ is here a smooth non-compact 4-manifold and we obviously have
\begin{equation}
\label{eq6.30}
\partial S'_u (M(\Gamma)-H)_{\rm I} = \Theta^4 \times S^{N-1} \underset{\overbrace{\mbox{\footnotesize$\partial \, \Theta^4 \times S^{N-1}$}}}{\cup} \partial \, \Theta^4 \times B^N \, ,
\end{equation}
from which we can perceive the arrow $\lambda$ below, lifting to $i_0$ from (6.29.1) and which is a $\pi_1$-isomorphism, since $N \gg 1$ and $\pi_1 \, S^{N-1} = 0$:
$$
\Theta^3 (fX^2 - H)' \ \overset{\lambda}{\!\!\!-\!\!\!-\!\!\!\longrightarrow} \partial S'_u (M(\Gamma)-H)_{\rm I}  \, .
\eqno (6.30.1)
$$

The arrow $\lambda$ enters now into the following commutative diagram too
$$
\xymatrix{
\Theta^3 (fX^2-H)' \ar[rr] \ar[d]_-{\lambda}  &&[\Theta^3 (fX^2 - H({\rm normal}))] \ar[d]^-\beta \\
\partial S'_u (M(\Gamma)-H)_{\rm I} \ar[rr] &&{\mbox{$\partial S'_u (M(\Gamma)-H)_{\rm I} \cup \partial \, \sum (\infty)^{\wedge} \, \cup$} \atop \mbox{$\cup \, \underset{p_{\infty\infty} ({\rm all})}{\overset{ \ }{\sum}} p_{\infty\infty} ({\rm all}) \times [-\varepsilon,\varepsilon] = \partial S'_u (M(\Gamma)-H)_{\rm I}^{\wedge} - \overset{\circ}{\sigma} (\infty)$} \, .}
} \eqno (6.30.2)
$$

In this last diagram, when moving from $\lambda$ to $\beta$, the relevant homotopical fact is that, at the level of both of the horizontal lines, the full $p_{\infty\infty} ({\rm all})$ contribution is being restored. It follows from these things that, via $\beta_*$, the $\pi_1 [\Theta^3 (fX^2 - H({\rm normal}))]$ completely catches the $\pi_1 (\partial S'_u (M(\Gamma)-H)_{\rm I}^{\wedge} - \overset{\circ}{\sigma} (\infty))$. \hfill $\Box$

\bigskip

[{\bf Remark.} Notice that $i_0$ (6.29.1) is a homotopy equivalence.]

\bigskip

We go back to the $\Lambda_n = \Lambda (H_n^-)$, $n=1,2,\ldots$ from the lemma~4.8 considered now inside $\partial S'_u (M(\Gamma)-H)_{\rm I} \subset S'_u (M(\Gamma)-H)_{\rm I} \subset \partial S'_u (M(\Gamma)-H)_{\rm II}$ where they are {\ibf null-homotopic}, because of (4.50.1).

\bigskip

\noindent {\bf The Main Lemma 6.5.} {\it For each $\Lambda_n$ there exists a singular disk $D_n^2 \subset \partial S'_u (M(\Gamma)-H)_{\rm II}$, cobounding $\Lambda_n$, such that}
$$
\lim_{n=\infty} \ D_n^2 = \infty \, .
$$

\bigskip

\noindent {\bf Proof.} To begin with, for further use, here is a more formal restatement of our lemma:

\medskip

\noindent (P$_1$) \quad $\forall \, \Lambda_n \ \exists$ a cobounding disk $D_n^2$ s.t. $\forall \, \{ K \, {\rm compact}\} \subset \partial S'_u (M(\Gamma)-H)_{\rm II}$ we have $\# \, \{ n \ {\rm s.t.} \ D_n^2 \cap K \ne \emptyset \} < \infty$.

\medskip

The next statement is the negation of (P$_1$):

\medskip

\noindent (NON P$_1$) \quad $\exists \, K \subset \partial S'_u (M(\Gamma)-H)_{\rm II}$ s.t. for every system of $D_n^2$'s cobounding the $\Lambda_n$'s, we have $\# \, \{ n \ {\rm s.t.} \ D_n^2 \cap K \ne \emptyset \} = \infty$.

\medskip

From now on, $K$ denotes a generic, arbitrary compact subset of $\partial S'_u (M(\Gamma)-H)_{\rm II}$.

\bigskip

\noindent {\bf Sublemma 6.5.1.} {\it The statement {\rm (NON P$_1$)} implies the following}

\medskip

\noindent (P$_2$) \quad {\it $\exists \, K$ and $\exists$ an infinite subsequence $\Lambda_{h_1} , \Lambda_{h_2} , \ldots \subset \{\Lambda_n\}$, s.t. $\forall$ cobounding systems $D_{h_1}^2 , D_{h_2}^2 , \ldots$, we have $D_{h_i}^2 \cap K \ne \emptyset$, $\forall \, i$.}

\bigskip

At the risk of being pedantic, we will prove this implication here. Define $\sum \subset \{ \Lambda_1 , \Lambda_2 , \Lambda_3 , \ldots\}$ by

\medskip

\noindent $(*)$ \quad $\Lambda_i \in \sum \Longleftrightarrow \exists$ a cobounding disk $D_i^2$ for $\Lambda_i$, s.t. $D_i^2 \cap K = \emptyset$, where $K$ is the compact space which has occured in (NON P$_1$).

\medskip

Next, let $\overline\sum \equiv \{ \Lambda_1 , \Lambda_2 , \ldots \} - \sum$, which means that $\Lambda_{h_i} \in \overline\sum \Longleftrightarrow \forall$ cobounding $D_{h_i}^2$ of $\Lambda_{h_i}$ are such that $D_{h_i}^2 \cap K \ne \emptyset$.

\smallskip

We have here two possible cases, either $\# \, \overline\sum < \infty$ or $\# \, \overline\sum = \infty$. Start by assuming that $\# \, \overline\sum < \infty$ and look at the disjoined partition 
$$
\{ \Lambda_1 , \Lambda_2 , \ldots \} = \sum + \overline\sum \, ,
$$
for which we choose the following system of cobounding $D_n^2$'s. If $\Lambda_i \in \sum$, then choose $D_i^2$ like in $(*)$ above and then choose any cobounding disk for the $\Lambda_i \in \overline \sum$. This system is in contradiction with (NON P$_1$) hence the implication
$$
\mbox{(NON P$_1$)} \Longrightarrow \left(\# \, \overline\sum = \infty\right) \, .
$$
So start now with the {\ibf infinite} system
$$
\overline\sum = \{ \Lambda_{h_1} , \Lambda_{h_2} , \ldots \}
$$
and the fact that $\overline\sum$ is infinite, clearly implies our (P$_2$). This proves the desired implication
$$
\mbox{(NON P$_1$)} \Longrightarrow (\mbox{P$_2$}) \, .
$$
\hfill $\Box$

\bigskip

So, we will assume now the (P$_2$) and from this we will deduce, eventually, an absurd conclusion. This will prove then our desired (P$_1$). The subsequence which is provided by (P$_2$), and to which all our attention will be devoted, from now on will be denoted by
$$
\{ \Lambda_{h_1} , \Lambda_{h_2} , \ldots \} \subset \{ \Lambda_1 , \Lambda_2 , \ldots \} \subset \partial S'_u (M(\Gamma)-H)_{\rm I} \, .
\eqno (6.31.1)
$$
We can apply the lemma~6.3, which poduces then for us sub-sub-sequence
$$
\{ \Lambda_{j_1} , \Lambda_{j_2} , \ldots \} \subset \{ \Lambda_{h_1} , \Lambda_{h_2} , \ldots \}
\eqno (6.31.2)
$$
which comes with
$$
\lim_{n=\infty} \, \Lambda_{j_n} = \Lambda_{\infty} \subset S'_u (M(\Gamma) - H)_{\rm I}^{\wedge} \subset \partial S'_u (M(\Gamma) - H)_{\rm II}^{\wedge} \, ,
$$ 
uniform convergence in $S'_u (M(\Gamma) - H)_{\rm II}^{\wedge}$, and hence in $\partial S'_u (M(\Gamma) - H)_{\rm II}$ too. Moreover, our (6.31.2) combined to (P$_2$) implies that any system $D_{j_1}^2 , D_{j_2}^2 , \ldots$ cobounding the $\{ \Lambda_{j_n} \}$ is such that $\forall \, D_{j_n}^2 \cap K \ne \emptyset$. What we will eventually show, and see here the sublemma~6.5.7 below, is that for $\Lambda_{j_n}$ high enough, there is a singular disk $\bar D_{j_n}^2$ which cobounds it inside $\partial S'_u (M(\Gamma)-H)_{\rm II}^{\wedge}$, and which avoids the $K$.

\smallskip

This will disprove the (P$_2$) and hence it will prove (P$_1$). The various discussions which will follow now will take place in the context of the diagram below
$$
\xymatrix{
&{\mbox{$\partial S'_u (M(\Gamma)-H)_{\rm I} \underset{ \ }{\subset} S'_u (M(\Gamma)-H)_{\rm I}$} \atop \mbox{$= \Theta^4 (\Theta^3 (\pi fX^2-H)',R) \times B^N \subset \partial S'_u (M(\Gamma)-H)_{\rm II}$}} \ar[d]^-{p_1} \\
\Theta^3 (\pi fX^2-H)' \ar[r]_-i &\Theta^4 (\Theta^3 (\pi fX^2-H)',R) \ar[d]^-{p_2} \\
&\Theta^3 (\pi fX^2-H)' \, .
}
\eqno{(6.31.3)}
$$

\noindent {\bf Sublemma 6.5.2.} {\it There is a sub-sub-sequence of the $\Lambda_{j_1} , \Lambda_{j_2} , \ldots$, in {\rm (6.31.2)}, which we will denote again $\Lambda_{j_1} , \Lambda_{j_2} , \ldots$, and for which we can find a sequence of positive numbers
$$
\widetilde\varepsilon_1 > \widetilde\varepsilon_2 > \ldots \, ,
$$
converging to zero, such that, for all $n$'s we have}
$$
\widetilde\varepsilon_{n+1} < d(\Lambda_{j_n} , \partial_{\infty}) \leq d (\Lambda_{j_n} , \Lambda_{\infty}) < \widetilde\varepsilon_n \, .
\eqno (6.32.1)
$$

The distances are computed here in the metric of $S'_u (M(\Gamma)-H)_{\rm II}^{\wedge}$ and, since we will have to be quite specific about it, this otherwise quite straightforward lemma will be proved only later on.

\smallskip

In the meanwhile, with an $\varepsilon_0$ like in (6.24.1), we consider a sequence of $C^{\infty}$ functions defined on $\sum (\infty)_*^{\wedge}$, and converging uniformly to zero
$$
\varepsilon_0 > \varepsilon_1 > \varepsilon_2 > \ldots > 0 \, .
\eqno (6.32.2)
$$
Once we have this sequence, then like in (6.25.1) we can define ``neighbourhoods''
\setcounter{equation}{32}
\begin{equation}
\label{eq6.33}
({\mathcal N}_n^3)' \equiv \sum (\infty)_*^{\wedge} \times (-\varepsilon_n , \varepsilon_n) - \partial \sum (\infty)_*^{\wedge} \times (\varepsilon_n = 0) \, ,
\end{equation}
so that we have $({\mathcal N}_0^3)' = \{$the $({\mathcal N}^3)'$ from (6.25.1)$\}$, and also
$$
({\mathcal N}_0^3)' \supset ({\mathcal N}_1^3)' \supset ({\mathcal N}_2^3)' \supset \ldots \, .
$$

By analogy with the retraction $r$ from (\ref{eq6.26}) we also have now an infinite sequence of such retractions, starting with $r_0 = r$,
$$
({\mathcal N}_n^3)' \overset{r_n}{-\!\!\!-\!\!\!-\!\!\!\longrightarrow} \sum (\infty)_*^{\wedge} \, .
\eqno (6.33.1)
$$
$$
\includegraphics[width=12cm]{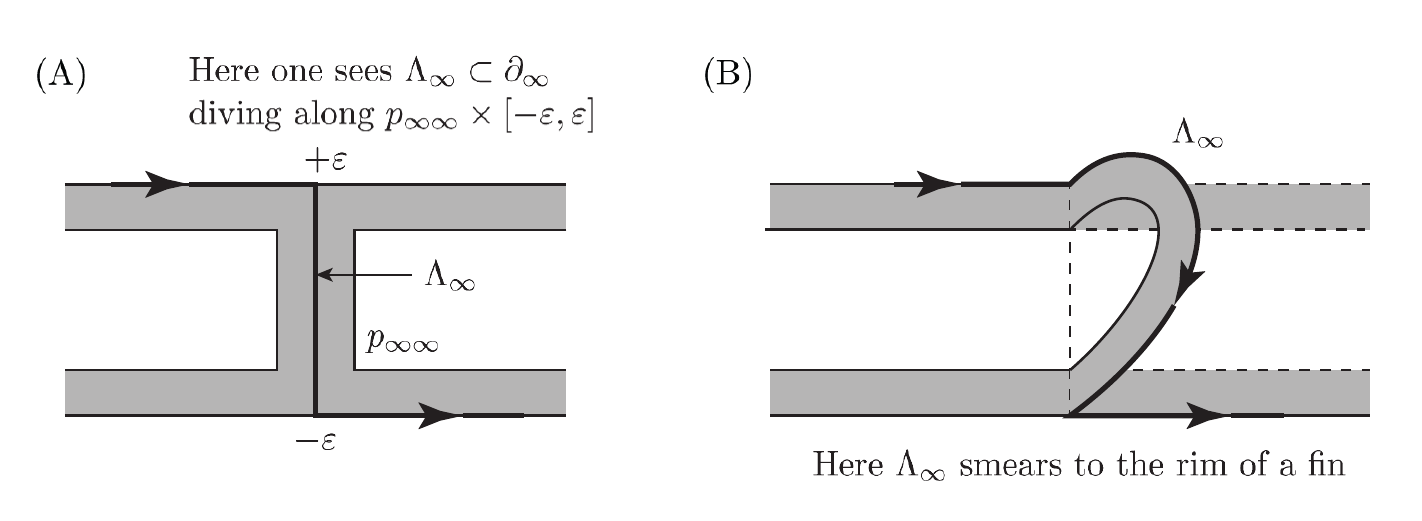}
$$
\label{fig6.4}
\vglue -8mm
\centerline {\bf Figure 6.4.} 
\begin{quote} 
We see here a typical detail of $\sum (\infty)_*^{\wedge}$ containing, in fat lines, some possible pieces of $\Lambda_{\infty}$.
\end{quote}

\bigskip

Also, for each of these objects we find that
$$
{\rm int} \left( \sum^0 (\infty)_* \cup {\rm fins} \right) = ({\mathcal N}_n^3)' \mid (\varepsilon_0 = 0) \, .
$$

In a different vein now, we can get a fairly clear picture of what $\partial_{\infty}$ looks like, by contemplating the figures~6.1, 6.1.bis, 6.1.ter, where its worst complicated details are suggested. We have the following decomposition of $\partial_{\infty}$ into connected components
$$
\partial_{\infty} = \partial_{\infty}^1 + \partial_{\infty}^2 + \ldots
$$
where $\partial_{\infty}^1 \supset \partial \sum (\infty)^{\wedge} \cup \{$arcs like $[x_{\infty} , y_{\infty}]$ which rest on $\partial \, \sum(\infty)^{\wedge}$, see figure~6.1$\}$, and we may well call this the main component, $\partial_{\infty}^2 \supset \{$the dotted trepod at $p_{\infty\infty} (\infty)$, in the left corner of the figure~6.1$\} \subset \sigma (\infty)$, a.s.o.

\smallskip

The figure~6.4 displays another kind of pieces of $\partial_{\infty}$, which are less obvious in the figure~6.1.

\smallskip

A belated COMMENT on the lemma~5.3.

\smallskip

The zipping flow, as described by lemma~5.3 travels, largely along the security walls $W_{\infty} ({\rm BLACK})_{H^{\varepsilon}}$. It is important here that, contrary to what had happened in \cite{31}, we have now $p_{\infty\infty} (W_{\infty}) \in {\rm int} \, W_{\infty}$. Most importantly, also, $W_{\infty} ({\rm BLACK})_{H^0}$ bites big into the infinite structure which comes with $\partial H_{\rm I}^{\wedge} (\gamma_n) \cap H^0$; see figure~1.2. \hfill $\Box$

\medskip

We finally move for good now, from context I to context II. We have here
$$
\left\{ {\rm int} \, \left(\Sigma (\infty)_* \cup {\rm fins}\right) \cup \sum_{R_0} {\rm int} \, R_0 \times [0,\infty) \ \mbox{which was added to} \ S'_u (M(\Gamma)-H)_{\rm I} \ \mbox{during (6.21.1)} \right\}
$$
\begin{equation}
\label{eq6.34}
\subset \partial S'_u (M(\Gamma)-H)_{\rm II} \supset \{ \mbox{the $K$ from (P$_2$)} \} \, .
\end{equation}

Our present problem is that, generically speaking $K \cap \overset{\circ}{\sum}  \, (\infty)_* \ne \emptyset$, which is bad, but then also, as a redeeming factor, there certainly is an $\eta > 0$ such that
\begin{equation}
\label{eq6.35}
0 < \eta < d(K,\partial_{\infty}) \ {\rm in} \ \partial S'_u (M(\Gamma)-H)_{\rm II}^{\wedge}  \, .
\end{equation}
In the context of (6.32.2) and (\ref{eq6.33}), it is only those $\varepsilon_n$'s which come with $2 \, \varepsilon_n < \eta$, which will be considered from now on.

\smallskip

The next lemma is a variation on the same themes as in the lemma~6.4, and hence some repetitions will be unavoidable.

\bigskip

\noindent {\bf Sublemma 6.5.3.} {\it The following inclusion is $\pi_1$-injective}
\begin{equation}
\label{eq6.36}
\underset{{\rm this \ is} \ \underset{R_0}{\bigcup} \, R_0 \times \{\infty\}}{\underbrace{\sum (\infty)_* \times \{\infty\}}} \ ({\rm cut}) \longrightarrow \partial S'_u (M(\Gamma)-H)_{\rm II}^{\wedge} - \overset{\circ}{\sigma} (\infty) \, .
\end{equation}
To begin with, we have the following commutative diagram of inclusion maps among smooth manifolds of high dimensions $\geq N+3$
$$
\xymatrix{
\partial S'_u (M(\Gamma)-H)_{\rm I} \ar[d] \ar[rr]_-a &&\partial S'_u (M(\Gamma)-H)_{\rm II} \ar[d] \\
S'_u (M(\Gamma)-H)_{\rm I} \ar[rr] &&S'_u (M(\Gamma)-H)_{\rm II} \, .
}
$$
The indices of the handles involved here are much smaller than $N+3$, and hence the vertical arrows are $\pi_1$-injective. The lower horizontal arrow, as described by (6.21.1) to (6.21.3) is clearly a homotopy equivalence. Here the upper horizontal arrow is $\pi_1$-injective; it also enters as the first vertical arrow in the next commutative diagram
$$
\xymatrix{
\partial S'_u (M(\Gamma)-H)_{\rm I} \ar[d]_-a \ar[rr]_-{c_{\rm I}} &&\partial S'_u (M(\Gamma)-H)_{\rm I}^{\wedge} - \overset{\circ}{\sigma} (\infty) \ar[d]^-b &&({\rm int} \, \sum(\infty)_*)({\rm cut}) \ar[ll]^-{i \, \circ \, \gamma} \ar[d]^-{\pi_1\mbox{\footnotesize -isomorphism}} \\
\partial S'_u (M(\Gamma)-H)_{\rm II} \ar[rr]_{c_{\rm II}} &&S'_u (M(\Gamma)-H)_{\rm II}^{\wedge} - \overset{\circ}{\sigma} (\infty) && \underset{R_0}{\bigcup} R_0 \times \{\infty\} ({\rm cut}) = \sum (\infty)_* ({\rm cut}) \ar[ll] \, ,
} \eqno (6.36.1)
$$
where $\gamma$ is like in (\ref{eq6.29}) and $i$ like in (\ref{eq6.26}). Also, we just saw that the $a$ is $\pi_1$-injective, while $i \circ \gamma$ has been shown to be so in the context of lemma~6.4. The $c_{\rm I}$ is the addition of the $\partial  \sum (\infty)_*^{\wedge}$ on top of the already existing ${\rm int} (\sum (\infty)_* \cup {\rm fins})$, while the $c_{\rm II}$ is the addition of the $\partial  \sum (\infty)_*^{\wedge} \cup \underset{R_0}{\bigcup} R_0 \times \{\infty\}$ on top of the already existing $\left\{ {\rm int} \left(\overset{\circ}{\sum} \, (\infty)_* \cup {\rm fins}\right) \cup \underset{R_0}{\bigcup} \, {\rm int} \, R_0 \times [0,\infty) \right\}$  (see the (\ref{eq6.34})). It follows that we have the implication
$$
a \ {\rm is} \ \pi_1\mbox{-injective} \Longrightarrow b \ {\rm is} \ \pi_1\mbox{-injective}.
$$

The second arrow on the last line of our diagram (6.36.1) is the one from (\ref{eq6.36}) and it follows from everything said that it is $\pi_1$-injective. \hfill $\Box$

\bigskip
We introduce now the following extension of the ${\mathcal N}^3$ from (\ref{eq6.25})
$$
{\mathcal X}^3 \equiv \underbrace{{\mathcal N}^3 \underset{\overbrace{\mbox{\footnotesize $\sum (\infty)_*^{\wedge}$}}}{\cup} \bigcup R_0 \times [0,\infty]}_{\mbox{\footnotesize this is actually equal to}} \supset \bigcup_{R_0} R_0 \times \{\infty\} \supset \partial \sum(\infty)_*^{\wedge} \, .
$$
\vglue -10mm
\begin{equation}
\label{eq6.37}
\qquad \qquad \qquad \qquad \qquad \qquad \qquad \qquad \qquad \qquad \qquad \qquad {\mathcal N}^3 \ \underset{\overbrace{\mbox{\footnotesize${\rm int} (\sum (\infty)_* \cup {\rm fins})$}}}{\cup} \ \underset{R_0}{\sum} \, {\rm int} \, R_0 \times [0,\infty]
\end{equation}

Of course, when one writes here $R_0 \times [0,\infty]$ we mean with every $x \times [0,1]$, for $x \in \partial R_0$, crushed to a point. With this, when ${\mathcal N}^3$ is replaced by the larger ${\mathcal X}^3$, there is an obvious extension of (\ref{eq6.26}), but only the following piece of it is written down explicitly here
\begin{equation}
\label{eq6.38}
\xymatrix{
&\partial S'_u (M(\Gamma)-H)_{\rm II}^{\wedge} \\
{\mathcal X}^3 \ar[ur]_-j \ar[rr]^-R &&\underset{R_0}{\bigcup} \, R_0 \times \{\infty\} \underset{{{\rm this \, is \, the \, obvious} \atop {\rm projection}} \atop {\rm (see \, here \, the \, (6.18.1) \, too)}}{-\!\!\!-\!\!\!-\!\!\!-\!\!\!-\!\!\!-\!\!\!-\!\!\!-\!\!\!-\!\!\!-\!\!\!-\!\!\!-\!\!\!\twoheadrightarrow} \sum (\infty)_*^{\wedge} \ar[ul] \\
\partial \sum (\infty)_*^{\wedge} \ar[u] \ar[rr]^-{\rm id} &&\partial \sum (\infty)_*^{\wedge} \ar[u]
}
\end{equation}
where all the vertical arrows are canonical injections and where $R$ is the following big RETRACTION. To begin with, proceeding on the same lines as in (6.21.4), we construct a retraction
\begin{equation}
\label{eq6.39}
\sum (\infty)_*^{\wedge} \cup \bigcup \, (R_0 \times [0,\infty]) \overset{R_{\infty}}{-\!\!\!-\!\!\!-\!\!\!-\!\!\!\longrightarrow} \, \bigcup_{R_0} \, R_0 \{\infty\} = \sum (\infty)_* \, ;
\end{equation}
this is suggested in the figure~6.5, and it has the virtue of combining continuously and seamlessly all the individual retractions
$$
R_0 \times [0,\infty] \longrightarrow R_0 \times \{\infty\} \, , \quad \forall \, R_0 \, ,
$$
into a single continuous map. Then, the big retraction $R$ appearing in (\ref{eq6.38}), is now the composition of the retraction $r$ from (\ref{eq6.26}), i.e. the ${\mathcal N}^3 \overset{r}{\longrightarrow} \overset{\wedge}{\sum} (\infty)_* \subset \{$source of $R_{\infty}\}$, with $R_{\infty}$.

\smallskip

From the combination of lemma~6.4 with the $\pi_1$-injectivity of the map $b$ in diagram (6.36.1), follows the following fact. The map
\begin{equation}
\label{eq6.40}
{\mathcal X}^3 ({\rm cut}) \overset{j}{-\!\!\!-\!\!\!-\!\!\!\longrightarrow} \partial S'_u (M(\Gamma)-H)_{\rm II}^{\wedge} - \overset{\circ}{\sigma} (\infty)
\end{equation}
is $\pi_1$-injective.

$$
\includegraphics[width=12cm]{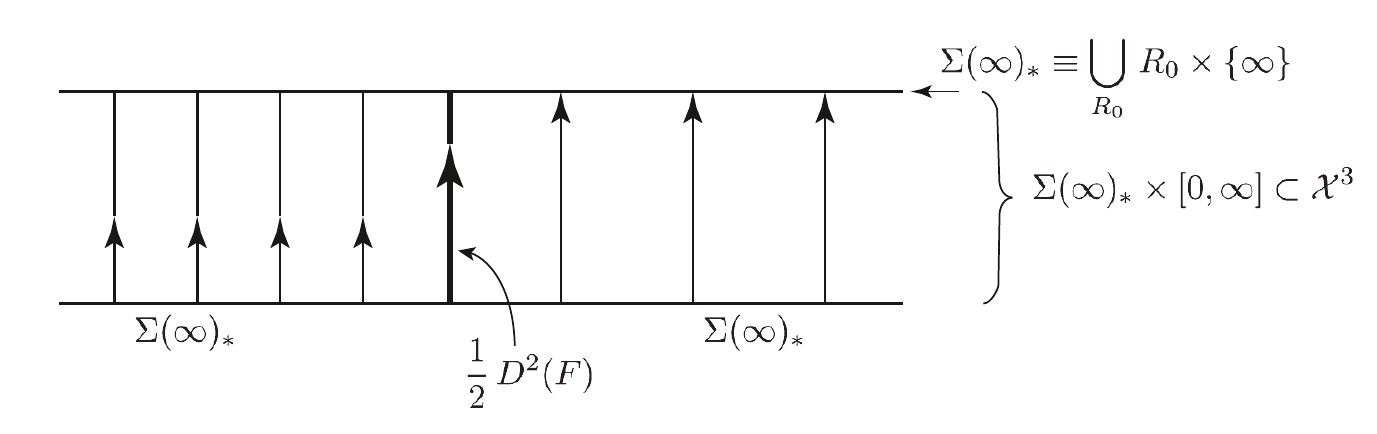}
$$
\label{fig6.5}
\centerline {\bf Figure 6.5.} 

\smallskip

\begin{quote} 
The retraction $R_{\infty}$ from (\ref{eq6.39}), the vertical arrows are its fibers. This figure should illustrate the (6.21.4). The big retraction $R$ in (\ref{eq6.38}) is $R = R_{\infty} \circ r$.
\end{quote}

\bigskip

For the convenience of the reader, we rewrite here the diagram (6.31.3), with some embellishments added, and this will be used in our next discussions

\bigskip

$$
\{\Lambda_n\} \subset \partial S'_u (M(\Gamma)-H)_{\rm I} \subset S'_u (M(\Gamma)-H)_{\rm I} = \underset{\left\downarrow_{{\quad \atop ^{p_1}} \atop \quad } \right.}{\Theta^4 (\Theta^3 (\pi f X^2-H)',R)} \times B^N \subset \underset{\underset{\mbox{\scriptsize $K$ (from NON P$_1$})}{\mbox{\large$\uparrow$}}}{\partial S'_u (M(\Gamma)-H)_{\rm II}}
$$ 
\vglue -4mm
\begin{equation}
\label{eq6.41}
\xymatrix{
\Theta^3 (\pi fX^2-H)' \ar[rr]^-i \ar@{<->}[dr]^-{\rm id} &&\Theta^4 (\Theta^3 (\pi f X^2-H)',R) \supset \{ p_1 \Lambda_n \} \ar[d]^-{p_2} \\
&&\!\!\!\!\!\!\!\!\!\!\!\!\!\!\!\!\!\!\!\!\!\!\!\!\!\!\!\!\Theta^3 (\pi f X^2 - H)' \supset \{ p_2 \circ p_1 \Lambda_n \} \, .
}
\end{equation}

Here, the lower commutative diagram stems from the fact that (as explained in detail in \cite{8}, \cite{18}) the $\Theta^4 (\Theta^3 , R)$ {\ibf is} a regular neighbourhood of $\Theta^3$. The $p_1$ is the obvious projection while $p_2$ is the retraction coming with $\Theta^4 (\ldots , R)$, when viewed as a regular neighbourhood. Also $(\Theta^3)'$ is here $(\Theta_{\rm I}^3)'$. 

\smallskip

As far as {\ibf precise} METRIC STRUCTURES are concerned, what we do have so far, is the ``riemannian'' metric $d$ on $M(\Gamma)$, which then lifts canonically to $\widetilde M (\Gamma)$, $Y(\infty)$, $X^2$. With this, the $\Theta^4$, $S_{\varepsilon}^{(')}$'s are clearly metrizable, and that was enough for everything  we did until now. At this point we stick for a short while at the level of the $\Theta^3 (fX^2 - H)'_{\rm I} \subset M(\Gamma)$ where in its normal isotopy class, this embedding will be made precise, so that the closure inside $M(\Gamma)$ should be exactly the following
\begin{equation}
\label{eq6.42}
\overline{\Theta^3 (fX^2-H)'} = \Theta^3 (fX^2-H)' \cup \partial_{\infty} \subset M(\Gamma) \, .
\end{equation}

We have a PROPER inclusion
$$
{\rm int} \left( \overset{\circ}{\sum} \, (\infty)_* \cup {\rm fins}\right) \subset \Theta^3 (\pi fX^2-H)' \eqno(6.42.1)
$$
and the only places where our $\Theta^3$ fails to be a smooth 3-manifold are the immortal singularities
\begin{eqnarray}
S(m,n) &= &W (m)({\rm BLACK}) \times [-\varepsilon_m , \varepsilon_m] \cap W^* (n)({\rm BLACK}) \times [-\varepsilon_n , +\varepsilon_n] \nonumber \\
&\cap &S_{\infty}^2 ({\rm BLUE}) \subset \overline S \subset {\rm Sing} \, M(\Gamma) \, , \ \mbox{with} \ (m,n) \in Z \times Z \, . \nonumber
\end{eqnarray}
The square $S(m,n)$ has exactly its four curves living at infinity, in $\partial \sum (\infty) \subset \partial_{\infty}$ and it sees two holes $H(m) , H(n)$, pertaining to $W(m) , W^* (n)$. They live above $S_{\infty}^2 ({\rm BLUE})$ at levels $x_{\infty} + x(m)$, $x_{\infty} + x(n)$, where $x(m) , x(n)$ are positive; the lower boundaries of $H(m) , H(n)$ are at $x_{\infty} + x(m)$, $x_{\infty} + x(n)$.

\smallskip

The $W(m) \times [-\varepsilon_m , \varepsilon_m]$ contains the infinitely many squares $S(m,n)$, $\vert n \vert \to \infty$, which accumulate on the two $(p_{\infty\infty} (S) (W(m)) \times [-\varepsilon_m , \varepsilon_m] \subset \partial_{\infty}$. The only specific metric fixing which our $\Theta^3$ needs, is to take $\underset{n=\infty}{\lim} \, x(m) = 0 = \underset{m=\infty}{\lim} \, x(n)$. With this our metric structure for $\Theta^3$ is fixed and (\ref{eq6.42}) is with us. But before we go to $\Theta^4$ and to $(S'_u)_{\rm I}$ we will continue for while with $\Theta^3 (\pi f X^2 - H)' \subset M(\Gamma)$. We define, for any $\alpha > 0$,
\begin{equation}
\label{eq6.43}
{\mathcal U}_{\alpha} \equiv \left\{\mbox{the $\alpha$-neighbourhood of $\bigcup \, \{$limit walls$\}$, in $M(\Gamma)$} \right\} \, .
\end{equation}
Consider next
\begin{equation}
\label{eq6.44}
\xymatrix{ 
{\rm int} \left( \overset{\circ}{\sum} (\infty)_* \cup {\rm fins} \right) \subset \underset{\equiv \, \overset{\circ}{\sum} (\infty)^{\wedge} \, ({\rm see} \, (6.24.2))}{\underbrace{{\rm int} \,{\mbox{\large$\Sigma$}} (\infty)_* \cup ({\rm fins} - {\rm rims})}} \subset \sum (\infty)_*^{\wedge} \cup \sigma (\infty) \supset \partial_{\infty} \ar@{_{(}->}[d] \\ 
\sum (\infty)_*^{\wedge} &{\mbox{$({\mathcal N}^3)' \supset ({\mathcal N}_n^3)' (6.33),$} \atop \mbox{${\rm with} \, r_n \equiv r \mid ({\mathcal N}_n^3)' \, .$}}  \ar[l]_-{r \, (6.26)} \\
}
\end{equation}
Everything in this diagram (\ref{eq6.44}) lives in $M(\Gamma)$, and in its context let us also define, for any $\beta > 0$, the sets
\begin{equation}
\label{eq6.45}
{\mathcal V}_{\beta} \equiv \left\{\mbox{the $\beta$-neighbourhood of $\partial_{\infty}$ inside ${\rm int} \left(\overset{\circ}{\sum} \, (\infty)_* \cup {\rm fins}\right)$}\right\} \, ,
\end{equation}
and 
\begin{equation}
\label{eq6.46}
U_{\beta} \equiv \{\mbox{the $\beta$-neighbourhood of $\partial_{\infty}$ in $M(\Gamma)$}\} \supset {\mathcal V}_{\beta} = U_{\beta} \cap {\rm int}  \left(\overset{\circ}{\sum} \, (\infty)_* \cup {\rm fins}\right) \, .
\end{equation}
It should be understood here that, with the $\eta > 0$ from the (\ref{eq6.35}), the $U_{\beta}$'s useful for us will be the ones coming with
$$
0 < \beta \ll \eta \, , \eqno (6.46.1)
$$
which will make sure that in the context of (6.31.3), $p_2 \, p_1 K \cap U_{\beta} = \emptyset$. In the rest of this section whenever things like $\beta$ and/or $\beta (\alpha)$ will be used, it should be understood that $\beta$ is always small enough such that (6.46.1) is verified. This will keep us safely away from things like
$$
K \cap \{W({\rm complementary}) \times [-\varepsilon , \varepsilon] - {\rm LIM} \, M_2(f) \times \{\pm \, \varepsilon\}\} \overset{?}{\ne} \emptyset \, ,
$$
which are bound to occur, possibly close to ${\rm LIM} \, M_2(f) \times \{\pm \, \varepsilon\}$ BUT, with our restriction (6.46.1), certainly NOT closer than $\eta$ (\ref{eq6.35}).

\smallskip

Some more notations will be necessary for what will follow next. In the context of (\ref{eq6.24}), we consider the decomposition
\begin{eqnarray}
\label{eq6.47}
\Theta^3 (\pi f X^2 - H)' &\equiv &\left[ \sum \, \mbox{thickened complementary walls} \ W_1 - H(p_{\infty\infty}) -  \{\mbox{the BLACK Holes} \} \right] \\
&\cup &\left[ \left(\sum \, \mbox{thickened non-complementary walls} \ W_2 \right) - \{ \mbox{completely normal Holes}\} \right] \, ; \nonumber
\end{eqnarray}
here $W_1, W_2$ are generic notations and the thickenings should be $W_1 \times [-\varepsilon , \varepsilon]$, $W_2^n \times [-\varepsilon_n , \varepsilon_n]$ where $\varepsilon_n \underset{n \to \infty}{-\!\!\!-\!\!\!-\!\!\!-\!\!\!\longrightarrow} 0$.

\bigskip

\noindent {\bf Sublemma 6.5.4.} {\it The sizes of the Holes can be chosen so that, without contradicting anything said so far, the following should happen too}

\medskip

\noindent (6.48.1) \quad {\it There is a continuous monotonically decreasing function $\alpha \in [0 , \infty) \overset{\beta}{\longrightarrow} [0,\infty)$, such that
$$
\lim_{\alpha = 0} \ \frac{\beta(\alpha)}{\alpha} = 1 \, , \ \mbox{and} \ \beta (\alpha) \leq \alpha \, ,
$$ 
with which come also the following features.

\medskip

We start by considering the following piece of $\Theta^3 (\pi f X^2 - H)'$, and we make use here of the notations from} (\ref{eq6.47})
$$
{\mathcal U}^3 (\alpha) \equiv \left(\bigcup_{W_1} N_1 (\beta(\alpha))_{W_1}\right) \cup \left[ \bigcup_{W_n \cap W_1 \subset (N_1)_{W_1}} ((W_n-H) \times [-\varepsilon_n , \varepsilon_n] ) - W_1 \times (-\varepsilon , \varepsilon)\right] \, ,
$$
{\it where
$$
N_1 (\beta(\alpha))_{W_1} \equiv \{\mbox{The $\beta(\alpha)$-neighbourhood of} \ \partial_{\infty} \cap (W_1 \times [-\varepsilon , \varepsilon]) = {\rm LIM} \, M_2 (f) \cap W_1 \times \{ \pm \, \varepsilon\} \} \, .
$$
Taking $W_n \cap W_1 \subset N_1$ in the formula above means that the $n$'s come with $n \geq n_0 (\alpha)$, where $Z_+ \overset{n_0}{-\!\!\!-\!\!\!\longrightarrow} Z_+$ is a function such that $\underset{\alpha = 0}{\lim} \ n_0 (\alpha) = \infty$. With all this, and without contradicting anything said so far, we can fix the $H \, ({\rm completely \, normal})$ so that we should have}
$$
{\mathcal U}^3 (\alpha) \subset \bigcup_{\beta(\alpha)} \cap \ \Theta^3 (\pi f X^2 - H)' \subset {\mathcal U}_{\alpha} \, . \eqno (6.48.2)
$$

\bigskip

\noindent {\bf Proof.} We let the $H = H (\mbox{completely normal})$ converge very fast towards the Ideal Holes. Then, asymptotically, the $\partial \, C(H)$'s are inside $U_{\beta(\alpha)}$ and the ${\mathcal U}^3 (\alpha)$ is now
$$
{\mathcal U}^3 (\alpha) = \{\mbox{the shaded region in figure 6.6}\} \cap \{\mbox{those regions marked $U_{\beta}$, in the same figure}\} \, .
$$
This takes care of the first inclusion in (6.48.2) and the second one follows from $\beta (\alpha) \leq \alpha$. \hfill $\Box$

\bigskip

The ${\mathcal V}_{\gamma}$'s occurring in figure~6.6 are part of the discussion which follows afterwards. In the context of (\ref{eq6.44}) $+$ (\ref{eq6.45}), we consider now, for
$$
({\mathcal N}_n^3)' \overset{r_n}{-\!\!\!-\!\!\!\longrightarrow} \sum (\infty)_*^{\wedge} \supset {\rm int} \left( \overset{\circ}{\sum} (\infty)_* \cup {\rm fins} \right) \, ,
$$
the $r_n^{-1} \, {\mathcal V}_{\beta} \subset ({\mathcal N}_n^3)'$, where it is understood that the contribution of
$$
\left( \sum_{p_{\infty\infty}({\rm all})} p_{\infty\infty} \times (-\varepsilon , \varepsilon) \right) \times (\varepsilon_n \ne 0) \subset \left( \sum p_{\infty\infty} \times (-\varepsilon , \varepsilon) \right) \times (-\varepsilon_n , \varepsilon_n)
$$
is also included in $({\mathcal N}_n^3)'$.
$$
\includegraphics[width=160mm]{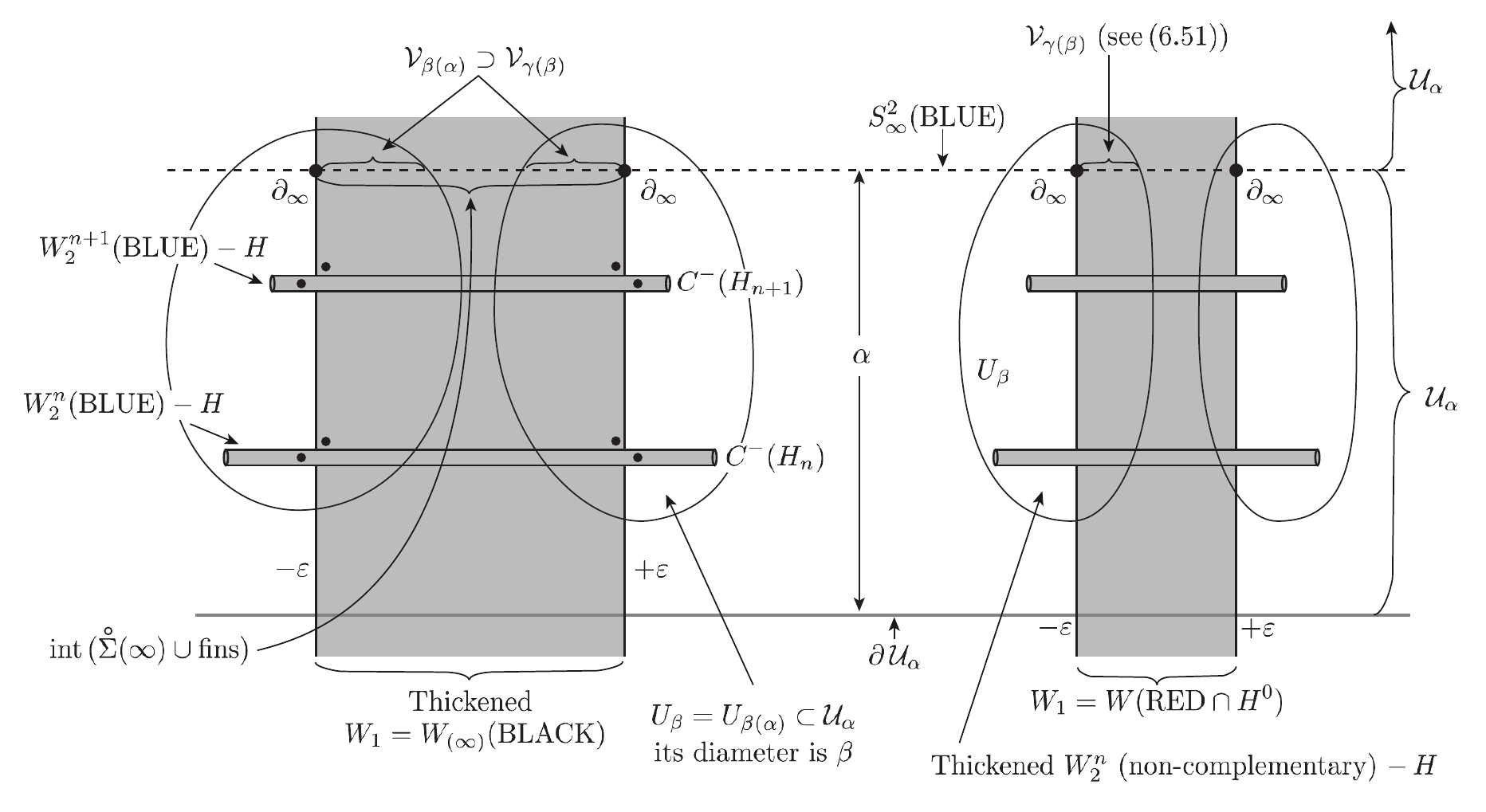}
$$
\label{fig6.6}
\centerline {\bf Figure 6.6.} 

\smallskip

\begin{quote} 
The present figure is supposed to live inside $M(\Gamma)$, and it should serve for vizualizing the various neighbourhoods ${\mathcal U}_{\alpha} \supset U_{\beta(\alpha)} \cap \Theta^3 (\pi f X^2-H)' \supset {\mathcal U}^3 (\alpha)$, $U_{\beta(\alpha)} \cap \Theta^3 (\pi f X^2-H)' \supset {\mathcal V}_{\beta(\alpha)} \subset {\rm int} \left(\overset{\circ}{\sum} \, (\infty)_* \cup {\rm fins}\right)$. The smaller ${\mathcal V}_{\gamma(\beta)} \subset {\mathcal V}_{\beta(\alpha)}$ mentioned in our drawing, is like in 4) from the Sublemma~6.5.5 and the formula (\ref{eq6.51}). The pair of points $(\cdot \cdot)$ drawn at the level of $W_{(\infty)}$ should suggest, very schematically, the trace of the zipping flow. Together with the $C^- (H_n)$'s they are the contribution of $\{ p_2 \, p_1 \Lambda_n \} \subset \Theta^3 (\pi f X^2 - H)'$, at the level of this figure. The piece of the dotted $S_{\infty}^2$ between a pair of red points in $\partial_{\infty}$ is in ${\rm int} \left(\sum \, (\infty)_* \cup {\rm fins}\right)$. Also ${\mathcal V}_{\beta (\alpha)} = U_{\beta} \cap {\rm int} \left(\sum \, (\infty)_* \cup {\rm fins}\right)$. The $\Theta^3 (\pi f X^2 - H)'$ is shaded.
\end{quote}

\bigskip

We fix now an $\alpha > 0$ and notice that, since the sequence $\{ \varepsilon_n \}$ is (6.32) converges uniformly to zero, we can certainly find an $N(\infty)$ high enough, so that we should have
\setcounter{equation}{48}
\begin{equation}
\label{eq6.49}
({\mathcal N}_{N(\alpha)}^3)' \subset {\mathcal U}_{\alpha} \cap \Theta^3 (\pi f X^2 - H)' \quad {\rm and} \quad \lim_{\alpha = 0} \, N(\alpha) = \infty \, ,
\end{equation}
and also
$$
{\rm dist} \left( ({\mathcal N}_{N(\alpha)}^3)' , \Sigma_1 (\infty) \ (1.14) \right) < \beta(\alpha) \, .
$$

\bigskip

\noindent {\bf Sublemma 6.5.5.} 1) {\it There is a homotopical retraction  $\rho_t$, $t \in [0,1]$ of ${\mathcal U}_{\alpha} \cap \Theta^3 (\pi fX^2-H)'$ into itself, such that $\rho_0 =$ identity, $\rho_1 (\Theta^3 (\pi fX^2-H') \cap {\mathcal U}_{\alpha}) \subset ({\mathcal N}^3_{N(\infty)})'$, with the features below.}

\medskip

\noindent 2) {\it In the context of} (6.48.2) {\it we have}
\begin{equation}
\label{eq6.50}
\rho_t (U_{\beta(\alpha)} \cap \Theta^3 (\pi fX^2-H)') \subset U_{\beta(\alpha)} \cap \Theta^3 (\pi fX^2-H)' \, , \quad \forall \, t \, .
\end{equation}

\noindent 3) {\it Points which are very close to $\partial_{\infty}$, when moved by $\rho_t$, stay sufficiently close to $\partial_{\infty}$.}

\medskip

\noindent 4) {\it There exists a continuous function $\beta \in [0,\infty) \overset{\gamma}{\longrightarrow} [0,\infty)$ coming with $\gamma(\beta) < \beta$ and $\underset{\beta = 0}{\lim} \, \gamma (\beta) = 0$, such that
\begin{equation}
\label{eq6.51}
\rho_1 (U_{\beta(\alpha)} \cap \Theta^3 (\pi fX^2-H)') \subset (r_{N(\alpha)}^{-1} \, {\mathcal V}_{\gamma(\beta)}) \cap ({\mathcal N}^3_{N(\infty)})' \, ,
\end{equation}
see here notations from} (\ref{eq6.33}), (6.33.1).

\bigskip

\noindent {\bf Proof.} Starting with (\ref{eq6.47}), we use a horizontal retraction which brings all the term $\sum$ [(thickened non-complementary walls $W_2) \, -$ (the completely normal Holes)] into $\sum$ (thickened complementary walls $W_1$); this part can be guessed from the figure~6.6. Next, we continue to retract vertically, towards $\sum (\infty)_*^{\wedge}$, staying close to the vertical counter-image $(r^{-1} \, \partial \sum (\infty)_*^{\wedge}) \subset {\mathcal N}^3$, see (\ref{eq6.26}). This again should be visualized on the figure~6.6. \hfill $\Box$

\bigskip

We fix now the metric structure on $\Theta^4 = \Theta^4 (\Theta^3 (\pi fX^2-H)' , {\mathcal R})$. Here, as we already know, $\Theta^3 = \Theta^3 (\pi fX^2-H)' \subset M(\Gamma)$, which induces a metric on $\Theta^3$, ${\rm Sing} \, \Theta^3 \subsetneqq {\rm Sing} \, M(\Gamma)$ and, as we also know, the ${\mathcal R}$-dependence of $\Theta^4$ is washed away when we multiply by $(\times \, B^N)$. So, we may as well assume that there is a desingularization ${\mathcal R}$ on $M(\Gamma)$ which induces {\ibf our} ${\mathcal R}$ in $\Theta^3$, so that
\begin{equation}
\label{eq6.52}
\Theta^4 (\Theta^3 , {\mathcal R}) \subset \{ \Theta^4 (M(\Gamma) , {\mathcal R}) , \ \mbox{which is compact}\} \, .
\end{equation}

The $\Theta^4 (M(\Gamma) , {\mathcal R})$ is a regular neighbourhood of $M(\Gamma)$, i.e. it comes with the standard diagram $M(\Gamma)$ $\underset{\rm inclusion}{-\!\!\!-\!\!\!-\!\!\!-\!\!\!-\!\!\!-\!\!\!-\!\!\!-\!\!\!\longrightarrow} \, \Theta^4 (M(\Gamma) , {\mathcal R}) \underset{\rm retraction}{-\!\!\!-\!\!\!-\!\!\!-\!\!\!-\!\!\!-\!\!\!-\!\!\!-\!\!\!\longrightarrow} M(\Gamma)$. We will endow $\Theta^4 (M(\Gamma) , {\mathcal R})$ with a metric structure, call it $d$, compatible with the inclusion of $(M (\Gamma) , \, {\rm its} \, d)$. We can define (see (\ref{eq6.42}))
$$
\Theta^4 (\Theta^3 , {\mathcal R}) \cup \partial_{\infty} = \Theta^4 (\Theta^3 , {\mathcal R}) \underset{\overbrace{\Theta^3}}{\cup} (\Theta^3 \cup \partial_{\infty}) \, , \quad {\rm where} \quad \Theta^3 = \Theta^3 (\pi fX^2-H)' \, ,
$$
and we will want now to fix the precise metric structure of $\Theta^4 (\Theta^3 , {\mathcal R})$ so that the closure inside $\Theta^4 (M(\Gamma) , {\mathcal R})$ should be, exactly,
\begin{equation}
\label{eq6.53}
\overline{\Theta^4 (\Theta^3 , {\mathcal R})} = \Theta^4 (\Theta^3 , {\mathcal R}) \cup \partial_{\infty} \, ,
\end{equation}
a formula which should be compared to (\ref{eq6.42}).

\smallskip

In terms of (\ref{eq6.41}), for each $x \in \Theta^3 (\pi fX^2-H)'$ there is a fiber $p_2^{-1} (x) \subset \Theta^4 (\Theta^3 , {\mathcal R})$, which topologically speaking, is either of the form $x \times [-1,+1]$, when $x$ is a smooth point, or $x \times \{$a finite tree with a single vertex at $x\}$, when $x \in {\rm Sing} \, \Theta^3$. We denote by $\Vert p_2^{-1} (x) \Vert$ the diameter of $p_2^{-1}$, measured in $(\Theta^4 (M(\Gamma) , {\mathcal R}) , d)$. We will impose that
\begin{equation}
\label{eq6.54}
\lim_{n=\infty} d(x_n , \partial_{\infty}) = 0 \Longrightarrow \lim \Vert p_2^{-1} (x_n) \Vert = 0.
\end{equation}
Our fixing of $\Vert p_2^{-1} (x) \Vert$ can also be thought of as a fixing of the embedding (\ref{eq6.52}), inside its allowed isotopy class.

\smallskip

With all this, $(\Theta^4 (M(\Gamma) , d))$ induces a metric, which we call again $d$, on the $\Theta^4 (\Theta^3 , {\mathcal R})$, coming with the correct compactification (\ref{eq6.53}).

\smallskip

We will move now to $S'_u (M(\Gamma) - H)_{\rm I} = \Theta^4 (\Theta^3 (M(\Gamma)-H)' , {\mathcal R}) \times B^N$, and for this object we have been working so far with a generic metric, which was defined only up to quasi-isometry.

\smallskip

We will think of this $S'_u$ as living inside
\begin{equation}
\label{eq6.55}
\Theta^4 (M(\Gamma) , {\mathcal R}) \times R^N \supset \Theta^4 (M(\Gamma) , {\mathcal R}) \times \{ 0 \}  \supset \Theta^4 (\Theta^3 (\pi f X^2-H)' , {\mathcal R}) \cup \partial_{\infty}  \times \{ 0 \} \, ,
\end{equation}
where the $0 \in R^N$ is the origin.

\smallskip

The $(S'_u)_{\rm I}$ needs a factor $B^N \subset R^N$ but we will {\ibf not} take it centered at $0$, but rather so that we should have $0 \in \partial B^N$, with $\{0\} = \{ s \} \equiv \{$the SOUTH POLE of $S^{N-1} = \partial B^N\}$ and $\{ n \} \equiv \{$the NORTH POLE of $S^{N-1}\}$ diametrally opposite to it.

\smallskip

With this, the canonical embedding $\Theta^4 (\Theta^3 (\pi f X^2-H)' , {\mathcal R}) \subset \partial S'_u (M(\Gamma)-H)_{\rm I}$, mentioned already several times, becomes
\begin{eqnarray}
\Theta^4 (\Theta^3 (\pi f X^2-H)' , {\mathcal R}) &= &\Theta^4 (\Theta^3 (\pi f X^2-H)' , {\mathcal R}) \times (\{ 0 \} = \{ s \}) \nonumber \\
&\subset &\Theta^4 (\Theta^3 (\pi f X^2-H)' , {\mathcal R}) \times S^{N-1} \subset \partial S'_u (M(\Gamma)-H)_{\rm I} \, . \nonumber
\end{eqnarray}
Here is now our PROGRAM for what will follow next. We want to fix an explicit metric on $S'_u (M(\Gamma)-H)_{\rm I}$, compatible with the embedding inside the $\Theta^4 (M(\Gamma) , {\mathcal R}) \times R^N$ from (\ref{eq6.55}), where the space in (\ref{eq6.55}) is endowed with its obvious product metric below
\begin{equation}
\label{eq6.56}
\{\mbox{the metric $d$ on $\Theta^4$, already mentioned}\} \times \{\mbox{euclidean metric on $R^N$}\} \, ,
\end{equation}
such that the following requirements should be fulfilled:

\medskip

\noindent (6.56.0) \quad In view of everything already said, what we want to do now is to fix, for every $x \in \Theta^4 (\Theta^3 (\pi f X^2-H)' , {\mathcal R})$, and its corresponding fiber $x \in B^N \subset S'_u (M(\Gamma) - H)_{\rm I}$, the $\Vert x \times B^N \Vert$.

\medskip

Our requests are the following.

\medskip

\noindent (6.56.1) \quad Reinforcing what was said before, in connection with the lemma~6.2, inside (\ref{eq6.55}) we should find that
$$
\overline{S'_u (M(\Gamma)-H)_{\rm I}} = S'_u (M(\Gamma)-H)_{\rm I} \cup \partial_{\infty} = S'_u (M(\Gamma) - H)_{\rm I}^{\wedge} \, .
$$

\medskip

\noindent (6.56.2) \quad Remember that our $\Lambda_{j_n}$'s live inside $\partial S'_u (M(\Gamma) - H)_{\rm I}$ and once the metric on $S'_u (M(\Gamma) - H)_{\rm I}$ is completely and explicitly fixed, any natural extension of it to $S'_u (M(\Gamma) - H)_{\rm II}$ will be OK. The (6.31.1 and 2) is a purely $\partial S'_u (M(\Gamma) - H)_{\rm I}$ affair, and we want that with our by now fixed $d$, not only the (6.31.1 and 2) should be with us but that, moreover, the sub-lemma~6.5.2 should be verified too. So, fulfilling our PROGRAM should provide a complete proof for the lemma in question.

\medskip

\noindent (6.56.3) \quad For the $\{ \Lambda_{j_n}$ from (6.30.2)$\} \subset \partial S'_u (M(\Gamma) - H)_{\rm I}$ and hence for $\{$our $\Lambda_{j_n}$ from the sub-lemma~6.5.2$\}$ too, we should have $\{\Lambda_{j_n} \} \cap [ \Theta^4 (\Theta^3 (\pi f X^2-H)' , {\mathcal R}) \times \{ n \}] = \emptyset$.

\medskip

We will implement now this program.

\smallskip

We start with the $\{\Lambda_{j_n}$ from (6.31.2)$\}$, which does converge uniformly to $\Lambda_{\infty} \subset \partial_{\infty}$, inside $\partial S'_u (M(\Gamma)-H)_{{\rm I} \, {\rm or} \, {\rm II}}$. The metric on $\Theta^4 (\Theta^3 (\pi f X^2-H)' , {\mathcal R})$ is by now already fixed and, whatever our choice of $\Vert x \times B^N \Vert$'s will be the $\{ p_1 \, \Lambda_{j_n}\}$ has to converge uniformly to $\Lambda_{\infty}$.

\smallskip

The $p_1 \, \Lambda_{j_n}$'s are disjoined from $\partial_{\infty}$, and hence
$$
0 < d (p_1 \, \Lambda_{j_n} , \partial_{\infty}) \leq d (p_1 \, \Lambda_{j_n} , \Lambda_{\infty}) \, .
$$

Putting together these things, it is not hard to find a sequence of positive numbers $\varepsilon'_1 > \varepsilon'_2 > \ldots$ converging to zero and also a sub-sub-sequence of the original $\{\Lambda_{j_n}\}$, which we will call again $\Lambda_{j_n}$, such that
\begin{equation}
\label{eq6.57}
\varepsilon'_{n+1} < d (p_1 \, \Lambda_{j_n} , \partial_{\infty}) \leq d (p_1 \, \Lambda_{j_n} , \Lambda_{\infty}) < \varepsilon'_n \, .
\end{equation}

This will be {\ibf the} $\Lambda_{j_n}$ from the sublemma~6.5.2.

\smallskip

Next, we will finally fix the metric structure on $S'_u (M(\Gamma)-H)_{\rm I}$.

\medskip

\noindent (6.58) \quad Modulo some smooth interpolations, we have the following. If $y \in \Theta^4 (\Theta^3 (\pi f X^2-H)' , {\mathcal R})$ is such that
$$
\varepsilon'_{n+1} < d(y, \partial_{\infty}) < \varepsilon'_n \, ,
$$
then we will take $\Vert y \times B^N \Vert = \varepsilon'_n$, and see here the figure~6.7 for an illustration. This takes care of (6.56.0). With this precise metric fixing, since we have that
$$
\lim \Vert y \times B^N \Vert = 0 \, , \quad {\rm when} \quad \lim (y , \partial_{\infty}) = 0 \, ,
$$
the condition (6.56.1), i.e. the correct compatification, is verified.

\medskip

$$
\includegraphics[width=13cm]{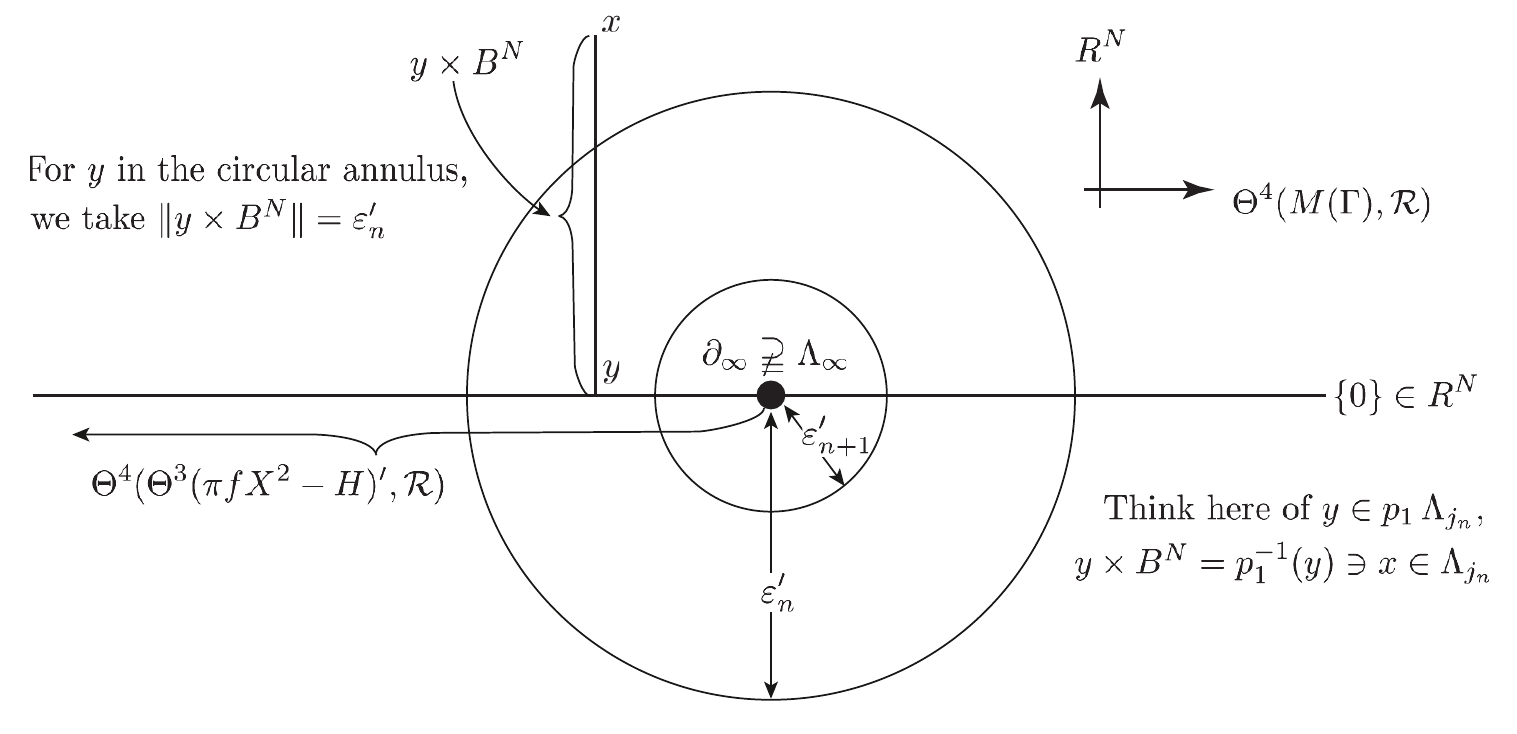}
$$
\label{fig6.7}

\centerline {\bf Figure 6.7.} 

\smallskip

\begin{quote} 
We see here some $\varepsilon'_p$-neighbourhoods of $\partial_{\infty}$ in $\Theta^4 (M(\Gamma) , {\mathcal R}) \times R^N$ and also a generic $y$ s.t. $\varepsilon'_{n+1} < d(y,\partial_{\infty}) < \varepsilon'_n$.
\end{quote}

\bigskip

We define now $\widetilde\varepsilon_n \equiv \sqrt2 \ \varepsilon'_n$, and these will be the actual $\widetilde\varepsilon_n$ which are occurring in the sublemma~6.5.2, i.e. in the (6.32.1).

\smallskip

Let now $x \in \Lambda_{j_n}$, with $y \equiv p_1 (x) \in p_1 \, \Lambda_{j_n}$. By (\ref{eq6.57}) we have that $\varepsilon'_{n+1} < d(y,\partial_{\infty}) < \varepsilon'_n$, which combined with (6.58) implies that $\Vert y \times B^N \Vert = \varepsilon'_n$. We also have $x \in y \times B^N$ and so, combining the second inequality in (\ref{eq6.57}) with (6.58), we have
$$
d(x,\Lambda_{\infty}) \leq \left[(d(y,\Lambda_{\infty})^2) + \Vert y \times B^N \Vert^2 \right]^{1/2} < \sqrt2 \ \varepsilon'_n = \widetilde\varepsilon_n \, .
$$

Next, by (\ref{eq6.57}), as already noticed, we have $d(y,\partial_{\infty}) > \varepsilon'_{n+1}$. Because $\Vert y \times B^N \Vert = \varepsilon'_n > \varepsilon'_{n+1}$, we have then, by Pythagora, that
$$
d(x,\partial_{\infty}) > \sqrt{2(\varepsilon'_{n+1})^2} = \sqrt2 \ \varepsilon'_{n+1} = \widetilde\varepsilon_{n+1} \, .
$$

Putting all these things together, we get finally that
$$
\widetilde\varepsilon_{n+1} < d(x,\partial_{\infty}) \leq d(x,\Lambda_{\infty}) < \widetilde\varepsilon_n \, .
$$

With this the item (6.56.2) in our program has been fulfilled. For the item (6.56.3) one notices that, inside the smooth $(N+3)$-manifold $\partial S'_u (M(\Gamma)-H)_{\rm I}$, the $\Theta^4 (\Theta^3 (\pi f X^2-H)' , {\mathcal R}) \times \{ n \}$ is a smooth submanifold of codimension $N-1 \gg 1$. We think here, of course, in terms of
$$
\partial S'_u (M(\Gamma)-H)_{\rm I} = \Theta^4 (\Theta^3 (\pi f X^2-H)' , {\mathcal R}) \times S^{N-1} \cup \partial \, \Theta^4 \times B^N \, .
$$
We are free to budge the curves $\Lambda_{j_n} \subset \partial S'_u (M(\Gamma)-H)_{\rm I}$ via a small isotopy, and so they can avoid the $\Theta^4 (\Theta^3 (\pi f X^2-H)' , {\mathcal R}) \times \{ n \}$.

\smallskip

This way, our PROGRAM stated above has been completely fulfilled.

\smallskip

Assume now that we are given a sequence of positive numbers $\alpha_n$, such that $\underset{n=\infty}{\lim} \, \alpha_n = 0$. The sublemmas~6.5.4, 6.5.5, and the formula (\ref{eq6.49}) provide us with quantities
\setcounter{equation}{58}
\begin{equation}
\label{eq6.59}
\beta_n = \beta (\alpha_n) \, , \ \gamma_n = \gamma (\beta_n) \, , \ N(\alpha_n) \ {\rm s.t.} \  \lim_{n=\infty} \, \beta_n = \lim_{n=\infty} \, \gamma_n = 0 \ {\rm and} \ \lim_{n=\infty} \, N(\alpha_n) = \infty \, .
\end{equation}

We have then also
\begin{equation}
\label{eq6.60}
d(U_{\beta_k} , \partial_{\infty}) \to 0 \, , \ d(p_2^{-1} \, U_{\beta_k} , \partial_{\infty}) \to 0 \, , \ d(p_1^{-1} \, p_2^{-1} \, U_{\beta_k} , \partial_{\infty}) \to 0 \, ,
\end{equation}
and this among other things, keeps the three objects occurring in formula (\ref{eq6.60}) safely away from $K$, provided that $k$ is high enough.

\smallskip

We go back now to (\ref{eq6.35}) and, for our same $K$, the compact set occurring in $(P_2)$, and coming with
$$
K \subset S'_u (M(\Gamma) - H)_{\rm II} \, ,
$$
we will strengthen (\ref{eq6.35}), with a possibly smaller $\eta > 0$ into the following inequality
\begin{equation}
\label{eq6.61}
d(K , \partial_{\infty}) \geq d \left( K , \bigcup_{R_0} R_0 \times \{\infty\} \right) > \eta > 0 \, .
\end{equation}
For this new smaller $\eta$, we continue with a $\beta$ abiding to (6.46.1). We introduce now the subset of $M(\Gamma)$, see figure~6.6 too,
$$
U_{\beta} \cap \Theta^3 \equiv U_{\beta} \cap \Theta^3 (\pi f X^2-H)' \, .
$$

With $\beta_n$ small enough, we will have
$$
K\cap (p_2 \, p_1)^{-1} (U_{\beta_n} \cap \Theta^3) = \emptyset
$$
and also
$$
K\cap (p_2 \, p_1)^{-1} \, r_{N(\alpha_n)}^{-1} \, {\mathcal V}_{\gamma (\beta_n)} = \emptyset \, .
$$

The next lemma is a direct consequence of the various things said above.

\bigskip

\noindent {\bf Sub-sub-lemma 6.5.6.} 1) {\it In the context of} (\ref{eq6.59}), {\it we can fix the sequence $\alpha_n \to 0$ so that, for $n$ high enough, the following should happen
\begin{equation}
\label{eq6.62}
\Lambda_{j_n} \subset p_1^{-1} \, p_2^{-1} \, (U_{\beta_n} \cap \Theta^3) \, , \ \mbox{i.e. also} \ p_1 \, \Lambda_{j_n} \subset p_2^{-1} (U_{\beta_n} \cap \Theta^3) \, .
\end{equation}
Using here {\rm (6.56.3)} and taking $\partial B^N = \{ n \} = R^{N-1}$, the map $p_1 \mid \Lambda_{j_n}$ is then the obvious canonical projection}
$$
\Theta^4 (\Theta^3 (\pi f X^2-H)' , {\mathcal R}) \times R^{N-1} \longrightarrow \Theta^4 (\Theta^3 (\pi f X^2-H)' , {\mathcal R}) \, .
$$

\medskip

\noindent (6.63) \quad {\it Making use of this last fact, one can find a homotopy $[\Lambda_{j_n} , p_1 \, \Lambda_{j_n}]$ which avoids the compact $K$.}

\medskip

\noindent (6.64) \quad {\it We move now to the $p_2^{-1} \, U_{\beta_n} \subset \Theta^4 (\Theta^3 (\pi f X^2-H)' , {\mathcal R})$ and to the $p_2 \, p_1 \, \Lambda_{j_n} \subset \Theta^3 (\pi f X^2-H)' ) \cap U_{\beta_n}$. There is a homotopy $[p_1 \, \Lambda_{j_n} , p_2 \, p_1 \, \Lambda_{j_n}] \subset p_2^{-1} (\Theta^3 \cap U_{\beta_n})$, which makes use of the union of fibers
$$
\bigcup_{\overbrace{z \, \in \, U_{\beta_n}}} \ p_2^{-1} (z) \, .
$$

By now, of course, we have $p_2 \, p_1 \, \Lambda_{j_n} \subset U_{\beta_n} \cap \Theta^3$. Also, with our $j_n$ which is supposed to be very high, the homotopy $[p_1 \, \Lambda_{j_n} , p_2 \, p_1 \, \Lambda_{j_n}]$ avoids $K$, by which we mean the various incarnations of $K$, via $p_1$ or $p_2 \, p_1$, at the level of }(6.31.3).

\medskip

\noindent (6.65) \quad {\it With all this, the next homotopy
$$
[p_2 \, p_1 \, \Lambda_{j_n} , \rho_1 \, p_2 \, p_1 \, \Lambda_{j_n}] \subset U_{\beta_n} \cap \Theta^3 (\pi f X^2-H)' \, ,
$$
which uses the $\rho_t$ from} (\ref{eq6.50}), {\it also avoids the compact $K$. At the end of this homotopy, we find
$$
\rho_1 \, p_2 \, p_1 \, \Lambda_{j_n} \subset ({\mathcal N}^3_{N(\alpha_n)})' \, ,
$$
with the $({\mathcal N}^3_{N(\alpha_n)})'$ like in} (\ref{eq6.49}), {\it with the $\alpha$ from} (\ref{eq6.49}) {\it changed to our very high $\alpha_n$.}

\medskip

\noindent 2) {\it In the context of } (\ref{eq6.37}) $+$ (\ref{eq6.38}) {\it we move now to
$$
\rho_1 \, p_2 \, p_1 \, \Lambda_{j_n} \subset r_{N(\alpha_n)}^{-1} \, {\mathcal V}_{\gamma_n} \subset ({\mathcal N}^3_{N(\alpha_n)})' \subset {\mathcal X}^3 \overset{R}{-\!\!\!-\!\!\!\longrightarrow} \bigcup_{R_0} \, R_0 \times \{ \infty \} \underset{\rm TOP}{=} \sum (\infty)_* \, .
$$
The first inclusion in the formula above, makes use of formula} (\ref{eq6.51}) {\it in the sublemma~{\rm 6.5.5.} The big retraction $R$ provides us with a homotopy
\setcounter{equation}{65}
\begin{equation}
\label{eq6.66}
[\rho_1 \, p_2 \, p_1 \, \Lambda_{j_n} , R \, \rho_1 \, p_2 \, p_1 \, \Lambda_{j_n}] \subset {\mathcal X}^3 \, ,
\end{equation}
localized in a close enough neighbourhood of $\partial_{\infty}$ so as to avoid again $K$, and at the end of this last homotopy, we have}
$$
R \, \rho_1 \, p_2 \, p_1 \, \Lambda_{j_n} \subset \bigcup_{R_0} \, R_0 \times \{\infty\} = \sum(\infty)_* \times \{\infty\} \, .
$$

\bigskip

Notice that all our homotopies in the sub-sub-lemma above are taking place close enough to $\partial_{\infty}$, so as to avoid $K$.

\bigskip

\noindent {\bf Sub-sub-lemma 6.5.7.} {\it For $n$ high enough, the
$$
\Lambda_{j_n} \subset \partial S'_u (M(\Gamma)-H)_{\rm I} \subset \partial S'_u (M(\Gamma)-H)_{\rm II} \subset \partial S'_u (M(\Gamma)-H)^{\wedge}_{\rm II}
$$
cobounds a disk $\overline D_{j_n}^2 \subset \partial S'_u (M(\Gamma)-H)^{\wedge}_{\rm II}$ which does not touch the compact set $K \subset \partial S'_u (M(\Gamma)-H)_{\rm II} \subset \partial S'_u (M(\Gamma)-H)^{\wedge}_{\rm II}$.

\smallskip

On can perturb $\overline D_{j_n}^2$ rel $\Lambda_{j_n} = \partial \, \overline D_{j_n}^2$ so as to change it into a disk 
$$
\overset{=}{D}_{j_n}^2 \subset \partial S'_u (M(\Gamma)-H)_{\rm II} - K \, .
$$
This contradicts the statement $(P_2)$, and hence it proves the MAIN LEMMA}~6.5.

\bigskip

\noindent {\bf Proof.} By (4.50.1), we know already that $\Lambda_{j_n}$ is null-homotopic inside $\partial S'_u (M(\Gamma)-H)_{\rm I}$. Hence, it is also null-homotopic inside the $\partial S'_u (M(\Gamma)-H)_{\rm II} \subset \partial S'_u (M(\Gamma)-H)_{\rm II}^{\wedge}$. But the problem, of course, is that the generic singular disk cobounding $\Lambda_{j_n}$ may happily cut through $K$.

\smallskip

So, here is how we produce our desired $\overline D_{j_n}^2$. Putting together the homotopies in (6.63), (6.64), (6.65), (\ref{eq6.66}), we get a log homotopy $[\Lambda_{j_n} , R \, \rho_1 \, p_2 \, p_1 \, \Lambda_{j_n}] \subset \partial S'_u (M(\Gamma)-H)_{\rm II}^{\wedge}$ {\ibf avoiding} $K$. Here, we have
$$
R \, \rho_1 \, P_2 \, P_1 \, \Lambda_{j_n} \subset \sum (\infty) \times \{\infty\} \subset \partial S'_u (M(\Gamma)-H)_{\rm II}^{\wedge} - \partial S'_u (M(\Gamma)-H)_{\rm II} \, .
$$

Remember at this point that through our $\sum (\infty)$ go the cut arcs $[\alpha (\infty) , \beta)$, which are avoided by the zipping flow. Without any loss of generality, it may be assumed that we also have
\begin{equation}
\label{eq6.67}
R \, \rho_1 \, P_2 \, P_1 \, \Lambda_{j_n} \subset \left(\sum (\infty)_* \times \{\infty \}) ({\rm cut}\right) \, .
\end{equation}
Invoking now the $\pi_1$-injectivity from the sub-lemma~6.5.3, we can find a singular disk $d_n^2 \to (\sum (\infty)_* \times \{\infty \}) ({\rm cut})$, such that $\partial d_n^2 = R \, \rho_1 \, P_2 \, P_1 \, \Lambda_{j_n}$. With this, we may take now
$$
\overline D_{j_n}^2 = [\Lambda_{j_n} , R \, \rho_1 \, P_2 \, P_1 \, \Lambda_{j_n}] \cup d_n^2 \, .
$$

This ENDS the proof of the sub-sub-lemma~6.5.7. But then, the compactness lemma~4.7 is by now proved too; see below. And, as already said earlier, this implies the GSC theorem~2.3, i.e. that
$$
S_u \, \widetilde M (\Gamma)_{\rm II} \in {\rm GSC} \, .
$$

\bigskip

\noindent FINAL COMMENTS CONCERNING THE MAIN LEMMA.

\smallskip

Our MAIN LEMMA~6.5 certainly has a certain $\pi_1^{\infty} = 0$ flavour, but this is deceptive. The lemma in question has nothing to do with the issue of $\pi_1^{\infty} \, \Gamma$, as I will explain, and certainly NO CLAIM concerning $\pi_1^{\infty} \, \Gamma$ has been made in this paper.

\smallskip

The whole argument in this section, in particular in connection with the MAIN LEMMA~6.5, had to do with
$$
\partial S'_u (M(\Gamma) - H)_{\rm II} \subset \{\mbox{the smooth {\ibf compact} manifold} \ \Theta^4 (M(\Gamma) , {\mathcal R}) \times B^L \, , \ L \ {\rm high}\} \, ,
$$
and these objects do not know, neither about the $\Gamma$-actions nor about $\widetilde M (\Gamma)$ and its asymptotic topology, where one has to look for the mysterious $\pi_1^{\infty} \, \Gamma$. For getting to those items one would have to go to the universal covering spaces, and that is NOT what our discussion has been all about in the present section.

\smallskip

Next, in a more technical vein, the spaces $({\mathcal N}^3)'$ and $({\mathcal N}_n^3)'$, which live inside the non compact $\Theta^3 (\pi f X^2 - H)'_{\rm I}$ are neither compact, nor with compact complements, so they are not the kind of neighbourhoods of infinity which the $\pi_1^{\infty}$-issue is all about. Anyway, what we have dealt with in this section is the man-created infinity of $S'_u (M(\Gamma)-H)$ and NOT the God-created infinity of $\widetilde M (\Gamma)$.

\smallskip

Now that we have our MAIN LEMMA we shall give

\bigskip

\noindent THE PROOF OF THE COMPACTNESS LEMMA~4.7. What we need to show in the context of (4.47.II) is that
$$
\eta \circ \beta \mid \sum_n \, C^- (H_n) \quad {\rm and} \quad \alpha \mid \sum_n C^- (H_n)
$$
are connected by a PROPER homotopy. The MAIN LEMMA tells us that
$$
\Lambda_n = \{ \alpha \, C^- (H_n) \underset{\overbrace{\gamma_n}}{\bullet} \eta \, \beta \, C^- (H_n) \} \subset \partial S'_u (M(\Gamma) - H)_{\rm II} \, ,
$$
which are such that $\underset{n=\infty}{\lim} \, \Lambda_n = \infty$ are, at the same time, cobounded by singular disks $D_n^2$ also coming with $\underset{n=\infty}{\lim} \, D_n^2 = \infty$. In view of the dumb-ball geometry of $\Lambda_n$, these $D_n^2$'s can be viewed as singular cylinders $(S^1 \times I)_n$, which connect the two closed curves, and which are such that $\underset{n=\infty}{\lim} \, (S^1 \times I)_n =\infty$. The little drawing below should illustrate this.

$$
\includegraphics[width=10cm]{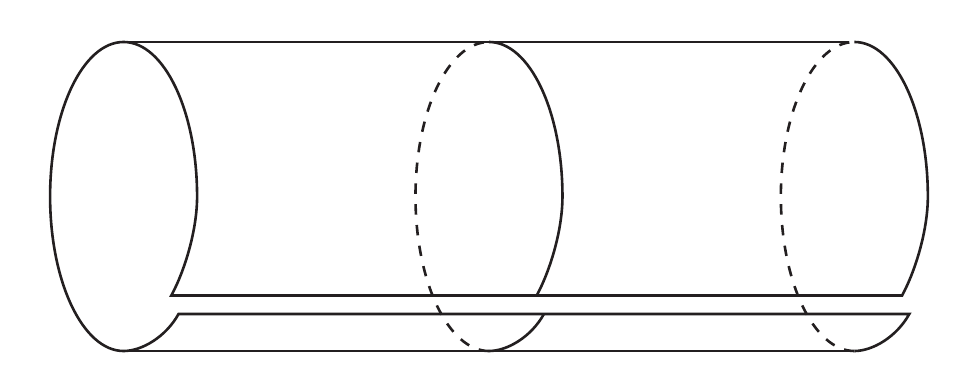}
$$
\label{figEND}


\end{document}